\documentclass[sn-mathphys,Numbered]{sn-jnl}

\usepackage{graphicx}
\usepackage{multirow}
\usepackage{amsmath,amssymb,amsfonts}
\usepackage{amsthm}
\usepackage{mathrsfs}
\usepackage[title]{appendix}
\usepackage{xcolor}
\usepackage{textcomp}
\usepackage{manyfoot}
\usepackage{booktabs}
\usepackage{algorithm}
\usepackage{algorithmicx}
\usepackage{algpseudocode}
\usepackage{listings}

\usepackage{subcaption}
\usepackage{url}
\usepackage{color}

\graphicspath{{figs/}}

\theoremstyle{thmstyleone}

\theoremstyle{thmstyletwo}

\theoremstyle{thmstylethree}

\begin{document}

\title[ADER-DG method with local DG predictor for IVP ODE systems]{%
	Arbitrary high order ADER-DG method with local DG predictor for 
	solutions of initial value problems for systems of first-order ordinary differential equations
}

\author*[1]{\fnm{Ivan S.} \sur{Popov}}\email{diphosgen@mail.ru, popovis@omsu.ru}

\affil*[1]{%
	\orgdiv{Department of Theoretical Physics}, 
	\orgname{Dostoevsky Omsk State University}, 
	\orgaddress{%
		\street{Mira prospekt}, 
		\city{Omsk}, 
		\postcode{644077}, 
		\country{Russia}
	}
}

\abstract{%
An adaptation of the arbitrary high order ADER-DG numerical method with local DG predictor for solving the IVP for a first-order non-linear ODE system is proposed. The proposed numerical method is a completely one-step ODE solver with uniform steps, and is simple in algorithmic and software implementations. It was shown that the proposed version of the ADER-DG numerical method is $A$-stable and $L$-stable. The ADER-DG numerical method demonstrates superconvergence with convergence order $2N+1$ for the solution at grid nodes, while the local solution obtained using the local DG predictor has convergence order $N+1$. It was demonstrated that an important applied feature of this implementation of the numerical method is the possibility of using the local solution as a solution with a subgrid resolution, which makes it possible to obtain a detailed solution even on very coarse coordinate grids. The scale of the error of the local solution, when calculating using standard representations of single or double precision floating point numbers, using large values of the degree $N$, practically does not differ from the error of the solution at the grid nodes. The capabilities of the ADER-DG method for solving stiff ODE systems characterized by extreme stiffness are demonstrated. Estimates of the computational costs of the ADER-DG numerical method are obtained.
}

\keywords{%
	discontinuous Galerkin method,
	ADER-DG method,
	local DG predictor,
	first-order ODE systems,
	superconvergence
}

\pacs[MSC Classification]{65L05, 65L60, 65L20}

\maketitle

\section*{Introduction}
\label{intro}

In this paper, a study of the arbitrary high order ADER discontinuous Galerkin (DG) method with local DG predictor, which is frequently used for solving problems for partial differential equations, is based on solving the initial value problems (IVP) for the first-order non-linear ordinary differential equation (ODE) system chosen in the following form
\begin{equation}\label{eq:ivp_chosen_form}
\frac{d\mathbf{u}}{dt} = \mathbf{F}\left(\mathbf{u},\, t\right),\quad t\in\Omega = \left\{t\, |\, t \in [t_{0},\ t_{0} + T]\right\},\quad \mathbf{u}(t_{0}) = \mathbf{u}_{0},
\end{equation}
where $\mathbf{u}: \Omega \rightarrow \mathbb{R}^{K}$ is a desired function; the function $\mathbf{F}: \mathbb{R}^{K}\times\Omega \rightarrow \mathbb{R}^{K}$ is a right side function, which is given. The initial condition $\mathbf{u}_{0}$ on the desired function $\mathbf{u}$ was chosen at the point $t_{0}$. The classical ODE theory shows that in case $\mathbf{F} \in C_{1}(\mathbb{R}^{K}\times\Omega)$ the solution of the problem exists and is unique. It is clear that the IVP for equations of the second and higher orders, uniquely solvable with respect to higher derivatives, as well as for systems of such equations, can be represented in the chosen form of the first-order non-linear ODE system (\ref{eq:ivp_chosen_form}).

The problem of obtaining a numerical IVP for an ODE system occurs in a wide range of scientific and technical problems. At the moment, there are many numerical methods for solving the initial value problem for ODE systems~\cite{Butcher_book_2016, Hairer_book_1, Hairer_book_2}. In recent years, the study of DG methods and their superconvergence has been an active research field in numerical analysis, see books~\cite{Babuska_book_2001, Wahlbin_lectures_1995} and recent articles~\cite{dg_ivp_ode_4, dg_ivp_ode_5, dg_ivp_ode_6}. DG methods occupy a special place among the numerical methods for solving differential equations, since they allow obtaining a numerical solution of very high accuracy, while being sufficiently simple and scalable to solve complex and large problems. DG methods were created in the work~\cite{lasl_rep_dg_1973} for numerical solution of the neutron transport equations. In 1981 Delfour \textit{et al} \cite{Delfour_1981} constructed a Runge-Kutta-like DG method for solving the IVP for the ODE system that demonstrates superconvergence up to order $2N+2$ for polynomials with a degree $N$. Cockburn, Shu \textit{et al} in a series of works~\cite{Cockburn_base_1, Cockburn_base_2, Cockburn_base_3, Cockburn_base_4, Cockburn_base_5} created an accurate and thoroughly developed mathematical basis of DG methods, which stimulated their further development and use for solving a wide class of problems. DG methods are used to solve the IVP for ODE systems~\cite{dg_ivp_ode_4, dg_ivp_ode_5, dg_ivp_ode_6, dg_ivp_ode_1, dg_ivp_ode_2, dg_ivp_ode_3}, the boundary value problem for ODE systems~\cite{dg_bvp_ode_1, dg_bvp_ode_2}, to solve stochastic differential equations~\cite{dg_stoch_ode_1, dg_stoch_ode_2}, and especially widely to solve PDEs~\cite{ader_dg_ideal_flows, ader_dg_diss_flows, ader_dg_ale, ader_dg_grmhd, ader_dg_gr_prd, ader_dg_gr_z4_2024, ader_dg_simple_mod, ader_dg_PNPM, fron_phys, ader_dg_axioms, exahype, PNPM_DG_2009}.

In this work, an adaptation of the arbitrary high order ADER-DG numerical method~\cite{ader_dg_ideal_flows, ader_dg_diss_flows, ader_dg_PNPM, fron_phys, ader_dg_axioms} with local DG predictor for solving the IVP for the first-order non-linear ODE system is proposed. In 2008 Dumbser \textit{et al}~\cite{ader_stiff_1} proposed a paradigm of local solution and DG predictor and showed the possibility of using the numerical method to obtain a sufficiently accurate local solution of IVP for a linear scalar ODE in comparison with using the Cauchy-Kovalewski procedure (using Taylor series)~\cite{ader_init_1, ader_init_2}, which was further used in~\cite{ader_stiff_1, ader_stiff_2} and in subsequent works related to the use of a local DG predictor. In 2010 Dumbser~\cite{PNPM_DG_2010} proposed an effective way to use the ADER-DG method to solve initial value and boundary value problems for ODE systems, which was demonstrated by solving the boundary layer equation. The implementation of the method is based on the use of the nodal basis of Lagrange interpolation polynomials with nodal points at the nodes of the Gauss-Legendre quadrature formula. This made it possible to use the properties of the nodal basis associated with the lack of calculation of the integrals of the function $\mathbf{F}$ on the right side of the ODE system (\ref{eq:ivp_chosen_form}). An important feature of this variant of the DG method is the significant simplicity of its algorithmic and software implementations, which in no way reduces the accuracy and efficiency of the method. It was shown that the proposed version of the ADER-DG numerical method is $A$-stable and $L$-stable, and implements superconvergence with convergence order $p = 2N+1$ for the solution at grid nodes, while the local solution obtained using the local DG predictor has convergence order $p = N+1$, where $N$ is the basis polynomials degree, which is usually expected from classical implementations of DG methods for solving ODE. Baccouch in~\cite{dg_ivp_ode_1, dg_ivp_ode_3} obtained rigorous proofs for the convergence orders $2N+1$ and $N+1$ for DG methods presented in a modal functional representation based on Legendre polynomials, which is expected to be correct for the ADER-DG numerical method studied in this paper. It should be noted that this numerical method is fully one-step ODE solver, with uniform steps in terms of formula implementation.

It should be noted that in 2021 Han Veiga \textit{et al}~\cite{dec_vs_ader} showed that the numerical methods of the ADER family are significantly interconnected with numerical methods based on the deferred correction (DeC) paradigm. DeC methods has a long history and its application to initial value problems for ODE systems goes back to~\cite{dec_src_1968}, and are effectively used to solve both ODE systems~\cite{dec_dutt_2000, dec_dutt_2000, dec_minion_2003, dec_shu_2008} and partial differential equations~\cite{dec_abgrall_2017, dec_abgrall_2019}. DeC methods, like the methods of the ADER family, allow one to obtain an arbitrarily high order and are characterized by a high accuracy of the numerical solution. The current state of research on DeC methods~\cite{dec_vs_ader} shows that they compete with methods of the ADER family. The work notes~\cite{dec_vs_ader} that the explicit ADER method can be considered as a special interpretation of the DeC paradigm presented in the work; it is also noted that the implementation of numerical methods of DeC is easier to implement than ADER methods.

The numerical methods of ADER-DG, based on the use of the ADER paradigm, allow the creation of numerical methods of an arbitrarily high order. Among the existing numerical methods, we can distinguish high-order numerical methods based on the Taylor expansion for solving ODE systems~\cite{ivp_ode_taylor_series_2017, ivp_ode_taylor_series_soft_2005}, which allow obtaining an arbitrary high order. High-order numerical Taylor methods, in this context, in their general structure are close to the use of the Cauchy-Kovalevskaya procedure, which was used in the original versions of the development of the ADER paradigm~\cite{ader_init_1, ader_init_2}, where this procedure was used to obtain a local solution. Modern implementations of high-precision methods based on the ADER paradigm use a local solution obtained by a local DG predictor~\cite{ader_stiff_1, ader_stiff_2, ader_weno_lstdg_ideal, ader_weno_lstdg_diss, ader_dg_dev_1, ader_dg_dev_2, ader_dg_ideal_flows, ader_dg_diss_flows}. Approximate Taylor methods allow arbitrarily high order, and there are efficient high-performance software implementations for using them~\cite{ivp_ode_taylor_series_soft_2005}. However, in comparison with numerical methods of ADER-DG, approximate Taylor methods have their own disadvantages, which led to the transition from using the Cauchy-Kovalevskaya procedure to using a local DG predictor to obtain a local solution in numerical methods based on the use of the ADER paradigm. Also among the existing numerical methods, high-order classical extrapolation methods for solving ODE~\cite{Hairer_book_1, Hairer_book_2, comp_rk_extr_dec} can be distinguished, which also make it possible to obtain an arbitrary high order. In contrast to the numerical methods of ADER-DG and the approximate Taylor methods, the use of extrapolation methods for solving ODE usually does not allow obtaining solutions only at grid nodes, while the numerical methods of ADER-DG and approximate Taylor methods allow obtaining a high-order numerical solution and in domains between grid nodes.

An important applied feature of this implementation of the ADER-DG numerical method is the possibility of using the local solution as a solution with a subgrid resolution, which makes it possible to obtain a detailed solution even on very coarse coordinate grids. The scale of the error of the local solution, when calculating using standard representations of single or double precision floating point numbers, using large values of the degree $N$, practically does not differ from the error of the solution at the grid nodes. It should be noted that the local solution is not some kind of reconstruction of the solution from the set of values at the nodes, it is precisely the representation of the solution in the form of a finite element. Further in the text of this paper, a detailed description of the numerical method, analysis of its stability, a demonstration example of the possibility of using a local solution as output numerical solution, and calculation of the convergence orders $p$ are presented.

\section{General description of the numerical method}
\label{sec:1}

\subsection{Formulation of the numerical method}
This paragraph presents the mathematical apparatus of the arbitrary high order ADER-DG numerical method with local DG predictor for solutions of IVP for first-order non-linear ODE systems studied in this work. The method proposed below is an adaptation of the arbitrary high order ADER-DG numerical method~\cite{ader_dg_ideal_flows, ader_dg_diss_flows, ader_dg_PNPM, fron_phys, ader_dg_axioms} with local DG predictor for solving the IVP for a first-order non-linear ODE system.

The numerical solution to the IVP for a first-order non-linear ODE system will be found on partitioning the domain of definition $\Omega$ of the desired function by a finite number of non-overlapping discretization domains $\Omega_{n} = \{t\, |\, t \in [t_{n},\, t_{n+1}]\}$, where $t_{n+1} = t_{n} + \Delta t_{n}$, and $\Delta t_{n}$ is the discretization step, which is not assumed to be constant for different discretization domains $\Omega_{n}$; therefore, the proposed numerical method allows implementation for a variable step. The set of discretization domains $\Omega_{n}$ covering the domain of definition $\Omega$ of the desired function represents a one-dimensional grid. Further points $t_{n}$ will be denoted by grid nodes, the space between nodes $t \in [t_{n},\, t_{n+1}]$ will be denoted by the domain between nodes.

The ADER-DG numerical method for solving the IVP for a first-order non-linear ODE system chosen in the form (\ref{eq:ivp_chosen_form}) is based on using the integral form of writing a system of differential equations in discretization domains $\Omega_{n}$:
\begin{equation}\label{eq:int_form_tn_tnp1}
\mathbf{u}_{n+1} = \mathbf{u}_{n} + \int\limits_{t_{n}}^{t_{n+1}} \mathbf{F}\left(\mathbf{u}(t),\, t\right) dt,
\end{equation}
where $\mathbf{u}_{n}$ and $\mathbf{u}_{n+1}$ denote the numerical solution at the discretization node $t_{n}$ and $t_{n+1}$, respectively (further referred to simply as the solution at the nodes); in the integrand, the function $\mathbf{u} = \mathbf{u}(t)$ is defined in the space between nodes $t \in [t_{n},\, t_{n+1}]$. The transition was made to the local coordinate system of the discretization domain $\Omega_{n}$:
\begin{equation}\label{eq:rescale}
t(\tau) = t_{n} + \Delta t \cdot \tau,
\end{equation}
where $\tau$ is a local coordinate that allows one to encapsulate the coordinates of a particular discretization domain and carry out the transformation to the reference domain $\tau \in [0, 1]$. In the case of a reference domain, the integral form (\ref{eq:int_form_tn_tnp1}) of a system of differential equations can be represented in the following form
\begin{equation}\label{eq:int_form_ref}
\mathbf{u}_{n+1} = \mathbf{u}_{n} + \int\limits_{0}^{1} \mathbf{f}\left(\mathbf{u}\left(t(\tau)\right),\, t(\tau)\right) d\tau,
\end{equation}
where $\mathbf{f} = \Delta t_{n} \cdot \mathbf{F}$ is the rescaled, to the reference domain, function of the right side of the ODE system (\ref{eq:ivp_chosen_form}), expression $t = t(\tau)$ is defined by the relation (\ref{eq:rescale}). 

In the space between nodes $t \in [t_{n},\ t_{n+1}]$, a local discrete solution $\mathbf{q}$ is introduced:
\begin{equation}\label{eq:local_sol}
\mathbf{u} = \mathbf{u}(t(\tau))\ \mapsto\ \mathbf{q} = \mathbf{q}(\tau),\ \tau\in[0, 1],
\end{equation}
which will replace the function $\mathbf{u}$ in the integrand of the integral form of the ODE system (\ref{eq:int_form_ref}); $\mathbf{q}$ is a $K$-component vector function, such as $\mathbf{u}$. The local discrete solution $\mathbf{q}$ will satisfy IVP for the ODE system in the weak form:
\begin{equation}\label{eq_ode_local_weak}
\frac{d\mathbf{q}}{d\tau} = \mathbf{f}\left(\mathbf{q},\, t(\tau)\right),\quad \mathbf{q}(0) = \mathbf{u}_{n},
\quad\Rightarrow\quad
\int\limits_{0}^{1} \varphi_{p}(\tau) \left[\frac{d\mathbf{q}(\tau)}{d\tau} - \mathbf{f}\left(\mathbf{q}(\tau),\, t(\tau)\right)\right] d\tau = 0,
\end{equation}
where $\varphi_{p}(\tau)$ is the set of basic functions of the DG representation, and $p$ is the index of the basis function in the set. The weak form represents the condition of orthogonality of the residual of the local solution $\mathbf{q}$ of the ODE system (\ref{eq:int_form_tn_tnp1}) to all functions of the basis set in domain $\Omega_{n}$. The resulting expression (\ref{eq_ode_local_weak}), prepared for computing a local discrete solution, represents a local DG predictor.

In this work, the Legendre interpolation polynomials were chosen as the nodal basis functions $\varphi_{p}(\tau)$, with nodal points at the nodes of the Gauss-Legendre quadrature formula in the range $\tau\in[0, 1]$. The interpolation condition was the expression $\varphi_{p}(\tau_{l}) = \delta_{p, l}$, where $\delta_{p, l}$ is the delta symbol, $\tau_{l}$ are the roots of classical shifted Legendre polynomials $L_{N+1}(\tau)$ of degree $N+1$; indices $(p, l)$ are in integer range $0 \leqslant p,\, l \leqslant N$.

The representation of the local discrete solution $\mathbf{q}$ was chosen in the form of expansion in basis polynomials $\varphi_{p}(\tau)$:
\begin{equation}\label{eq:q_repr}
\mathbf{q}(\tau) = \sum\limits_{p = 0}^{N} \hat{\mathbf{q}}_{p}\cdot\varphi_{p}(\tau),
\end{equation}
where the expansion coefficients $\hat{\mathbf{q}}_{p}$ are the $K$-component vectors; the full set of coefficients $\hat{\mathbf{q}}_{p}$, $0 \leqslant p \leqslant N$, in the representation of a local discrete solution $\mathbf{q}$ can be represented by a matrix with sizes $(N+1) \times K$. The function $\mathbf{f}$ of the weak form of the ODE system (\ref{eq_ode_local_weak}) was also presented in the form of an expansion in terms of a set of polynomials:
\begin{equation}
\mathbf{f}\left(\mathbf{q}(\tau),\, t(\tau)\right) = \sum\limits_{p = 0}^{N} \hat{\mathbf{f}}_{p}\cdot\varphi_{p}(\tau)\ \mapsto\
\sum\limits_{p = 0}^{N} \hat{\mathbf{f}}\left(\hat{\mathbf{q}}_{p},\, t\left(\tau_{p}\right)\right)\cdot\varphi_{p}(\tau),
\end{equation}
where the point-wise representation was used~\cite{Jackson_2017} for the expansion coefficients $\hat{\mathbf{f}}_{p}$ for the subsequent substitution into the integral relations (\ref{eq_ode_local_weak}), which is an effective property of nodal bases.

The weak form (\ref{eq_ode_local_weak}) of the ODE system, after integration by parts and substitution of representation (\ref{eq:q_repr}), was rewritten as follows:
\begin{equation}\label{dg_pred_src}
\sum\limits_{q = 0}^{N} \Big[\mathrm{K}_{pq} \hat{\mathbf{q}}_{q} - \mathrm{M}_{pq} \mathbf{f}\left(\hat{\mathbf{q}}_{q},\, t\left(\tau_{q}\right)\right)\Big] = 
\varphi_{p}(0)\, \mathbf{u}_{n},
\end{equation}
where the matrices $\mathrm{K} = ||\mathrm{K}_{pq}||$ and $\mathrm{M} = ||\mathrm{M}_{pq}||$ were expressed and calculated in terms of the coefficients of the basis polynomials:
\begin{equation}
\mathrm{K}_{pq} = \varphi_{p}(1)\varphi_{q}(1) - \int\limits_{0}^{1} \frac{d\varphi_{p}(\tau)}{d\tau}\, \varphi_{q}(\tau)\, d\tau,\quad
\mathrm{M}_{pq} = \int\limits_{0}^{1} \varphi_{p}(\tau)\, \varphi_{q}(\tau)\, d\tau;
\end{equation}
it should be noted that the functional basis $\varphi_{p}(\tau)$ is orthogonal in $L_{2}$, so the mass matrix $\mathrm{M}$ is diagonal: $\mathrm{M}_{pq} = \mathrm{M}_{p}\delta_{pq}$.
The final expression, convenient for the implementation of the local DG predictor, was obtained in the following form:
\begin{equation}\label{dg_pred}
\hat{\mathbf{q}}_{p} - 
\sum\limits_{q = 0}^{N} \mathrm{B}_{pq} \cdot \mathbf{f}\left(\hat{\mathbf{q}}_{q},\, t\left(\tau_{q}\right)\right) = \mathbf{u}_{n},
\end{equation}
where the notation for the matrix $\mathrm{B} = \mathrm{K}^{-1} \cdot \mathrm{M}$ was introduced, and the property
\begin{equation}
\sum\limits_{q = 0}^{N} \left[\mathrm{K}^{-1}\right]_{pq} \varphi_{q}(0) = 1,\quad 0 \leqslant p \leqslant N,
\end{equation}
was also used. The expression (\ref{dg_pred}) implements a local DG predictor and represents a system of equations for the coefficients $\hat{\mathbf{q}}_{p}$ of a local discrete solution $\mathbf{q}$. In the case of a linear in $\mathbf{u}$ function $\mathbf{F} = \mathbf{F}(\mathbf{u}, t)$ on the right side of the ODE system (\ref{eq:ivp_chosen_form}), the system of algebraic equations (\ref{dg_pred}) is also linear; in the general case of a non-linear function $F$, the system of equations (\ref{dg_pred}) is also non-linear. It is convenient to solve the system of equations in the non-linear case by Newton's method or by the Picard iteration method. The initial iteration can be chosen in the form of an initial condition $\mathbf{q}_{p} = \mathbf{u}_{n}$, $0 \leqslant p \leqslant N$, or more efficient, in terms of convergence speed, methods for choosing the initial iteration can be used~\cite{ader_dg_axioms}.

The local solution $\mathbf{q}$ obtained in the discrete domain $\Omega_{n}$ using the local DG predictor (\ref{dg_pred}) was used to obtain a numerical solution $\mathbf{u}_{n+1}$ at the grid node $t_{n+1}$ using the integral form of the ODE system (\ref{eq:int_form_ref}) and the Gauss-Legendre quadrature formula:
\begin{equation}\label{eq:sol_in_node}
\mathbf{u}_{n+1} = \mathbf{u}_{n} + \sum\limits_{p = 0}^{N} w_{p}\, \mathbf{f}\left(\hat{\mathbf{q}}_{p},\, t\left(\tau_{p}\right)\right),
\end{equation}
where $w_{p}$ are the weights of the Gauss-Legendre quadrature formula in the range $[0, 1]$ (shifted Legendre polynomials) and also the defining property $\varphi_{p}(\tau_{l}) = \delta_{p, l}$ of the interpolation basis polynomials $\varphi_{p}(\tau)$ was used.

\subsection{Stability of the numerical method}

The linear stability of the ADER-DG numerical method with a local DG prector was investigated. Linear stability was studied based on the analysis of the Dahlquist's test equation $\dot{u} = \lambda u$. The numerical solution of the equation has the form $u_{n+1} = R(\lambda\cdot\Delta t_{n})\, u_{n}$, where $R = R(z)$ is called the stability function of the numerical method. The stability function argument $z = \lambda\cdot\Delta t_{n}$ is a complex number: $z = \Re(z) + i\cdot\Im(z) \in \mathbb{C}$. Using the equation (\ref{dg_pred}) and formula (\ref{eq:sol_in_node}), the following result was obtained:
\begin{equation}
\sum\limits_{q = 0}^{N}\left[\delta_{pq} - \left(\lambda\cdot\Delta t_{n}\right)\mathrm{B}_{pq}\right]\cdot\hat{q}_{q} = u_{n},\qquad
u_{n+1} = u_{n} + \left(\lambda\cdot\Delta t_{n}\right) \sum\limits_{p = 0}^{N} w_{p}\, \hat{q}_{p},
\end{equation}
which allows using the well-known expression used to analyze the stability of implicit Runge-Kutta methods using matrix coefficients $\mathrm{B} = ||\mathrm{B}_{pq}||$ and vector weights $\mathbf{w} = [w_{0} \ldots w_{p}]^{T}$ of the Gauss-Legendre quadrature formula:
\begin{equation}\label{eq:stab_func_expr}
R(z) = 1 + \mathbf{w}^{T} \cdot \left[\mathrm{E} - z\, \mathrm{B}\right]^{-1} \cdot \mathbf{1} 
=\frac{\left|\mathrm{E} - z\, \mathrm{B} + z\, \mathbf{1}\otimes\mathbf{w}\right|}{\left|\mathrm{E} - z\mathrm{B}\right|},
\end{equation}
where $\mathrm{E}$ is the identity matrix with sizes $(N+1)\times(N+1)$, $\mathbf{1} = [1 \ldots 1]^{T}$ is a vector of ones with size $N+1$. 

\begin{figure}[h!]
\captionsetup[subfigure]{%
position=bottom,
textfont=normalfont,
singlelinecheck=off,
justification=raggedright,
font=tiny
}
\centering
\begin{subfigure}{0.48\textwidth}
\includegraphics[width=\textwidth]{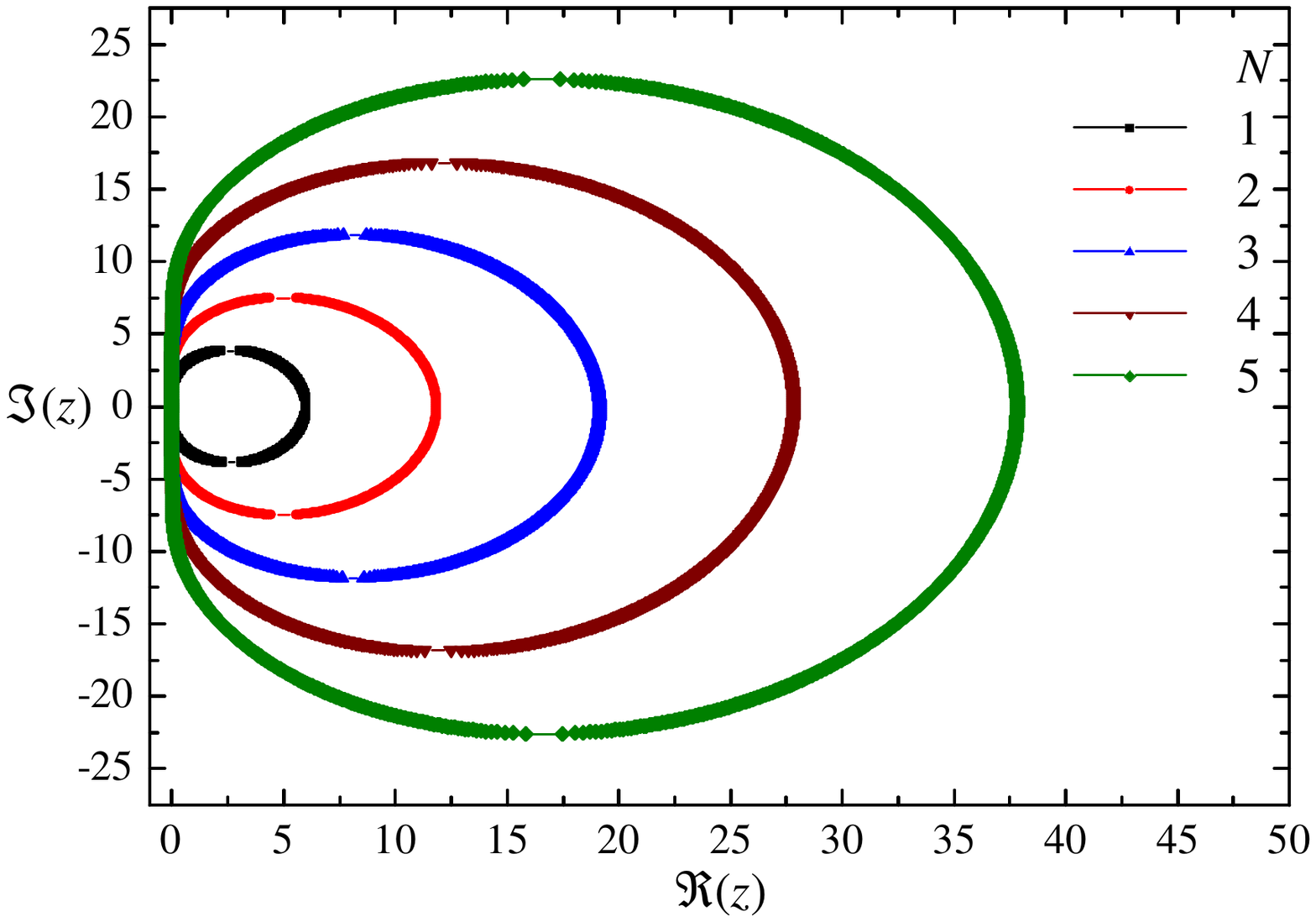}
\vspace{-8mm}\caption{\label{fig:stab_domain:a}\vspace{2mm}}
\end{subfigure}
\begin{subfigure}{0.48\textwidth}
\includegraphics[width=\textwidth]{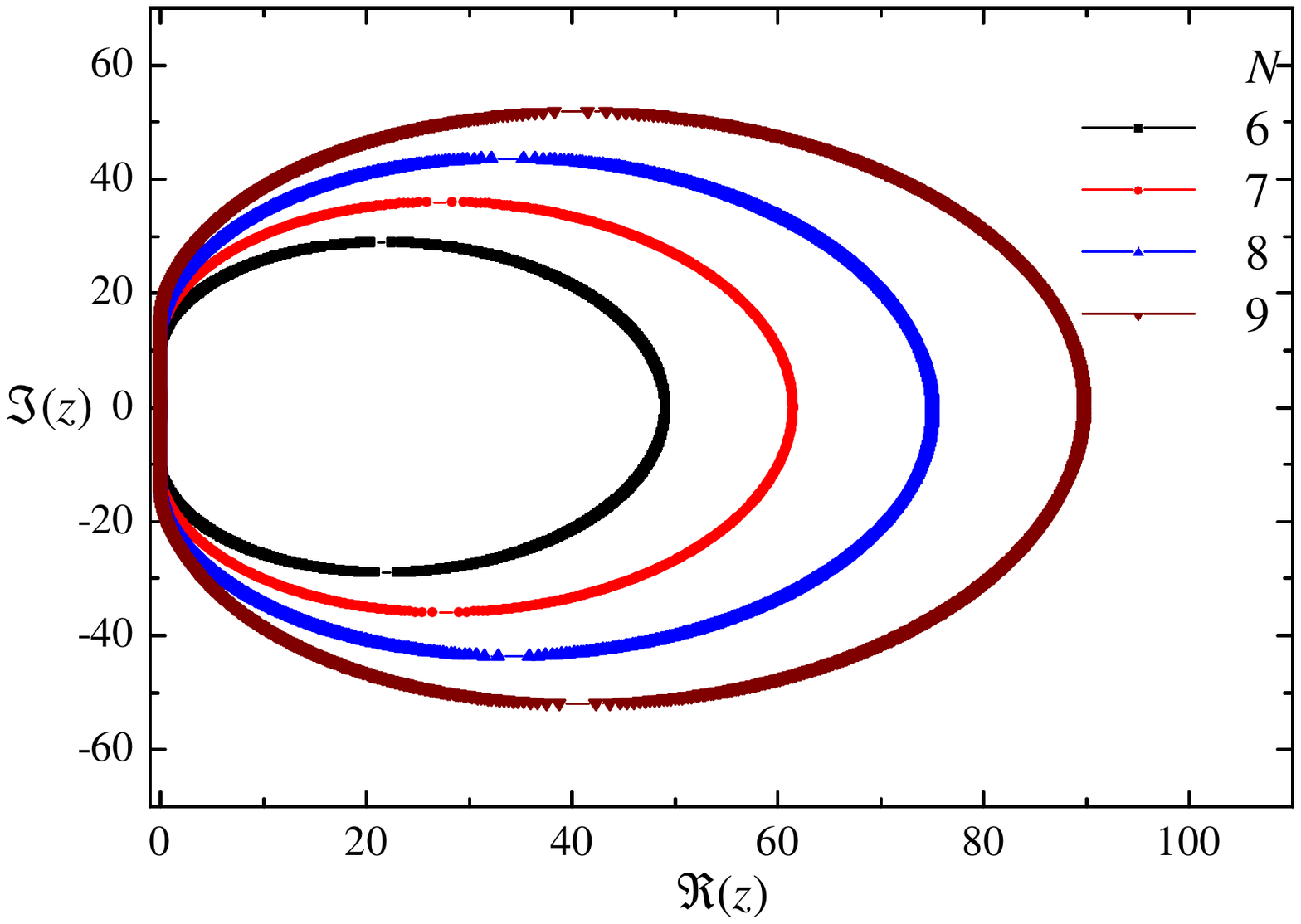}
\vspace{-8mm}\caption{\label{fig:stab_domain:b}\vspace{2mm}}
\end{subfigure}\\
\begin{subfigure}{0.48\textwidth}
\includegraphics[width=\textwidth]{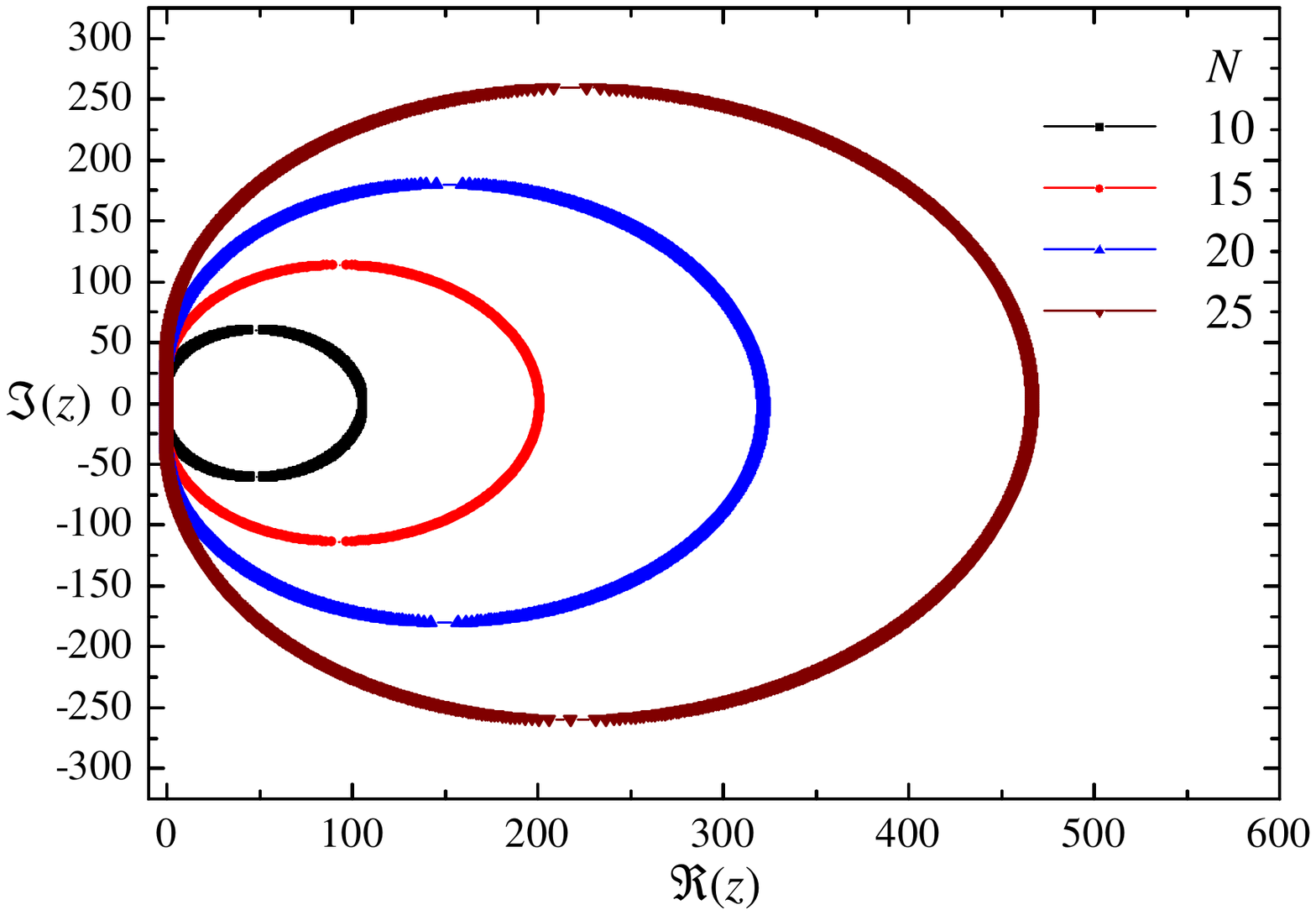}
\vspace{-8mm}\caption{\label{fig:stab_domain:c}\vspace{2mm}}
\end{subfigure}
\begin{subfigure}{0.48\textwidth}
\includegraphics[width=\textwidth]{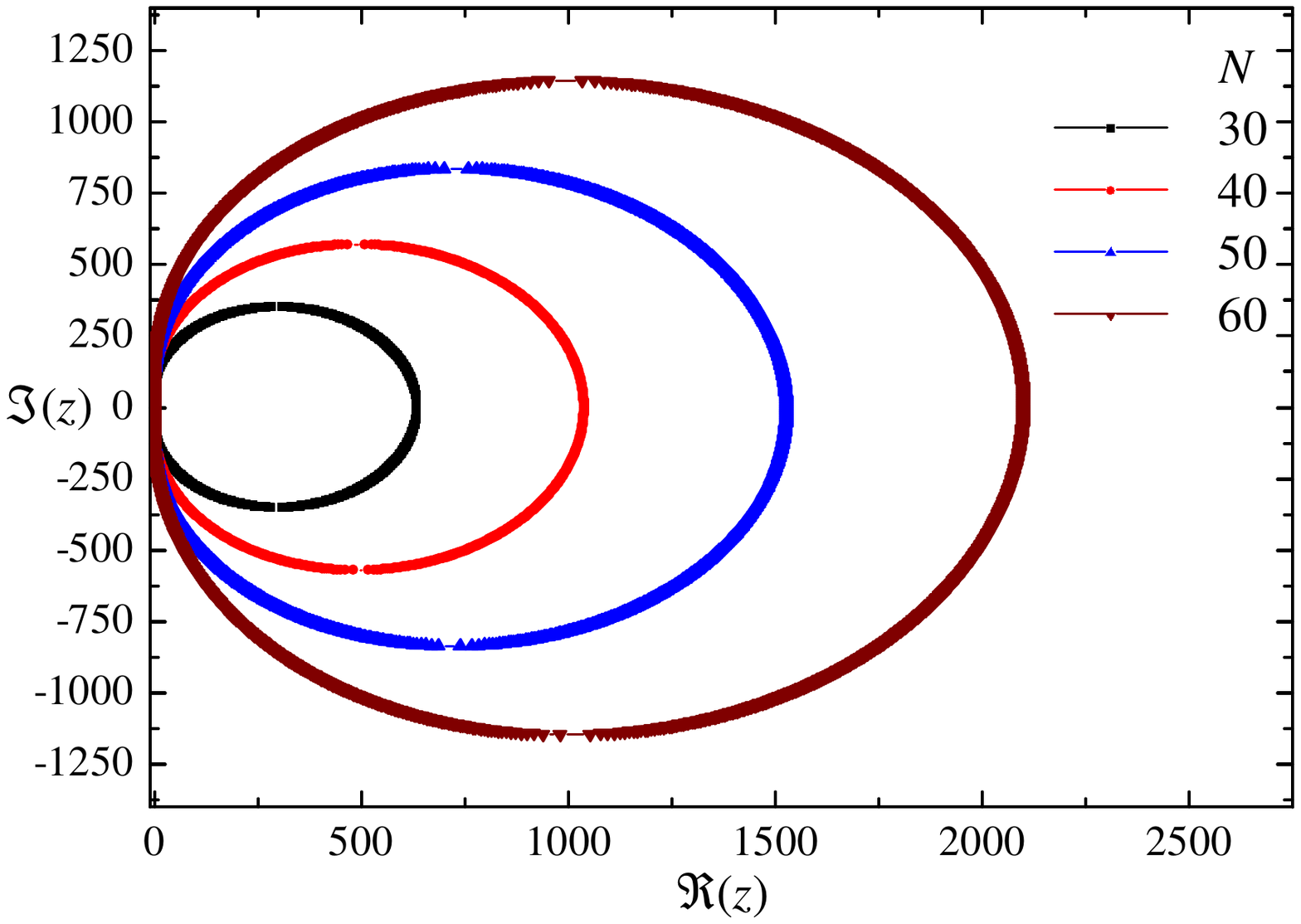}
\vspace{-8mm}\caption{\label{fig:stab_domain:d}\vspace{2mm}}
\end{subfigure}
\caption{%
Regions of absolute stability $|R(z)| < 1$ (outer part of closed curves) in the complex plane $z$ of the ADER-DG numerical method with a local DG predictor for the degrees $1 \leqslant N \leqslant 60$ of polynomials in the DG representation: $N = 1,\, 2,\, 3,\, 4,\, 5$ (a), $N = 6,\, 7,\, 8,\, 9$ (b), $N = 10,\, 15,\, 20,\, 25$ (c), $N = 30,\, 40,\, 50,\, 60$ (d). Horizontal and vertical axes represent real $\Re(z)$ and imaginary $\Im(z)$ parts of complex number $z = \lambda\cdot\mathrm{\Delta}t_{n}$.
}
\label{fig:stab_domain}
\end{figure}

The calculated regions of absolute stability $|R(z)| < 1$ are presented in Fig.~\ref{fig:stab_domain} for the degrees $1 \leqslant N \leqslant 60$ of polynomials in the DG representation. In the range $\Re(z) < 0$ the numerical method is absolutely stable. Therefore, it was concluded that the ADER-DG numerical method with a local DG predictor is $A$-stable.

\begin{figure}[h!]
\captionsetup[subfigure]{%
position=bottom,
textfont=normalfont,
singlelinecheck=off,
justification=raggedright,
font=tiny
}
\centering
\begin{subfigure}{0.48\textwidth}
\includegraphics[width=\textwidth]{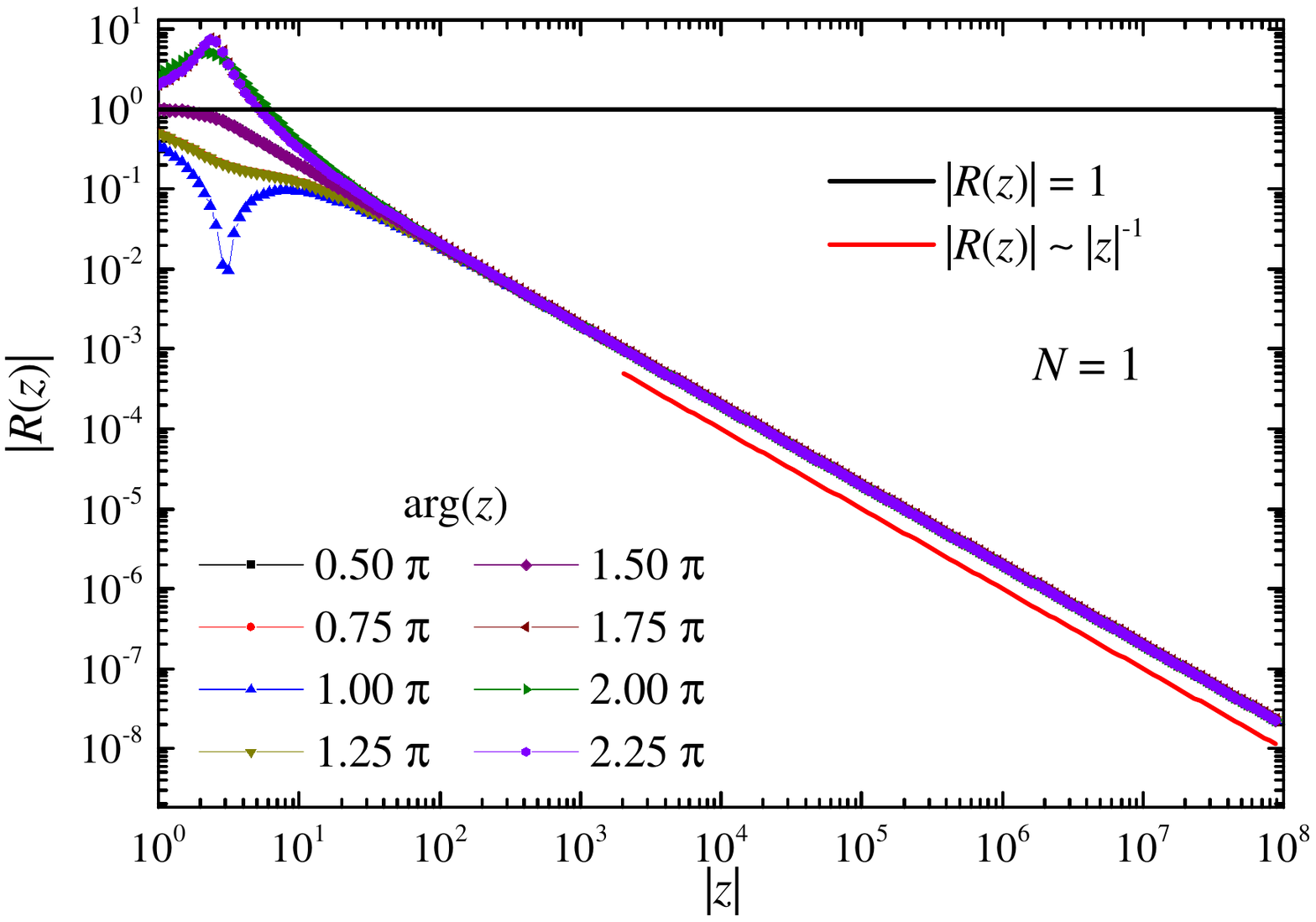}
\vspace{-8mm}\caption{\label{fig:stab_function_rays:a}\vspace{2mm}}
\end{subfigure}
\begin{subfigure}{0.48\textwidth}
\includegraphics[width=\textwidth]{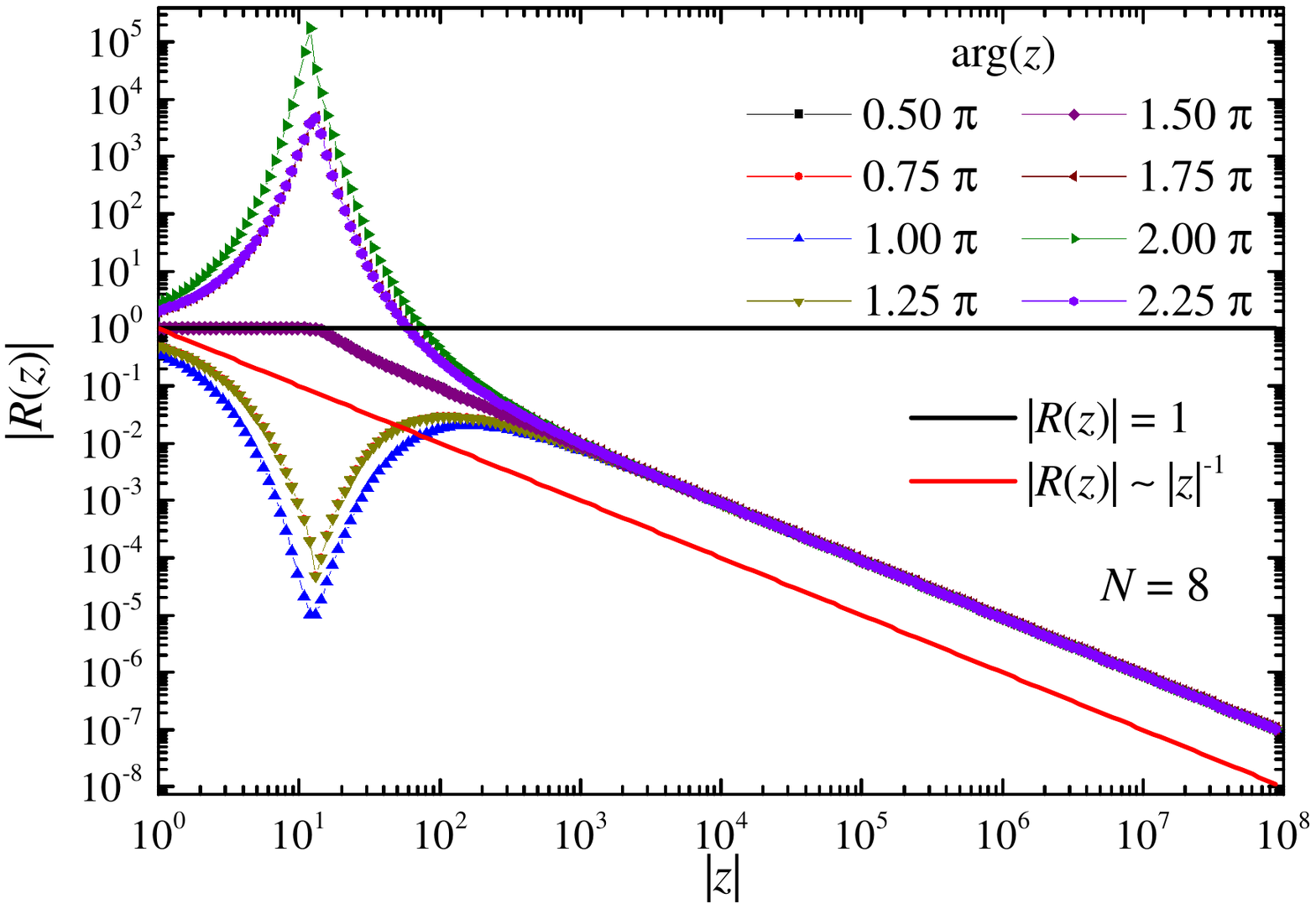}
\vspace{-8mm}\caption{\label{fig:stab_function_rays:c}\vspace{2mm}}
\end{subfigure}\\
\begin{subfigure}{0.48\textwidth}
\includegraphics[width=\textwidth]{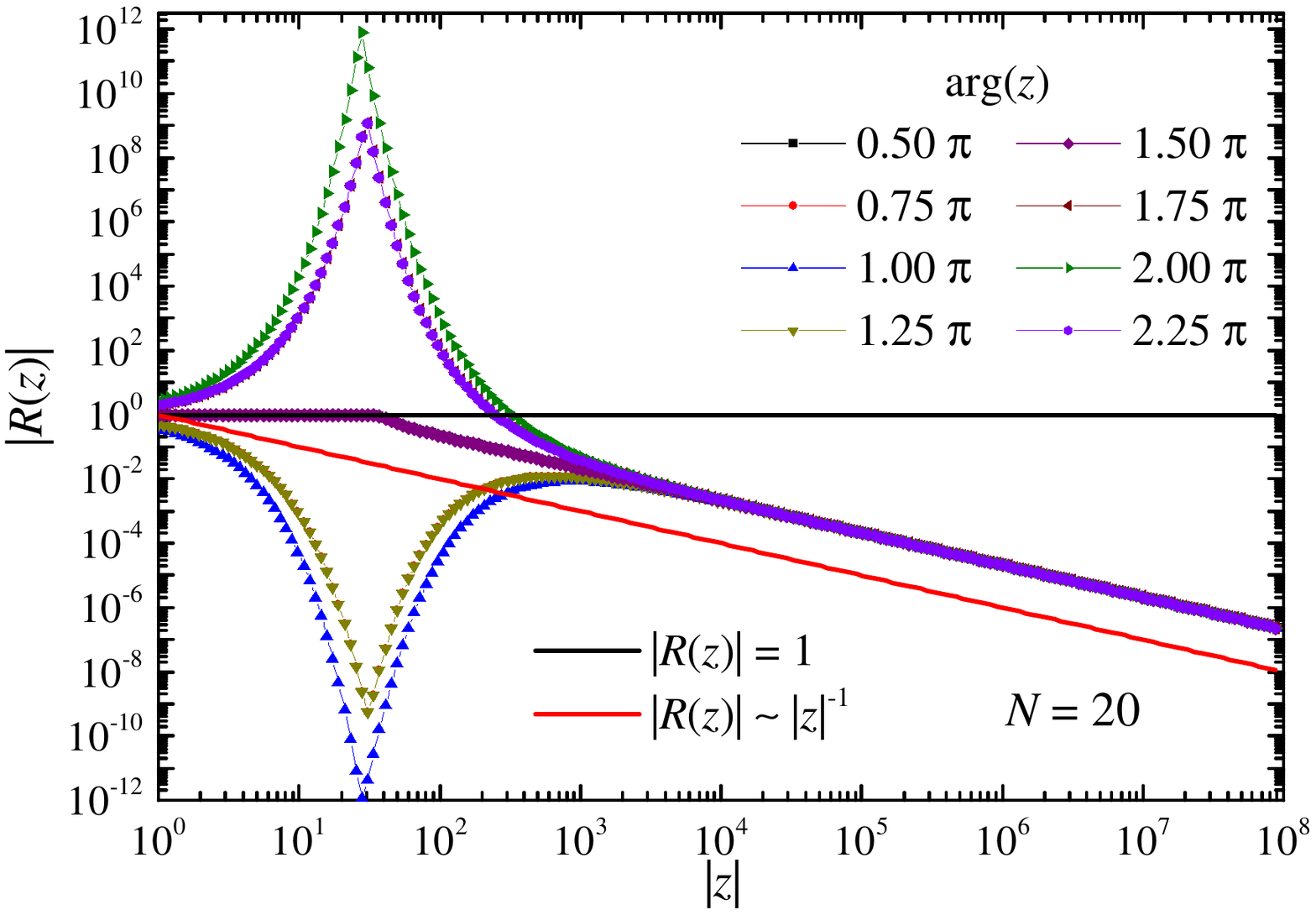}
\vspace{-8mm}\caption{\label{fig:stab_function_rays:d}\vspace{2mm}}
\end{subfigure}
\begin{subfigure}{0.48\textwidth}
\includegraphics[width=\textwidth]{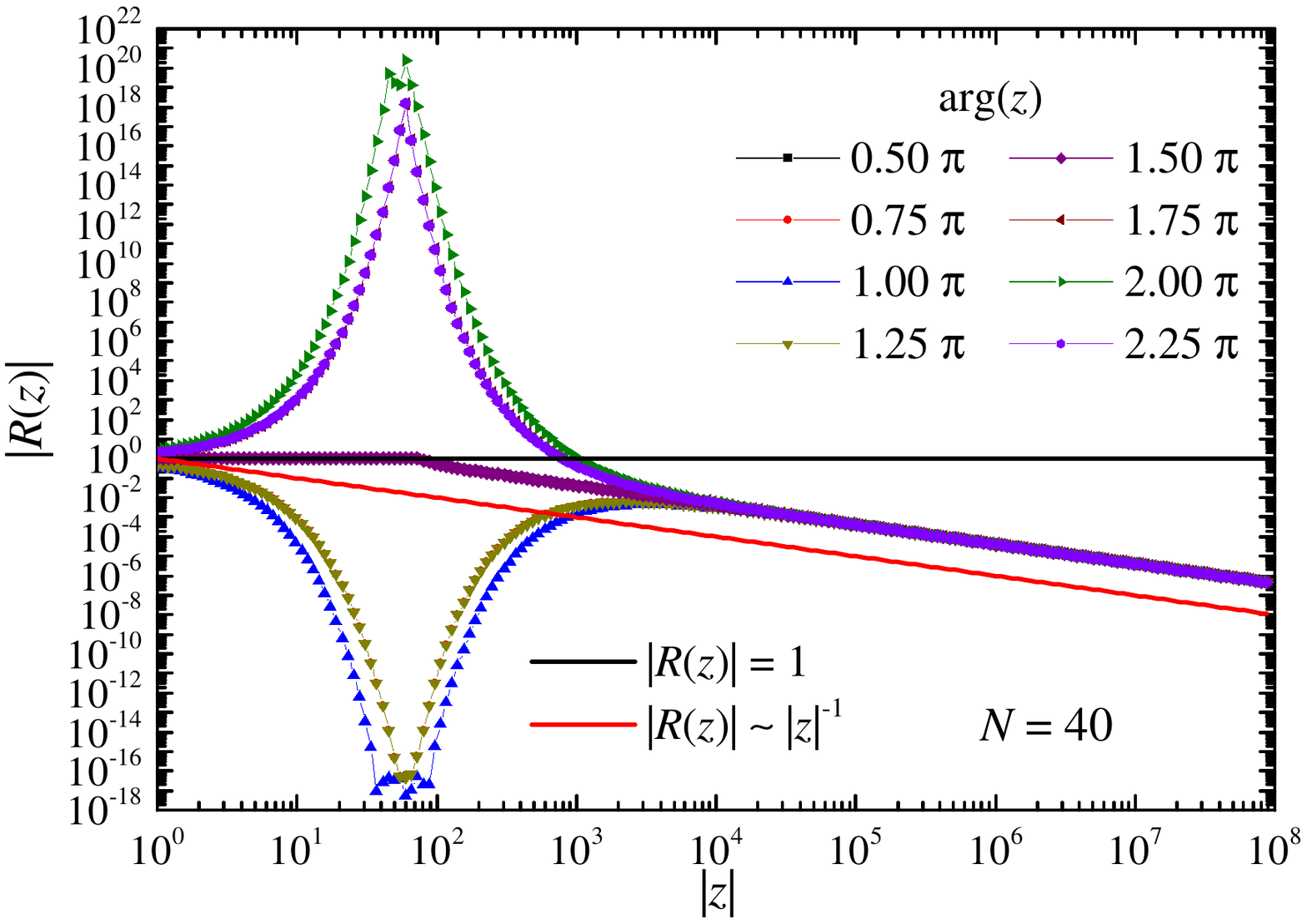}
\vspace{-8mm}\caption{\label{fig:stab_function_rays:f}\vspace{2mm}}
\end{subfigure}\\
\caption{%
The absolute values of the stability function $|R(z)|$ of the ADER-DG numerical method with a local DG predictor for the degrees $N$ on a set of radial rays $z = |z|\exp(i\,\arg(z))$ in the complex plane $z$, with $0.5\pi \leqslant \arg(z) < 2.5\pi$: $N = 1$ (a), $N = 8$ (b), $N = 20$ (c), $N = 40$ (d). The line $|R(z)| \sim |z|^{-1}$ demonstrates the asymptotic behavior appropriate in the $L$-stability. The line $|R(z)| = 1$ shows the region of absolute stability; the rays $\arg(z) = 1.75\pi,\, 2.0\pi,\, 2.25\pi$ have unstable regions with $|R(z)| > 1$ corresponding to the regions inside the closed curves in Fig.~\ref{fig:stab_domain}. 
}
\label{fig:stab_function_rays}
\end{figure} 

It should be noted that the determinants in the numerator and in the denominator of the final expression (\ref{eq:stab_func_expr}) for the stability function $R$ represent characteristic polynomials for matrices $[\mathrm{B} - \mathbf{1}\otimes\mathbf{w}]$ and $\mathrm{B}$, in which numbers $z$ that are reciprocals $z = 1/\nu$ of the classical definition of a characteristic polynomial $|\mathrm{A} - \nu\, \mathrm{E}|$ are presented as arguments (taking out the value $z^{N+1}$ in the numerator and denominator does not change the final expression).

Direct calculations showed that the rank of the matrix $[\mathrm{B} - \mathbf{1}\otimes\mathbf{w}]$ equals $N$, and the rank of the matrix $\mathrm{B}$ equals $N+1$ for the values of the degrees $1 \leqslant N \leqslant 80$, so the rank of the matrix $[\mathrm{B} - \mathbf{1}\otimes\mathbf{w}]$ is one less than the rank of the matrix $\mathrm{B}$. The expansion of the determinant in the numerator of the expression (\ref{eq:stab_func_expr}) always has a factor of the form $z = 1/\nu$. Therefore, after reducing the final fraction by a factor $z^{N+1}$, the numerator represents a polynomial of degree $N$ from the argument $z$, and the denominator represents a polynomial of degree $N+1$ from the argument $z$. Therefore, it is expected that the absolute value of the stability function $R(z)$ will have an asymptotic behavior of the form $|R(z)| \sim |z|^{-1}$ in the domain $|z|\rightarrow\infty$ for any directions in the complex plane. Therefore, it was concluded that the ADER-DG numerical method with a local DG predictor is $L$-stable; more precisely --- $L_{1}$-stable.

In Fig.~\ref{fig:stab_function_rays} additionally presented are the absolute values of the stability function $|R(z)|$ for several basic values of degrees $N$ of polynomials in the DG representation on a set of radial rays $z = |z|\exp(i\,\arg(z))$ in the complex plane $z$, with $0.5\pi \leqslant \arg(z) < 2.5\pi$; the reference point $0.5\pi$ is chosen for the boundary ray completely passing in the region of absolute stability (see Fig.~\ref{fig:stab_domain}). The line $|R(z)| \sim |z|^{-1}$ demonstrates the asymptotic behavior appropriate in the $L$-stability. The values on the rays with $\arg(z) = 0.5\pi,\, 0.75\pi, 1.75\pi$ are the same as the values on the rays with $\arg(z) = 1.5\pi,\, 1.25\pi,\, 2.25\pi$, respectively, which is due to complex conjugation. The rays $\arg(z) = 1.75\pi,\, 2.0\pi,\, 2.25\pi$ have unstable regions with $|R(z)| > 1$ corresponding to the regions inside the closed curves in Fig.~\ref{fig:stab_domain}.

\subsection{Numerical solution at nodes and local solution}

As a result of using the ADER-DG numerical method with a local DG predictor, two types of solutions are obtained: a solution at the grid nodes $\mathbf{u}_{n}$ and a local solution $\mathbf{q}$ given in the space between nodes. In this paper, we propose the use of both types of solutions, and show that the numerical method has a high order of accuracy both for the solution at the nodes and for the local solution. The local solution can be used both in the form of expansion coefficients in polynomials of the nodal basis, and directly in the form of a tabulation of the local solution in the grid sub-nodes specified in the spaces between grid nodes. The local solution at the grid sub-nodes can be obtained on the basis of the standard matrix-matrix multiplication, for which it is necessary to pre-compute the values of the basis polynomials at the grid sub-nodes at the beginning of the numerical calculation.

\begin{figure}[h!]
\captionsetup[subfigure]{%
position=bottom,
textfont=normalfont,
singlelinecheck=off,
justification=raggedright,
font=tiny
}
\centering
\begin{subfigure}{0.32\textwidth}
\includegraphics[width=\textwidth]{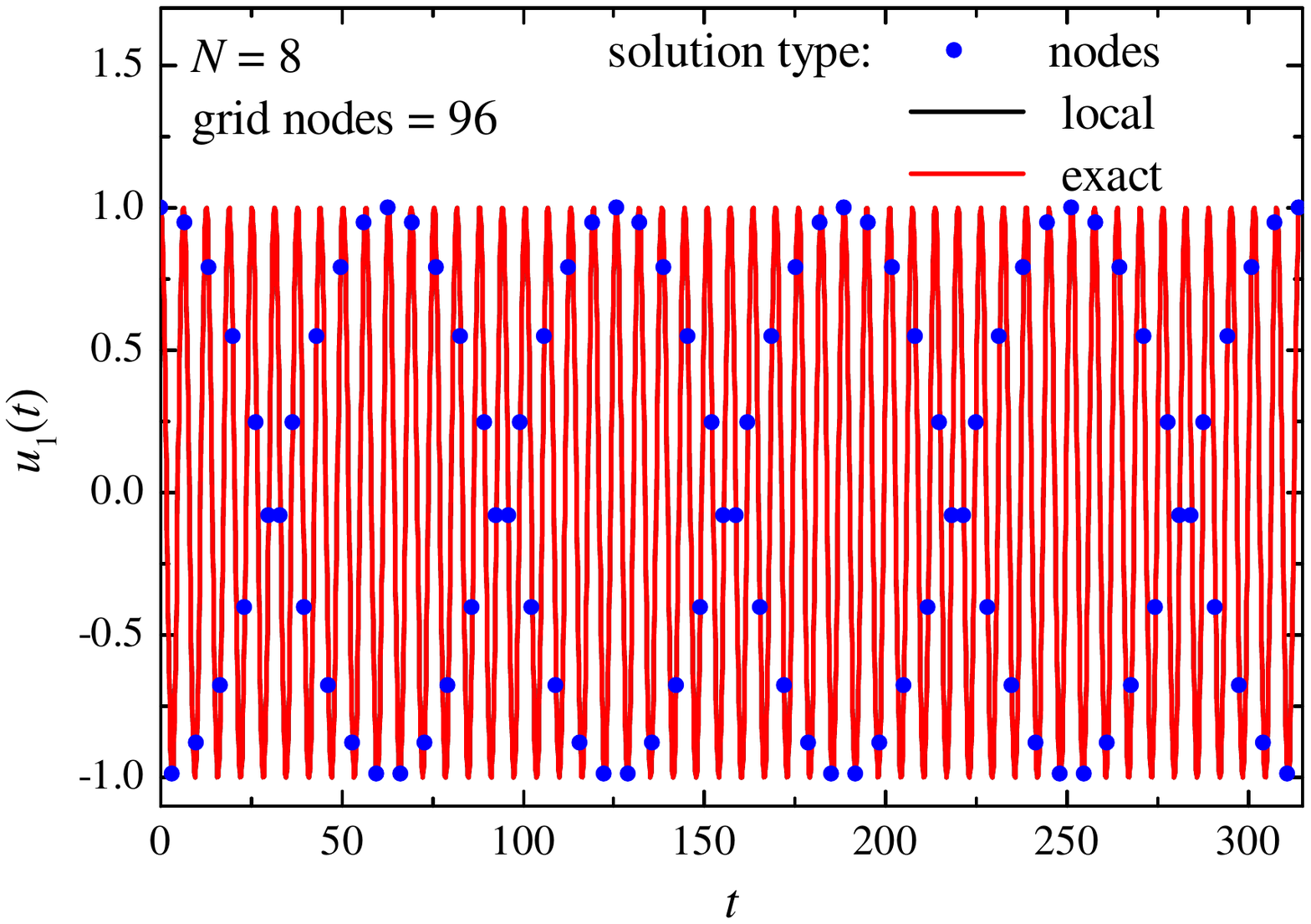}
\vspace{-8mm}\caption{\label{fig:harm_osc_sols_prec:a1}}
\end{subfigure}
\begin{subfigure}{0.32\textwidth}
\includegraphics[width=\textwidth]{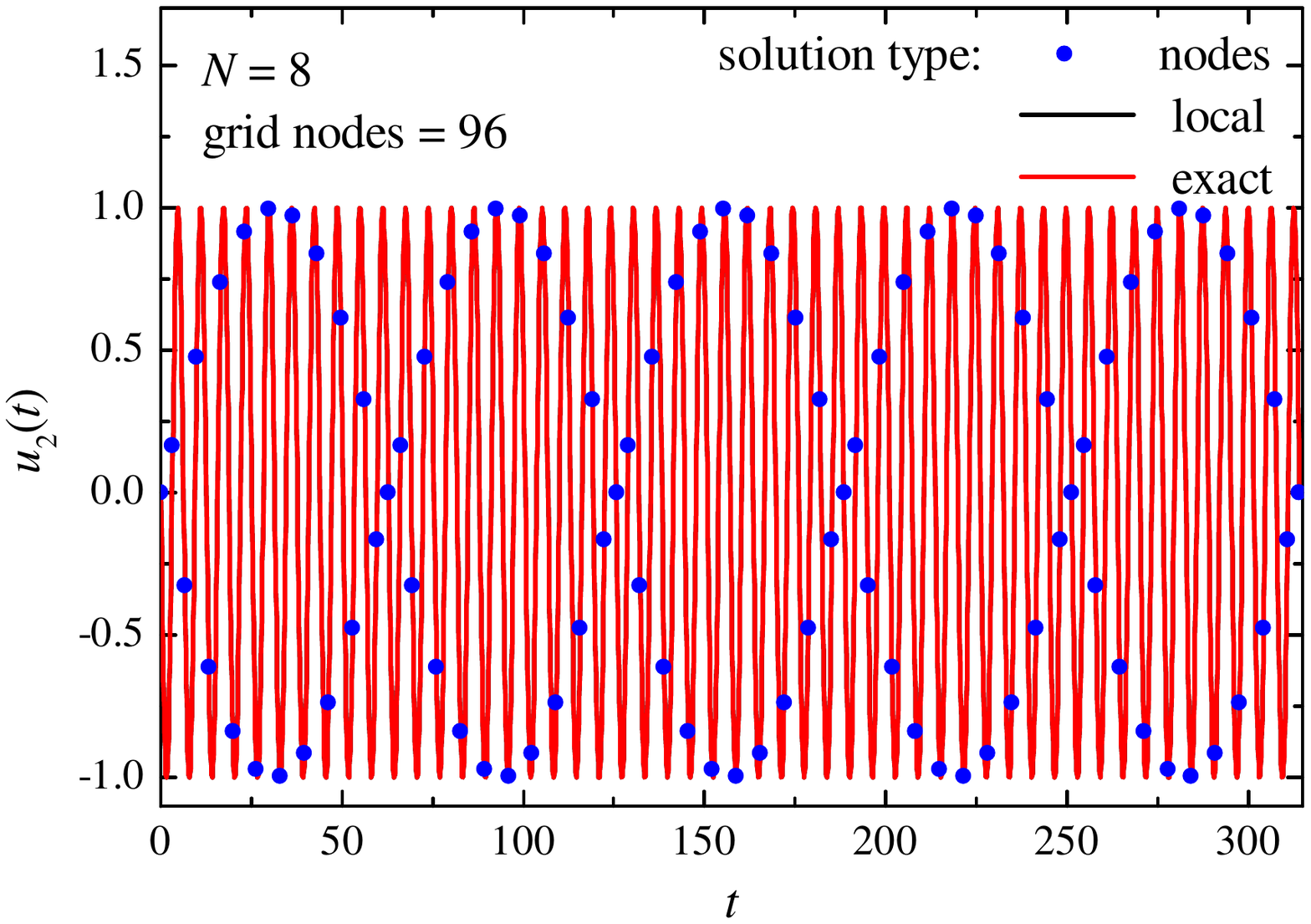}
\vspace{-8mm}\caption{\label{fig:harm_osc_sols_prec:a2}}
\end{subfigure}
\begin{subfigure}{0.32\textwidth}
\includegraphics[width=\textwidth]{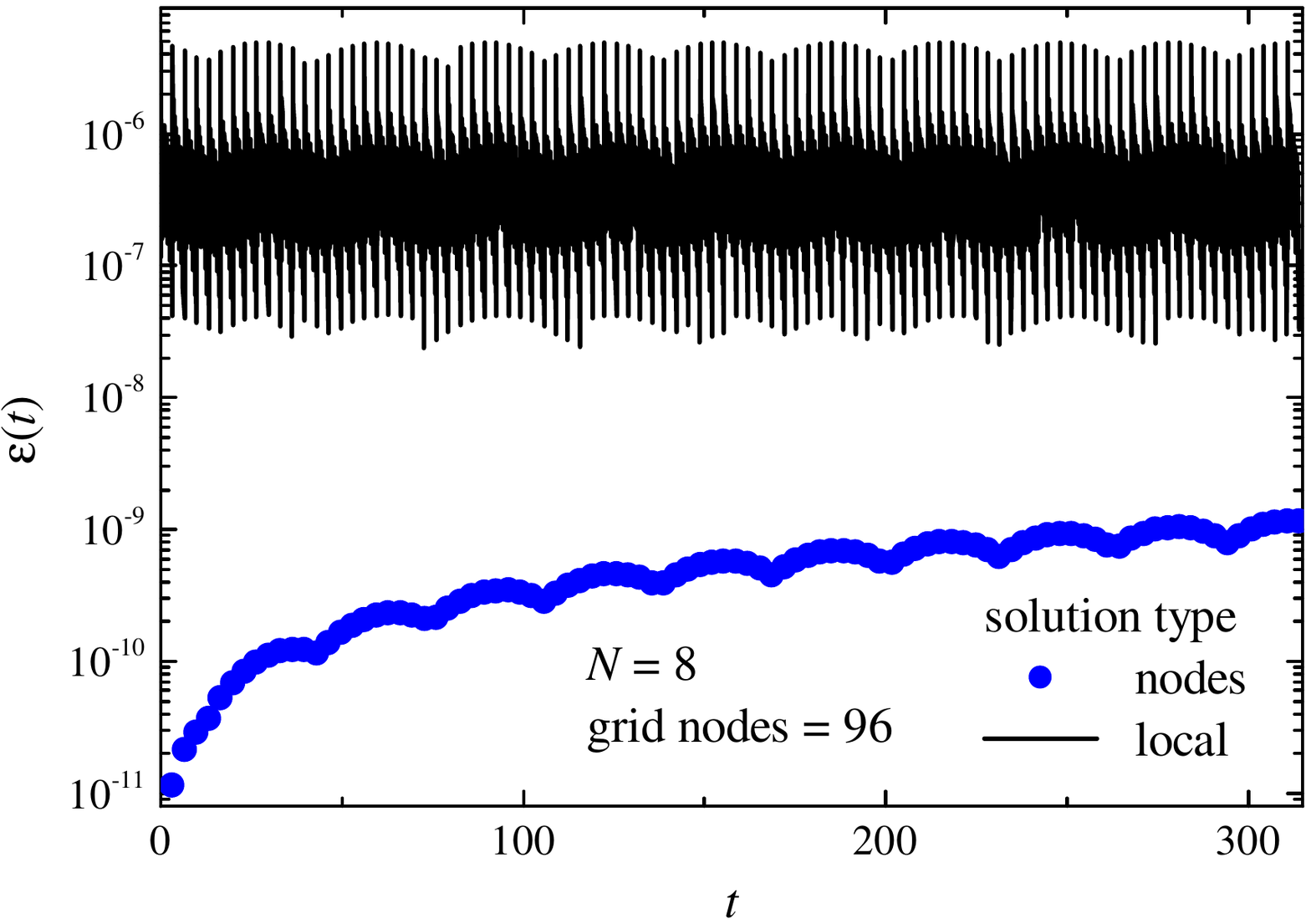}
\vspace{-8mm}\caption{\label{fig:harm_osc_sols_prec:a3}}
\end{subfigure}\\
\begin{subfigure}{0.32\textwidth}
\includegraphics[width=\textwidth]{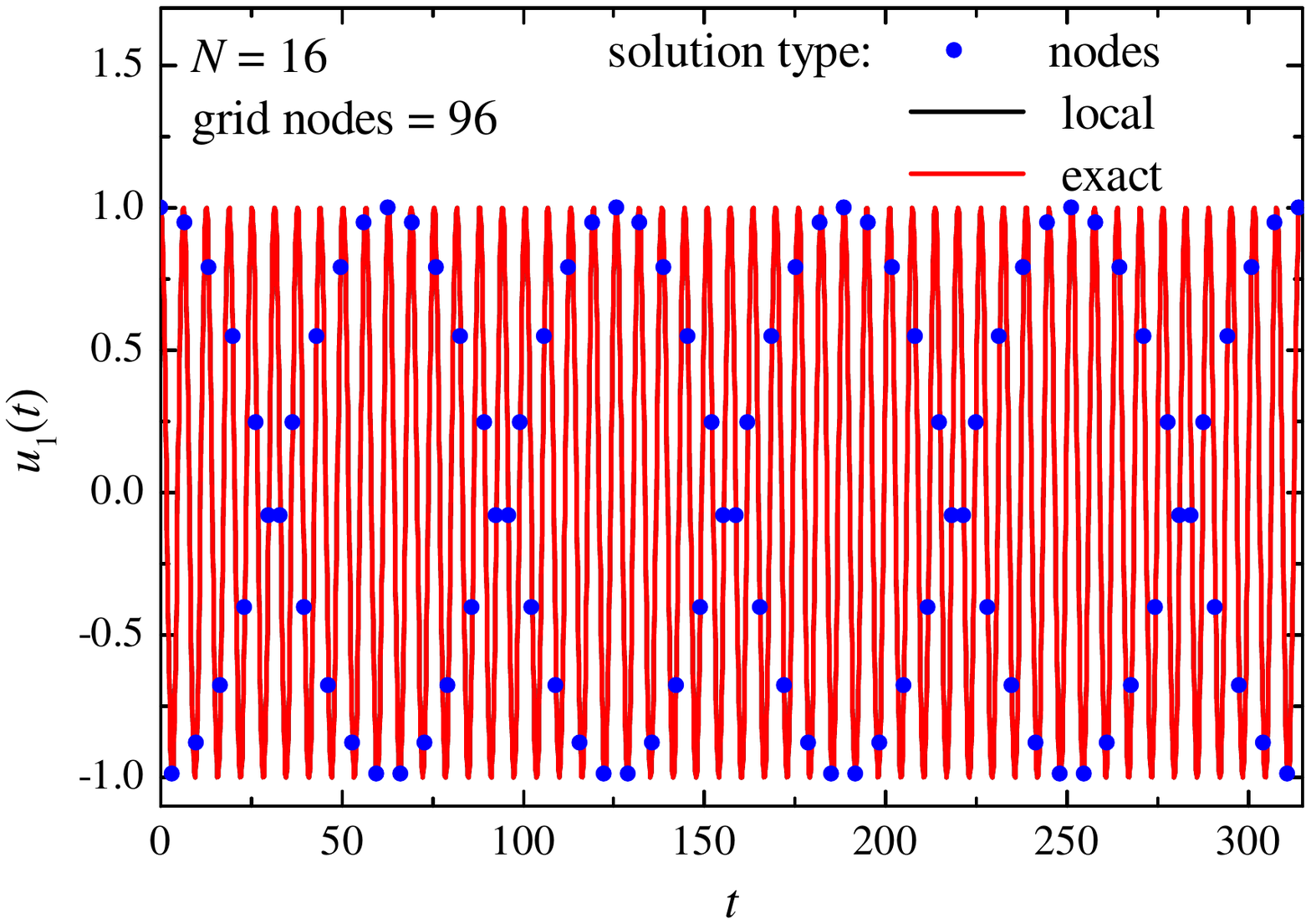}
\vspace{-8mm}\caption{\label{fig:harm_osc_sols_prec:b1}}
\end{subfigure}
\begin{subfigure}{0.32\textwidth}
\includegraphics[width=\textwidth]{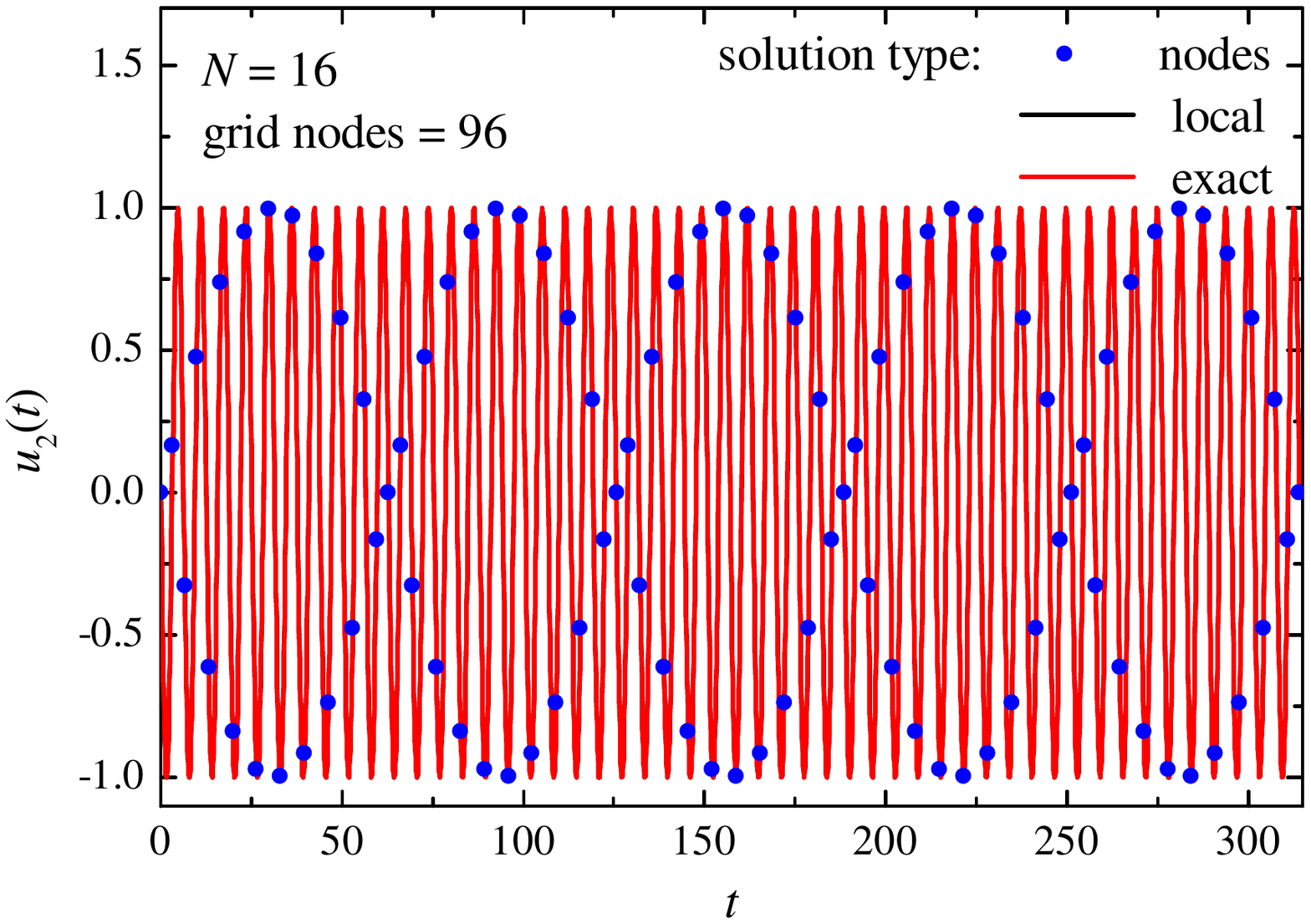}
\vspace{-8mm}\caption{\label{fig:harm_osc_sols_prec:b2}}
\end{subfigure}
\begin{subfigure}{0.32\textwidth}
\includegraphics[width=\textwidth]{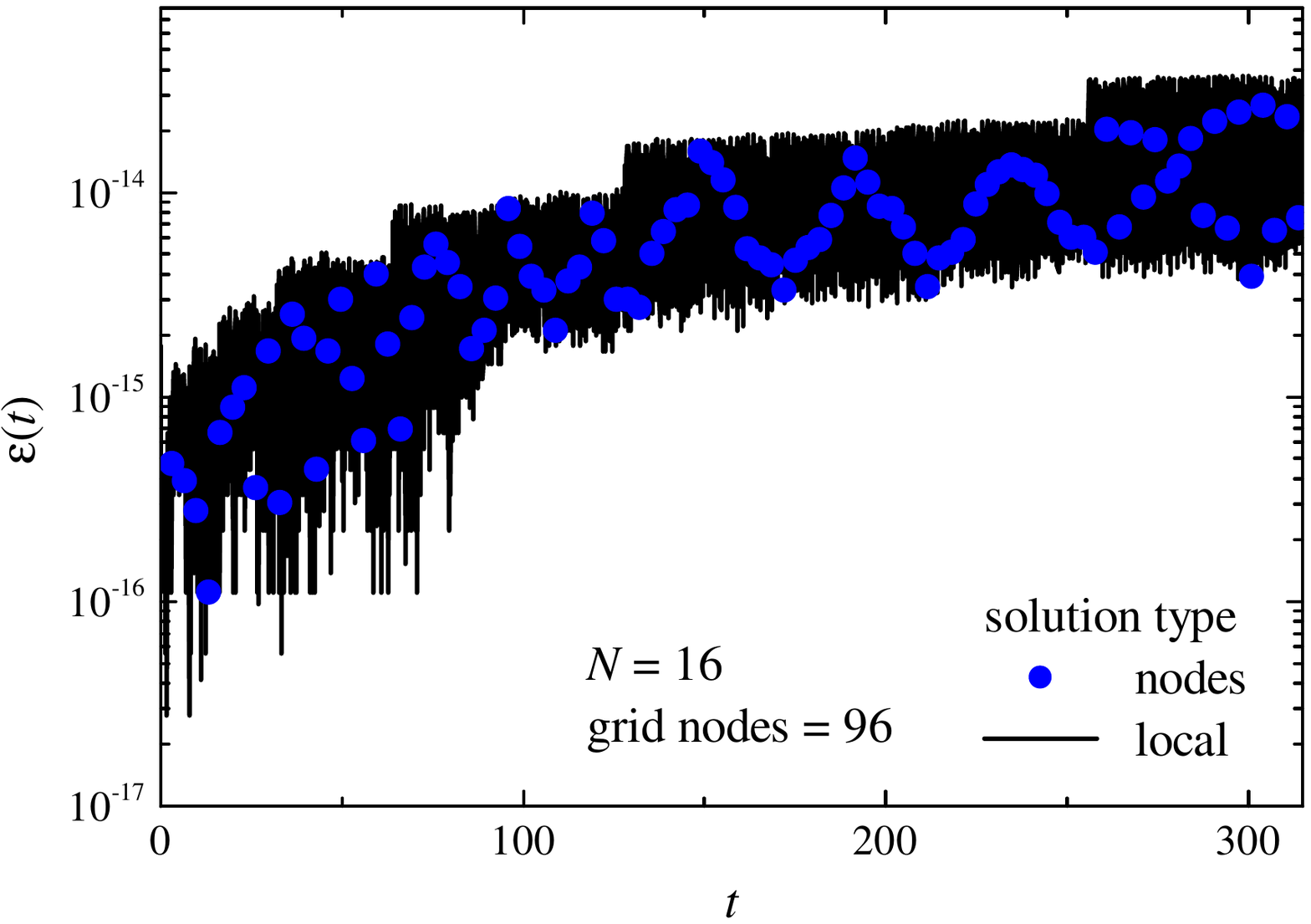}
\vspace{-8mm}\caption{\label{fig:harm_osc_sols_prec:b3}}
\end{subfigure}\\
\begin{subfigure}{0.32\textwidth}
\includegraphics[width=\textwidth]{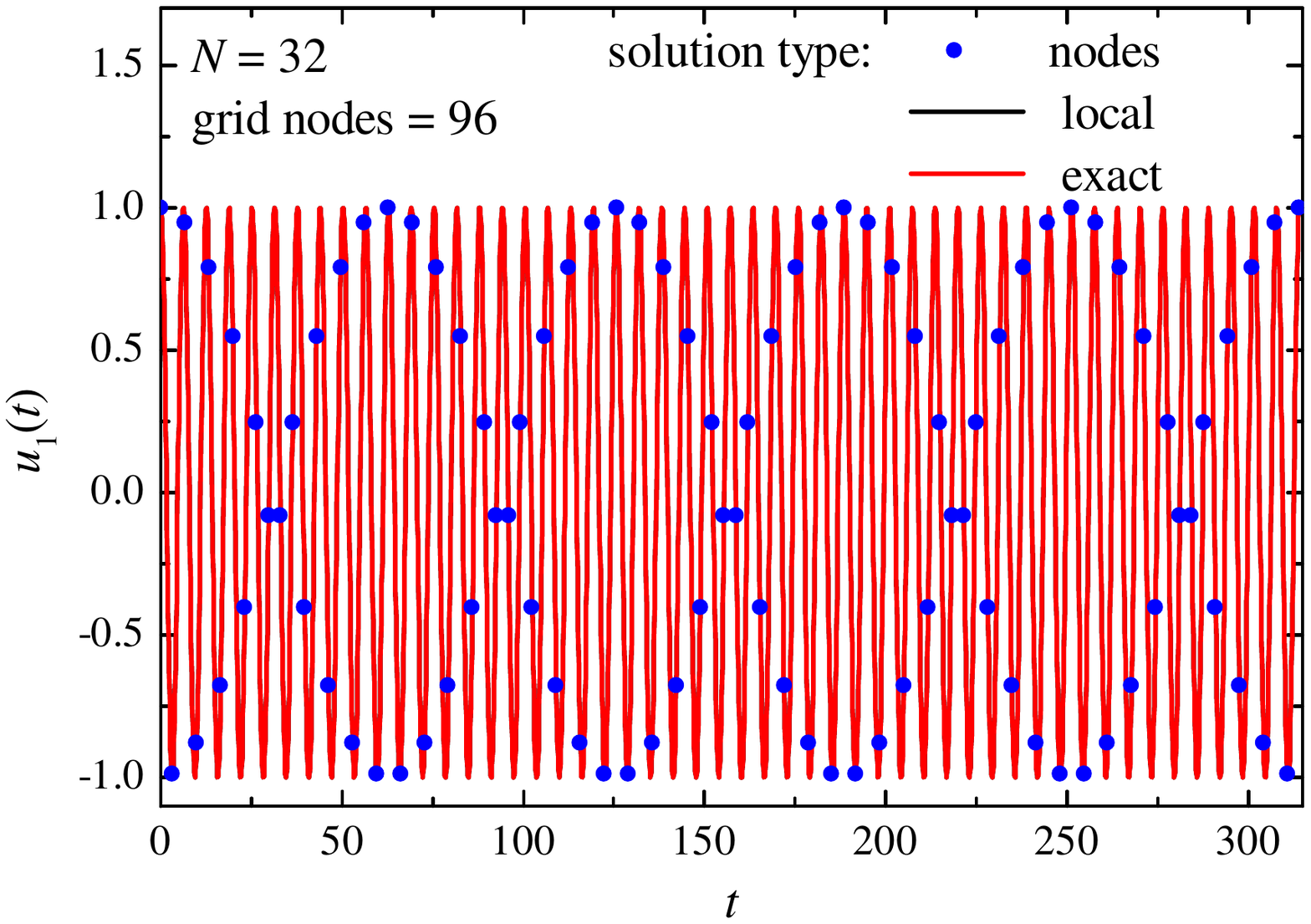}
\vspace{-8mm}\caption{\label{fig:harm_osc_sols_prec:c1}}
\end{subfigure}
\begin{subfigure}{0.32\textwidth}
\includegraphics[width=\textwidth]{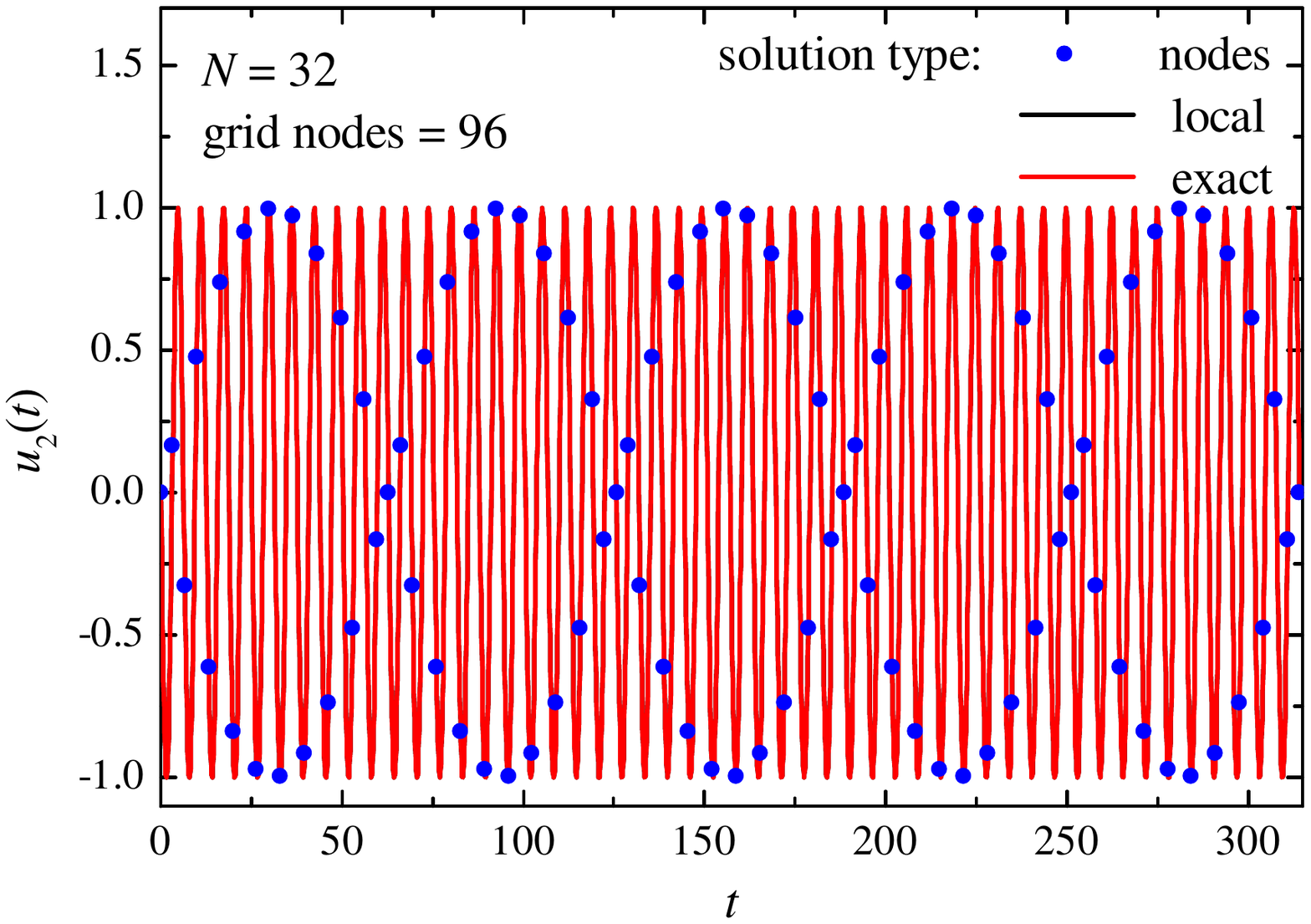}
\vspace{-8mm}\caption{\label{fig:harm_osc_sols_prec:c2}}
\end{subfigure}
\begin{subfigure}{0.32\textwidth}
\includegraphics[width=\textwidth]{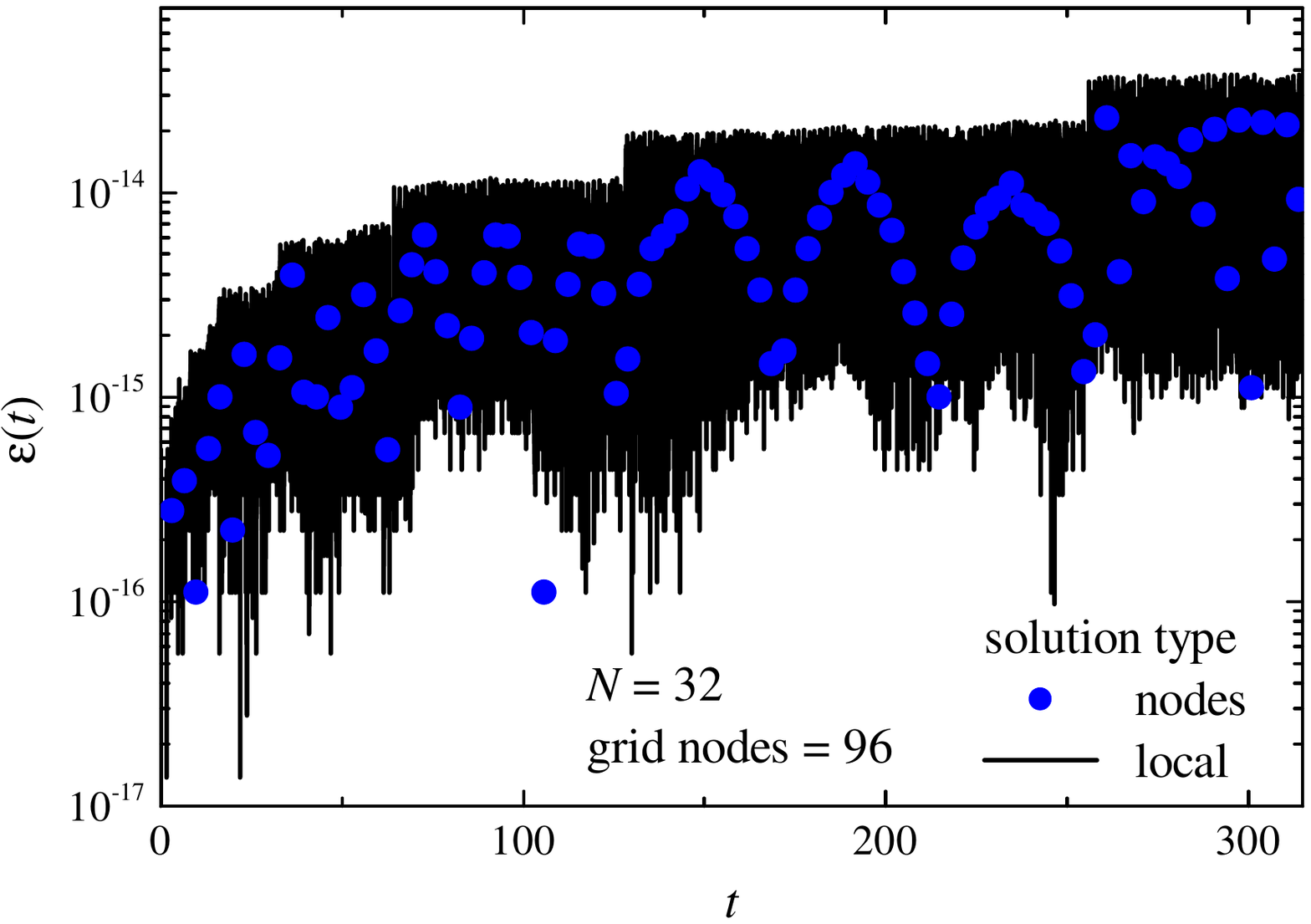}
\vspace{-8mm}\caption{\label{fig:harm_osc_sols_prec:c3}}
\end{subfigure}\\
\caption{%
Numerical solution of the problem (\ref{eq:harm_osc_prec}), in the domain $0 \leqslant t \leqslant 100\pi$ with step $\mathrm{\Delta}t = \frac{100}{96}\pi$,
using the ADER-DG numerical method with a local DG predictor with degrees of polynomials $N = 8$ (\subref{fig:harm_osc_sols_prec:a1}, \subref{fig:harm_osc_sols_prec:a2}, \subref{fig:harm_osc_sols_prec:a3}), $16$ (\subref{fig:harm_osc_sols_prec:b1}, \subref{fig:harm_osc_sols_prec:b2}, \subref{fig:harm_osc_sols_prec:b3}) and $32$ (\subref{fig:harm_osc_sols_prec:c1}, \subref{fig:harm_osc_sols_prec:c2}, \subref{fig:harm_osc_sols_prec:c3}). Presented data: numerical solution at nodes, local solution and exact solution separately for components $u_{1}$ (\subref{fig:harm_osc_sols_prec:a1}, \subref{fig:harm_osc_sols_prec:b1}, \subref{fig:harm_osc_sols_prec:c1}) and $u_{2}$ (\subref{fig:harm_osc_sols_prec:a2}, \subref{fig:harm_osc_sols_prec:b2}, \subref{fig:harm_osc_sols_prec:c2}); dependence of the point-wise error $\varepsilon(t)$ of the numerical solution for the solution at the nodes and the local solution (\subref{fig:harm_osc_sols_prec:a3}, \subref{fig:harm_osc_sols_prec:b3}, \subref{fig:harm_osc_sols_prec:c3}).
}
\label{fig:harm_osc_sols_prec}
\end{figure} 

The notation $\mathbf{u}_{L}(t)$ for a local solution in the domain $\Omega$ of the desired function has been introduced:
\begin{equation}\label{eq:def_local_sol_t}
\mathbf{u}_{L}(t) = \mathbf{q}_{n}\left(\frac{t - t_{n}}{\Delta t_{n}}\right) = 
\sum\limits_{p = 0}^{N} \mathbf{q}_{n, p} \varphi_{p}\left(\frac{t - t_{n}}{\Delta t_{n}}\right),\quad 
t_{n} \leqslant t < t_{n+1},
\end{equation}
where $\mathbf{q}_{n} = \mathbf{q}_{n}(\tau)$ denotes the local solution obtained at the discretization domain $\Omega_{n}$, $\mathbf{q}_{n, p}$ denotes the expansion coefficients of the local solution $\mathbf{q}_{n} = \mathbf{q}_{n}(\tau)$ in the DG representation. Discretization $\xi_{m}$, $0 \leqslant m < M$, in units of the local coordinate $\tau$, was introduced by grid sub-nodes located in the space between grid nodes. The values of the basis polynomials $\varphi_{p}$ were calculated in sub-nodes $\xi_{m}$ of discretization, which made it possible to calculate the tabulated local solution $\mathbf{u}_{L}(t_{n, m})$ in the form of multiplication the polynomials values matrix, with sizes $M \times (N+1)$, on the matrix of coefficients $\mathbf{q}_{p}$, with sizes $(N+1) \times K$.

An important feature of the ADER-DG numerical method with a local DG predictor is the possibility of using the local solution $\mathbf{u}_{L}(t)$ as the final solution of the problem even in the case of very coarse grids with a large discretization step $\Delta t$. As a demonstration example a numerical solution of the IVP for a linear one-dimensional harmonic oscillator $\ddot{x} + x = 0$, $x(0) = 1, \dot{x}(0) = 0$:
\begin{equation}\label{eq:harm_osc_prec}
\dfrac{du_{1}}{dt} = u_{2};\quad \dfrac{du_{2}}{dt} = -u_{1};\quad u_{1}(0) = 1;\quad u_{2}(0) = 0;
\end{equation}
with an analytical exact solution $\mathbf{u}^{\rm ex} = [\cos(t),\, -\sin(t)]^{T}$ is presented in Fig.~\ref{fig:harm_osc_sols_prec}. The solution was obtained in the domain $0 \leqslant t \leqslant 100\pi$ with a discretization step $\Delta t = \frac{100}{96}\pi$ equal to more than a half of harmonic oscillation period. The calculations for this demonstration example were carried out using standard double precision floating point arithmetic. To obtain a tabular local solution in all spaces between grid nodes uniform in the local coordinate $\tau$, grids with $M = 10^{3}$ sub-nodes were additionally included. 

The presented results show that even in the case of polynomial degree $N = 8$, the local solution $\mathbf{u}_{L}$ has a very high accuracy, and the solution at the nodes $\mathbf{u}_{n}$ reaches the round-off error of double-precision floating-point numbers. In the case of using polynomials of degree $N = 16$, the errors $\varepsilon(t)$ of both the node solution $\mathbf{u}_{n}$ and the local solution $\mathbf{u}_{L}$ are compared in order, resulting in round-off double-precision floating-point errors. In the case of using polynomials with a degree of $N = 32$, the approximate order of accuracy no longer improves, however, there are more tabulation points $t_{n, m}$ of the local solution $\mathbf{u}_{L}$ with a very small error $\sim 10^{-15}\div10^{-16}$ than in the case $N = 8$.

This example is a demonstration example, and it shows that the real accuracy of the numerical solution obtained using the predictor can already reach the highest applied calculation accuracy.

\subsection{Accuracy, convergence and norms}

A quantitative study of the accuracy and convergence of the ADER-DG numerical method with a local DG predictor in this work was carried out on the basis of solving a set of test problems for the IVP for ODE systems for which there is an exact analytical solution: three problems for second-order equations and three problems for third-order equations, which created systems of two and three first-order equations, respectively. The results of the study are presented in the following Section~\ref{sec:2}.

This paragraph presents the mathematical apparatus of the norms used in this work and methods for their calculation, on the basis of which the errors of the numerical solution and the convergence orders of the numerical method were calculated. Accuracy and convergence were separately investigated for the solution at the nodes $\mathbf{u}_{n}$ and the local solution $\mathbf{q}_{n} = \mathbf{q}_{n}(\tau)$ in the space between nodes $[t_{n},\, t_{n+1}]$. The solution in nodes represents the values $\mathbf{u}_{n}$ of the desired function $\mathbf{u}$ at individual grid points $t_{n}$. The local solution $\mathbf{u}_{L}(t) = \mathbf{q}_{n}(\tau(t))$ is a continuous function (see formula (\ref{eq:def_local_sol_t})) of local coordinate $\tau$ determined by the expansion in basis polynomials $\varphi_{p}(\tau)$, the coefficients $\hat{\mathbf{q}}_{n, p}$ of which are determined for each separate space $[t_{n},\, t_{n+1}]$ between adjacent grid nodes. In this paragraph, the local solution $\mathbf{q}_{n}$ will also be matched with the index $n$ of the grid node, which is the left grid node $t_{n}$ for the discretization domain $\Omega_{n} = \{t\, |\, t \in [t_{n},\ t_{n+1}]\}$ in which the local solution is defined.

The accuracy of the numerical solution was determined point-wise, based on the values of the local error of the numerical solution:
\begin{equation}\label{eq:local_eps_def}
\varepsilon(t) = \max\limits_{1 \leqslant k \leqslant K}|u_{k}(t) - u^{\rm ex}_{k}(t)|,
\end{equation}
where $u_{k}(t)$ and $u^{\rm ex}_{k}(t)$ are the $k$-th components of the vector function of the numerical solution $\mathbf{u}(t)$ and the exact solution $\mathbf{u}^{\rm ex}(t)$ at a point $t$ respectively. Among the results presented below are graphs of the dependence of the error on $t$. The error $\varepsilon(t_{n})$ of the numerical solution in nodes $\mathbf{u}_{n}$ will be presented only in the nodes $t_{n}$ of the grid. The error $\varepsilon(t_{n, m})$ of the local solution $\mathbf{u}_{L}(t)$ will be represented in sub-nodes $t_{n, m}$ located between the nodes $t_{n}$ of the grid. In this paper, $M = 10^{3}$ points were uniformly added in the spaces $\Omega_{n}$ between the nodes, at which estimates of the local numerical solution were made.

The calculation of the convergence orders $p$ of the numerical solution was carried out on the basis of the analysis of the global error $e$ of the numerical solution, integrally defined for the entire domain $\Omega$ of the definition of the desired function $\mathbf{u}$. The global error $e$ was calculated in three classical functional norms:
\begin{equation}\label{eq:global_eps_def}
\begin{split}
&e_{L_{1}} = \int\limits_{\Omega} \left| \mathbf{u}(t) - \mathbf{u}^{\rm ex}(t) \right| dt,\\
&e_{L_{2}}^{2} = \int\limits_{\Omega} \left| \mathbf{u}(t) - \mathbf{u}^{\rm ex}(t) \right|^{2} dt,\\
&e_{L_{\infty}} = \sup\limits_{t\in\Omega} \left| \mathbf{u}(t) - \mathbf{u}^{\rm ex}(t) \right| dt,
\end{split}
\end{equation}
where under the symbols $|\ldots|$ we denote the maximum absolute value taken over all components of the vector function: $|\mathbf{u}| = \max_{k} |u_{k}|$, as in (\ref{eq:local_eps_def}). The calculation of global errors (\ref{eq:global_eps_def}) for the solution at the nodes $\mathbf{u}_{n}$ was carried out on the basis of the replacement of integrals by finite sums multiplied by sampling steps; the operation $\sup$ by $t\in\Omega$ has been replaced by the operation $\max$ by nodes $t_{n}$. The calculation of global errors (\ref{eq:global_eps_def}) for the local solution $\mathbf{u}_{L}(t)$ in the space between nodes $\Omega_{n}$ was carried out on the basis of the replacement of integrals by finite sums over sub-nodes $t_{n, m}$, multiplied by the discretization steps between sub-nodes $\Delta t_{n}/M$; the operation $\sup$ by $t\in\Omega$ has been replaced by the operation $\max$ by sub-nodes $t_{n, m}$. In the following Figs. the global error $e$ of the solution at the nodes $\mathbf{u}_{n}$ is denoted $e_{G}$, the global error of the local solution $\mathbf{u}_{L}(t)$ in the space between nodes is denoted $e_{L}$.

The convergence orders $p$ were calculated based on the approximation of the dependence of the global error $e$ of the numerical solution on the discretization step $\Delta t$ by the power dependence $e(\Delta t) \sim (\Delta t)^{p}$. The calculation of the convergence order was carried out on the basis of the analysis of the numerical solution, performed with a constant discretization step $\Delta t$. The convergence orders $p$ in this work were calculated separately for the solution at the nodes and the local solution in the space between the nodes, separately for the functional norms: $p_{L_{1}}$, $p_{L_{2}}$, $p_{L_{\infty}}$.

The ADER-DG numerical method with a local DG predictor have high accuracy and sufficiently large convergence orders. Therefore, already on coarse grids, the error of the numerical solution reaches the round-off error of double-precision floating-point numbers. Therefore, calculations where the values of local errors reached values $10^{-14}\div10^{-15}$ were carried out using arbitrary precision floating-point arithmetic. The implementation of arbitrary precision floating-point arithmetic \texttt{mpmath} was used (math package for the \texttt{python} programming language), with the parameter \texttt{mpmath.mp.dps = 500}.

According to the results of the works~\cite{dg_ivp_ode_1, dg_ivp_ode_3}, it is expected that the ADER-DG numerical method with a local DG predictor with degrees of polynomials $N$ considered in this work can demonstrate superconvergence with the convergence order $p_{G}$ for the solution at grid nodes $\mathbf{u}_{n}$ and the convergence order $p_{L}$ for the local solution $\mathbf{u}_{L}(t)$:
\begin{equation}\label{eq:expect_orders}
p_{\rm (G)} = 2N+1;\quad
p_{\rm (L)} = N+1;
\end{equation}
which is typical for DG methods for solving the IVP for ODE systems.

\section{Applications of the numerical method}
\label{sec:2}

This Section presents applications of using the ADER-DG numerical method. The main applications are examples of initial value problems for linear and nonlinear ODE systems and the test for the stiff ODE, designed to demonstrate the capabilities of this numerical method for solving stiff problems. The last Subsection presents estimates of the computational costs of the ADER-DG numerical method and their comparison with the computational costs of classical Runge-Kutta numerical methods.

\subsection{Classical tests}

This Subsection considers six classical examples of solving the IVP for a first-order ODE system. The numerical solution was obtained using the ADER-DG numerical method with a local DG predictor with degrees of polynomials $N = 1, \ldots, 60$. The local errors $\varepsilon(t)$ of the numerical solution were calculated based on the analysis of the dependence of the global error $e$ on the discretization step $\Delta t$ for various functional norms (\ref{eq:global_eps_def}), both for the solution at the nodes $\mathbf{u}_{n}$ and for the local solution $\mathbf{u}_{L}$. To obtain a tabular local solution in all spaces between grid nodes uniformly in the local coordinate $\tau$, grids with $M = 10^{3}$ sub-nodes were additionally included. 

Empirical convergence orders $p$ were calculated based on the approximation of the dependence of the global error $e$ on the discretization step $\Delta t$ by a power law, and technically, on the basis of the approximation of the dependence of the logarithm of the global error $\log(e)$ on the logarithm of the step $\log(\Delta t)$ by a linear function $\log(e) \sim p\log(\Delta t)$. Empirical convergence orders of the numerical solution were calculated and compared with the theoretical values of the convergence order $p_{G}$ for the solution at grid nodes and the convergence order $p_{L}$ for the local solution defined by expression (\ref{eq:expect_orders}).

The results were obtained in terms of the accuracy of the numerical solution when compared with the exact analytical solution. In order to avoid repetition in each individual Example, it should be immediately noted, that the numerical solution (the results are presented in Figures~\ref{fig:harm_osc_sols}, \ref{fig:exp_diss_sols},~\ref{fig:bratu_problem_sols}, \ref{fig:third_order_ivp_ex_1_sols}, \ref{fig:third_order_ivp_ex_2_sols}, \ref{fig:third_order_ivp_ex_3_sols} and will also be discussed in detail below) obtained using the ADER-DG numerical method with a local DG predictor accurately reproduces all the main features of the problem solution even on coarse grids. In the case of degree $N = 1$, the discontinuity of the local solution $\mathbf{u}_{L}$ at the nodes is clearly observed. A similar behavior is also observed in cases of other values of the degrees $N > 1$, but this is not clearly visible in these figures. This aspect is manifested in a point-wise comparison of the results for the errors $\varepsilon(t)$ of the local and global solutions. The error $\varepsilon(t)$ scales for the local solution and the solution at the nodes in the case of degree $N = 1$ practically do not differ, except for the case of the domain $\Omega_{0}$, where $t_{0} = 0$ is the left point. An analysis of the error $\varepsilon(t)$ for $N > 1$ shows that the error scales for the local solution and the solution at nodes differ by several orders of magnitude. It can be noted that for the error of the numerical solution at the nodes $\mathbf{u}_{n}$, there is an approximately linear increase in the error $\varepsilon(t)$ with an increase in the argument $t$.

\subsubsection{Example 1}

\begin{figure}[h!]
\captionsetup[subfigure]{%
position=bottom,
textfont=normalfont,
singlelinecheck=off,
justification=raggedright,
font=tiny
}
\centering
\begin{subfigure}{0.275\textwidth}
\includegraphics[width=\textwidth]{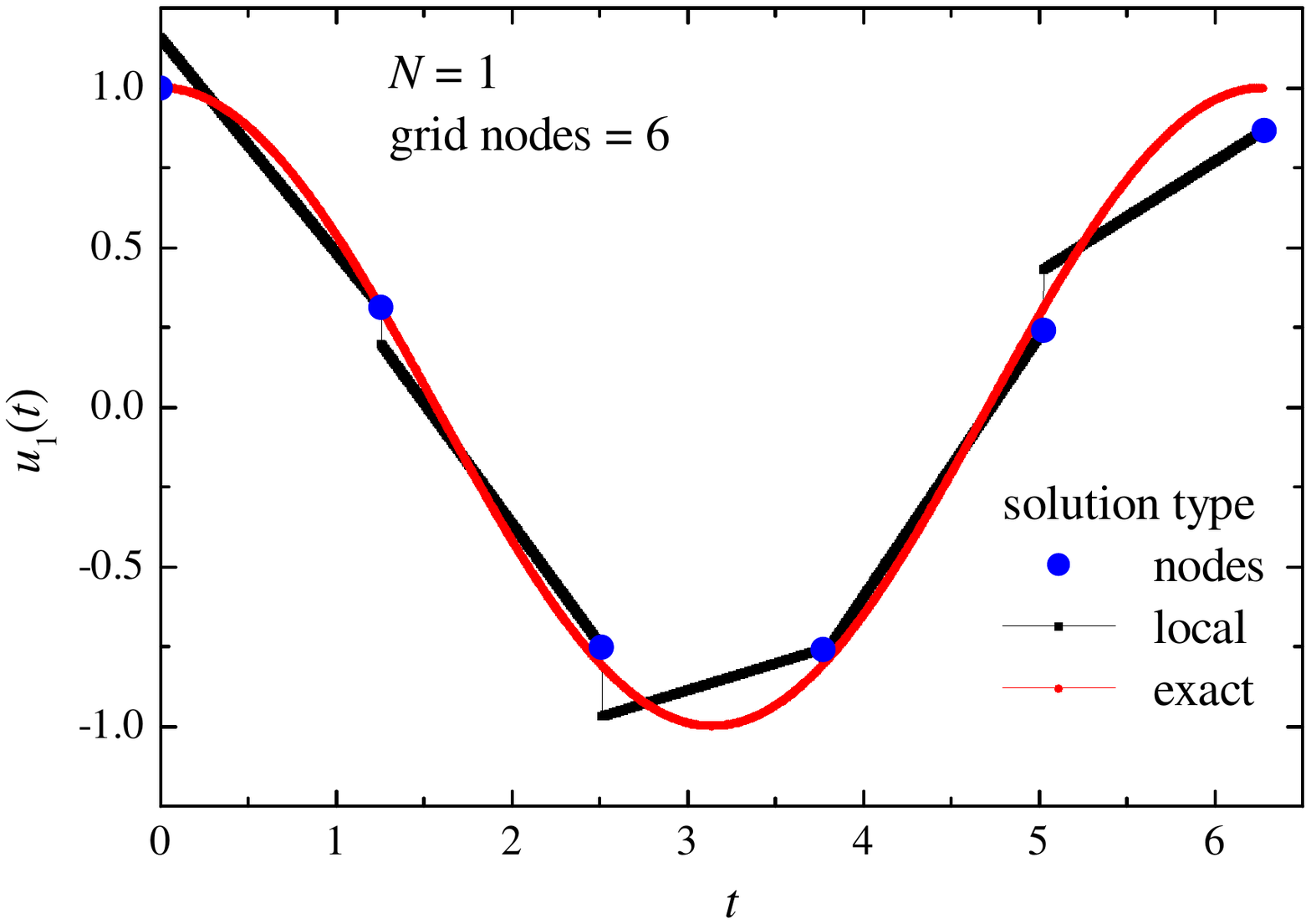}
\vspace{-8mm}\caption{\label{fig:harm_osc_sols:a1}}
\end{subfigure}
\begin{subfigure}{0.275\textwidth}
\includegraphics[width=\textwidth]{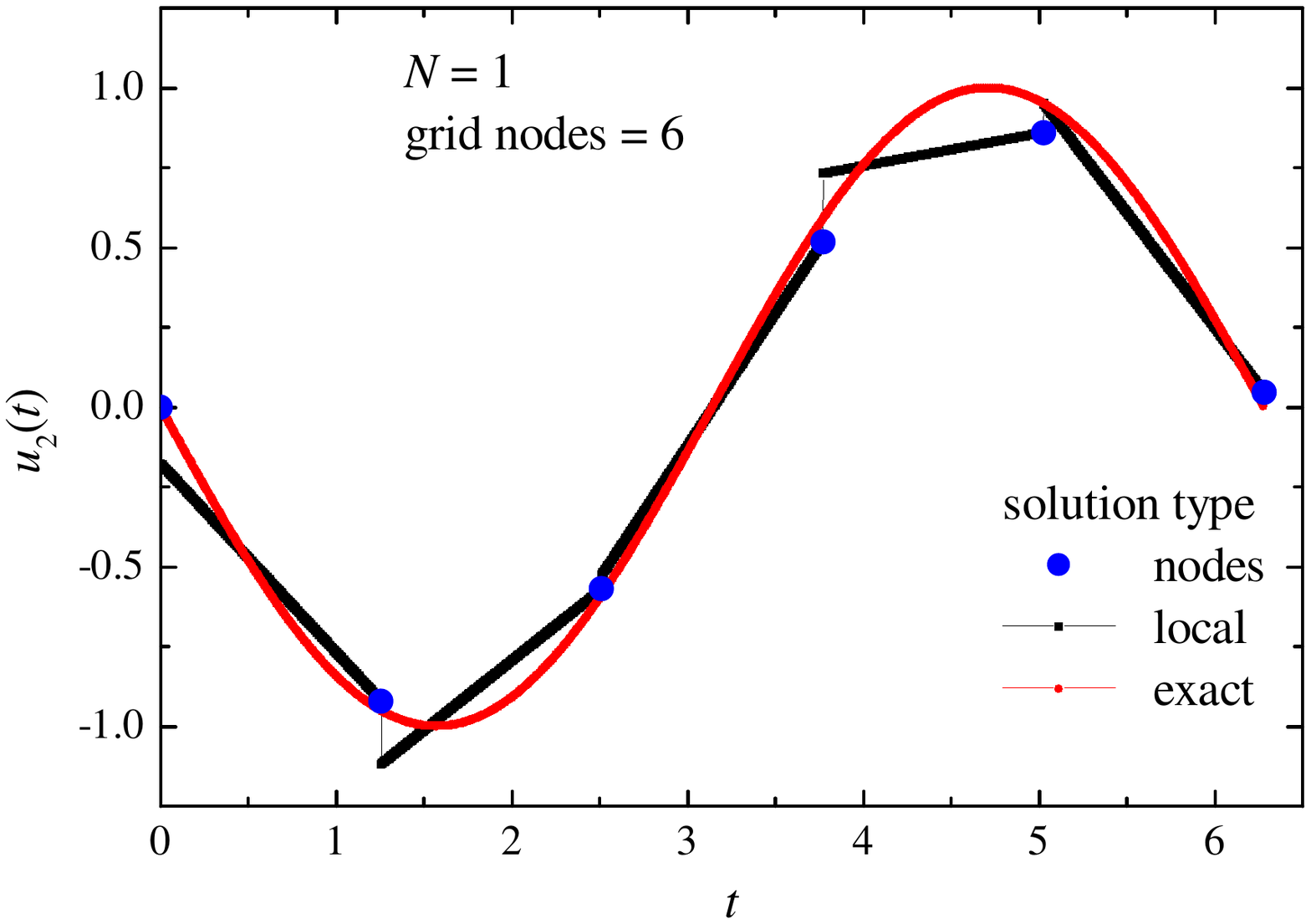}
\vspace{-8mm}\caption{\label{fig:harm_osc_sols:a2}}
\end{subfigure}
\begin{subfigure}{0.275\textwidth}
\includegraphics[width=\textwidth]{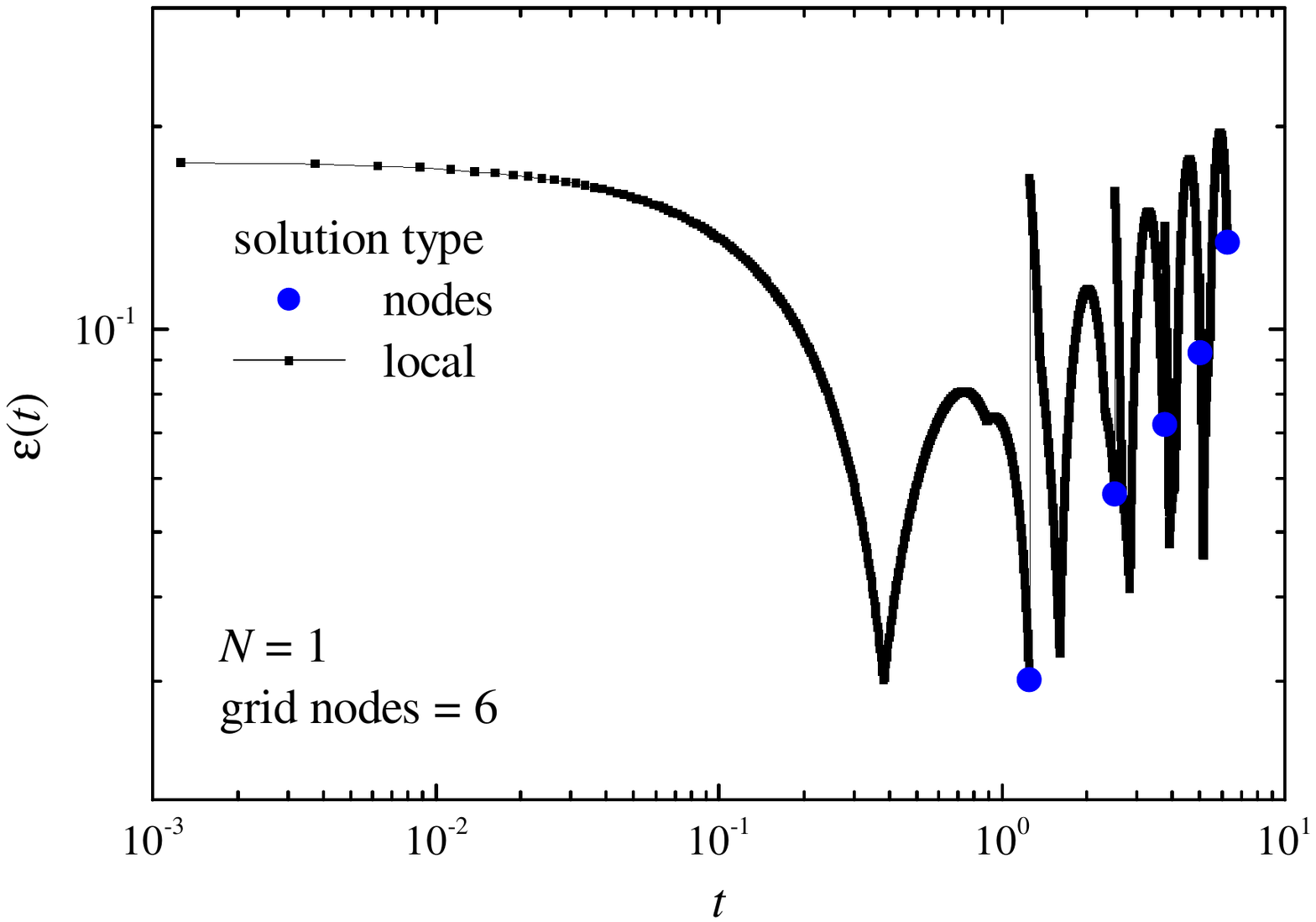}
\vspace{-8mm}\caption{\label{fig:harm_osc_sols:a3}}
\end{subfigure}\\
\begin{subfigure}{0.275\textwidth}
\includegraphics[width=\textwidth]{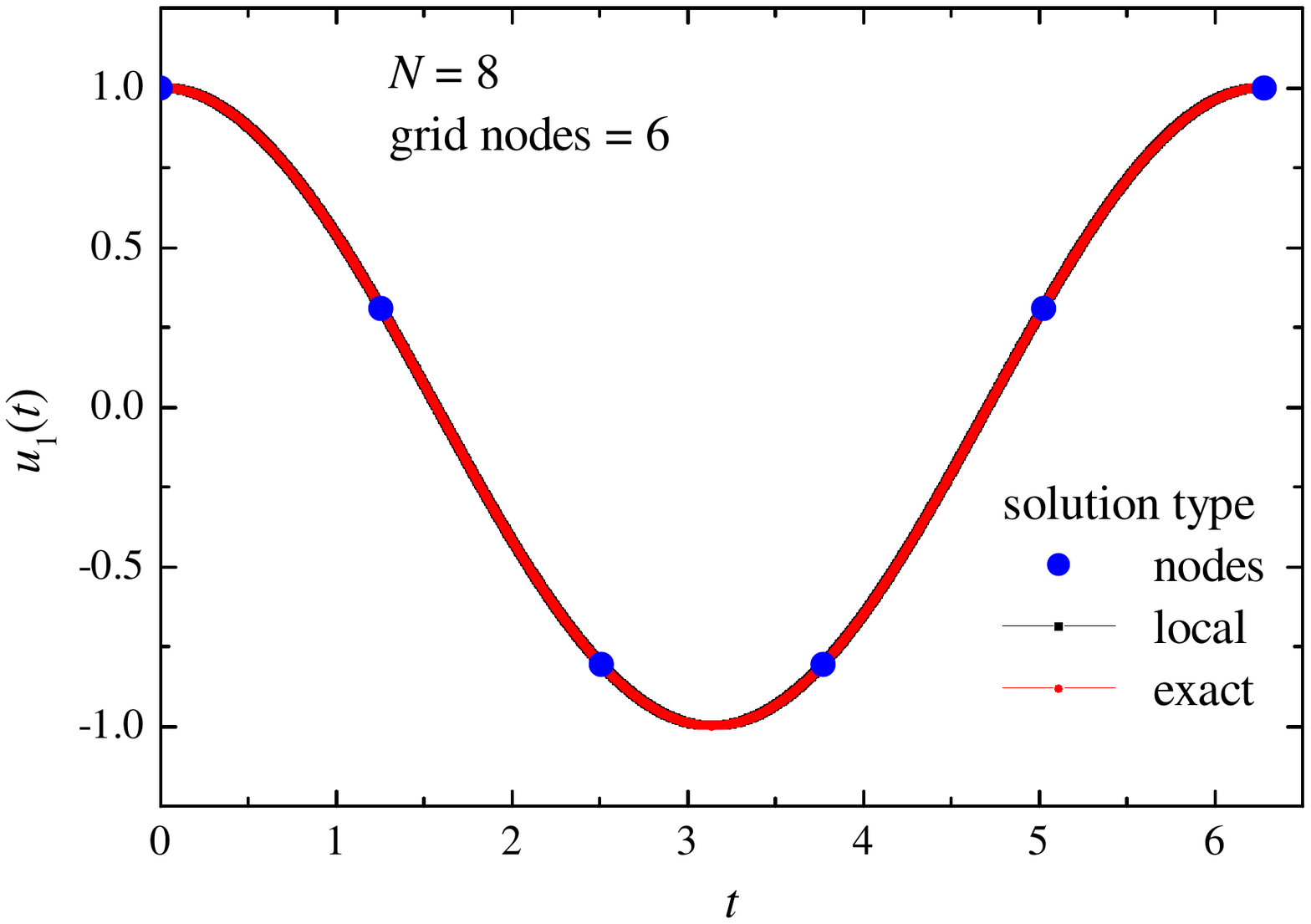}
\vspace{-8mm}\caption{\label{fig:harm_osc_sols:b1}}
\end{subfigure}
\begin{subfigure}{0.275\textwidth}
\includegraphics[width=\textwidth]{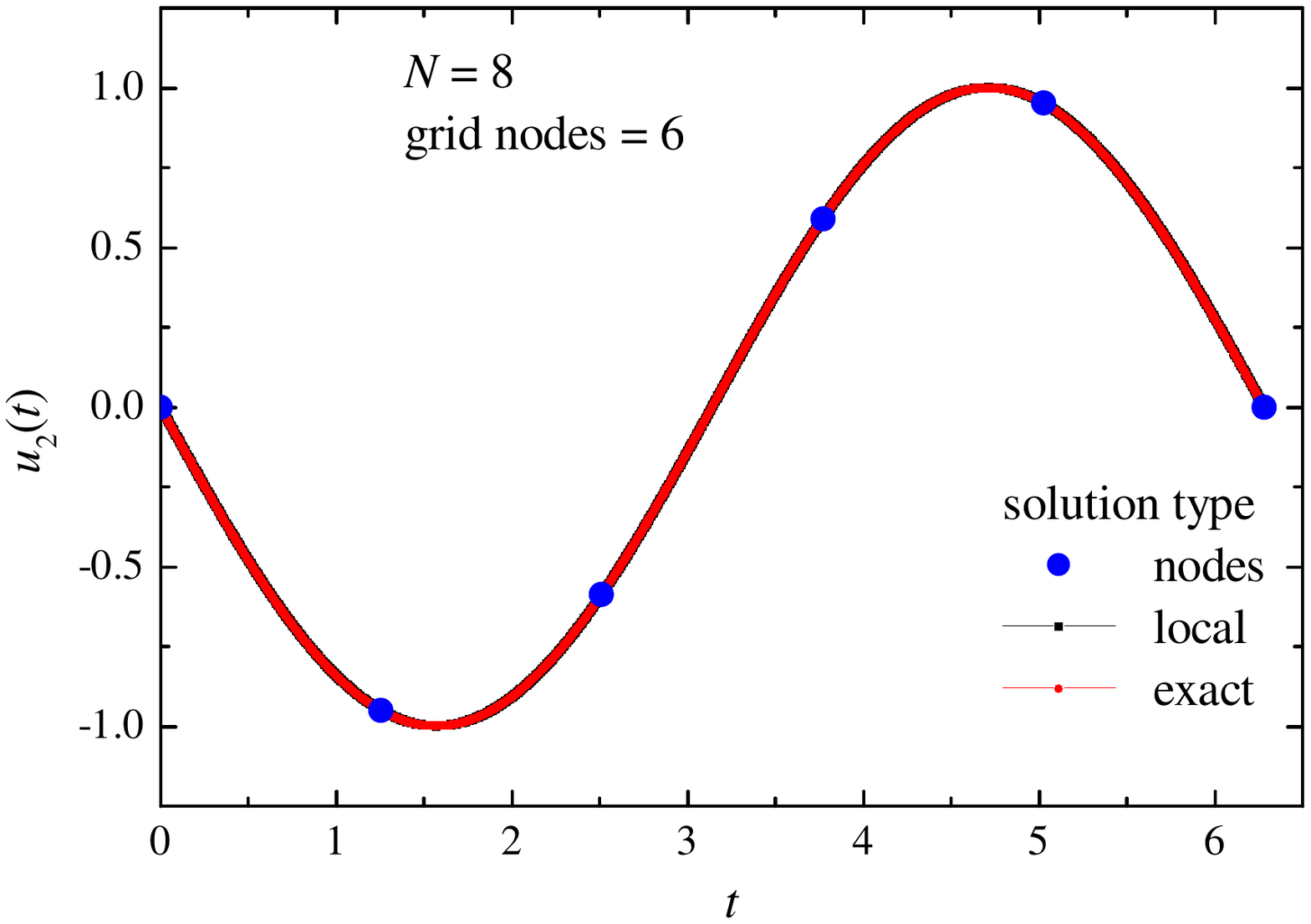}
\vspace{-8mm}\caption{\label{fig:harm_osc_sols:b2}}
\end{subfigure}
\begin{subfigure}{0.275\textwidth}
\includegraphics[width=\textwidth]{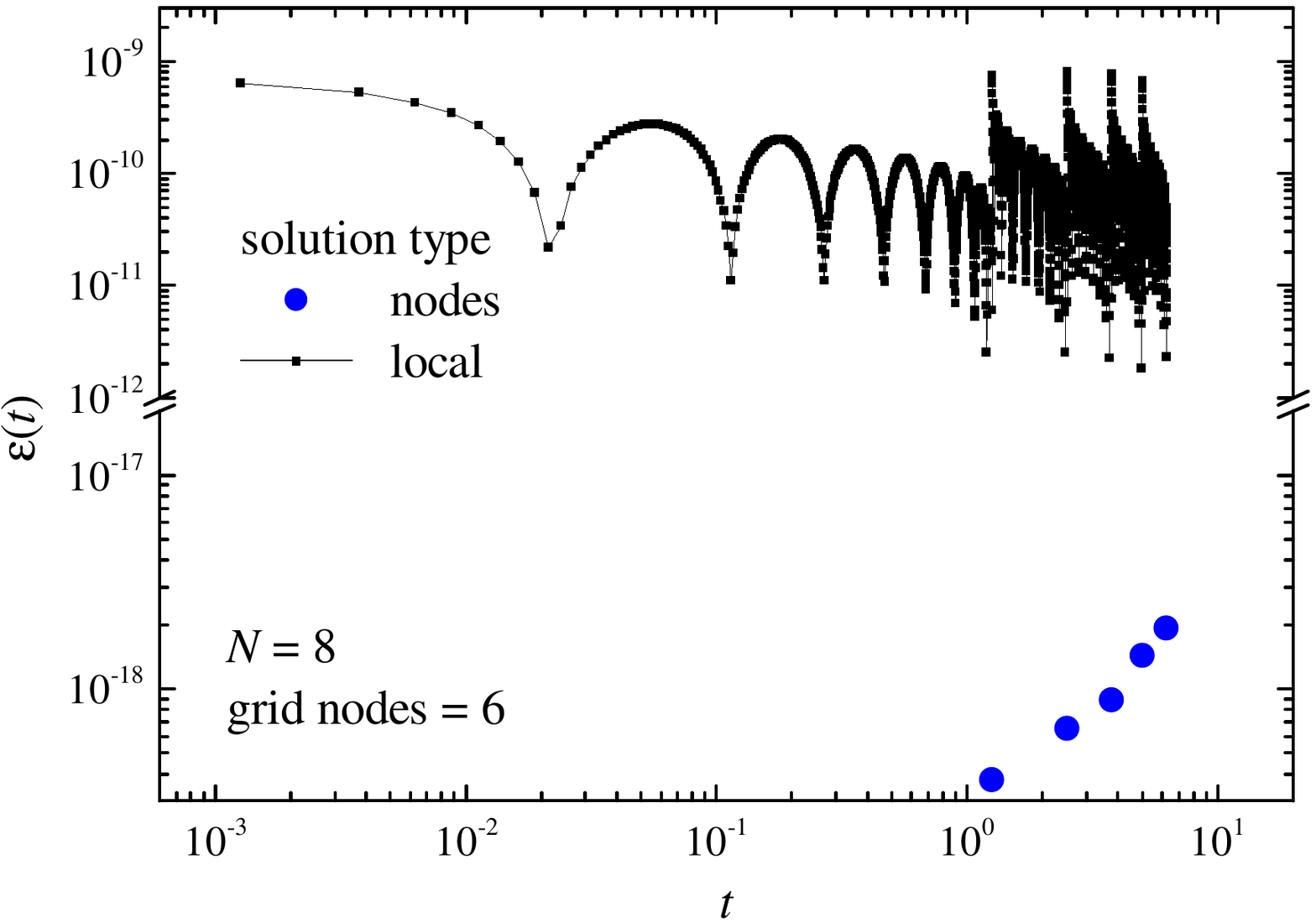}
\vspace{-8mm}\caption{\label{fig:harm_osc_sols:b3}}
\end{subfigure}\\
\begin{subfigure}{0.275\textwidth}
\includegraphics[width=\textwidth]{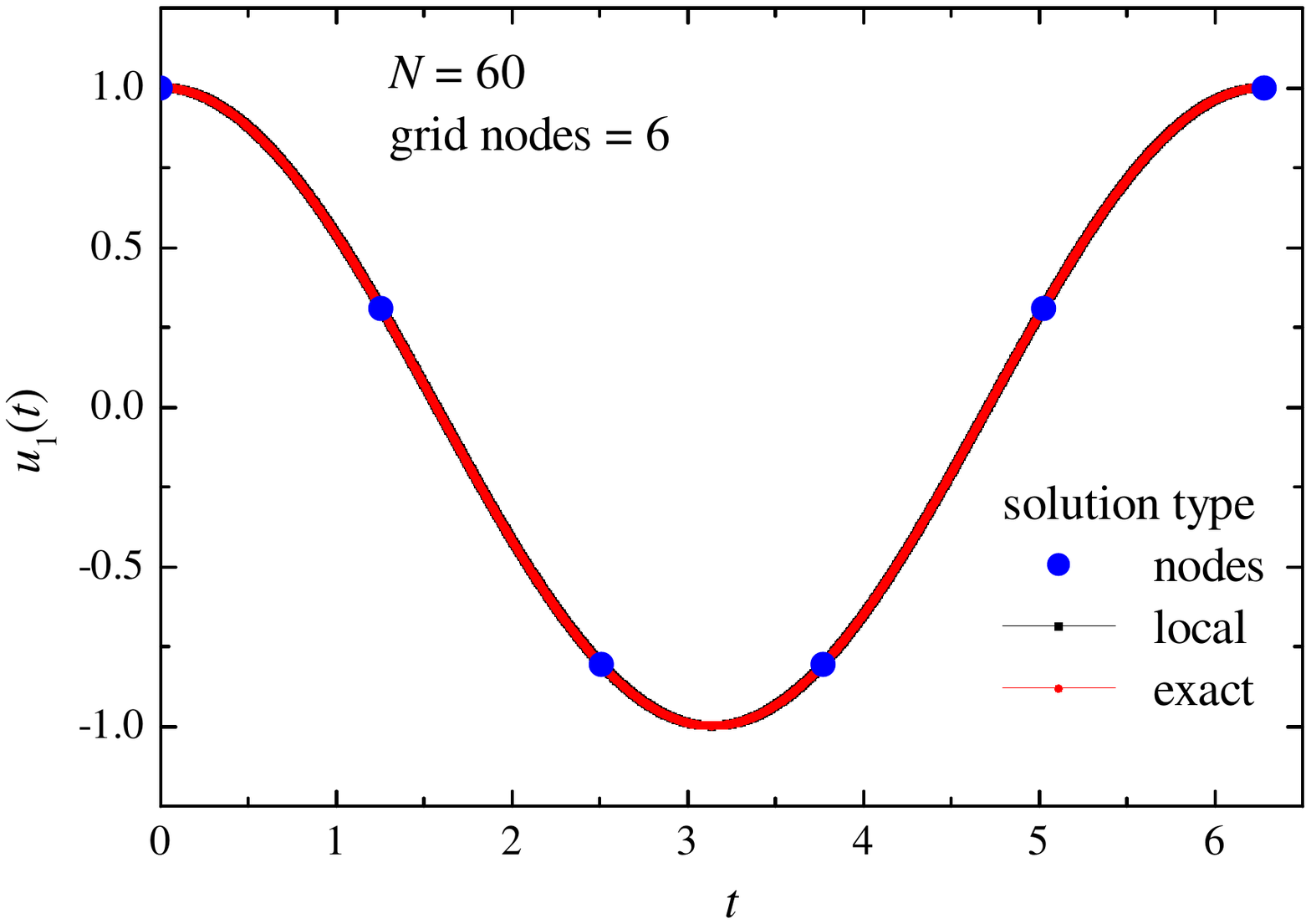}
\vspace{-8mm}\caption{\label{fig:harm_osc_sols:c1}}
\end{subfigure}
\begin{subfigure}{0.275\textwidth}
\includegraphics[width=\textwidth]{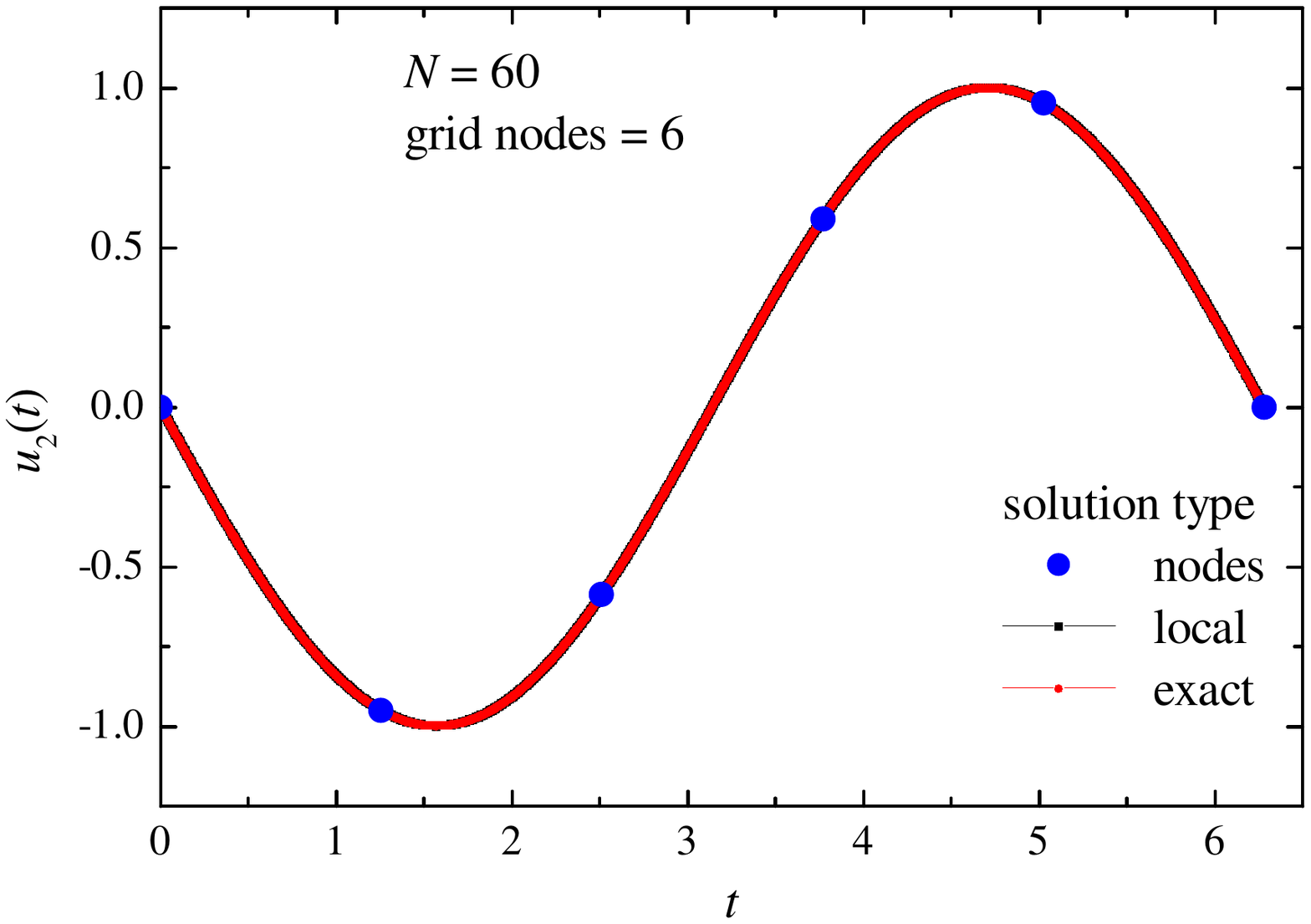}
\vspace{-8mm}\caption{\label{fig:harm_osc_sols:c2}}
\end{subfigure}
\begin{subfigure}{0.275\textwidth}
\includegraphics[width=\textwidth]{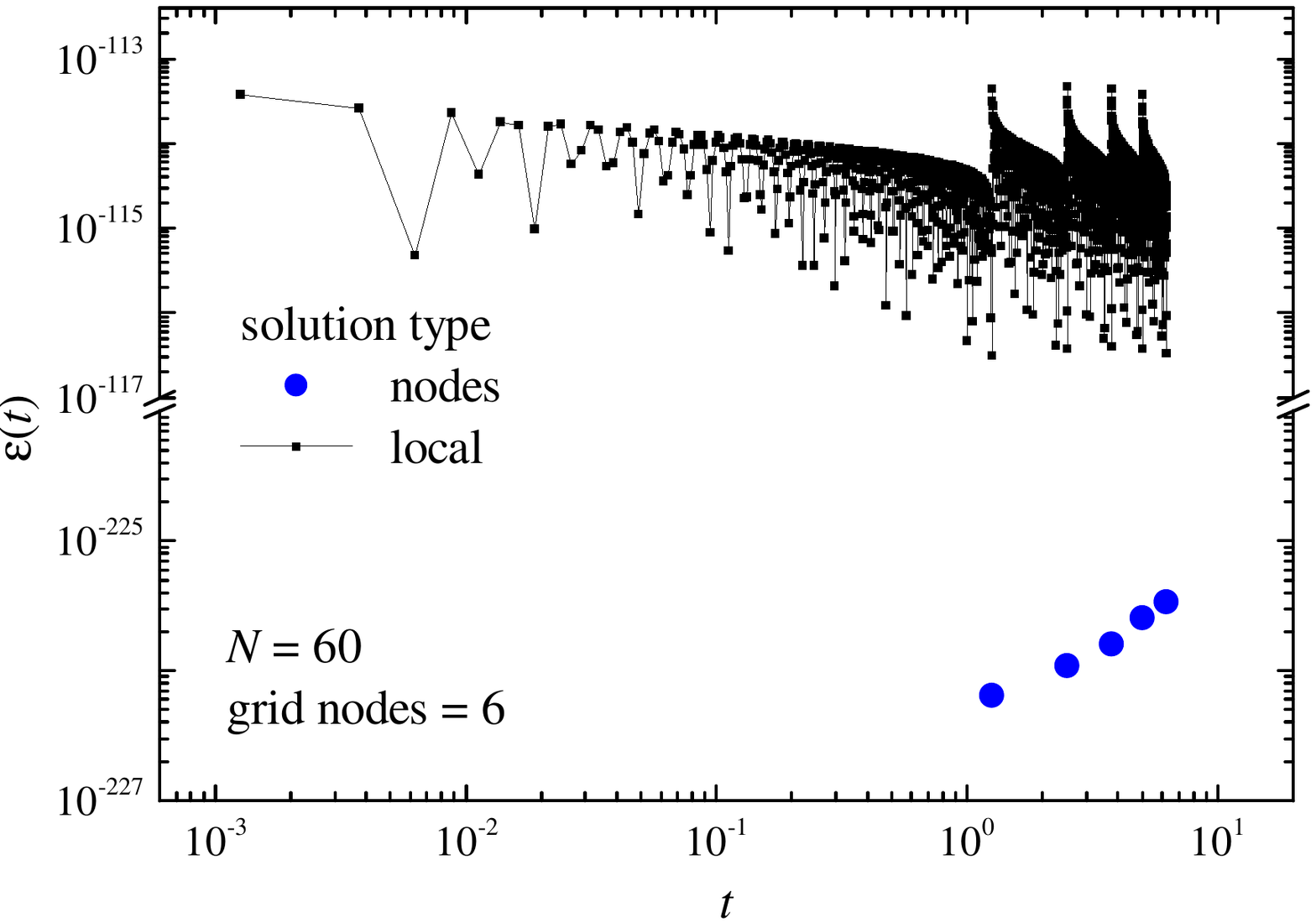}
\vspace{-8mm}\caption{\label{fig:harm_osc_sols:c3}}
\end{subfigure}\\
\begin{subfigure}{0.275\textwidth}
\includegraphics[width=\textwidth]{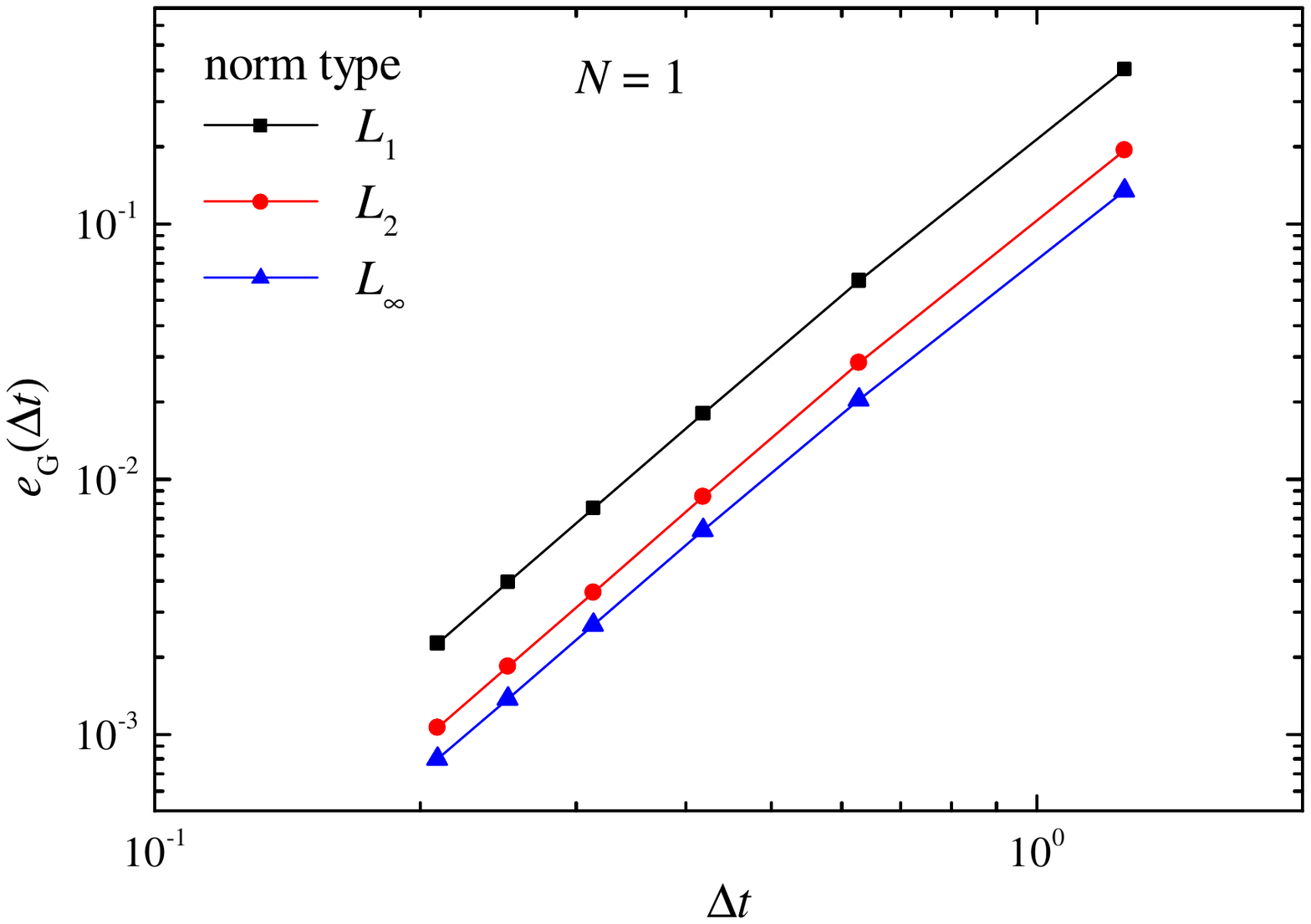}
\vspace{-8mm}\caption{\label{fig:harm_osc_sols:d1}}
\end{subfigure}
\begin{subfigure}{0.275\textwidth}
\includegraphics[width=\textwidth]{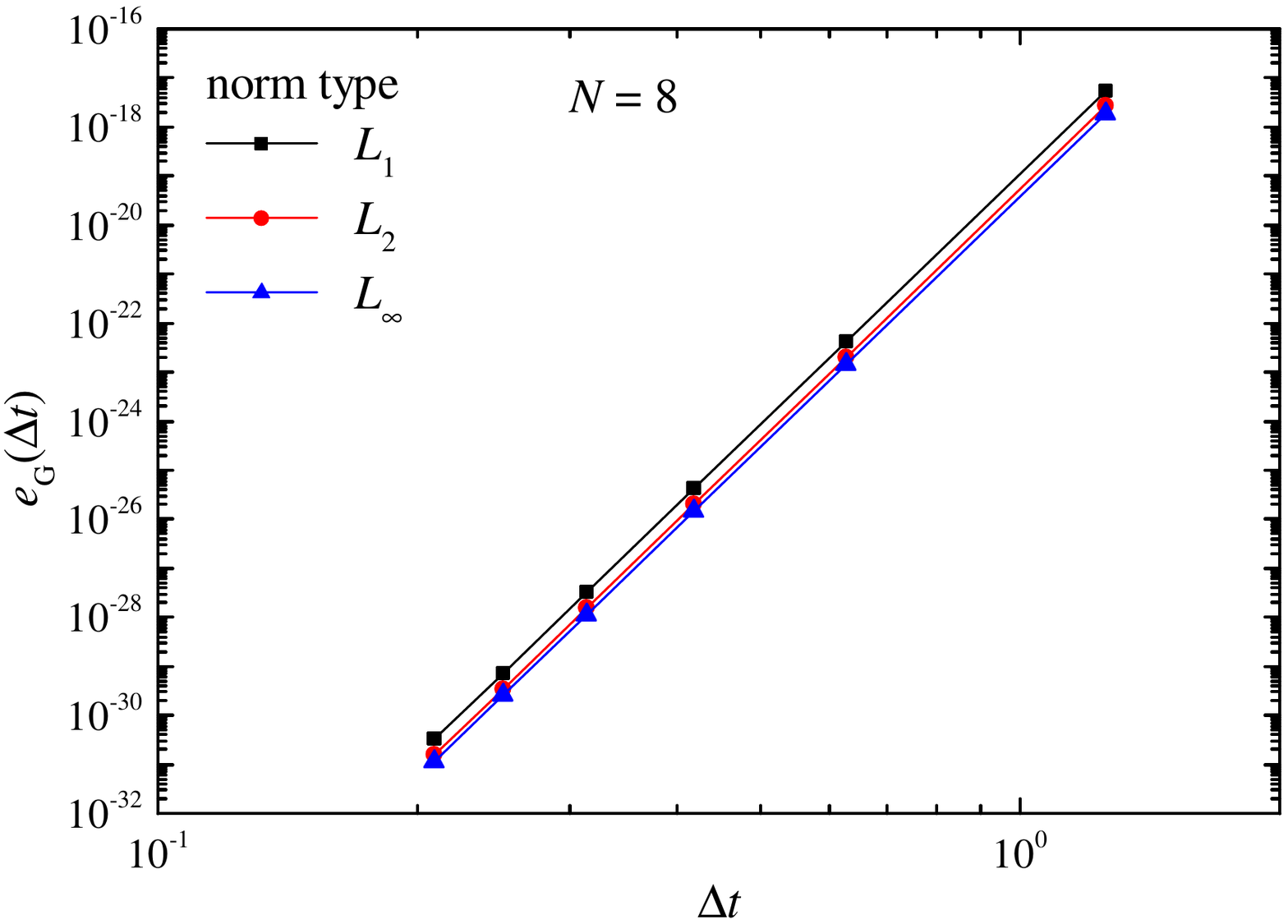}
\vspace{-8mm}\caption{\label{fig:harm_osc_sols:d2}}
\end{subfigure}
\begin{subfigure}{0.275\textwidth}
\includegraphics[width=\textwidth]{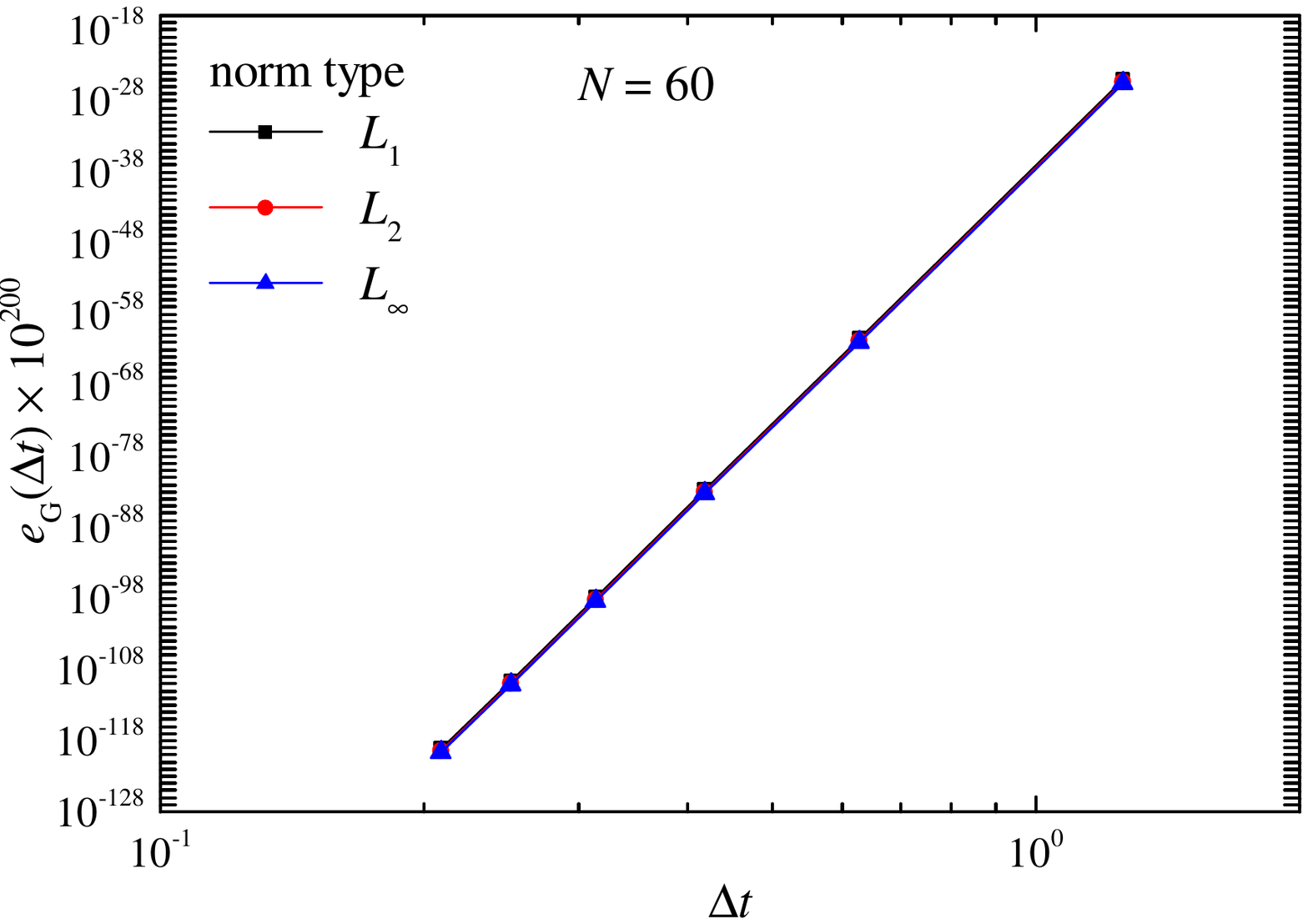}
\vspace{-8mm}\caption{\label{fig:harm_osc_sols:d3}}
\end{subfigure}\\
\begin{subfigure}{0.275\textwidth}
\includegraphics[width=\textwidth]{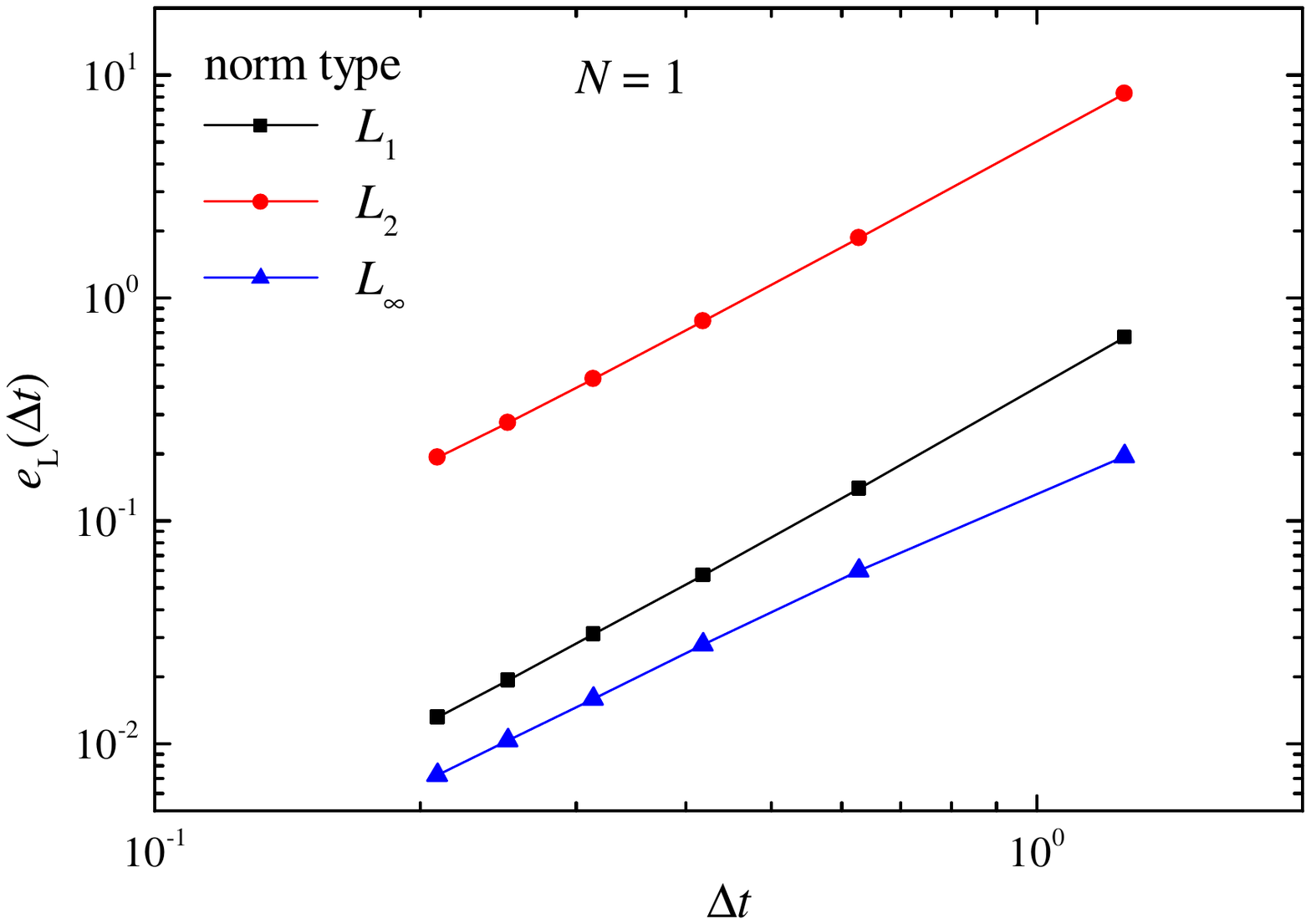}
\vspace{-8mm}\caption{\label{fig:harm_osc_sols:e1}}
\end{subfigure}
\begin{subfigure}{0.275\textwidth}
\includegraphics[width=\textwidth]{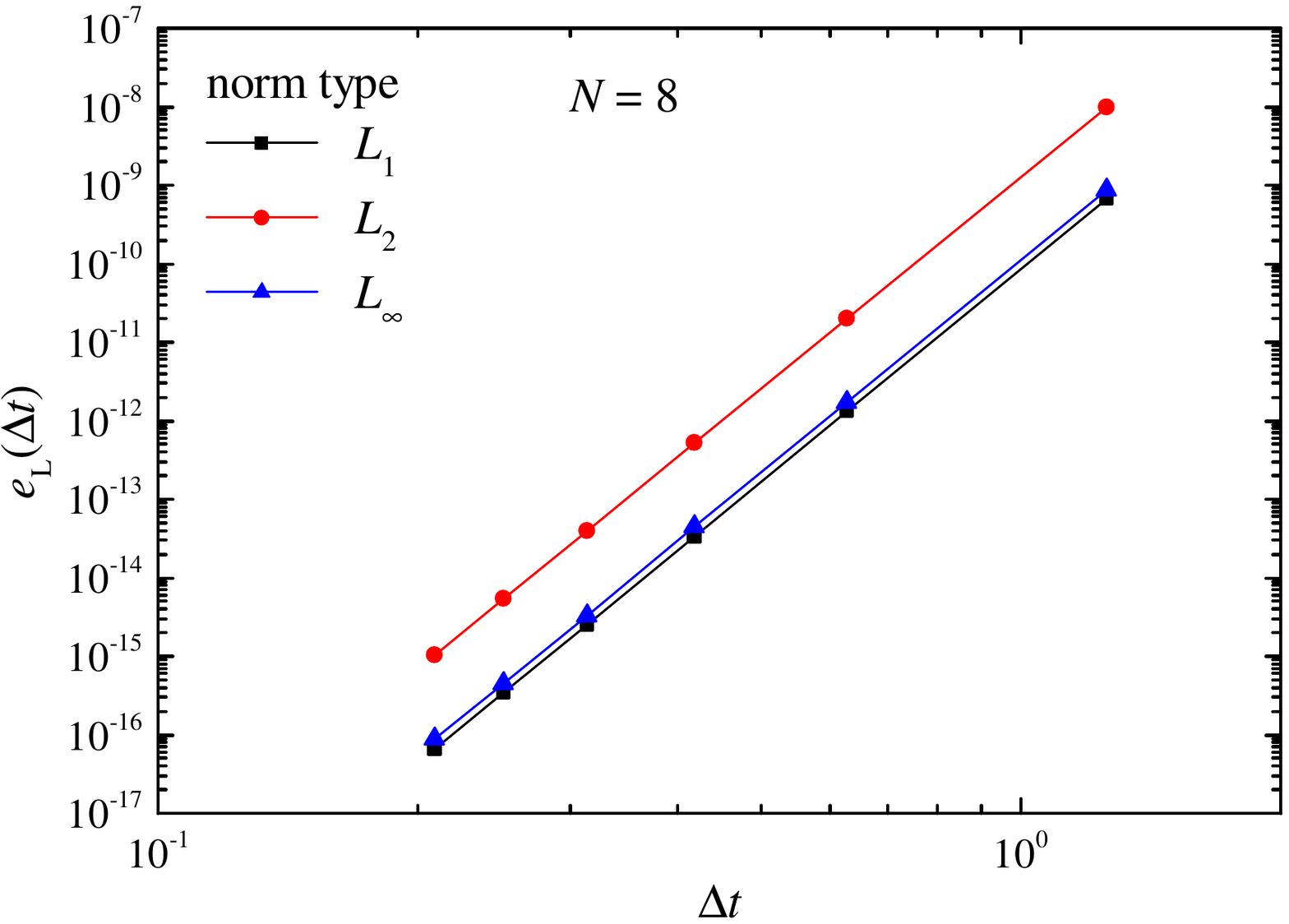}
\vspace{-8mm}\caption{\label{fig:harm_osc_sols:e2}}
\end{subfigure}
\begin{subfigure}{0.275\textwidth}
\includegraphics[width=\textwidth]{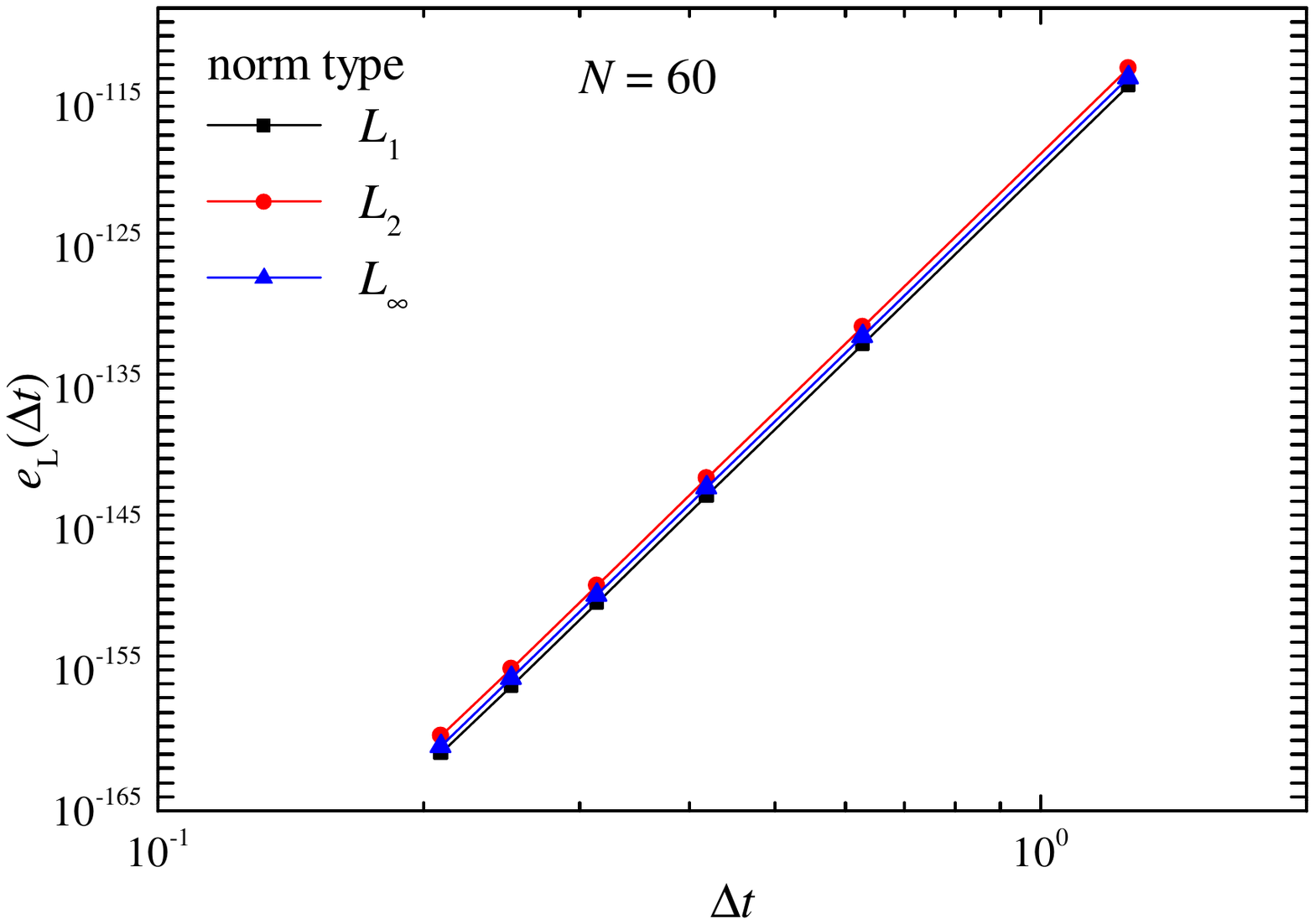}
\vspace{-8mm}\caption{\label{fig:harm_osc_sols:e3}}
\end{subfigure}\\
\caption{%
Numerical solution of the problem (\ref{eq:harm_osc}). Comparison of the solution at nodes $\mathbf{u}_{n}$, the local solution $\mathbf{u}_{L}(t)$ and the exact solution $\mathbf{u}^{\rm ex}(t)$ for components $u_{1}$ (\subref{fig:harm_osc_sols:a1}, \subref{fig:harm_osc_sols:b1}, \subref{fig:harm_osc_sols:c1}) and $u_{2}$ (\subref{fig:harm_osc_sols:a2}, \subref{fig:harm_osc_sols:b2}, \subref{fig:harm_osc_sols:c2}), the error $\varepsilon(t)$ (\subref{fig:harm_osc_sols:a3}, \subref{fig:harm_osc_sols:b3}, \subref{fig:harm_osc_sols:c3}), obtained using polynomials with degrees $N = 1$ (\subref{fig:harm_osc_sols:a1}, \subref{fig:harm_osc_sols:a2}, \subref{fig:harm_osc_sols:a3}), $N = 8$ (\subref{fig:harm_osc_sols:b1}, \subref{fig:harm_osc_sols:b2}, \subref{fig:harm_osc_sols:b3}) and $N = 60$ (\subref{fig:harm_osc_sols:c1}, \subref{fig:harm_osc_sols:c2}, \subref{fig:harm_osc_sols:c3}). Dependence of the global error for the solution at nodes $e_{G}$ (\subref{fig:harm_osc_sols:d1}, \subref{fig:harm_osc_sols:d2}, \subref{fig:harm_osc_sols:d3}) and the local solution $e_{L}$ (\subref{fig:harm_osc_sols:e1}, \subref{fig:harm_osc_sols:e2}, \subref{fig:harm_osc_sols:e3}) on the discretization step $\mathrm{\Delta}t$, obtained in the norms $L_{1}$, $L_{2}$ and $L_{\infty}$, obtained using polynomials with degrees $N = 1$ (\subref{fig:harm_osc_sols:d1}, \subref{fig:harm_osc_sols:e1}), $N = 8$ (\subref{fig:harm_osc_sols:d2}, \subref{fig:harm_osc_sols:e2}) and $N = 60$ (\subref{fig:harm_osc_sols:d3}, \subref{fig:harm_osc_sols:e3}).
}
\label{fig:harm_osc_sols}
\end{figure} 
\begin{table}[h!]
\centering
\caption{%
Orders of convergence $p_{L_{1}}$, $p_{L_{2}}$, $p_{L_{\infty}}$ of the ADER-DG numerical method with local DG predictor calculated in norms $L_{1}$, $L_{2}$, $L_{\infty}$ for the problem (\ref{eq:harm_osc}); $N$ is the degree of the basis polynomials $\varphi_{p}$ used for the DG representations. The orders of convergence $p$ for the solution at nodes $\mathbf{u}_{n}$ are compared to the theoretical value $p_{\rm (G)} = 2N+1$, the orders of convergence $p$ for the local solution $\mathbf{u}_{L}$ represented between the nodes are compared to the theoretical value $p_{\rm (L)} = N+1$.
}
\label{tab:harm_osc_sols}
\begin{tabular}{@{}lllllllll@{}}
\toprule
& \multicolumn{4}{c}{order $p$ for solution at nodes} & \multicolumn{4}{c}{order $p$ for local solution} \\
\midrule
$N$ & $p_{L_{1}}$ & $p_{L_{2}}$ & $p_{L_{\infty}}$ & theor. & $p_{L_{1}}$ & $p_{L_{2}}$ & $p_{L_{\infty}}$ & theor. \\
\midrule
$1$	&	$2.90$	&	$2.91$	&	$2.87$		&	$3$	&	$2.19$	&	$2.10$	&	$1.84$		&	$2$	\\
$2$	&	$4.96$	&	$4.98$	&	$4.95$		&	$5$	&	$3.04$	&	$2.96$	&	$2.96$		&	$3$	\\
$3$	&	$6.97$	&	$7.00$	&	$6.97$		&	$7$	&	$4.01$	&	$3.95$	&	$3.98$		&	$4$	\\
$4$	&	$8.98$	&	$9.01$	&	$8.97$		&	$9$	&	$5.00$	&	$4.95$	&	$4.99$		&	$5$	\\
$5$	&	$10.98$	&	$11.01$	&	$10.98$		&	$11$	&	$6.00$	&	$5.95$	&	$5.99$		&	$6$	\\
$6$	&	$12.99$	&	$13.02$	&	$12.98$		&	$13$	&	$7.00$	&	$6.96$	&	$6.99$		&	$7$	\\
$7$	&	$14.99$	&	$15.02$	&	$14.99$		&	$15$	&	$8.00$	&	$7.96$	&	$7.99$		&	$8$	\\
$8$	&	$16.99$	&	$17.02$	&	$16.99$		&	$17$	&	$9.00$	&	$8.96$	&	$8.99$		&	$9$	\\
$9$	&	$18.99$	&	$19.02$	&	$18.99$		&	$19$	&	$10.00$	&	$9.96$	&	$10.00$		&	$10$	\\
$10$	&	$21.00$	&	$21.03$	&	$20.99$		&	$21$	&	$11.00$	&	$10.96$	&	$11.00$		&	$11$	\\
\midrule
$15$	&	$31.00$	&	$31.03$	&	$30.99$		&	$31$	&	$16.00$	&	$15.96$	&	$16.00$		&	$16$	\\
$20$	&	$41.00$	&	$41.03$	&	$40.99$		&	$41$	&	$21.00$	&	$20.96$	&	$21.00$		&	$21$	\\
$25$	&	$51.00$	&	$51.03$	&	$51.00$		&	$51$	&	$26.00$	&	$25.96$	&	$26.00$		&	$26$	\\
$30$	&	$61.01$	&	$61.04$	&	$61.00$		&	$61$	&	$31.00$	&	$30.96$	&	$31.00$		&	$31$	\\
$35$	&	$71.01$	&	$71.04$	&	$71.00$		&	$71$	&	$36.00$	&	$35.96$	&	$36.00$		&	$36$	\\
$40$	&	$81.01$	&	$81.04$	&	$81.00$		&	$81$	&	$41.00$	&	$40.96$	&	$41.00$		&	$41$	\\
$45$	&	$91.01$	&	$91.04$	&	$91.00$		&	$91$	&	$46.00$	&	$45.96$	&	$46.00$		&	$46$	\\
$50$	&	$101.01$	&	$101.04$	&	$101.00$		&	$101$	&	$51.00$	&	$50.96$	&	$51.00$		&	$51$	\\
$55$	&	$111.01$	&	$111.04$	&	$111.00$		&	$111$	&	$56.00$	&	$55.96$	&	$56.00$		&	$56$	\\
$60$	&	$121.01$	&	$121.04$	&	$121.00$		&	$121$	&	$61.00$	&	$60.96$	&	$61.00$		&	$61$	\\
\bottomrule
\end{tabular}
\end{table} 

For the first example the IVP was chosen for a linear one-dimensional harmonic oscillator described by a second-order differential equation $\ddot{x} + x = 0$, with initial conditions on the coordinate $x(0) = 1$ and velocity $\dot{x}(0) = 0$. The IVP for a first-order ODE system corresponding to this problem has the following form: 
\begin{equation}\label{eq:harm_osc}
\frac{du_{1}}{dt} = u_{2};\ \frac{du_{2}}{dt} = -u_{1};\qquad 
u_{1}(0) = 1;\ u_{2}(0) = 0;
\end{equation}
where $\mathbf{u} = [u_{1},\, u_{2}]^{T} = [x,\, \dot{x}]^{T}$ is a desired two component vector function. The exact analytical solution of this problem has the form $\mathbf{u}^{\rm ex} = [\cos(t),\, -\sin(t)]^{T}$. The solution definition domain $\Omega$ was chosen to be equal to one oscillation period as $0 \leqslant t \leqslant 2\pi$.

The study of the convergence of the numerical solution was carried out on the basis of the analysis of the solution for the set of $6$ different uniform discretizations of the domain of definition $\Omega$ of the desired function, with the number of grid nodes $L = 6$, $11$, $16$, $21$, $26$, $31$; note that the zero node, in which the initial condition was defined, also refers to grid nodes, so the discretization step was determined by the expression $\Delta t = 2\pi/(L-1)$. In contrast to the demonstration example shown in Fig.~\ref{fig:harm_osc_sols_prec}, in this case, arbitrary precision arithmetic for floating point numbers was used for high-precision calculations, which made it possible to correctly calculate errors that are much smaller than round-off errors for standard numbers of built-in float and double types.

The results of the study of the application of the numerical method for this problem are presented in Fig.~\ref{fig:harm_osc_sols} and in Table~\ref{tab:harm_osc_sols}. An analysis of the error $\varepsilon(t)$ for $N > 1$ shows that the error scales for the local solution and the solution at nodes differ by several orders of magnitude, reaching a value of $100$ orders of magnitude for a degree of $N = 60$. The calculated convergence orders $p_{L_{1}}$, $p_{L_{2}}$, $p_{L_{\infty}}$ for various functional norms are presented in Table~\ref{tab:harm_osc_sols}.

All obtained convergence orders $p$ sufficiently correspond to the expected theoretical values $p_{\rm (G)}$ and $p_{\rm (L)}$. There is a slight downward mismatch for low degrees of polynomials $N = 1,\, 2$. As a result, it becomes clear that the solution at nodes $\mathbf{u}_{n}$ obtained by the ADER-DG numerical method with a local DG predictor demonstrates the classical superconvergence $2N+1$ expected from such numerical methods, at least in problem (\ref{eq:harm_osc}). The local solution $\mathbf{u}_{L}$ obtained by the DG predictor demonstrates the classical convergence characteristic of DG methods.

\subsubsection{Example 2}

For the second example the IVP was chosen for a linear second-order differential equation $\ddot{x} - x = 0$ characterized by two exponential solutions (conveniently expressed in terms of $\sinh$ and $\cosh$ functions), with initial conditions on the $x(0) = 0$ and $\dot{x}(0) = 1$. The IVP for a first-order ODE system corresponding to this problem has the following form: 
\begin{equation}\label{eq:exp_diss}
\frac{du_{1}}{dt} = u_{2};\ \frac{du_{2}}{dt} = u_{1};\qquad
u_{1}(0) = 0;\ u_{2}(0) = 1;
\end{equation}
where $\mathbf{u} = [u_{1},\, u_{2}]^{T} = [x,\, \dot{x}]^{T}$ is a desired two component vector function. The exact analytical solution of this problem has the form $\mathbf{u}^{\rm ex} = [\sinh(t),\, \cosh(t)]^{T}$. The solution definition domain $\Omega$ was chosen as $0 \leqslant t \leqslant 2\pi$, similar to the previous problem (\ref{eq:harm_osc}).

The study of the convergence of the numerical solution was carried out on the basis of the analysis of the solution for the set of $6$ different uniform discretizations of the domain of definition $\Omega$ of the desired function, with the number of grid nodes $L = 6$, $11$, $16$, $21$, $26$, $31$; note that the zero node, in which the initial condition was defined, also refers to grid nodes, so the discretization step was determined by the expression $\Delta t = 2\pi/(L-1)$.

\begin{figure}[h!]
\captionsetup[subfigure]{%
position=bottom,
textfont=normalfont,
singlelinecheck=off,
justification=raggedright,
font=tiny
}
\centering
\begin{subfigure}{0.275\textwidth}
\includegraphics[width=\textwidth]{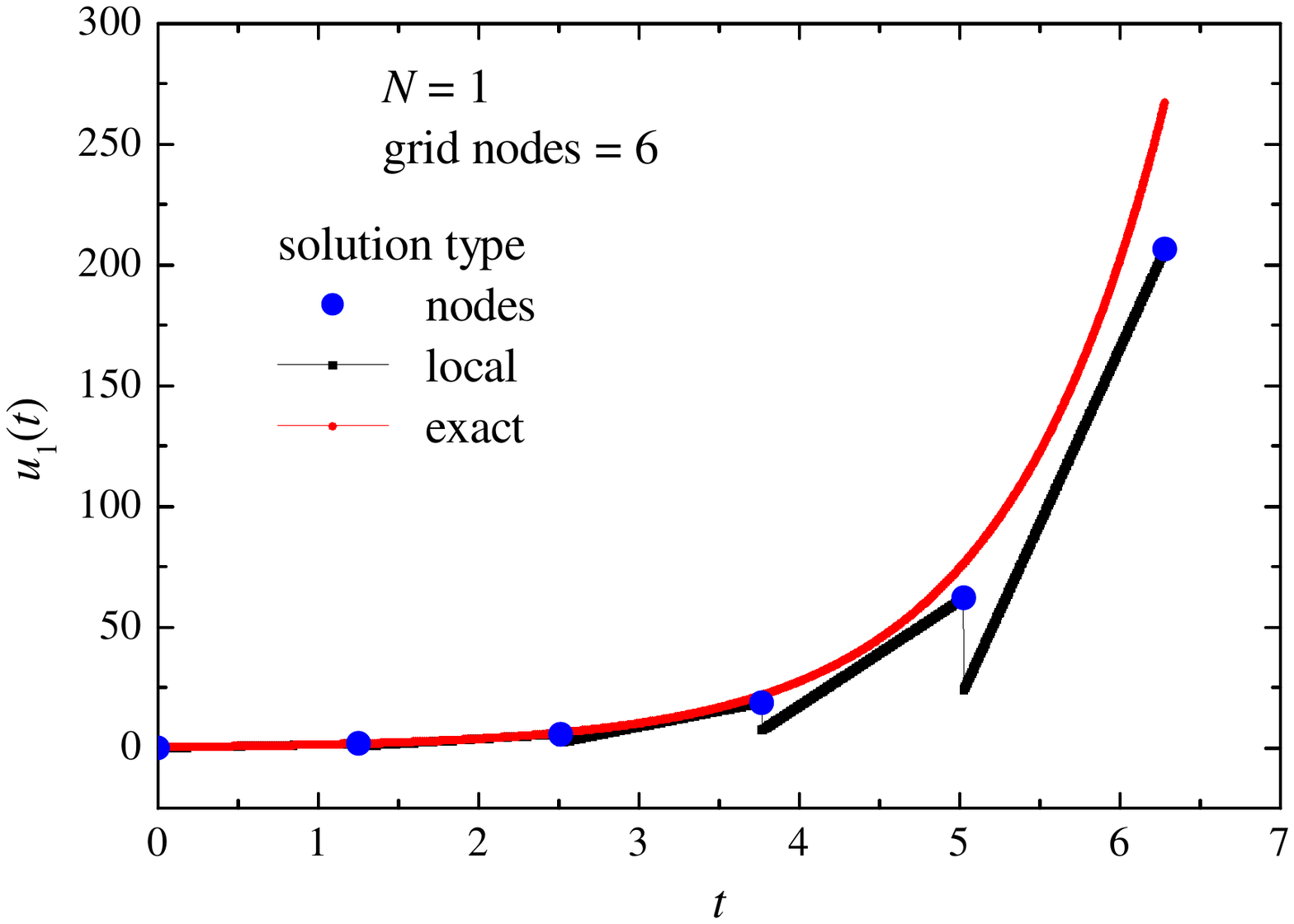}
\vspace{-8mm}\caption{\label{fig:exp_diss_sols:a1}}
\end{subfigure}
\begin{subfigure}{0.275\textwidth}
\includegraphics[width=\textwidth]{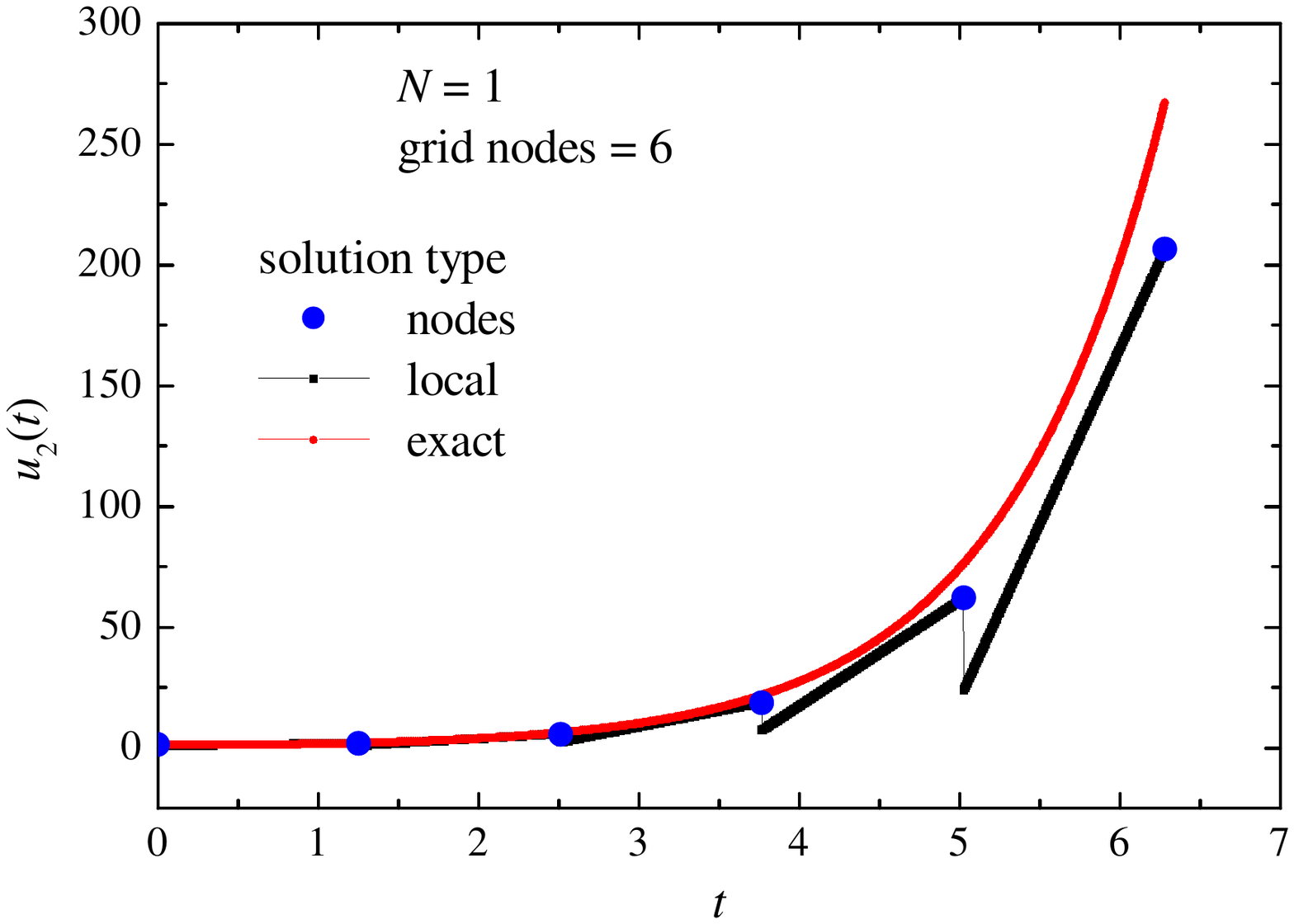}
\vspace{-8mm}\caption{\label{fig:exp_diss_sols:a2}}
\end{subfigure}
\begin{subfigure}{0.275\textwidth}
\includegraphics[width=\textwidth]{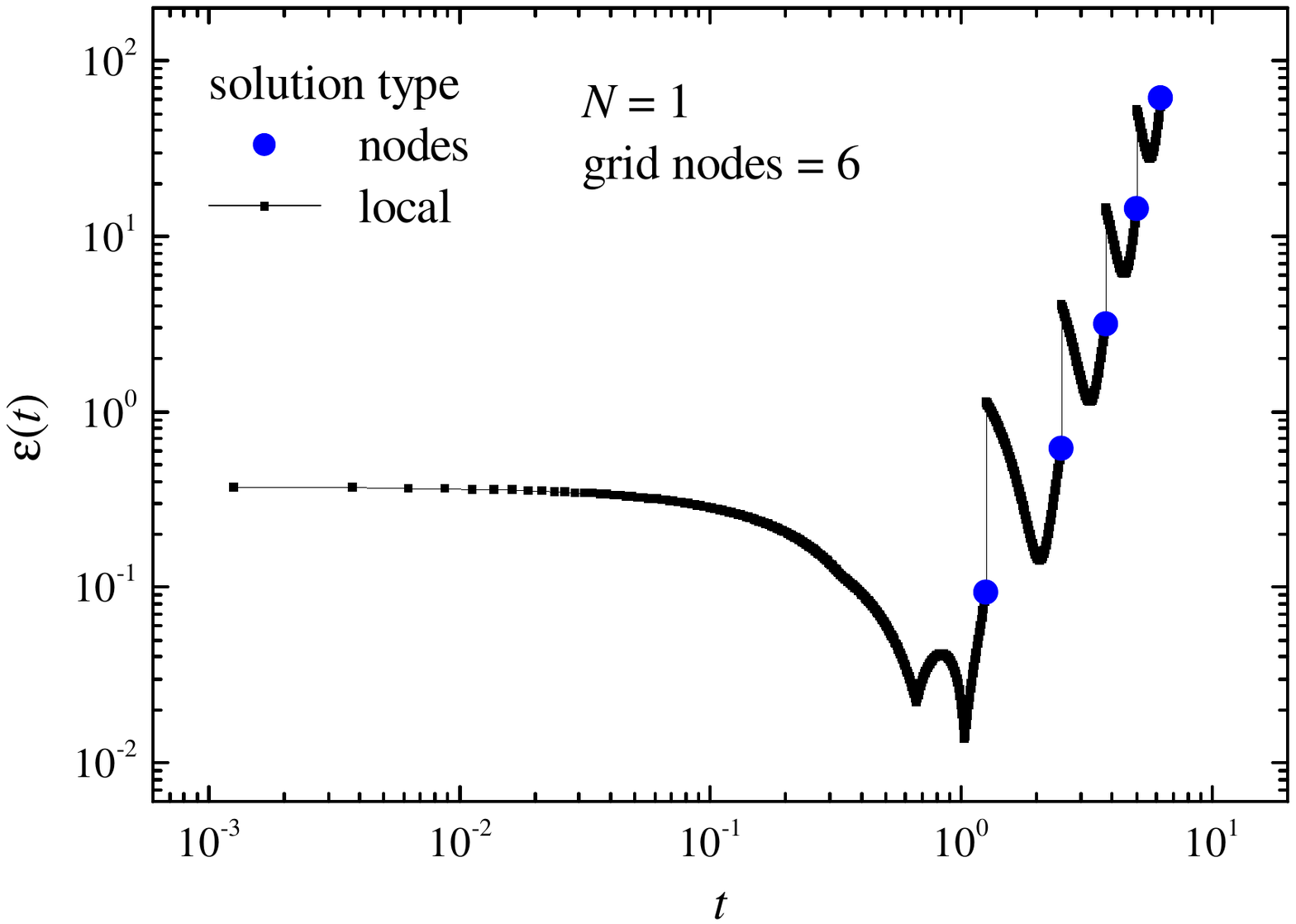}
\vspace{-8mm}\caption{\label{fig:exp_diss_sols:a3}}
\end{subfigure}\\
\begin{subfigure}{0.275\textwidth}
\includegraphics[width=\textwidth]{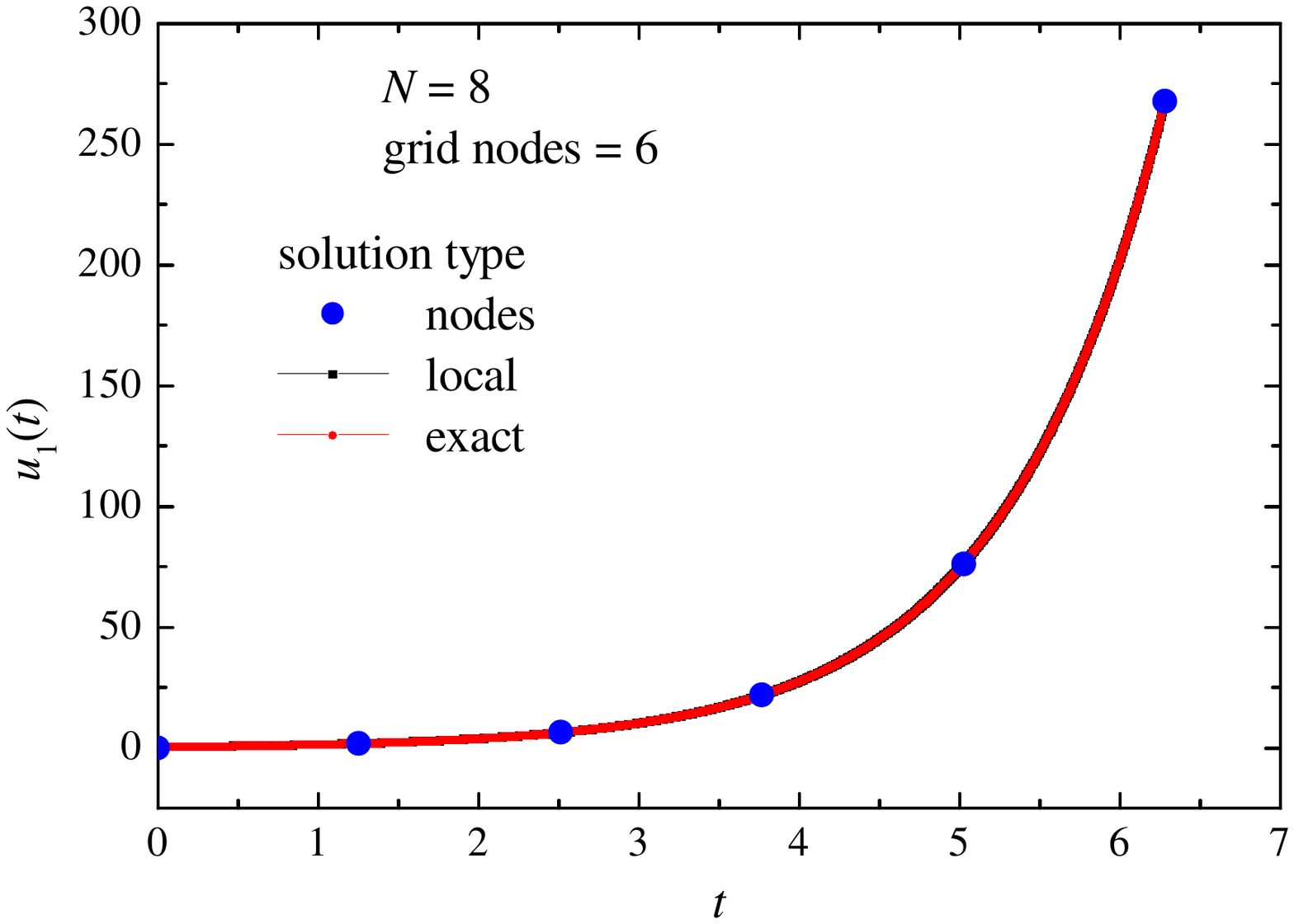}
\vspace{-8mm}\caption{\label{fig:exp_diss_sols:b1}}
\end{subfigure}
\begin{subfigure}{0.275\textwidth}
\includegraphics[width=\textwidth]{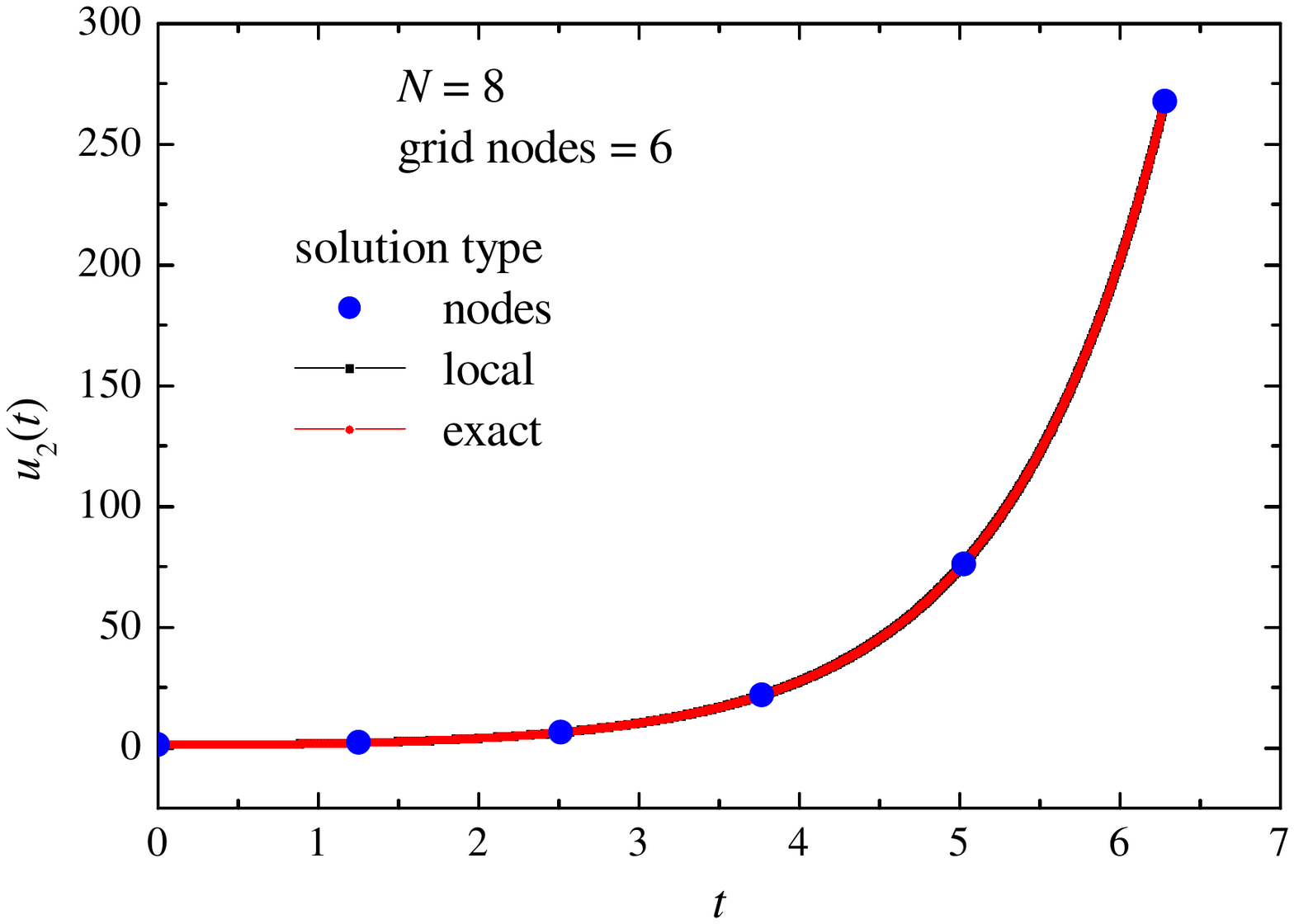}
\vspace{-8mm}\caption{\label{fig:exp_diss_sols:b2}}
\end{subfigure}
\begin{subfigure}{0.275\textwidth}
\includegraphics[width=\textwidth]{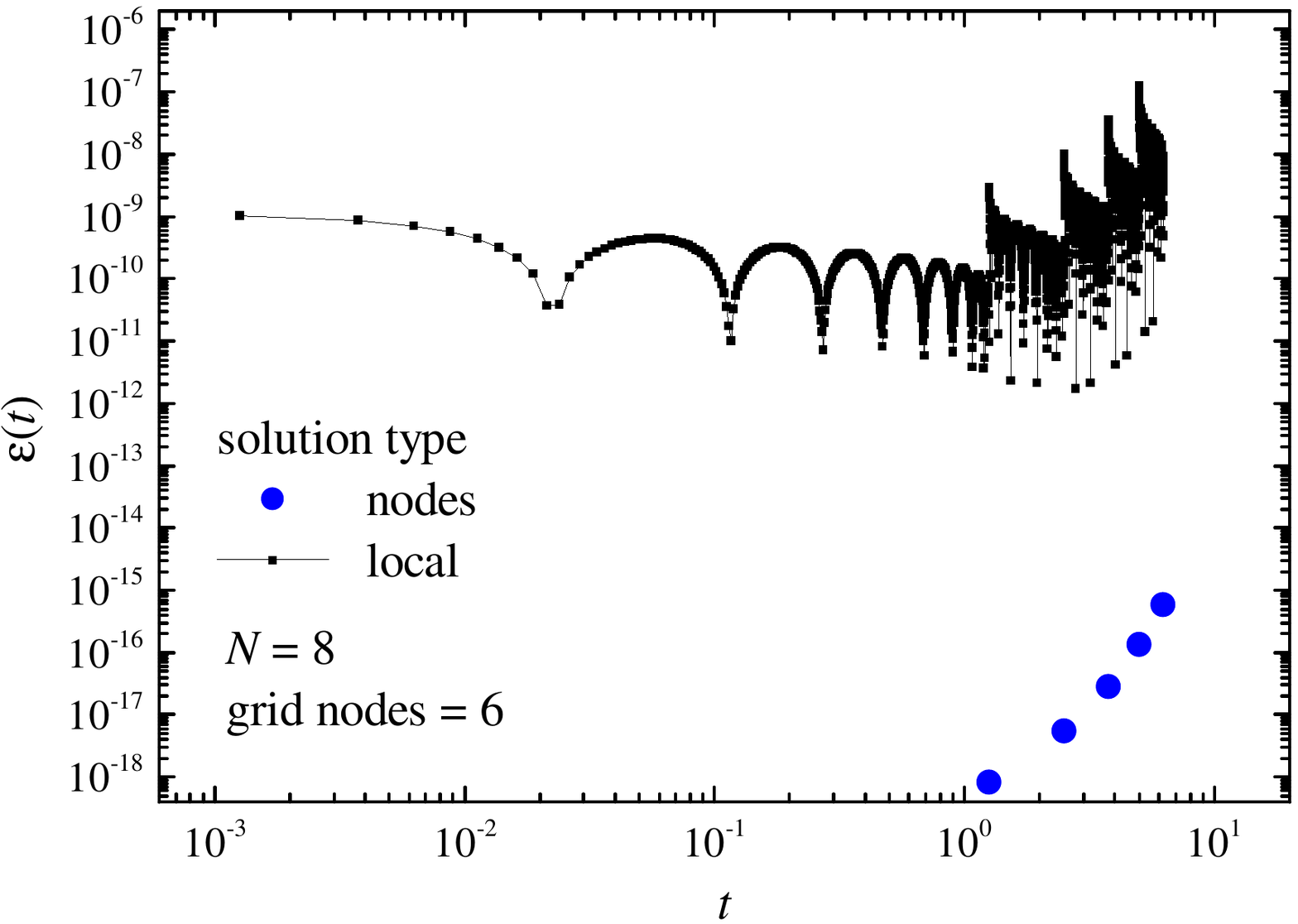}
\vspace{-8mm}\caption{\label{fig:exp_diss_sols:b3}}
\end{subfigure}\\
\begin{subfigure}{0.275\textwidth}
\includegraphics[width=\textwidth]{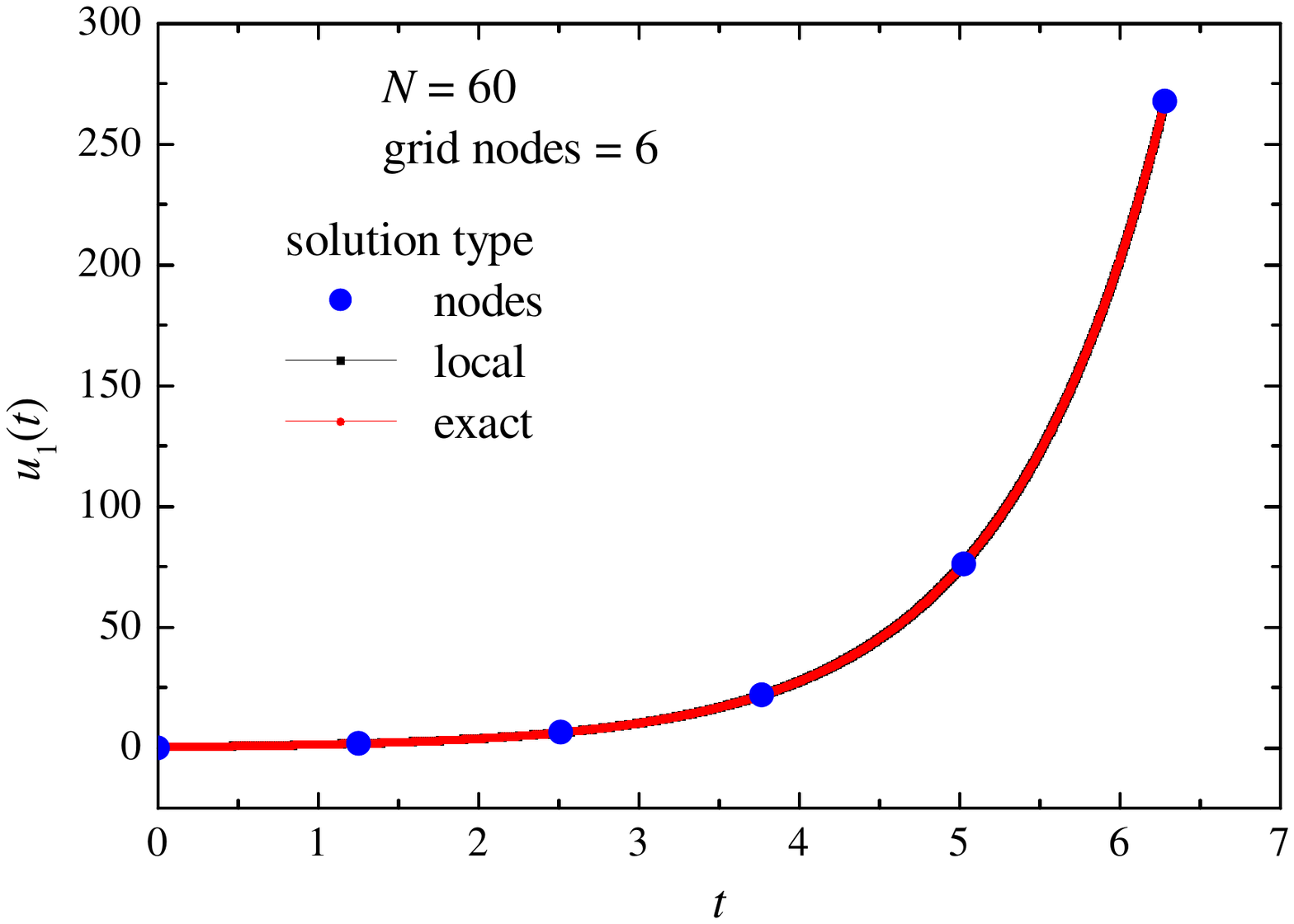}
\vspace{-8mm}\caption{\label{fig:exp_diss_sols:c1}}
\end{subfigure}
\begin{subfigure}{0.275\textwidth}
\includegraphics[width=\textwidth]{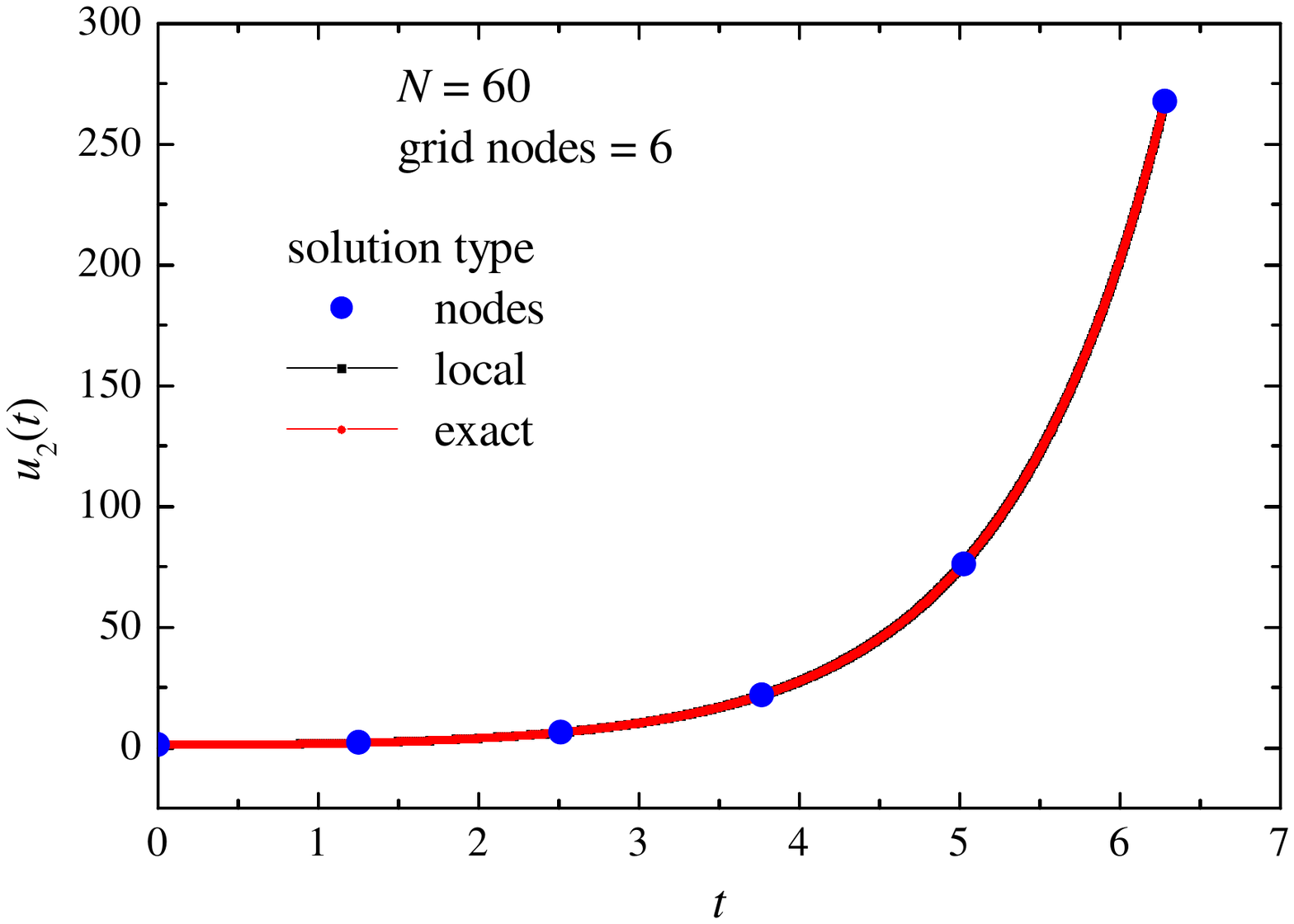}
\vspace{-8mm}\caption{\label{fig:exp_diss_sols:c2}}
\end{subfigure}
\begin{subfigure}{0.275\textwidth}
\includegraphics[width=\textwidth]{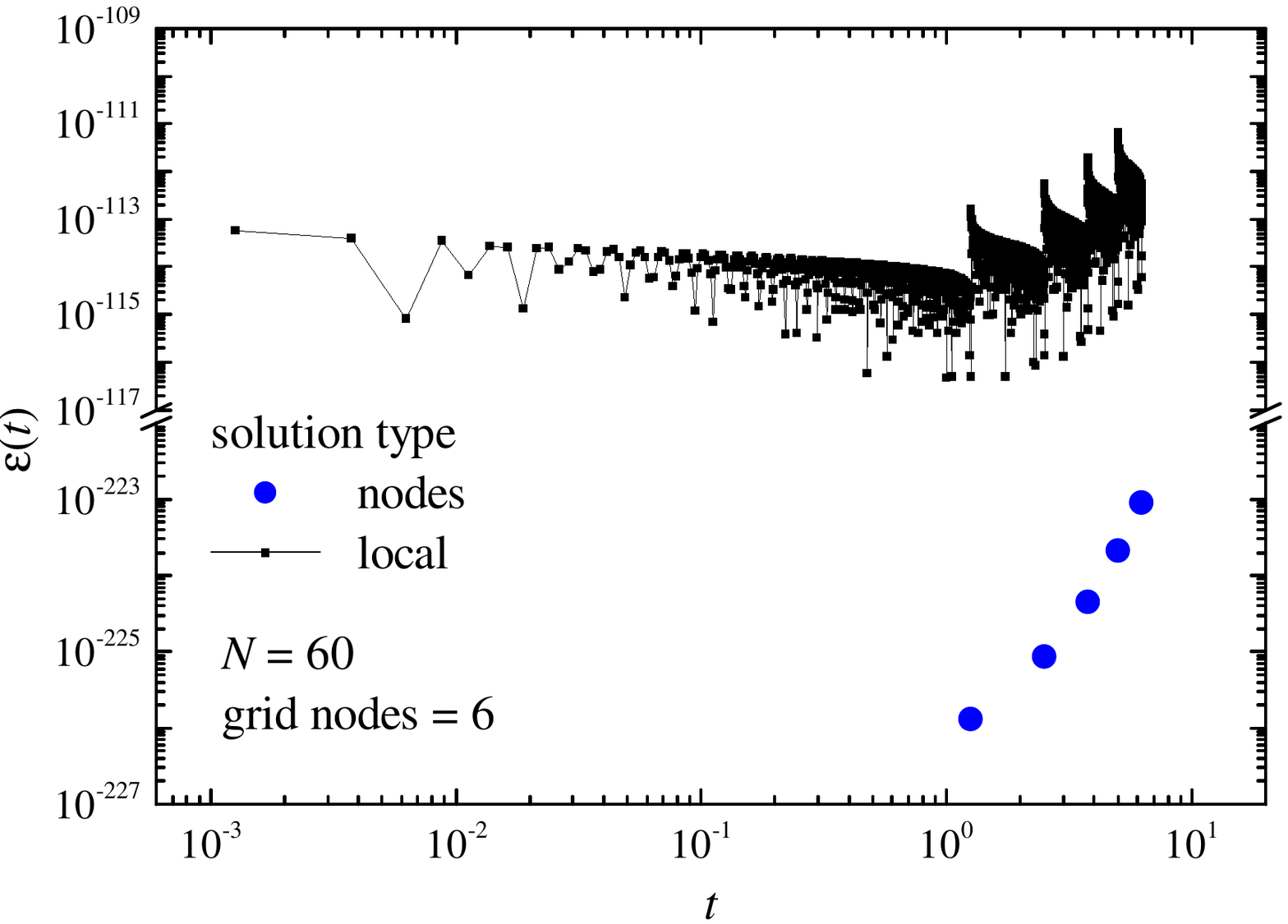}
\vspace{-8mm}\caption{\label{fig:exp_diss_sols:c3}}
\end{subfigure}\\
\begin{subfigure}{0.275\textwidth}
\includegraphics[width=\textwidth]{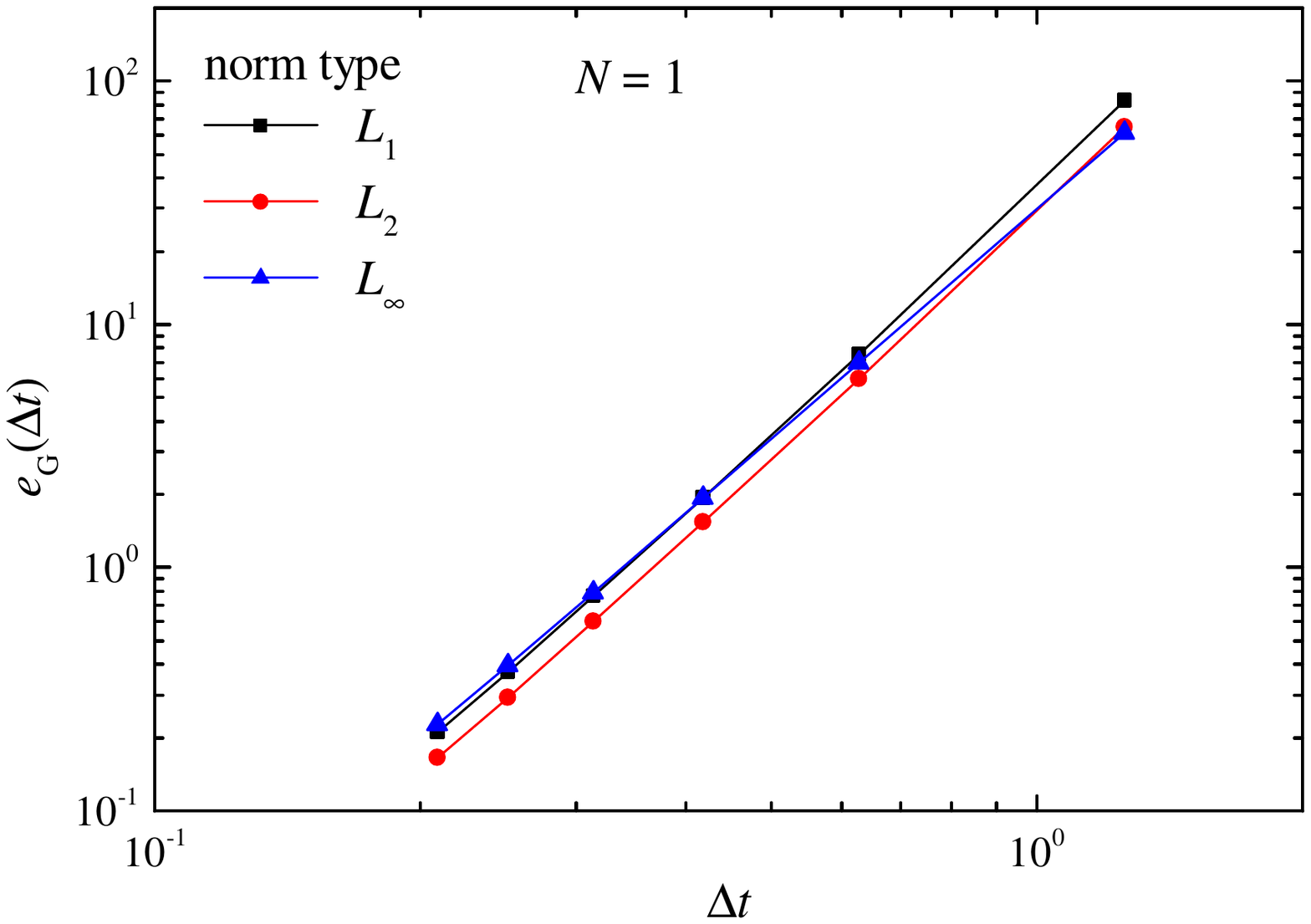}
\vspace{-8mm}\caption{\label{fig:exp_diss_sols:d1}}
\end{subfigure}
\begin{subfigure}{0.275\textwidth}
\includegraphics[width=\textwidth]{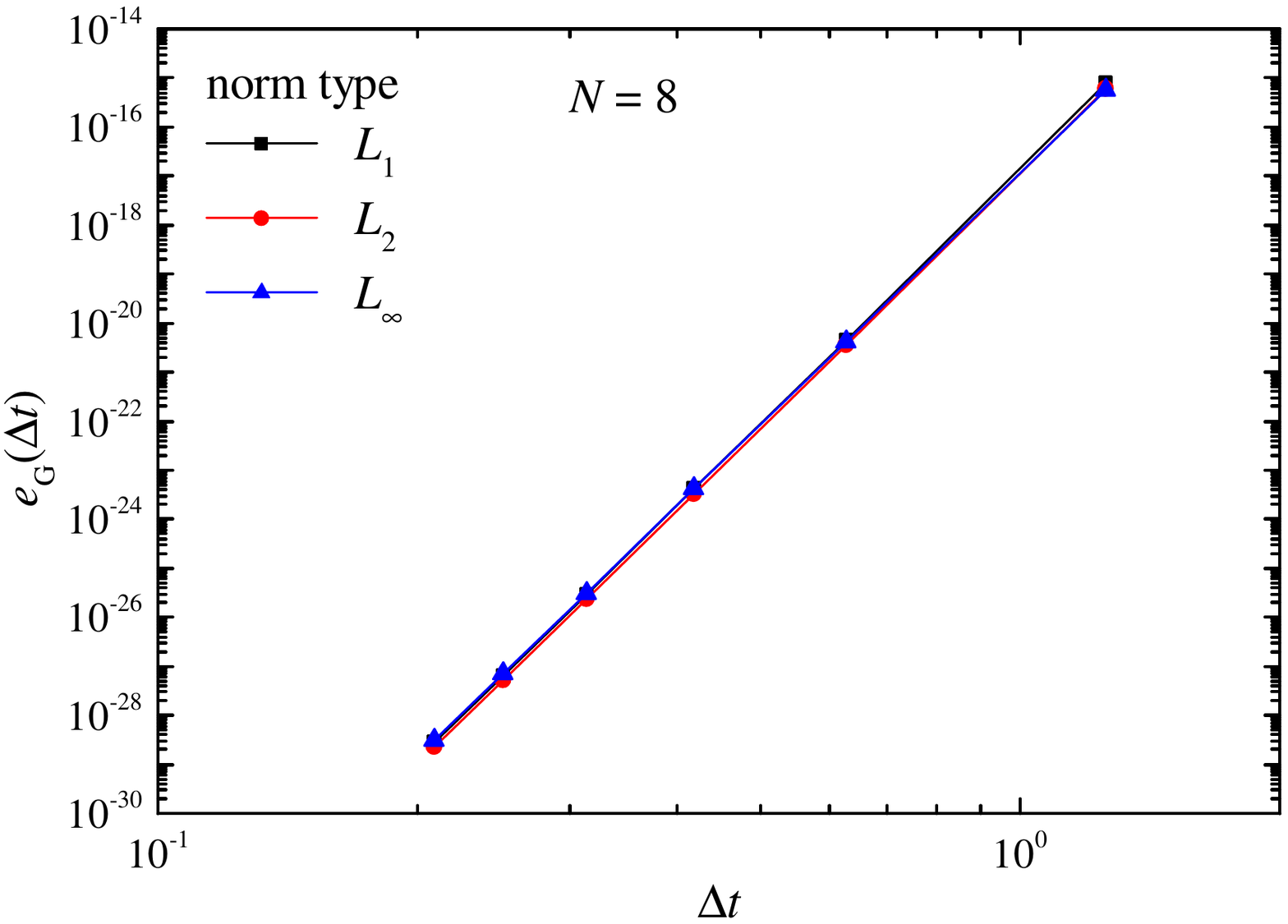}
\vspace{-8mm}\caption{\label{fig:exp_diss_sols:d2}}
\end{subfigure}
\begin{subfigure}{0.275\textwidth}
\includegraphics[width=\textwidth]{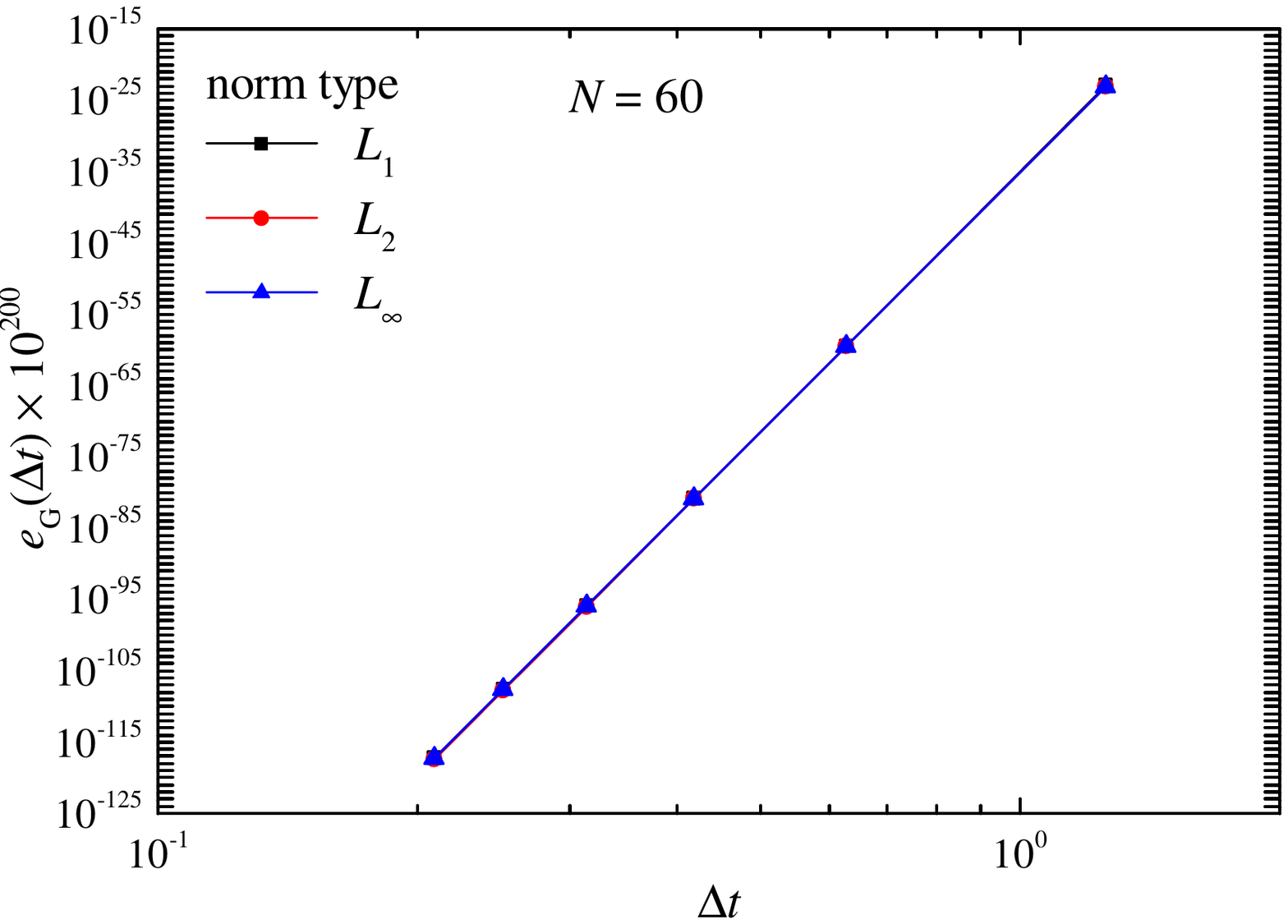}
\vspace{-8mm}\caption{\label{fig:exp_diss_sols:d3}}
\end{subfigure}\\
\begin{subfigure}{0.275\textwidth}
\includegraphics[width=\textwidth]{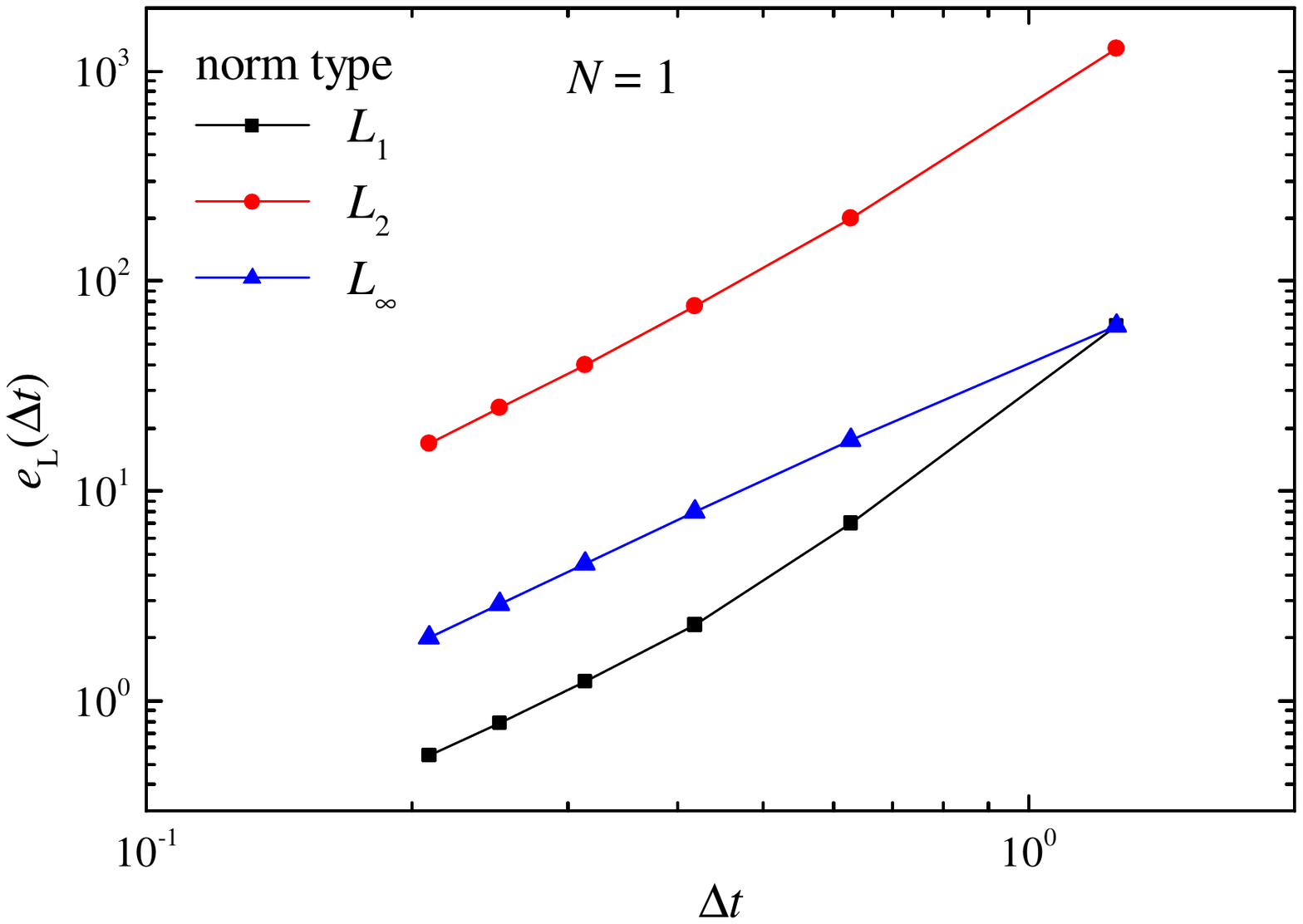}
\vspace{-8mm}\caption{\label{fig:exp_diss_sols:e1}}
\end{subfigure}
\begin{subfigure}{0.275\textwidth}
\includegraphics[width=\textwidth]{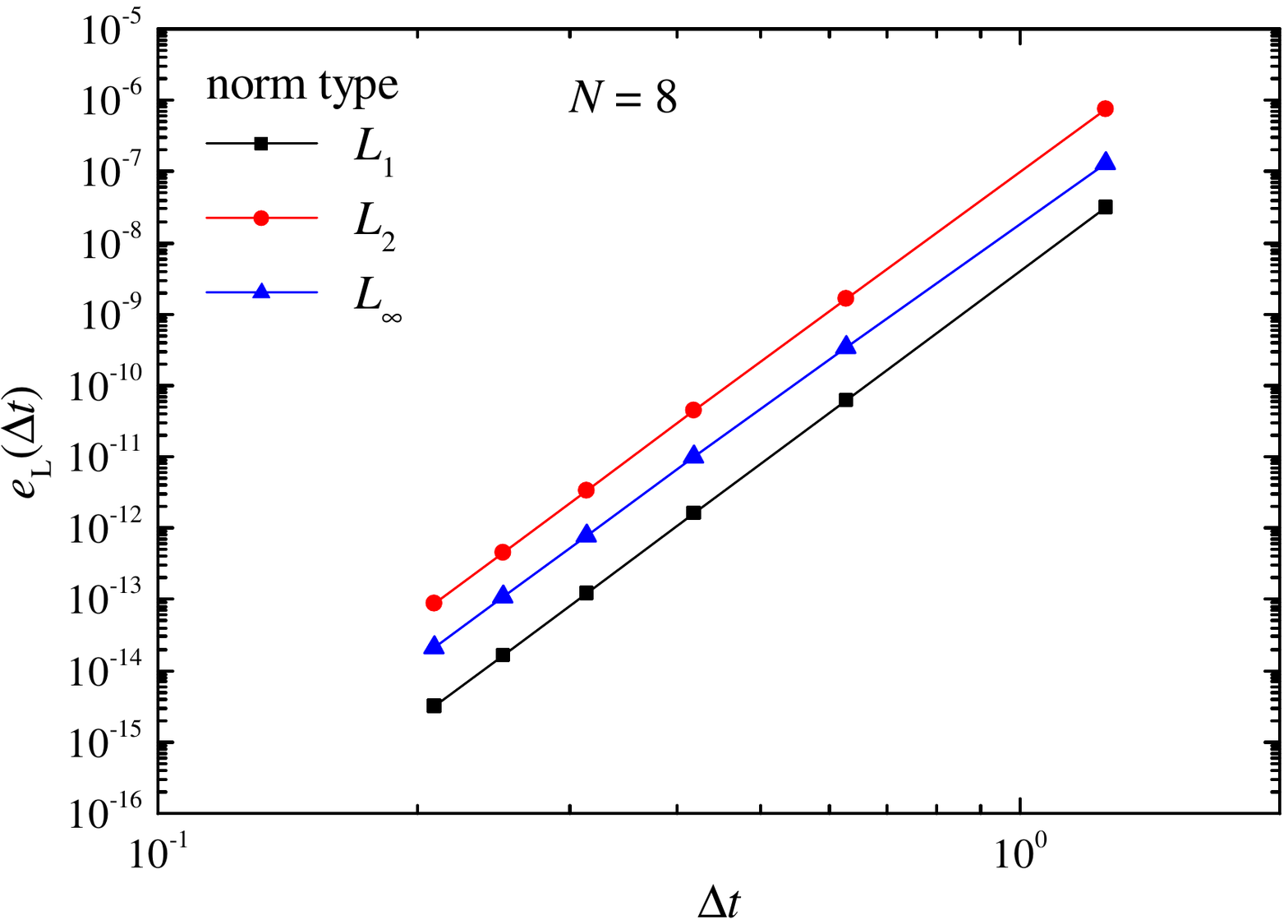}
\vspace{-8mm}\caption{\label{fig:exp_diss_sols:e2}}
\end{subfigure}
\begin{subfigure}{0.275\textwidth}
\includegraphics[width=\textwidth]{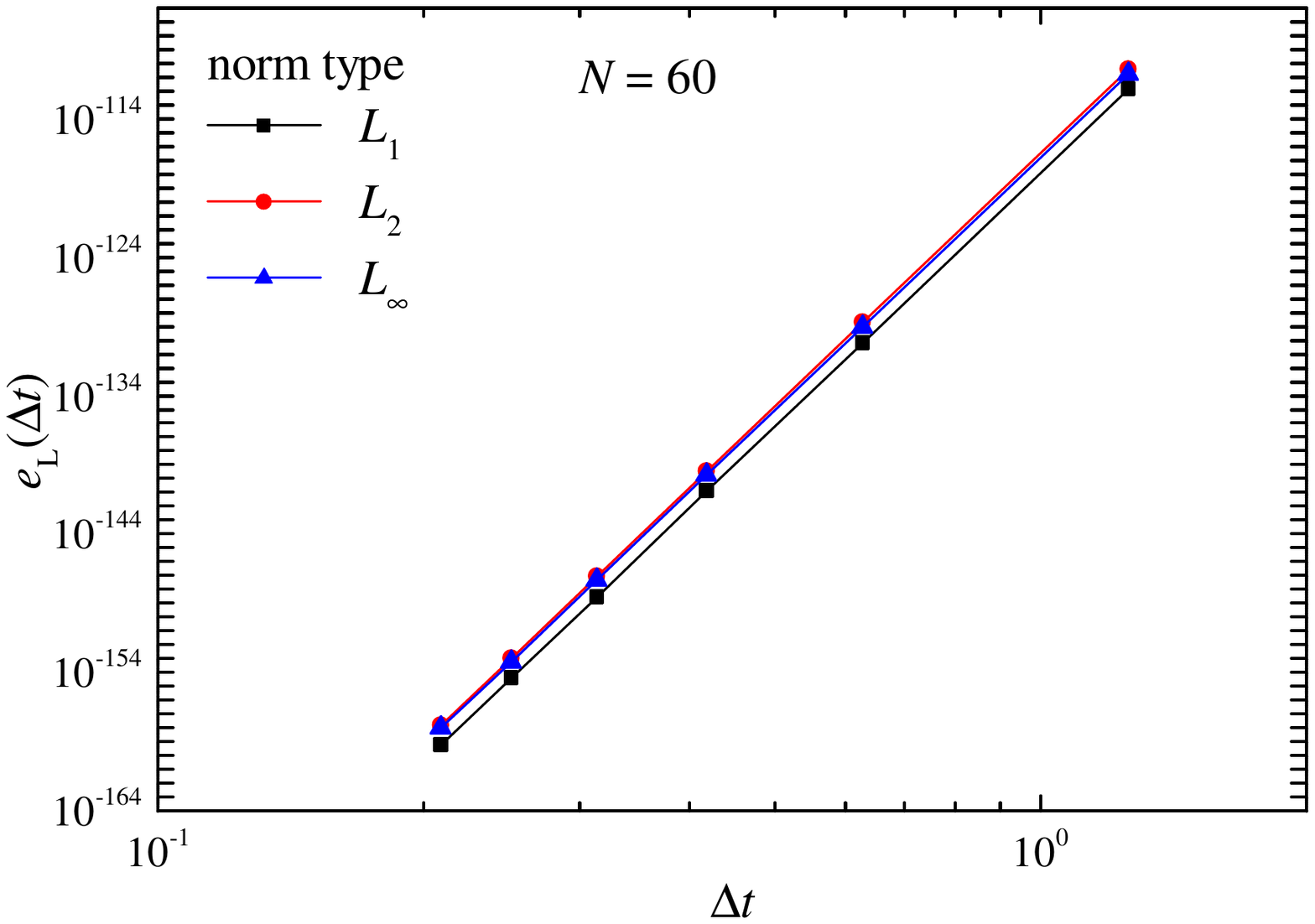}
\vspace{-8mm}\caption{\label{fig:exp_diss_sols:e3}}
\end{subfigure}\\
\caption{%
Numerical solution of the problem (\ref{eq:exp_diss}). Comparison of the solution at nodes $\mathbf{u}_{n}$, the local solution $\mathbf{u}_{L}(t)$ and the exact solution $\mathbf{u}^{\rm ex}(t)$ for components $u_{1}$ (\subref{fig:exp_diss_sols:a1}, \subref{fig:exp_diss_sols:b1}, \subref{fig:exp_diss_sols:c1}) and $u_{2}$ (\subref{fig:exp_diss_sols:a2}, \subref{fig:exp_diss_sols:b2}, \subref{fig:exp_diss_sols:c2}), the error $\varepsilon(t)$ (\subref{fig:exp_diss_sols:a3}, \subref{fig:exp_diss_sols:b3}, \subref{fig:exp_diss_sols:c3}), obtained using polynomials with degrees $N = 1$ (\subref{fig:exp_diss_sols:a1}, \subref{fig:exp_diss_sols:a2}, \subref{fig:exp_diss_sols:a3}), $N = 8$ (\subref{fig:exp_diss_sols:b1}, \subref{fig:exp_diss_sols:b2}, \subref{fig:exp_diss_sols:b3}) and $N = 60$ (\subref{fig:exp_diss_sols:c1}, \subref{fig:exp_diss_sols:c2}, \subref{fig:exp_diss_sols:c3}). Dependence of the global error for the solution at nodes $e_{G}$ (\subref{fig:exp_diss_sols:d1}, \subref{fig:exp_diss_sols:d2}, \subref{fig:exp_diss_sols:d3}) and the local solution $e_{L}$ (\subref{fig:exp_diss_sols:e1}, \subref{fig:exp_diss_sols:e2}, \subref{fig:exp_diss_sols:e3}) on the discretization step $\mathrm{\Delta}t$, obtained in the norms $L_{1}$, $L_{2}$ and $L_{\infty}$, obtained using polynomials with degrees $N = 1$ (\subref{fig:exp_diss_sols:d1}, \subref{fig:exp_diss_sols:e1}), $N = 8$ (\subref{fig:exp_diss_sols:d2}, \subref{fig:exp_diss_sols:e2}) and $N = 60$ (\subref{fig:exp_diss_sols:d3}, \subref{fig:exp_diss_sols:e3}).
}
\label{fig:exp_diss_sols}
\end{figure}

The results of the study of the application of the numerical method for this problem are presented in Fig.~\ref{fig:exp_diss_sols} and in Table~\ref{tab:exp_diss_sols}. An analysis of the error $\varepsilon(t)$ for $N > 1$ shows that the error scales for the local solution and the solution at nodes differ by several orders of magnitude, reaching a value of $100$ orders of magnitude for a degree of $N = 60$. The calculated convergence orders $p_{L_{1}}$, $p_{L_{2}}$, $p_{L_{\infty}}$ for various functional norms are presented in Table~\ref{tab:exp_diss_sols}.

\begin{table}[h!]
\centering
\caption{%
Orders of convergence $p_{L_{1}}$, $p_{L_{2}}$, $p_{L_{\infty}}$ of the ADER-DG numerical method with local DG predictor calculated in norms $L_{1}$, $L_{2}$, $L_{\infty}$ for the problem (\ref{eq:exp_diss}); $N$ is the degree of the basis polynomials $\varphi_{p}$ used for the DG representations. The orders of convergence $p$ for the solution at nodes $\mathbf{u}_{n}$ are compared to the theoretical value $p_{\rm (G)} = 2N+1$, the orders of convergence $p$ for the local solution $\mathbf{u}_{L}$ represented between the nodes are compared to the theoretical value $p_{\rm (L)} = N+1$.
}
\label{tab:exp_diss_sols}
\begin{tabular}{@{}lllllllll@{}}
\toprule
& \multicolumn{4}{c}{order $p$ for solution at nodes} & \multicolumn{4}{c}{order $p$ for local solution} \\
\midrule
$N$ & $p_{L_{1}}$ & $p_{L_{2}}$ & $p_{L_{\infty}}$ & theor. & $p_{L_{1}}$ & $p_{L_{2}}$ & $p_{L_{\infty}}$ & theor. \\
\midrule
$1$	&	$3.34$	&	$3.34$	&	$3.13$		&	$3$	&	$2.63$	&	$2.41$	&	$1.91$		&	$2$	\\
$2$	&	$5.35$	&	$5.34$	&	$5.14$		&	$5$	&	$3.07$	&	$3.00$	&	$2.81$		&	$3$	\\
$3$	&	$7.31$	&	$7.31$	&	$7.10$		&	$7$	&	$4.01$	&	$3.95$	&	$3.77$		&	$4$	\\
$4$	&	$9.29$	&	$9.29$	&	$9.08$		&	$9$	&	$5.00$	&	$4.94$	&	$4.76$		&	$5$	\\
$5$	&	$11.28$	&	$11.27$	&	$11.07$		&	$11$	&	$6.00$	&	$5.93$	&	$5.75$		&	$6$	\\
$6$	&	$13.27$	&	$13.26$	&	$13.06$		&	$13$	&	$6.99$	&	$6.92$	&	$6.74$		&	$7$	\\
$7$	&	$15.26$	&	$15.26$	&	$15.05$		&	$15$	&	$7.99$	&	$7.92$	&	$7.74$		&	$8$	\\
$8$	&	$17.25$	&	$17.25$	&	$17.04$		&	$17$	&	$8.99$	&	$8.92$	&	$8.74$		&	$9$	\\
$9$	&	$19.25$	&	$19.25$	&	$19.04$		&	$19$	&	$9.99$	&	$9.91$	&	$9.73$		&	$10$	\\
$10$	&	$21.25$	&	$21.24$	&	$21.04$		&	$21$	&	$10.98$	&	$10.91$	&	$10.73$		&	$11$	\\
\midrule
$15$	&	$31.24$	&	$31.23$	&	$31.02$		&	$31$	&	$15.98$	&	$15.91$	&	$15.72$		&	$16$	\\
$20$	&	$41.23$	&	$41.23$	&	$41.02$		&	$41$	&	$20.98$	&	$20.90$	&	$20.72$		&	$21$	\\
$25$	&	$51.23$	&	$51.22$	&	$51.02$		&	$51$	&	$25.97$	&	$25.90$	&	$25.72$		&	$26$	\\
$30$	&	$61.22$	&	$61.22$	&	$61.01$		&	$61$	&	$30.97$	&	$30.90$	&	$30.72$		&	$31$	\\
$35$	&	$71.22$	&	$71.22$	&	$71.01$		&	$71$	&	$35.97$	&	$35.90$	&	$35.72$		&	$36$	\\
$40$	&	$81.22$	&	$81.22$	&	$81.01$		&	$81$	&	$40.97$	&	$40.90$	&	$40.72$		&	$41$	\\
$45$	&	$91.22$	&	$91.22$	&	$91.01$		&	$91$	&	$45.97$	&	$45.90$	&	$45.72$		&	$46$	\\
$50$	&	$101.22$	&	$101.22$	&	$101.01$		&	$101$	&	$50.97$	&	$50.90$	&	$50.72$		&	$51$	\\
$55$	&	$111.22$	&	$111.22$	&	$111.01$		&	$111$	&	$55.97$	&	$55.90$	&	$55.72$		&	$56$	\\
$60$	&	$121.22$	&	$121.21$	&	$121.01$		&	$121$	&	$60.97$	&	$60.90$	&	$60.72$		&	$61$	\\
\bottomrule
\end{tabular}
\end{table}

All obtained convergence orders $p$ sufficiently correspond to the expected theoretical values $p_{\rm (G)}$ and $p_{\rm (L)}$, in accordance with (\ref{eq:expect_orders}). All convergence orders $p$ for the solution at nodes $\mathbf{u}_{n}$ demonstrate superconvergence in all considered functional norms (\ref{eq:global_eps_def}). For the values of the degrees of polynomials $N = 3\div5$, there is a significant excess of the expected values $2N+1$. The empirical values of the orders exceed the expected values $p_{\rm (G)}$ by no more than $1$. Orders $p_{L_{1}}$ and $p_{L_{2}}$ show values $0.2$-$0.3$ higher than expected values. The orders of $p_{L_{\infty}}$ are slightly higher than the expected values. In the case of a local solution $\mathbf{u}_{L}$, the convergence orders $p$ correspond to the expected values $p_{\rm (L)}$. 

As a result, it becomes clear that the solution at nodes $\mathbf{u}_{n}$ obtained by the ADER-DG numerical method with a local DG predictor demonstrates the classical superconvergence $2N+1$ expected from such numerical methods, at least in problem (\ref{eq:exp_diss}). The local solution $\mathbf{u}_{L}$ obtained by the DG predictor demonstrates the classical convergence characteristic of DG methods.

\subsubsection{Example 3}

\begin{figure}[h!]
\captionsetup[subfigure]{%
position=bottom,
textfont=normalfont,
singlelinecheck=off,
justification=raggedright,
font=tiny
}
\centering
\begin{subfigure}{0.275\textwidth}
\includegraphics[width=\textwidth]{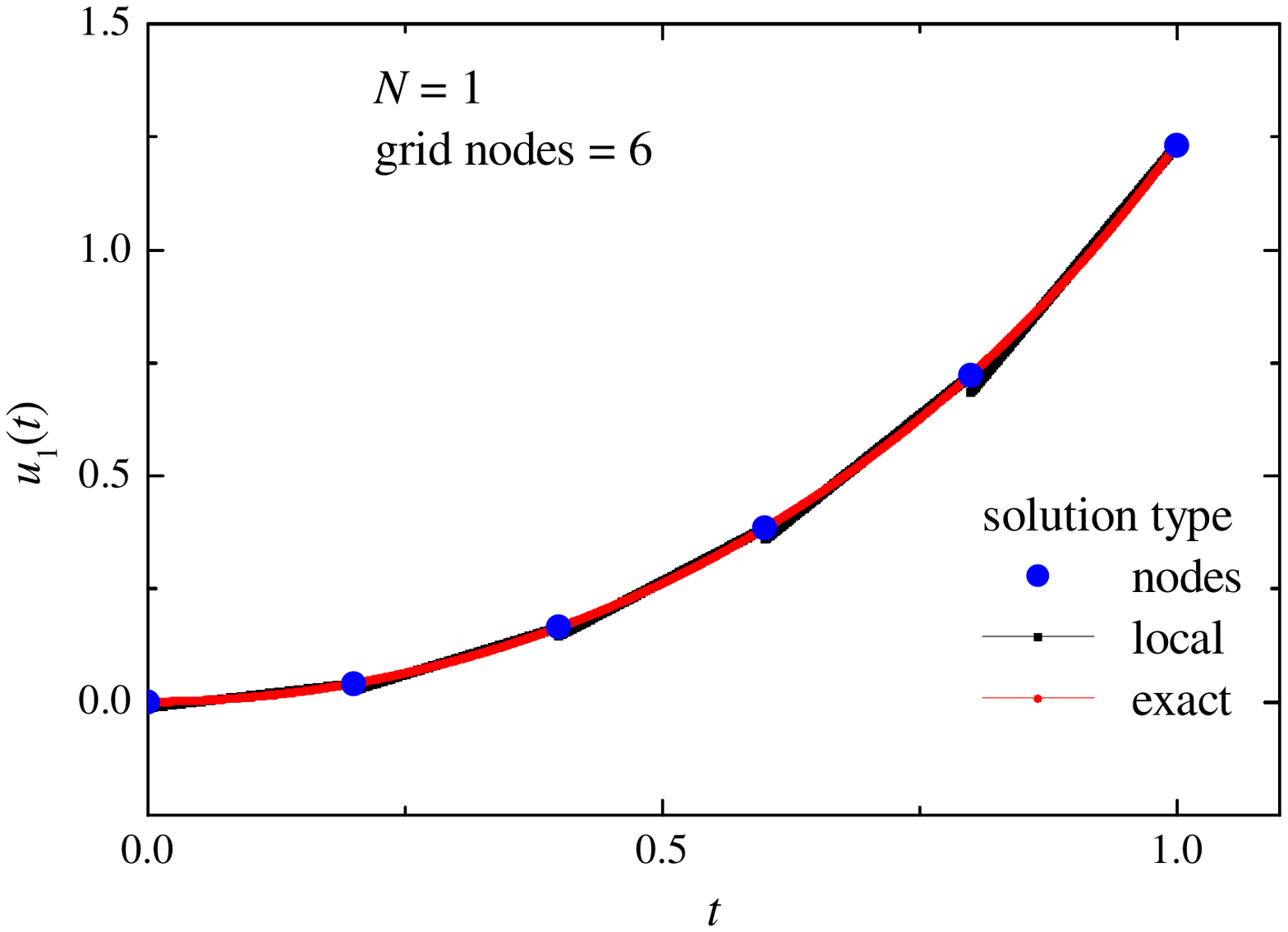}
\vspace{-8mm}\caption{\label{fig:bratu_problem_sols:a1}}
\end{subfigure}
\begin{subfigure}{0.275\textwidth}
\includegraphics[width=\textwidth]{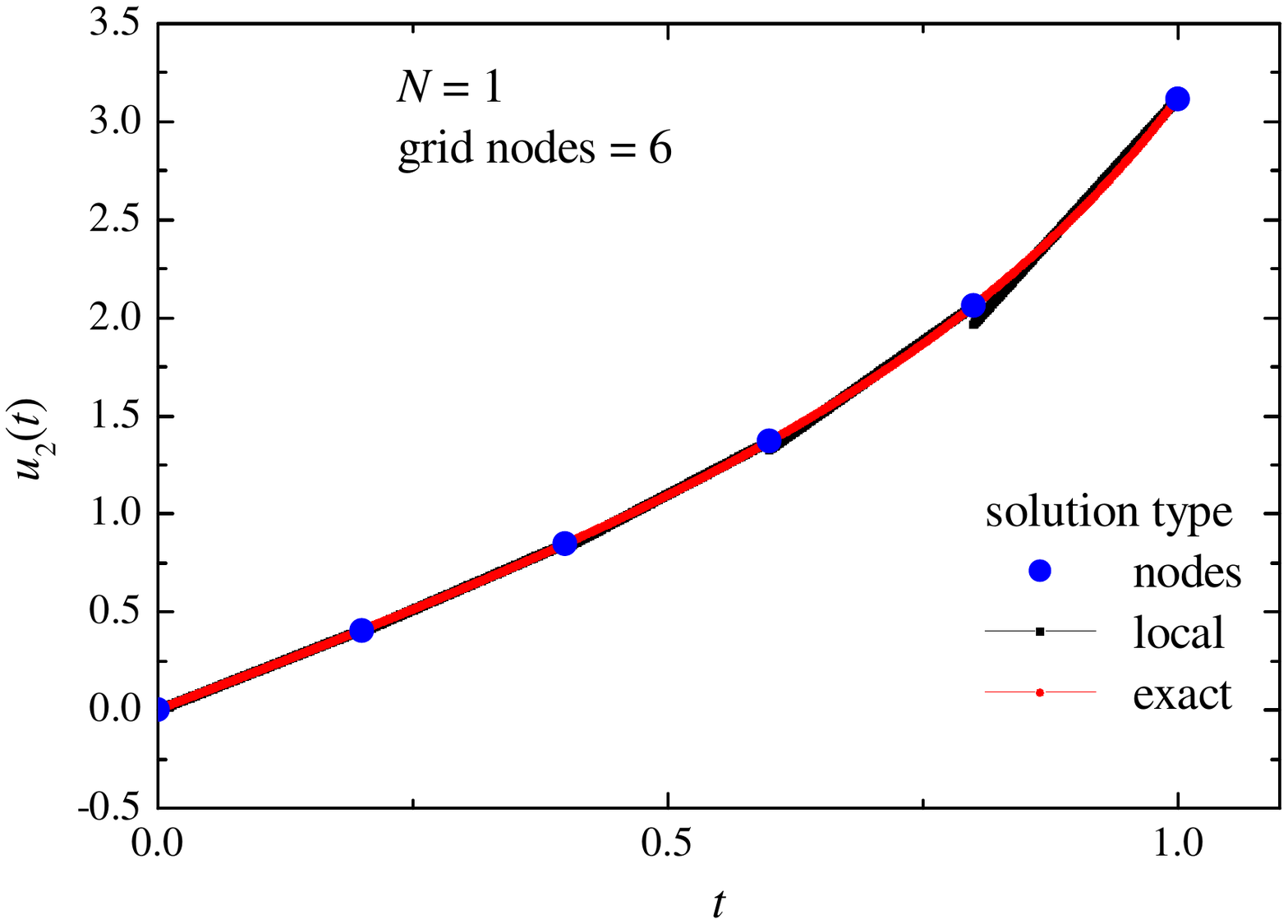}
\vspace{-8mm}\caption{\label{fig:bratu_problem_sols:a2}}
\end{subfigure}
\begin{subfigure}{0.275\textwidth}
\includegraphics[width=\textwidth]{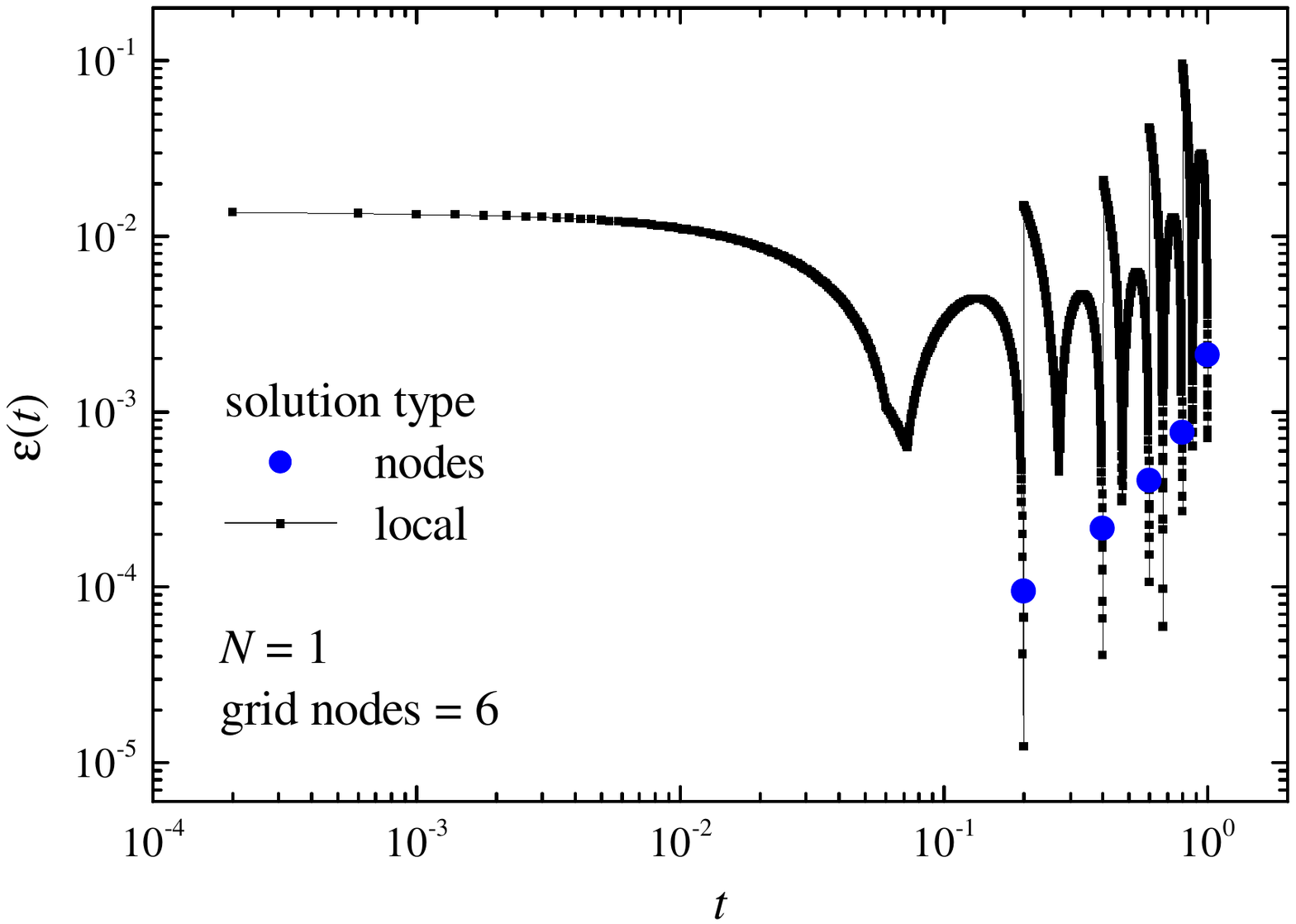}
\vspace{-8mm}\caption{\label{fig:bratu_problem_sols:a3}}
\end{subfigure}\\
\begin{subfigure}{0.275\textwidth}
\includegraphics[width=\textwidth]{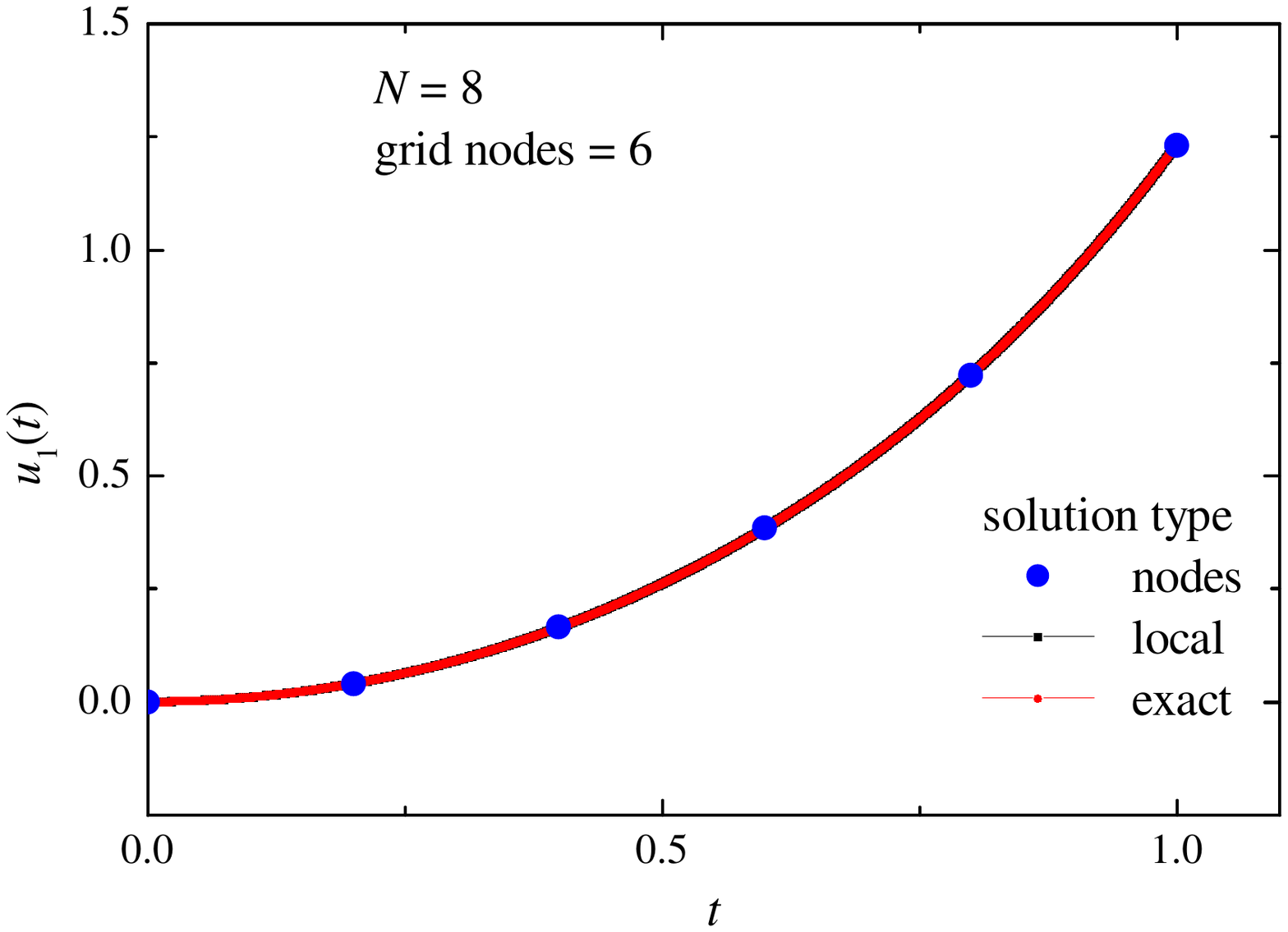}
\vspace{-8mm}\caption{\label{fig:bratu_problem_sols:b1}}
\end{subfigure}
\begin{subfigure}{0.275\textwidth}
\includegraphics[width=\textwidth]{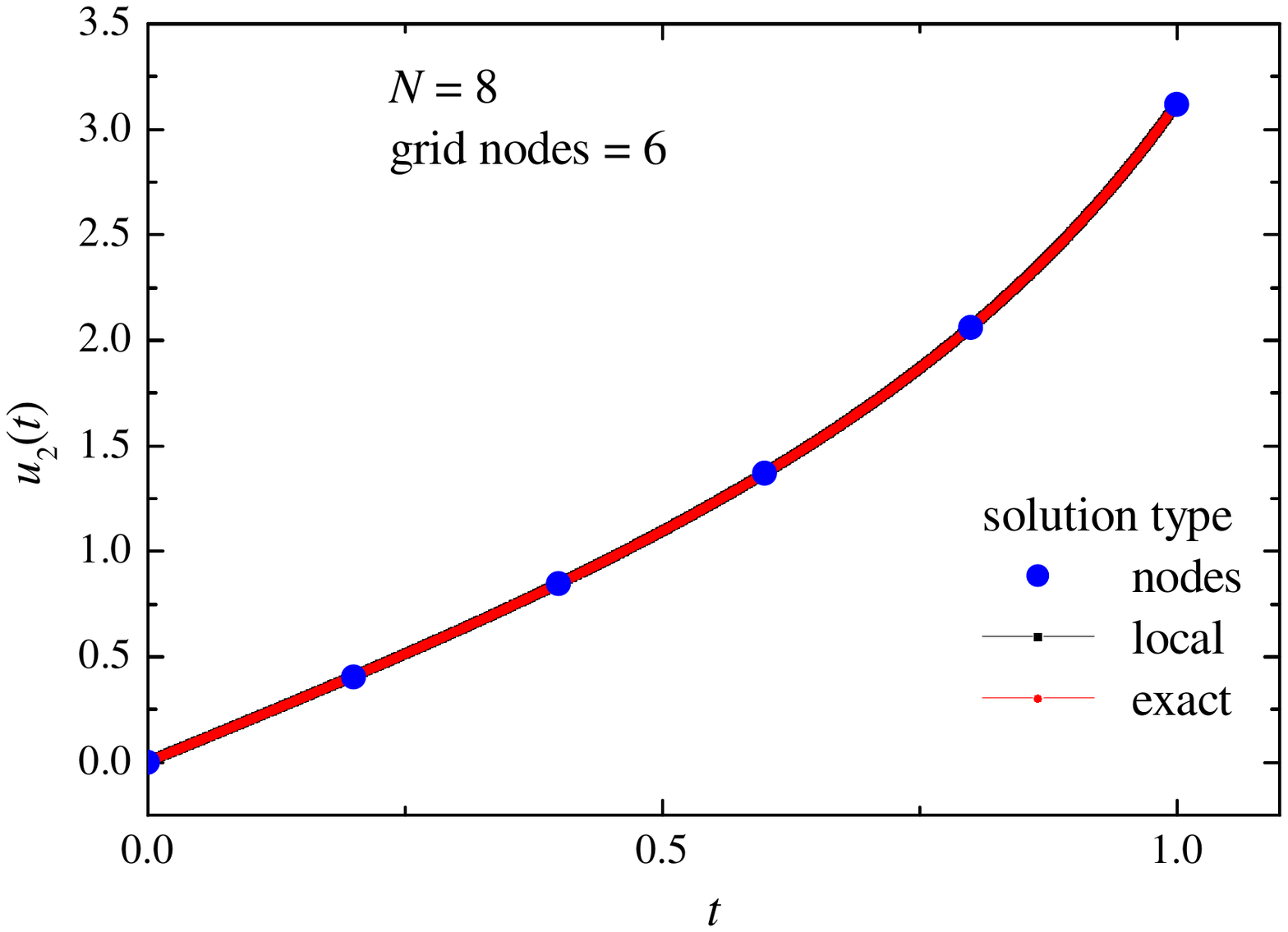}
\vspace{-8mm}\caption{\label{fig:bratu_problem_sols:b2}}
\end{subfigure}
\begin{subfigure}{0.275\textwidth}
\includegraphics[width=\textwidth]{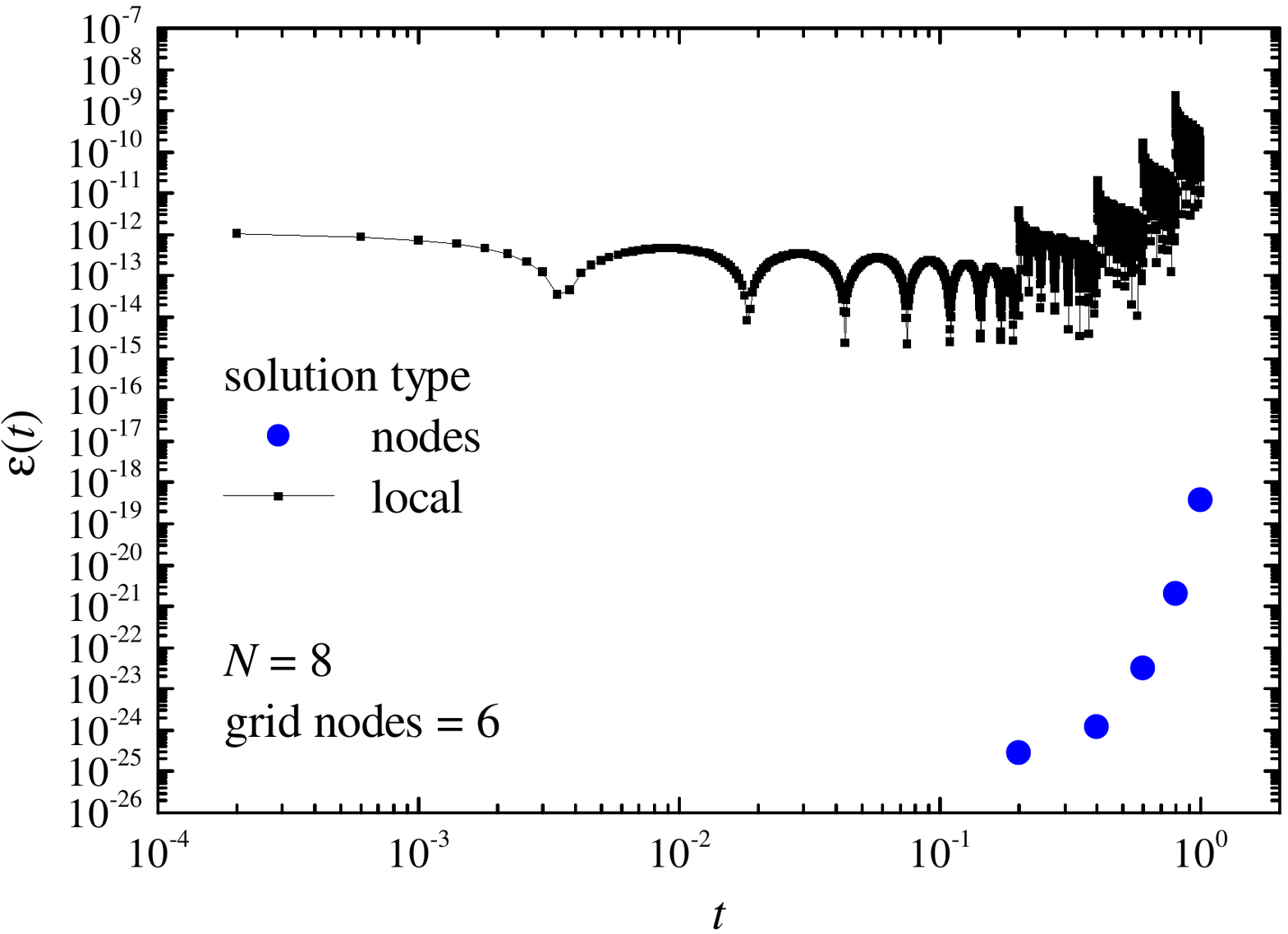}
\vspace{-8mm}\caption{\label{fig:bratu_problem_sols:b3}}
\end{subfigure}\\
\begin{subfigure}{0.275\textwidth}
\includegraphics[width=\textwidth]{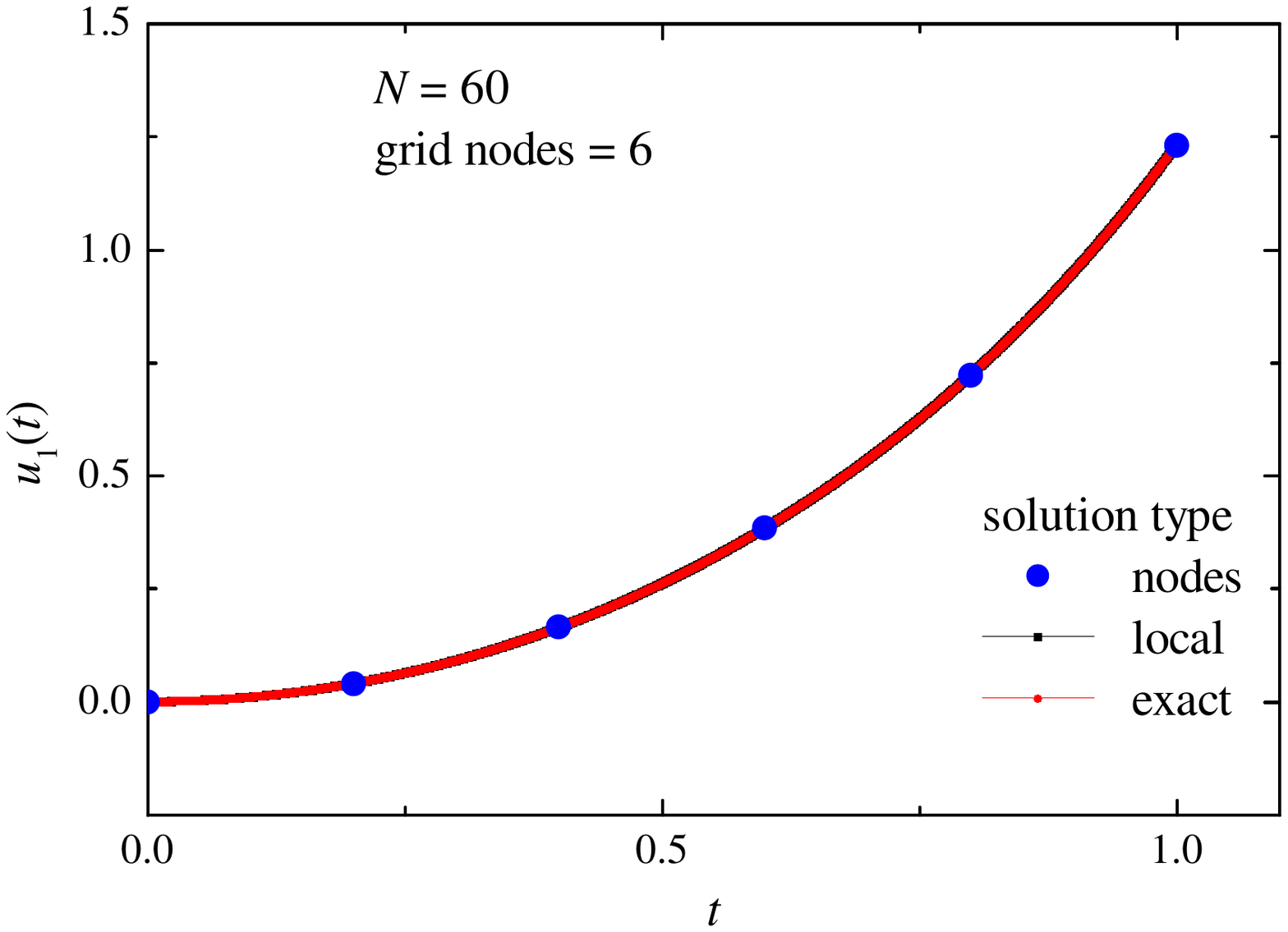}
\vspace{-8mm}\caption{\label{fig:bratu_problem_sols:c1}}
\end{subfigure}
\begin{subfigure}{0.275\textwidth}
\includegraphics[width=\textwidth]{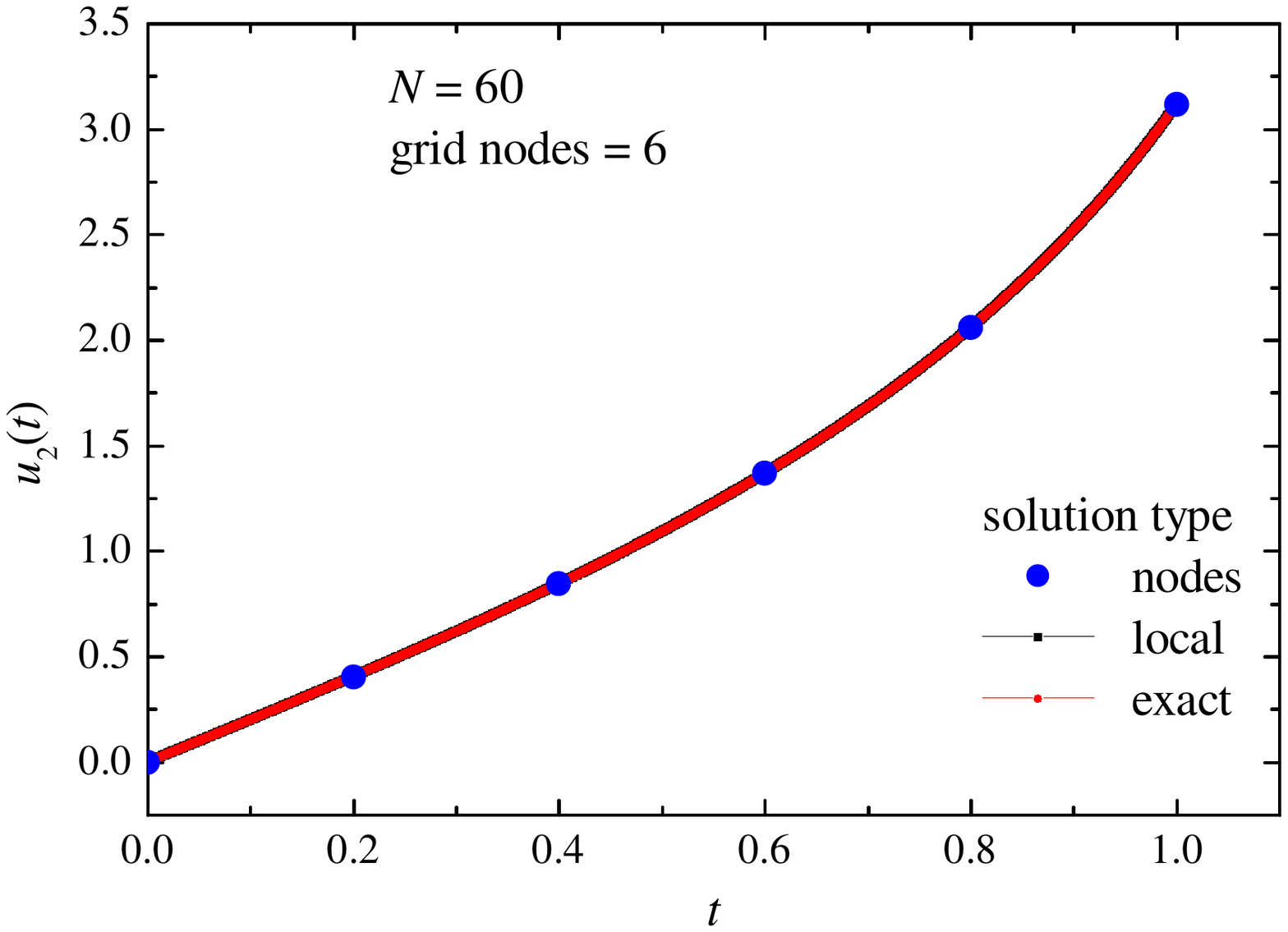}
\vspace{-8mm}\caption{\label{fig:bratu_problem_sols:c2}}
\end{subfigure}
\begin{subfigure}{0.275\textwidth}
\includegraphics[width=\textwidth]{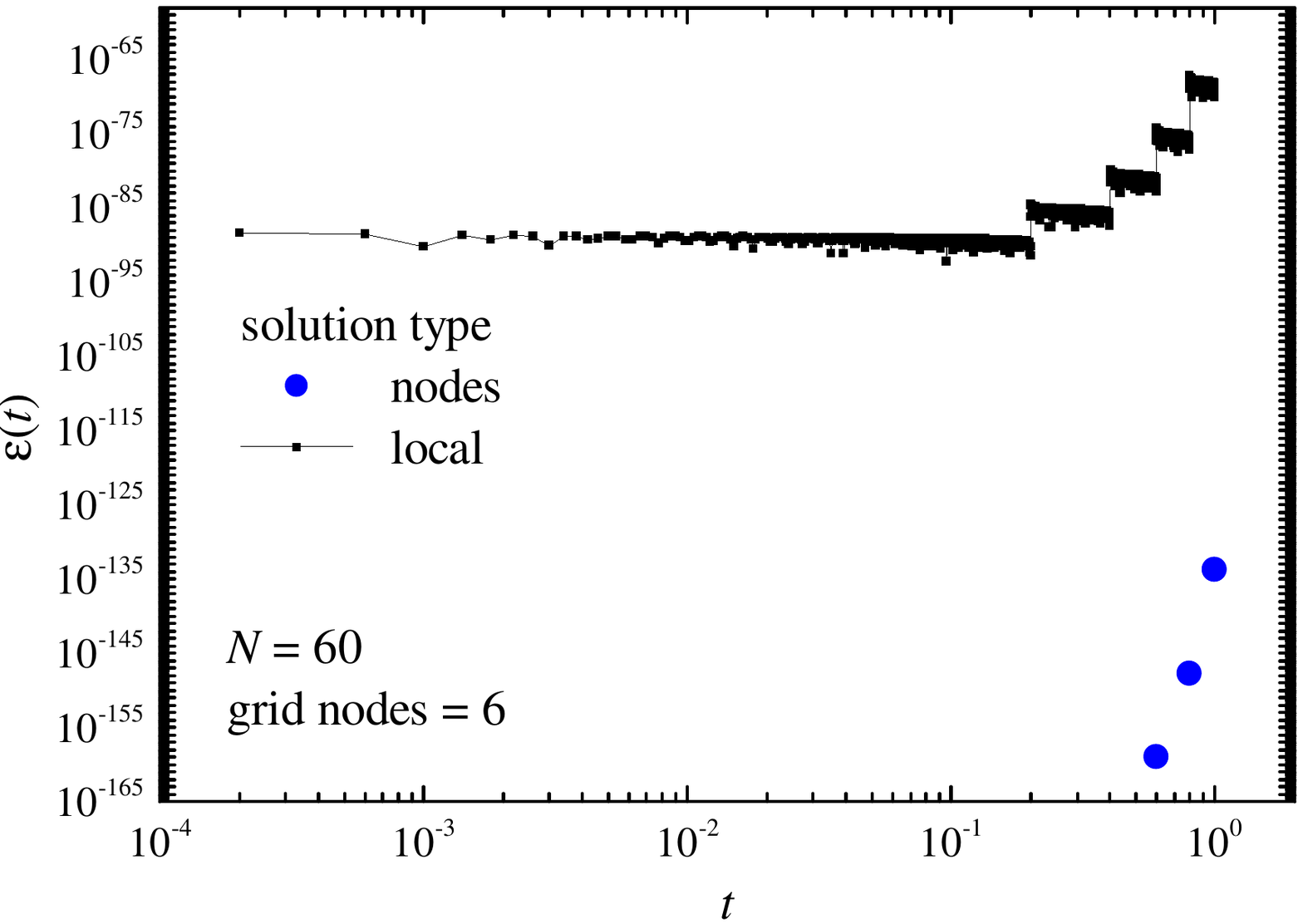}
\vspace{-8mm}\caption{\label{fig:bratu_problem_sols:c3}}
\end{subfigure}\\
\begin{subfigure}{0.275\textwidth}
\includegraphics[width=\textwidth]{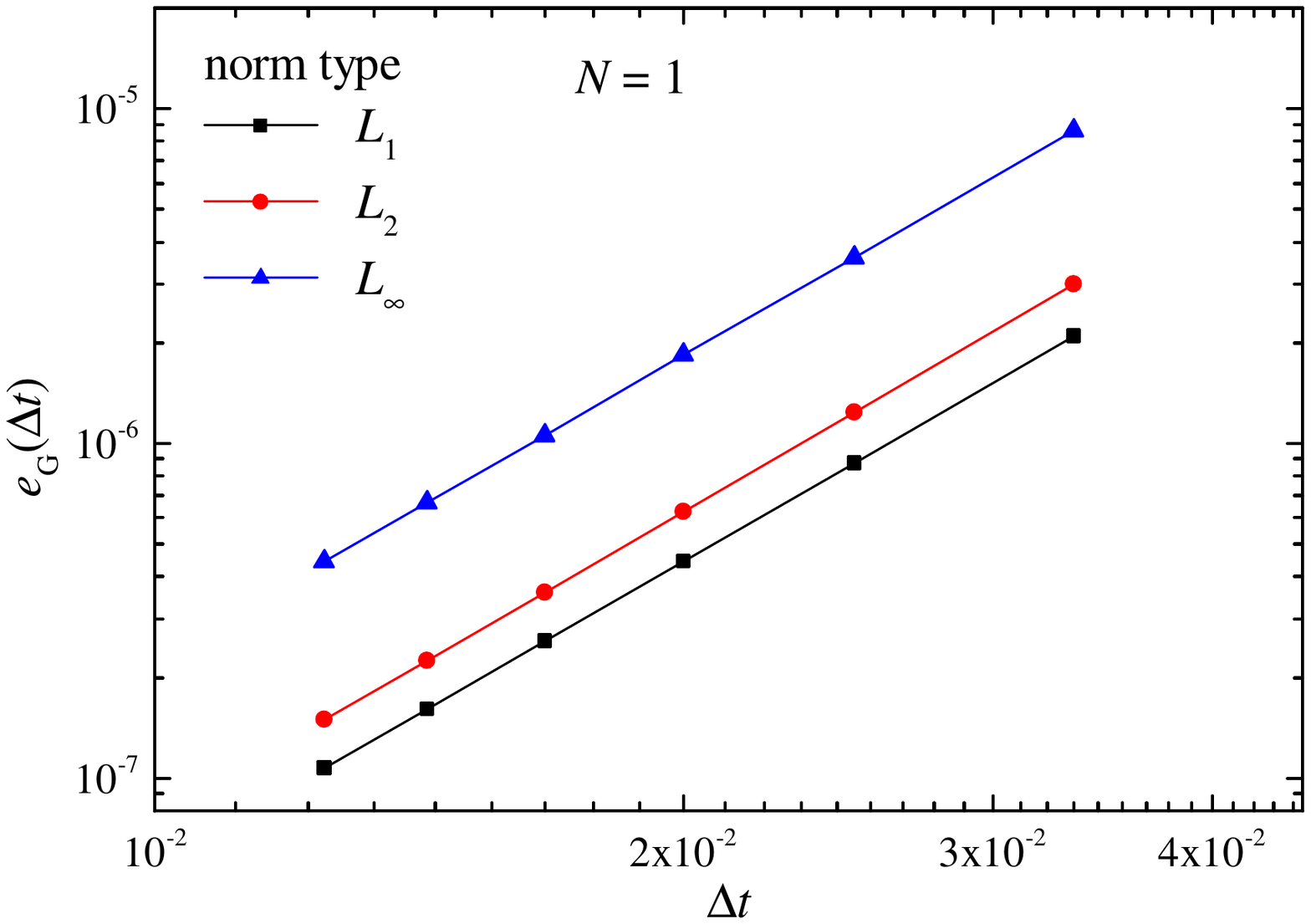}
\vspace{-8mm}\caption{\label{fig:bratu_problem_sols:d1}}
\end{subfigure}
\begin{subfigure}{0.275\textwidth}
\includegraphics[width=\textwidth]{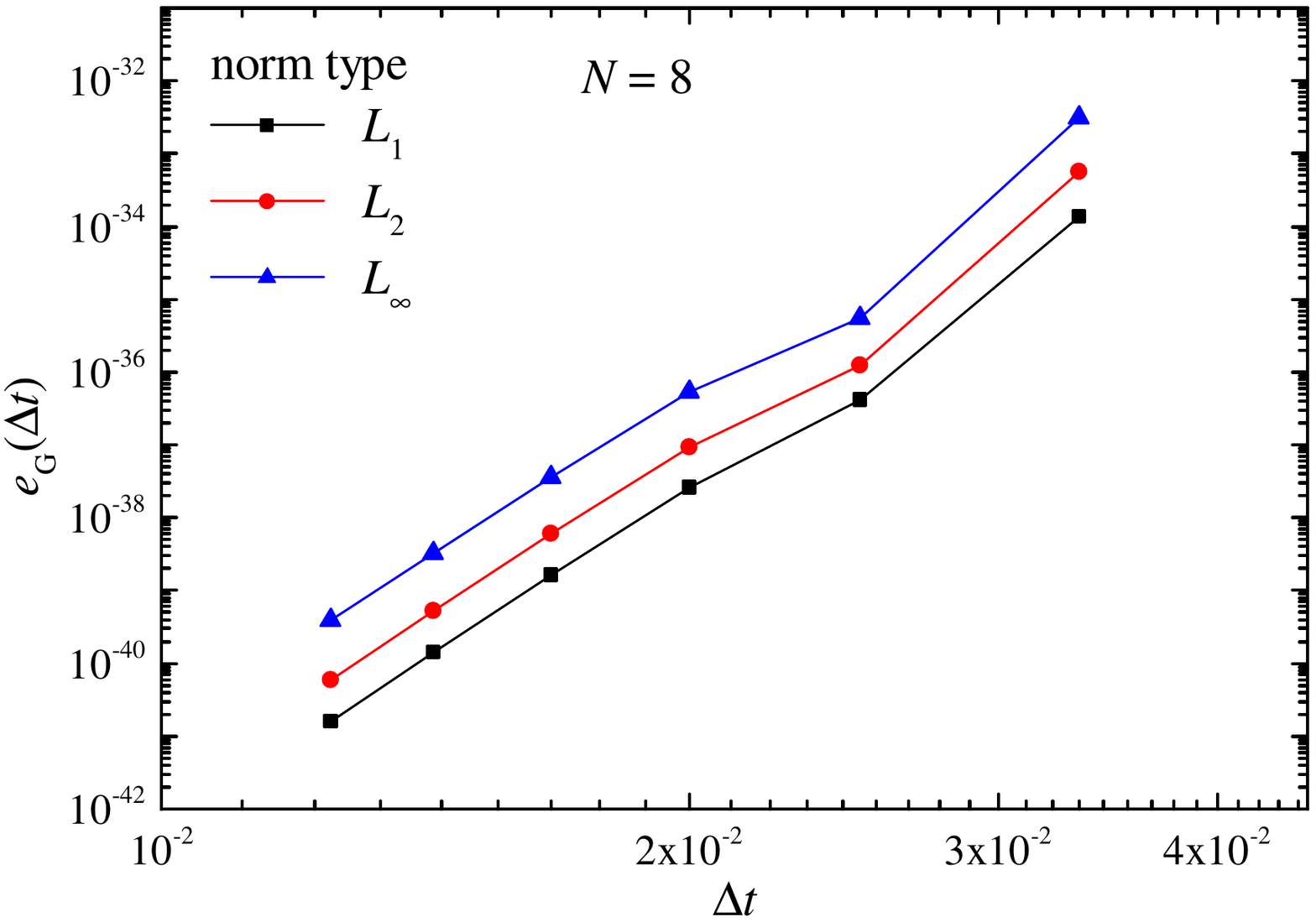}
\vspace{-8mm}\caption{\label{fig:bratu_problem_sols:d2}}
\end{subfigure}
\begin{subfigure}{0.275\textwidth}
\includegraphics[width=\textwidth]{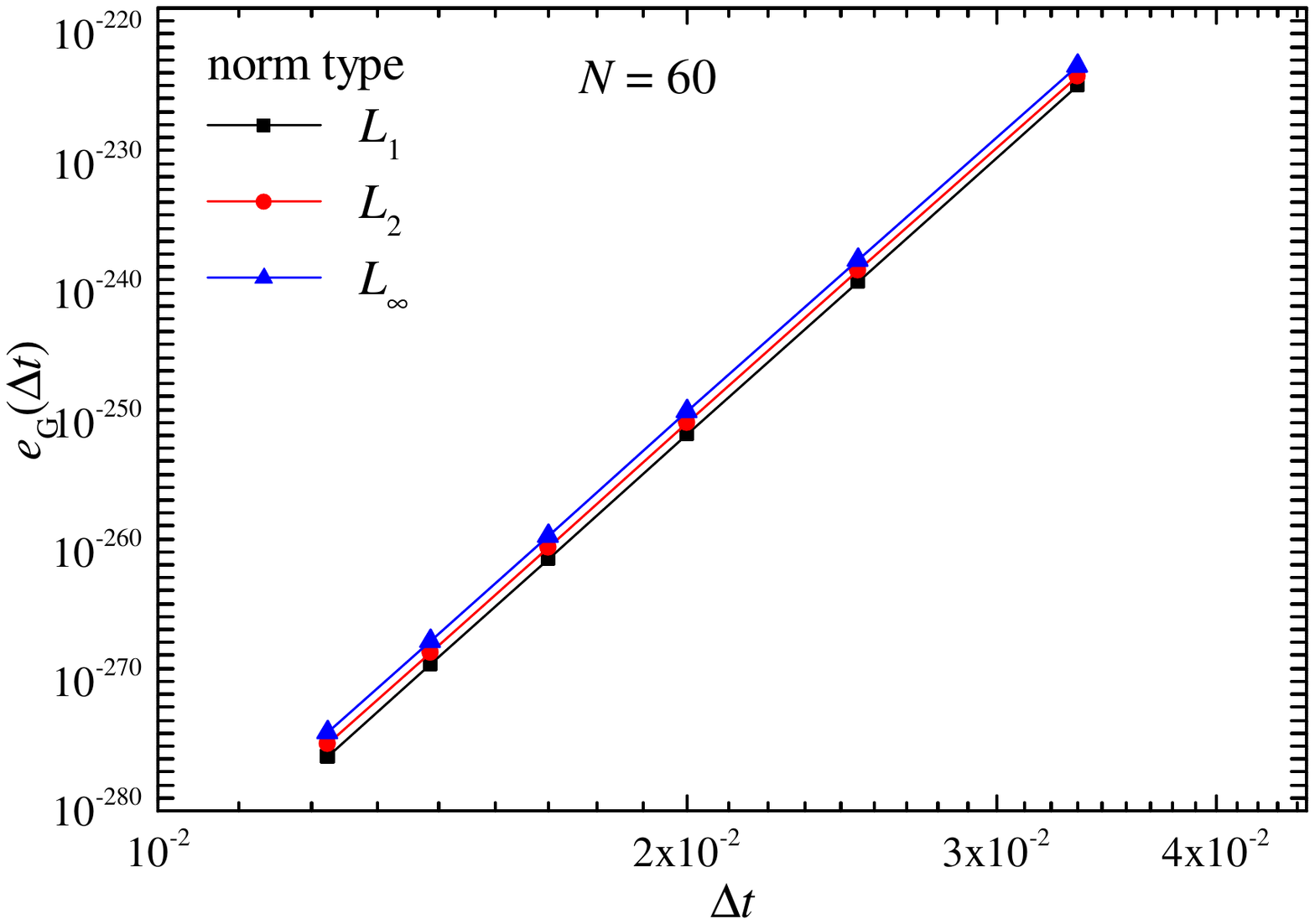}
\vspace{-8mm}\caption{\label{fig:bratu_problem_sols:d3}}
\end{subfigure}\\
\begin{subfigure}{0.275\textwidth}
\includegraphics[width=\textwidth]{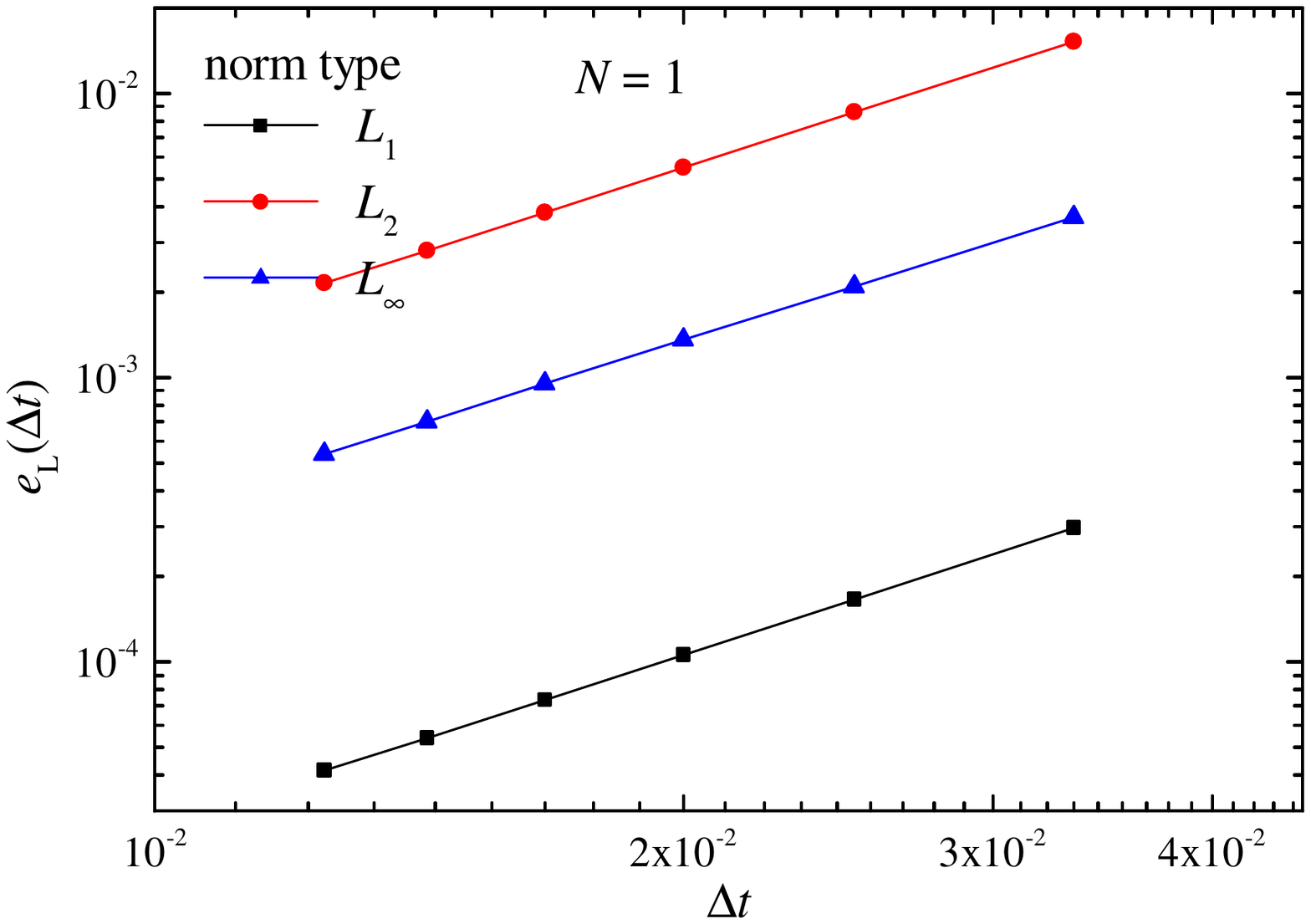}
\vspace{-8mm}\caption{\label{fig:bratu_problem_sols:e1}}
\end{subfigure}
\begin{subfigure}{0.275\textwidth}
\includegraphics[width=\textwidth]{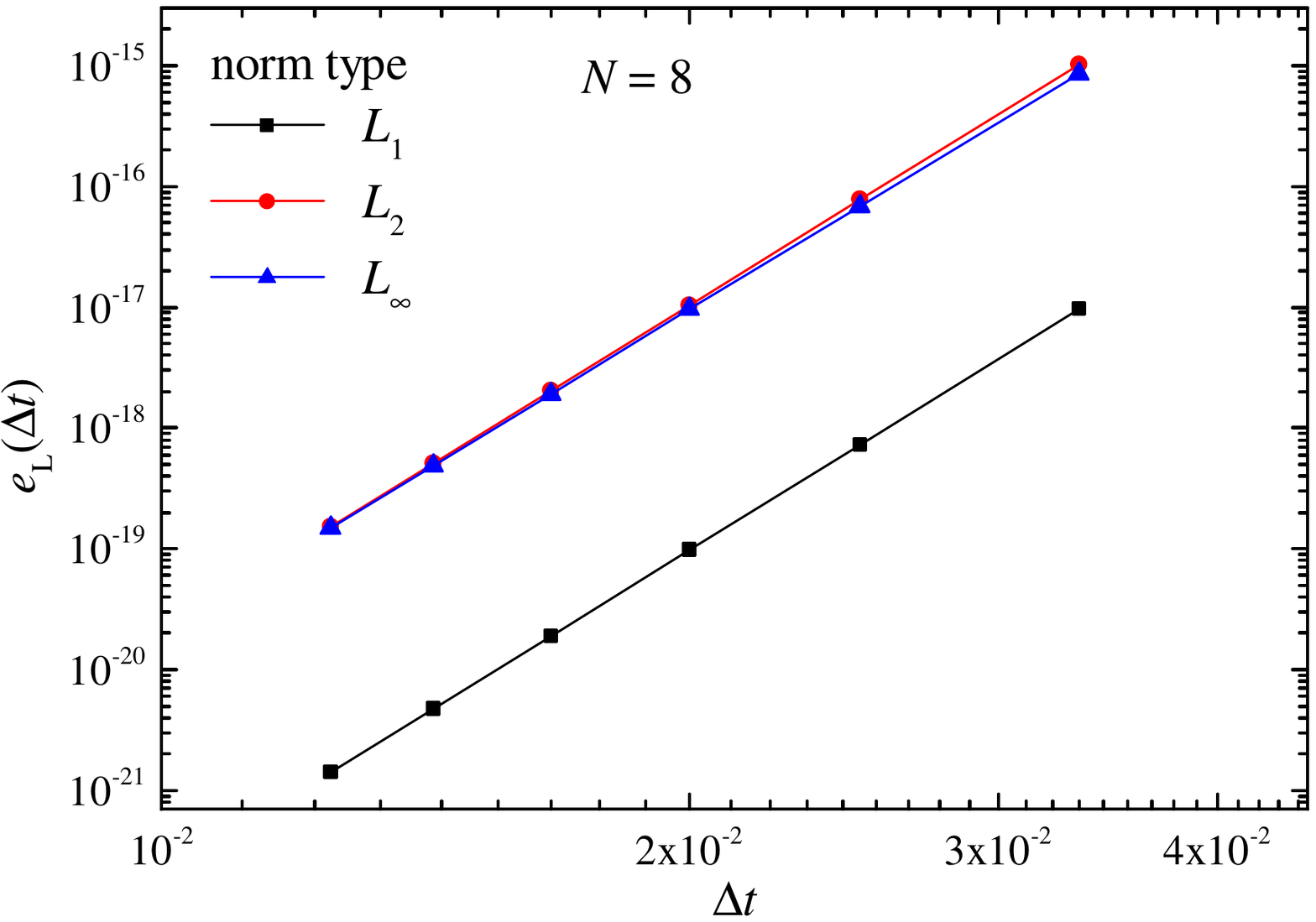}
\vspace{-8mm}\caption{\label{fig:bratu_problem_sols:e2}}
\end{subfigure}
\begin{subfigure}{0.275\textwidth}
\includegraphics[width=\textwidth]{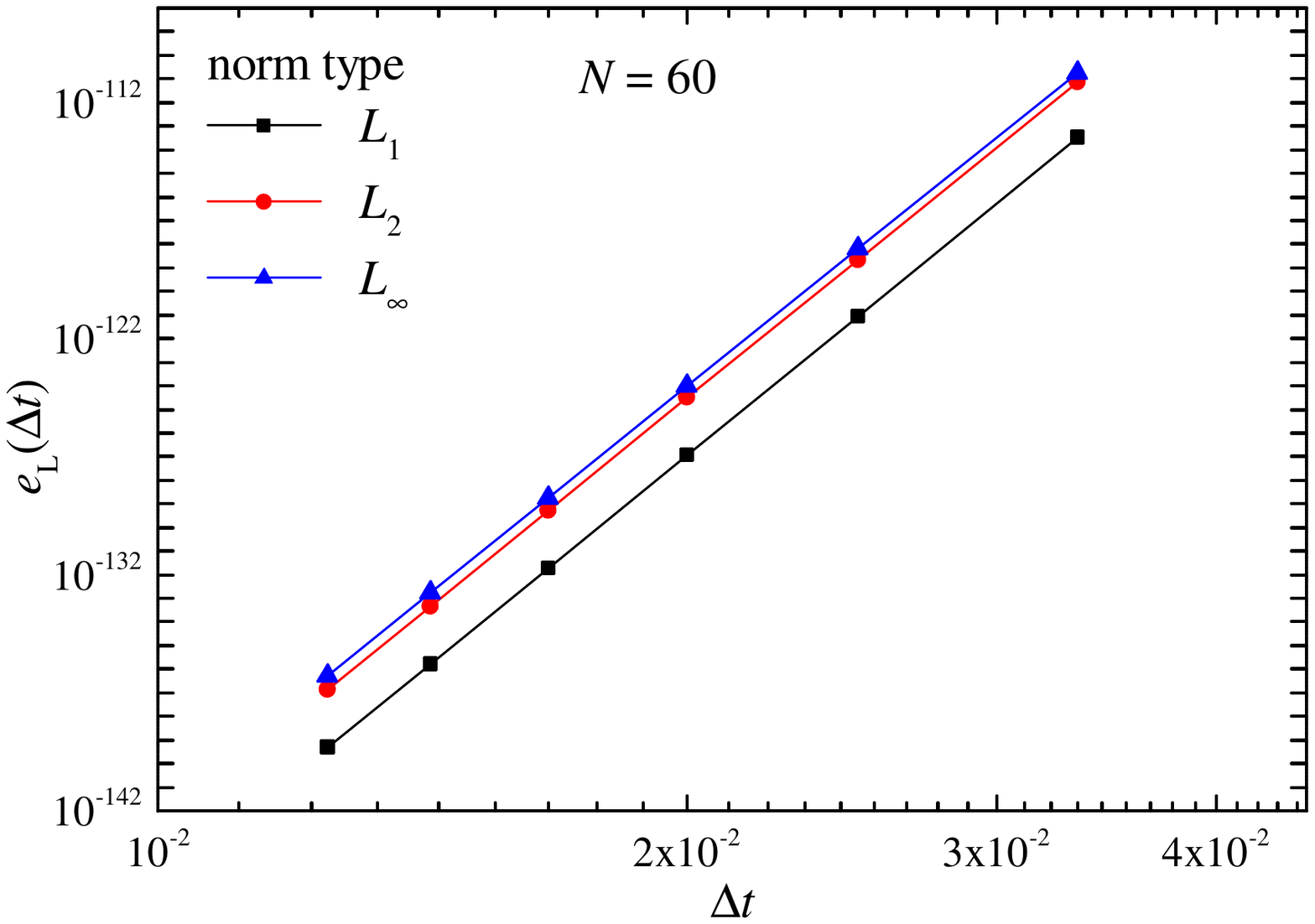}
\vspace{-8mm}\caption{\label{fig:bratu_problem_sols:e3}}
\end{subfigure}\\
\caption{%
Numerical solution of the problem (\ref{eq:bratu_problem}). Comparison of the solution at nodes $\mathbf{u}_{n}$, the local solution $\mathbf{u}_{L}(t)$ and the exact solution $\mathbf{u}^{\rm ex}(t)$ for components $u_{1}$ (\subref{fig:bratu_problem_sols:a1}, \subref{fig:bratu_problem_sols:b1}, \subref{fig:bratu_problem_sols:c1}) and $u_{2}$ (\subref{fig:bratu_problem_sols:a2}, \subref{fig:bratu_problem_sols:b2}, \subref{fig:bratu_problem_sols:c2}), the error $\varepsilon(t)$ (\subref{fig:bratu_problem_sols:a3}, \subref{fig:bratu_problem_sols:b3}, \subref{fig:bratu_problem_sols:c3}), obtained using polynomials with degrees $N = 1$ (\subref{fig:bratu_problem_sols:a1}, \subref{fig:bratu_problem_sols:a2}, \subref{fig:bratu_problem_sols:a3}), $N = 8$ (\subref{fig:bratu_problem_sols:b1}, \subref{fig:bratu_problem_sols:b2}, \subref{fig:bratu_problem_sols:b3}) and $N = 60$ (\subref{fig:bratu_problem_sols:c1}, \subref{fig:bratu_problem_sols:c2}, \subref{fig:bratu_problem_sols:c3}). Dependence of the global error for the solution at nodes $e_{G}$ (\subref{fig:bratu_problem_sols:d1}, \subref{fig:bratu_problem_sols:d2}, \subref{fig:bratu_problem_sols:d3}) and the local solution $e_{L}$ (\subref{fig:bratu_problem_sols:e1}, \subref{fig:bratu_problem_sols:e2}, \subref{fig:bratu_problem_sols:e3}) on the discretization step $\mathrm{\Delta}t$, obtained in the norms $L_{1}$, $L_{2}$ and $L_{\infty}$, obtained using polynomials with degrees $N = 1$ (\subref{fig:bratu_problem_sols:d1}, \subref{fig:bratu_problem_sols:e1}), $N = 8$ (\subref{fig:bratu_problem_sols:d2}, \subref{fig:bratu_problem_sols:e2}) and $N = 60$ (\subref{fig:bratu_problem_sols:d3}, \subref{fig:bratu_problem_sols:e3}).
}
\label{fig:bratu_problem_sols}
\end{figure} 
\begin{table}[h!]
\centering
\caption{%
Orders of convergence $p_{L_{1}}$, $p_{L_{2}}$, $p_{L_{\infty}}$ of the ADER-DG numerical method with local DG predictor calculated in norms $L_{1}$, $L_{2}$, $L_{\infty}$ for the problem (\ref{eq:bratu_problem}); $N$ is the degree of the basis polynomials $\varphi_{p}$ used for the DG representations. The orders of convergence $p$ for the solution at nodes $\mathbf{u}_{n}$ are compared to the theoretical value $p_{\rm (G)} = 2N+1$, the orders of convergence $p$ for the local solution $\mathbf{u}_{L}$ represented between the nodes are compared to the theoretical value $p_{\rm (L)} = N+1$.
}
\label{tab:bratu_problem_sols}
\begin{tabular}{@{}lllllllll@{}}
\toprule
& \multicolumn{4}{c}{order $p$ for solution at nodes} & \multicolumn{4}{c}{order $p$ for local solution} \\
\midrule
$N$ & $p_{L_{1}}$ & $p_{L_{2}}$ & $p_{L_{\infty}}$ & theor. & $p_{L_{1}}$ & $p_{L_{2}}$ & $p_{L_{\infty}}$ & theor. \\
\midrule
$1$	&	$3.03$	&	$3.05$	&	$3.01$		&	$3$	&	$2.01$	&	$2.00$	&	$1.95$		&	$2$	\\
$2$	&	$5.03$	&	$5.04$	&	$4.97$		&	$5$	&	$3.01$	&	$2.99$	&	$2.93$		&	$3$	\\
$3$	&	$6.97$	&	$6.94$	&	$6.77$		&	$7$	&	$4.01$	&	$3.99$	&	$3.91$		&	$4$	\\
$4$	&	$8.90$	&	$8.89$	&	$8.74$		&	$9$	&	$5.00$	&	$4.99$	&	$4.89$		&	$5$	\\
$5$	&	$10.68$	&	$10.65$	&	$10.46$		&	$11$	&	$6.00$	&	$5.99$	&	$5.88$		&	$6$	\\
$6$	&	$12.33$	&	$12.26$	&	$12.03$		&	$13$	&	$7.00$	&	$6.98$	&	$6.86$		&	$7$	\\
$7$	&	$13.06$	&	$12.76$	&	$12.24$		&	$15$	&	$8.00$	&	$7.98$	&	$7.84$		&	$8$	\\
$8$	&	$15.81$	&	$15.83$	&	$15.56$		&	$17$	&	$9.00$	&	$8.97$	&	$8.82$		&	$9$	\\
$9$	&	$19.08$	&	$19.21$	&	$19.01$		&	$19$	&	$10.00$	&	$9.97$	&	$9.80$		&	$10$	\\
$10$	&	$22.70$	&	$22.87$	&	$22.74$		&	$21$	&	$10.99$	&	$10.96$	&	$10.79$		&	$11$	\\
\midrule
$15$	&	$33.07$	&	$32.95$	&	$32.58$		&	$31$	&	$15.98$	&	$15.93$	&	$15.70$		&	$16$	\\
$20$	&	$42.65$	&	$42.49$	&	$42.08$		&	$41$	&	$20.96$	&	$20.89$	&	$20.61$		&	$21$	\\
$25$	&	$52.50$	&	$52.30$	&	$51.86$		&	$51$	&	$25.93$	&	$25.84$	&	$25.52$		&	$26$	\\
$30$	&	$62.39$	&	$62.14$	&	$61.68$		&	$61$	&	$30.90$	&	$30.78$	&	$30.43$		&	$31$	\\
$35$	&	$72.28$	&	$71.99$	&	$71.52$		&	$71$	&	$35.86$	&	$35.72$	&	$35.34$		&	$36$	\\
$40$	&	$82.17$	&	$81.84$	&	$81.37$		&	$81$	&	$40.82$	&	$40.65$	&	$40.25$		&	$41$	\\
$45$	&	$92.05$	&	$91.69$	&	$91.21$		&	$91$	&	$45.77$	&	$45.58$	&	$45.16$		&	$46$	\\
$50$	&	$101.91$	&	$101.53$	&	$101.05$		&	$101$	&	$50.72$	&	$50.51$	&	$50.07$		&	$51$	\\
$55$	&	$111.77$	&	$111.37$	&	$110.88$		&	$111$	&	$55.67$	&	$55.43$	&	$54.98$		&	$56$	\\
$60$	&	$121.62$	&	$121.20$	&	$120.72$		&	$121$	&	$60.61$	&	$60.35$	&	$59.89$		&	$61$	\\
\bottomrule
\end{tabular}
\end{table} 

For the third example the second-order non-linear Bratu IVP problem~\cite{dg_ivp_ode_3} was chosen:
\begin{equation}
\ddot{x} - 2\exp(x) = 0;\qquad
x(0) = 0;\ \dot{x}(0) = 0;\qquad
0 \leqslant t \leqslant 1;
\end{equation}
with exact analytical solution of this problem having the form $x(t) = -2\ln(\cos(t))$. The IVP for a first-order non-linear ODE system corresponding to this problem has the following form:
\begin{equation}\label{eq:bratu_problem}
\frac{du_{1}}{dt} = u_{2};\ \frac{du_{2}}{dt} = 2\exp(u_{1});\qquad
u_{1}(0) = 0;\ u_{2}(0) = 0;
\end{equation}
where $\mathbf{u} = [u_{1},\, u_{2}]^{T} = [x,\, \dot{x}]^{T}$ is a desired two component vector function. The full two component exact analytical solution of this problem has the form $\mathbf{u}^{\rm ex} = [-2\ln(\cos(t)),\, 2\sin(t)/\cos(t)]^{T}$. The solution definition domain $\Omega$ was chosen as $0 \leqslant t \leqslant 1$.

Unlike the previous considered problems (\ref{eq:harm_osc}) and (\ref{eq:exp_diss}) for linear second-order differential equations, this problem (\ref{eq:bratu_problem}) is related to the numerical solution of a non-linear equation. The study of the convergence of the numerical solution was carried out on the basis of the analysis of the solution for the set of $6$ different uniform discretizations of the domain of definition $\Omega$ of the desired function, with the number of grid nodes $L = 31$, $41$, $51$, $61$, $71$, $81$; note that the zero node, in which the initial condition was defined, also refers to grid nodes, so the discretization step was determined by the expression $\Delta t = 2\pi/(L-1)$. To demonstrate the possibilities of the numerical method for solving the problem, calculations were also performed for a coarse grid, with a number of steps $L = 6$.

The results of the study of the application of the numerical method for this problem are presented in Fig.~\ref{fig:bratu_problem_sols} and in Table~\ref{tab:bratu_problem_sols}. The results obtained, in terms of the accuracy of the numerical solution when compared with the exact analytical solution, demonstrate all the same properties that were revealed in the previous considered problems (\ref{eq:harm_osc}) and (\ref{eq:exp_diss}) for linear second-order differential equations. The error $\varepsilon(t)$ for $N > 1$ shows that the error scales for the local solution and the solution at nodes differ by several orders of magnitude, reaching a value of $50$-$80$ orders of magnitude for a degree of $N = 60$. The calculated convergence orders $p_{L_{1}}$, $p_{L_{2}}$, $p_{L_{\infty}}$ for various functional norms are presented in Table~\ref{tab:bratu_problem_sols}. The convergence orders $p$ for the solution at nodes are compared to the theoretical value $p_{\rm (G)} = 2N+1$, the convergence orders $p$ for the local solution represented between the nodes are compared to the theoretical value $p_{\rm (L)} = N+1$, in accordance with (\ref{eq:expect_orders}).

The main part of the obtained convergence orders $p$ sufficiently correspond to the expected theoretical values $p_{\rm (G)}$ and $p_{\rm (L)}$. However, for the degrees of polynomials $N = 7$-$8$, there are significantly lower values of the convergence orders $p_{\rm (G)}$ for the solution at nodes $\mathbf{u}_{n}$, which is associated with the asymptotic feature in the chosen values of the number of discretization nodes $L$. In the region of ``small'' values of the degrees of polynomials $N < 6$, there is a correspondence with the expected values $p_{\rm (G)}$ of the orders. In the region of large values of the degrees of polynomials $N \geqslant 10$, convergence orders $p$ are observed that exceed the expected values $p_{\rm (G)}$, which was an unexpected result. In the range of polynomial degrees $N = 10, \ldots, 45$, the empirical values of the convergence orders $p$ for the solution at nodes exceed the expected values $p_{\rm (G)}$ by more than $1$. In the case of a local solution $\mathbf{u}_{L}$, the convergence orders $p$ correspond to the expected values $p_{\rm (L)}$.

As a result, it becomes clear that the solution at nodes $\mathbf{u}_{n}$ obtained by the ADER-DG numerical method with a local DG predictor demonstrates the classical superconvergence $2N+1$ expected from such numerical methods, at least in problem (\ref{eq:bratu_problem}), same as in previous considered problems (\ref{eq:harm_osc}) and (\ref{eq:exp_diss}). The local solution $\mathbf{u}_{L}$ obtained by the DG predictor demonstrates the classical convergence characteristic of DG methods. 

\subsubsection{Example 4}

\begin{figure}[h!]
\captionsetup[subfigure]{%
position=bottom,
textfont=normalfont,
singlelinecheck=off,
justification=raggedright,
font=tiny
}
\centering
\begin{subfigure}{0.23\textwidth}
\includegraphics[width=\textwidth]{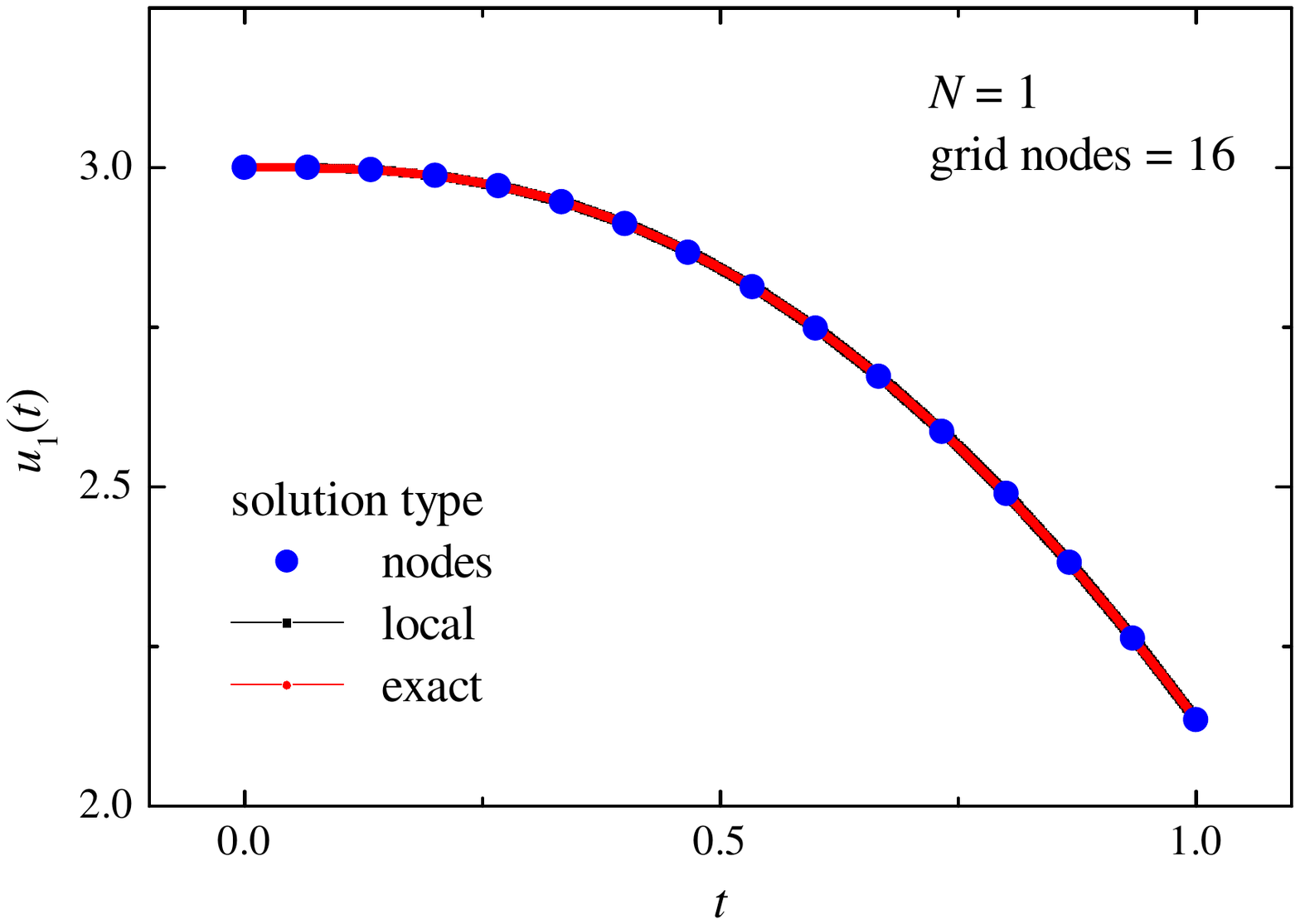}
\vspace{-8mm}\caption{\label{fig:third_order_ivp_ex_1_sols:a1}}
\end{subfigure}
\begin{subfigure}{0.23\textwidth}
\includegraphics[width=\textwidth]{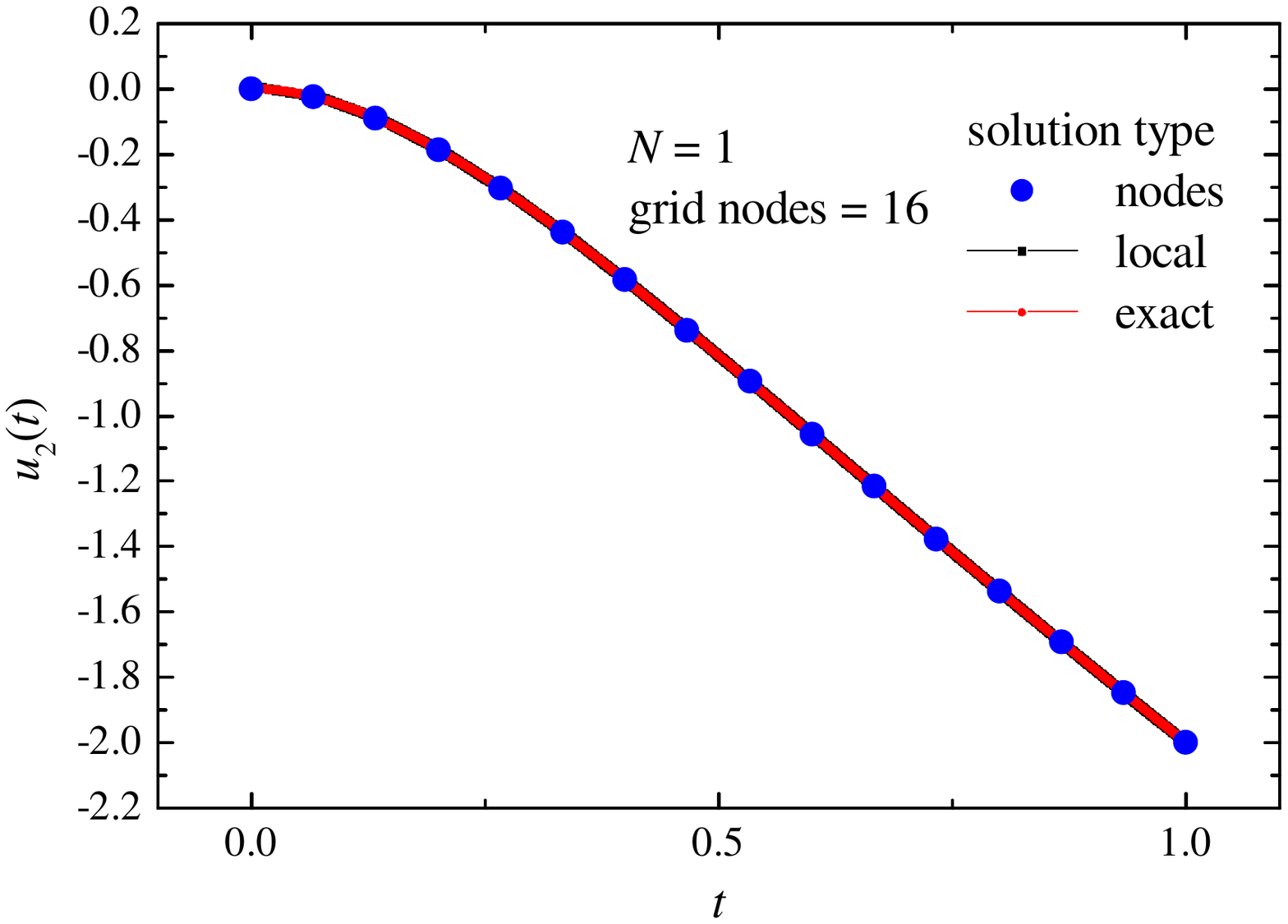}
\vspace{-8mm}\caption{\label{fig:third_order_ivp_ex_1_sols:a2}}
\end{subfigure}
\begin{subfigure}{0.23\textwidth}
\includegraphics[width=\textwidth]{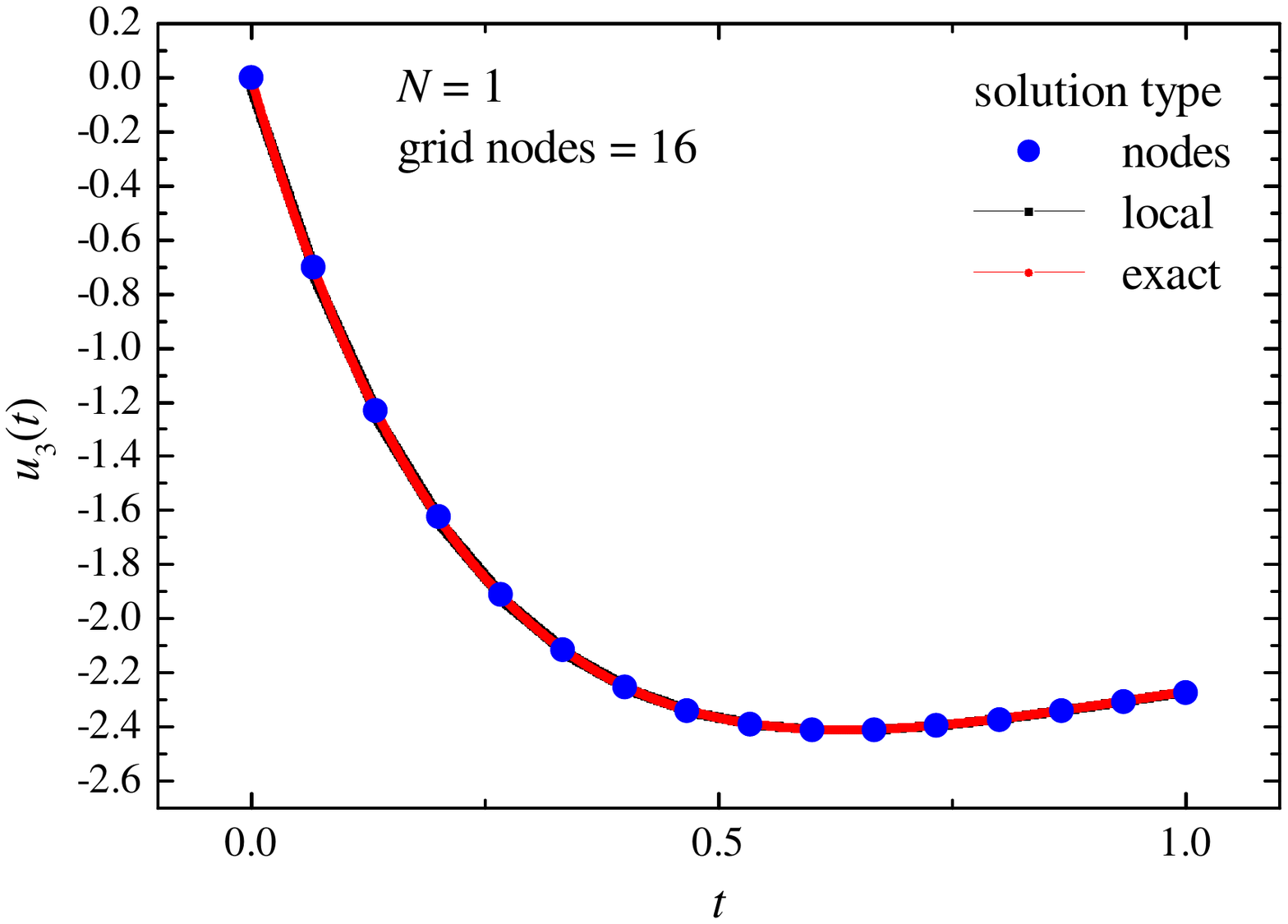}
\vspace{-8mm}\caption{\label{fig:third_order_ivp_ex_1_sols:a3}}
\end{subfigure}
\begin{subfigure}{0.23\textwidth}
\includegraphics[width=\textwidth]{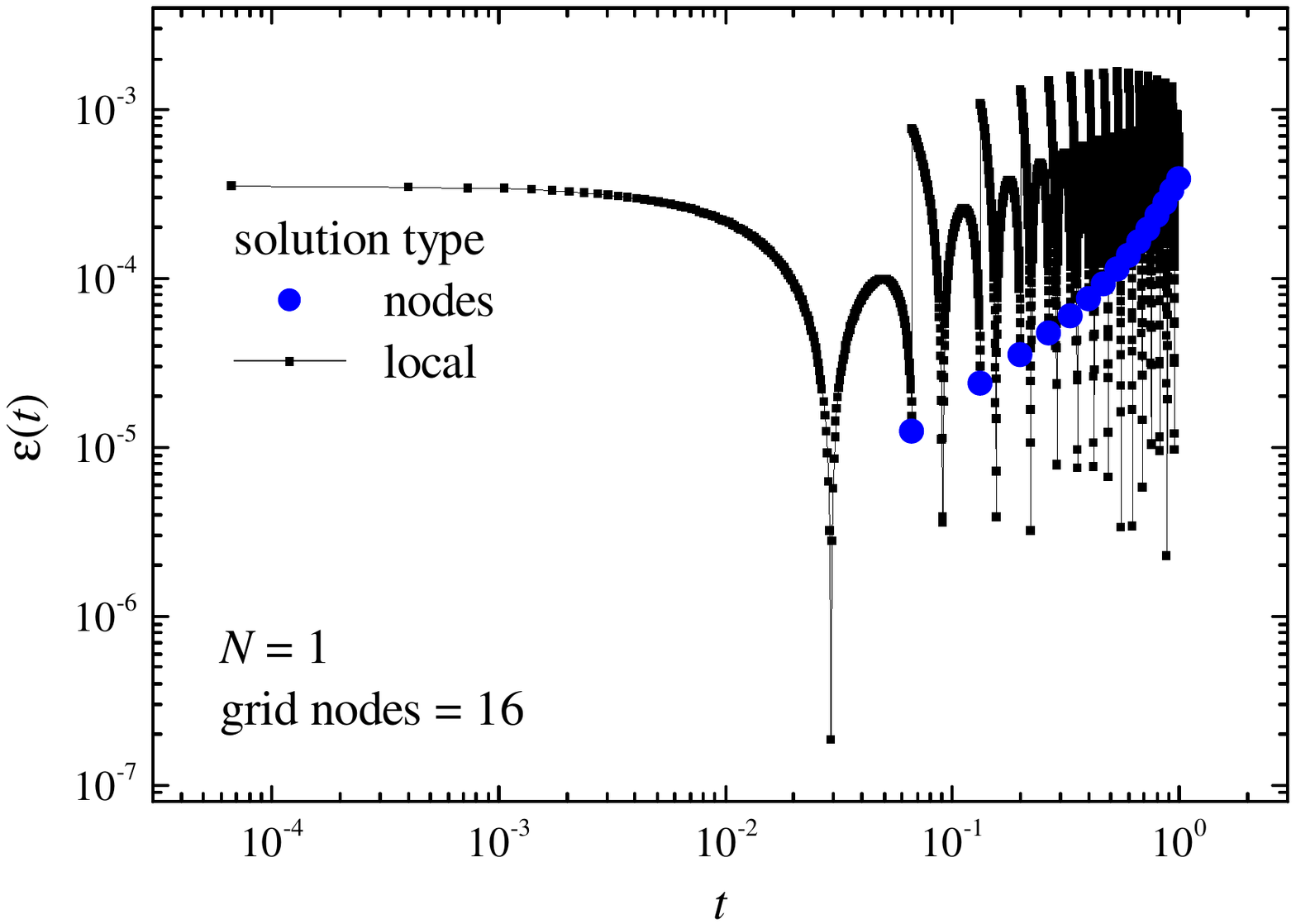}
\vspace{-8mm}\caption{\label{fig:third_order_ivp_ex_1_sols:a4}}
\end{subfigure}\\
\begin{subfigure}{0.23\textwidth}
\includegraphics[width=\textwidth]{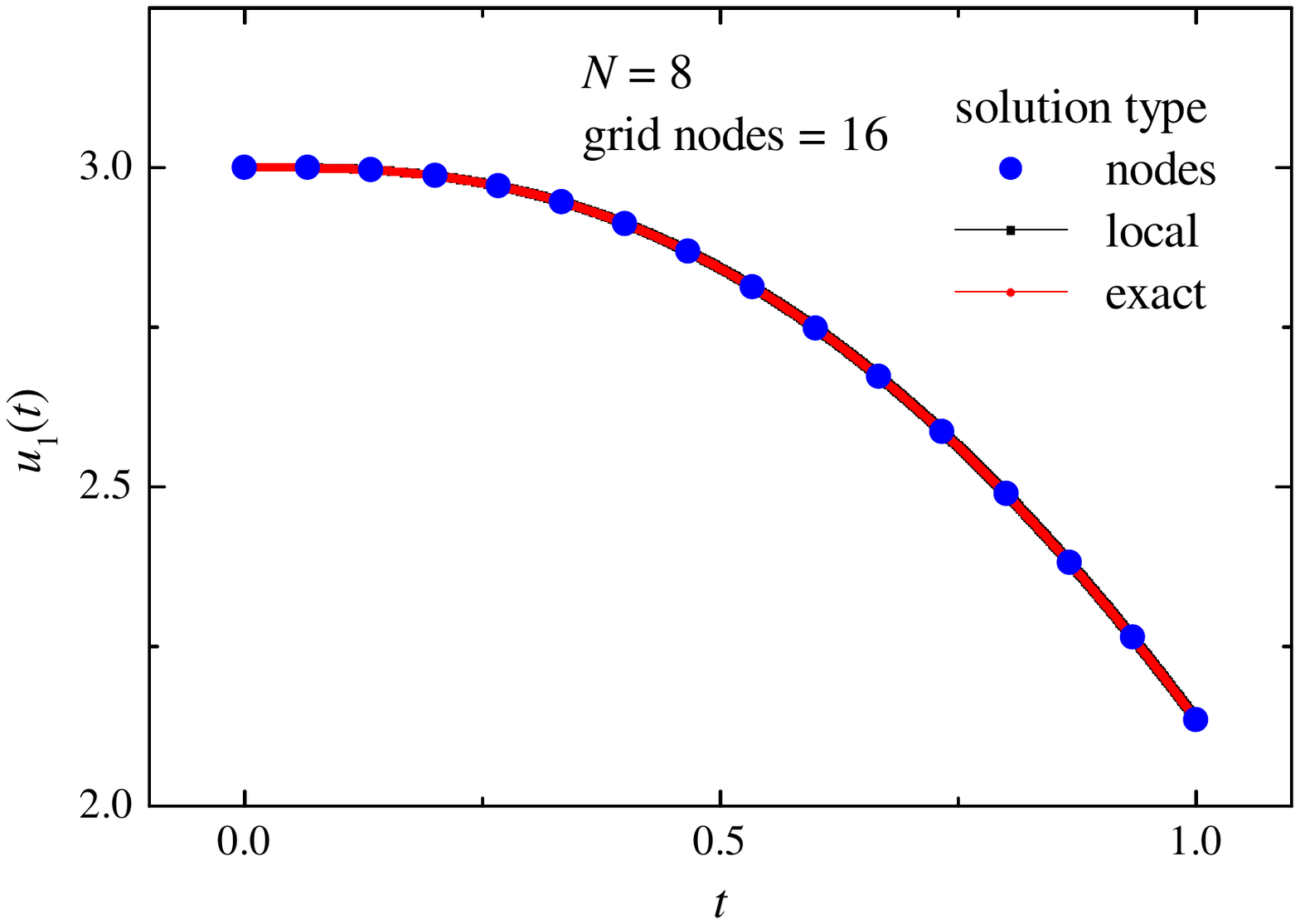}
\vspace{-8mm}\caption{\label{fig:third_order_ivp_ex_1_sols:b1}}
\end{subfigure}
\begin{subfigure}{0.23\textwidth}
\includegraphics[width=\textwidth]{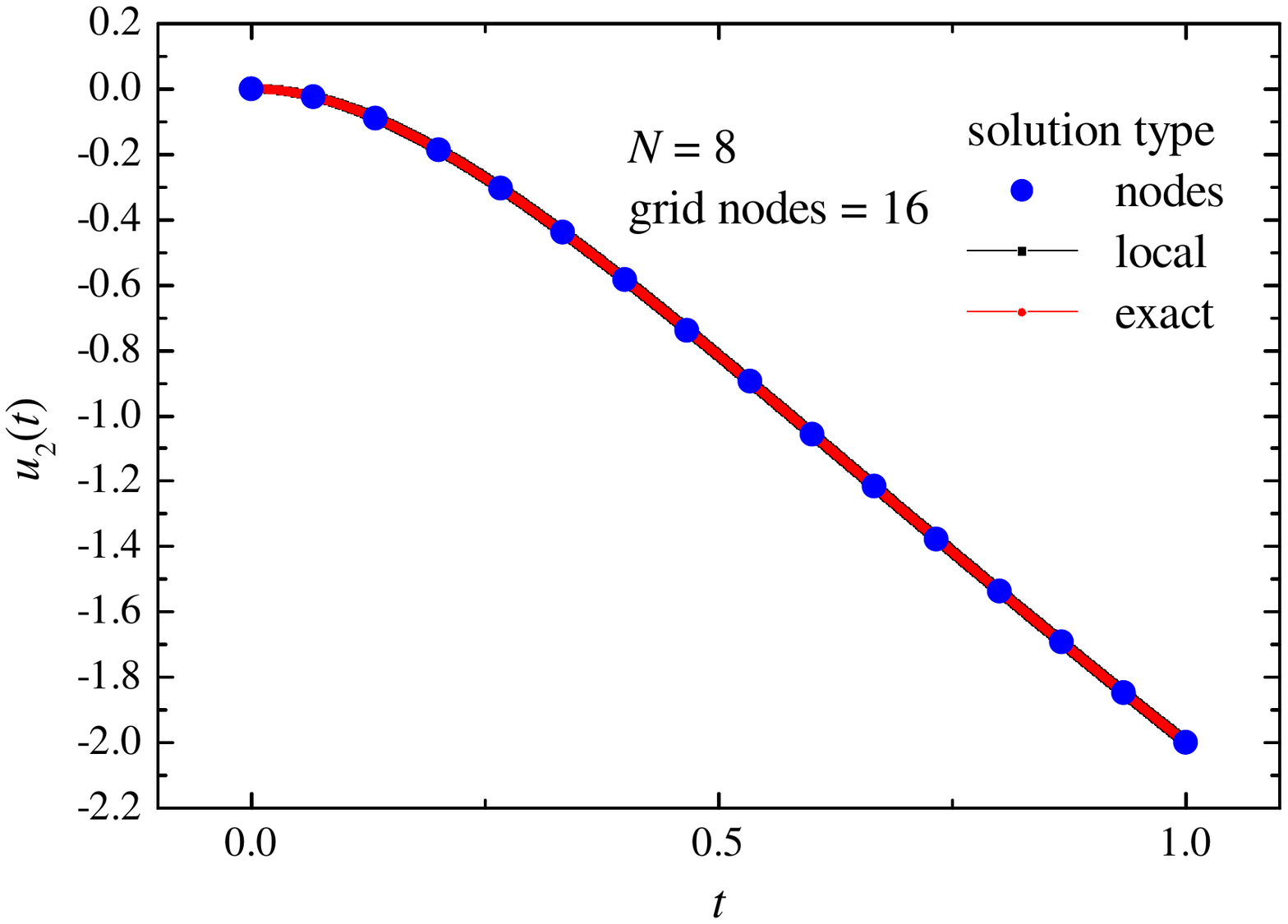}
\vspace{-8mm}\caption{\label{fig:third_order_ivp_ex_1_sols:b2}}
\end{subfigure}
\begin{subfigure}{0.23\textwidth}
\includegraphics[width=\textwidth]{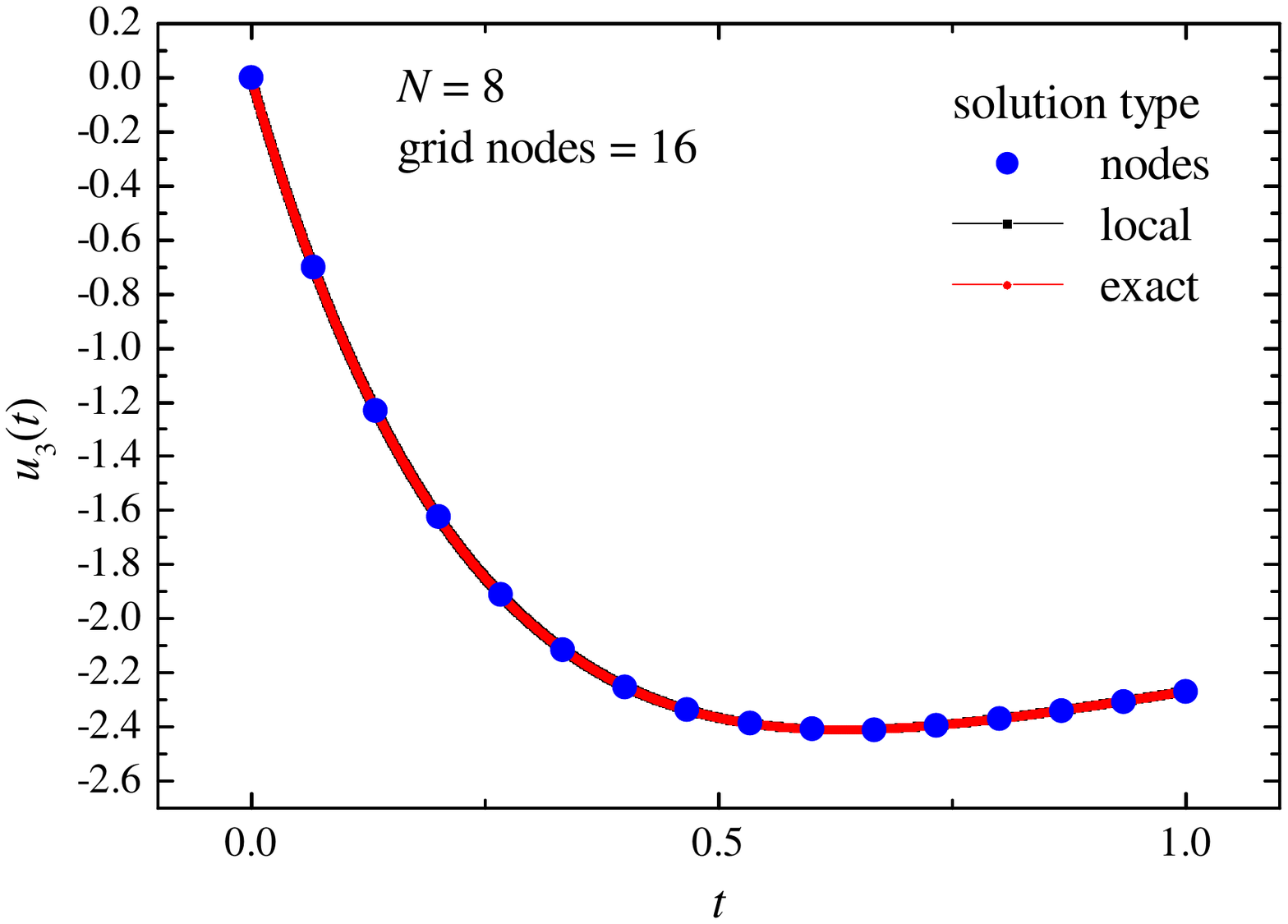}
\vspace{-8mm}\caption{\label{fig:third_order_ivp_ex_1_sols:b3}}
\end{subfigure}
\begin{subfigure}{0.23\textwidth}
\includegraphics[width=\textwidth]{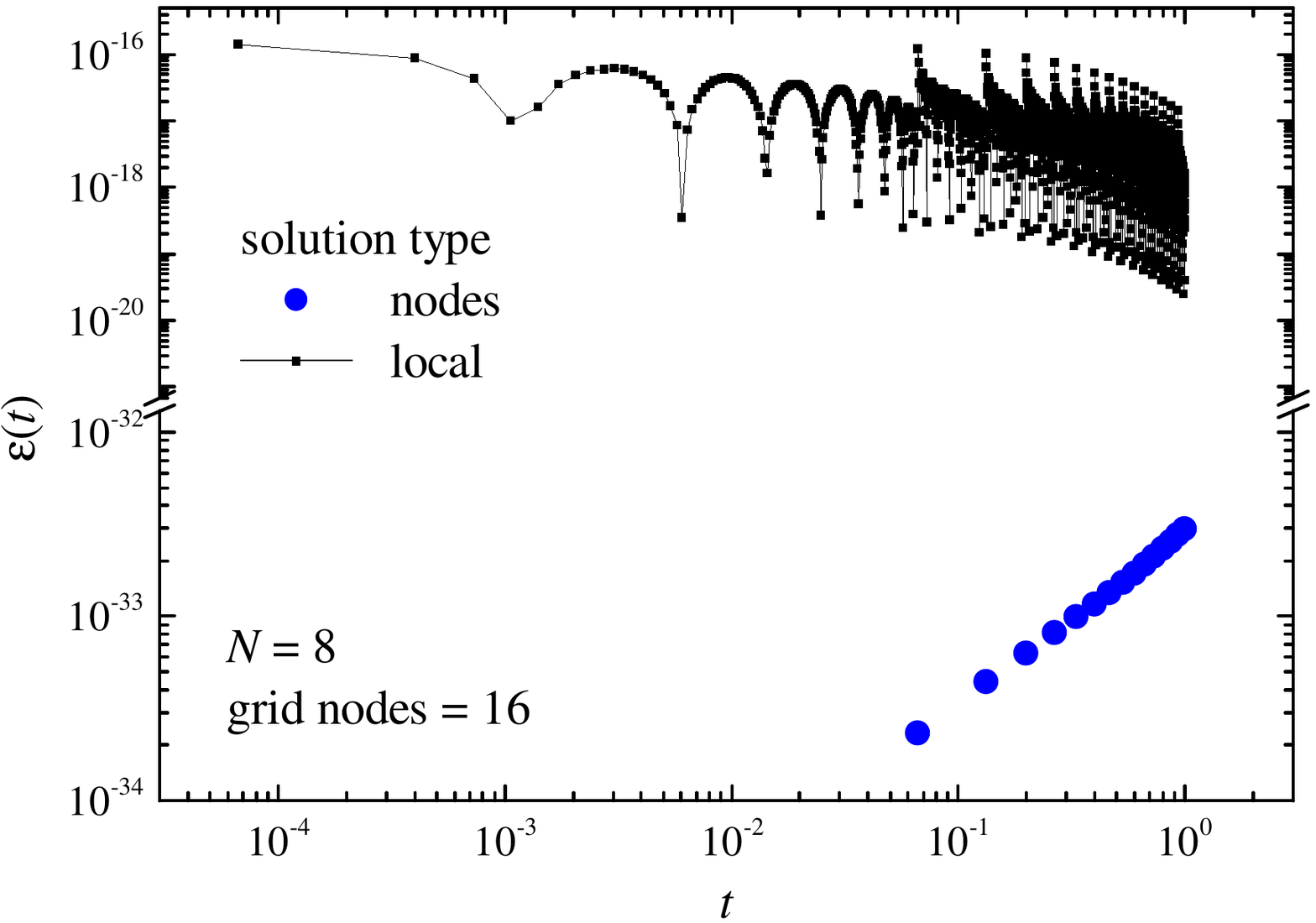}
\vspace{-8mm}\caption{\label{fig:third_order_ivp_ex_1_sols:b4}}
\end{subfigure}\\
\begin{subfigure}{0.23\textwidth}
\includegraphics[width=\textwidth]{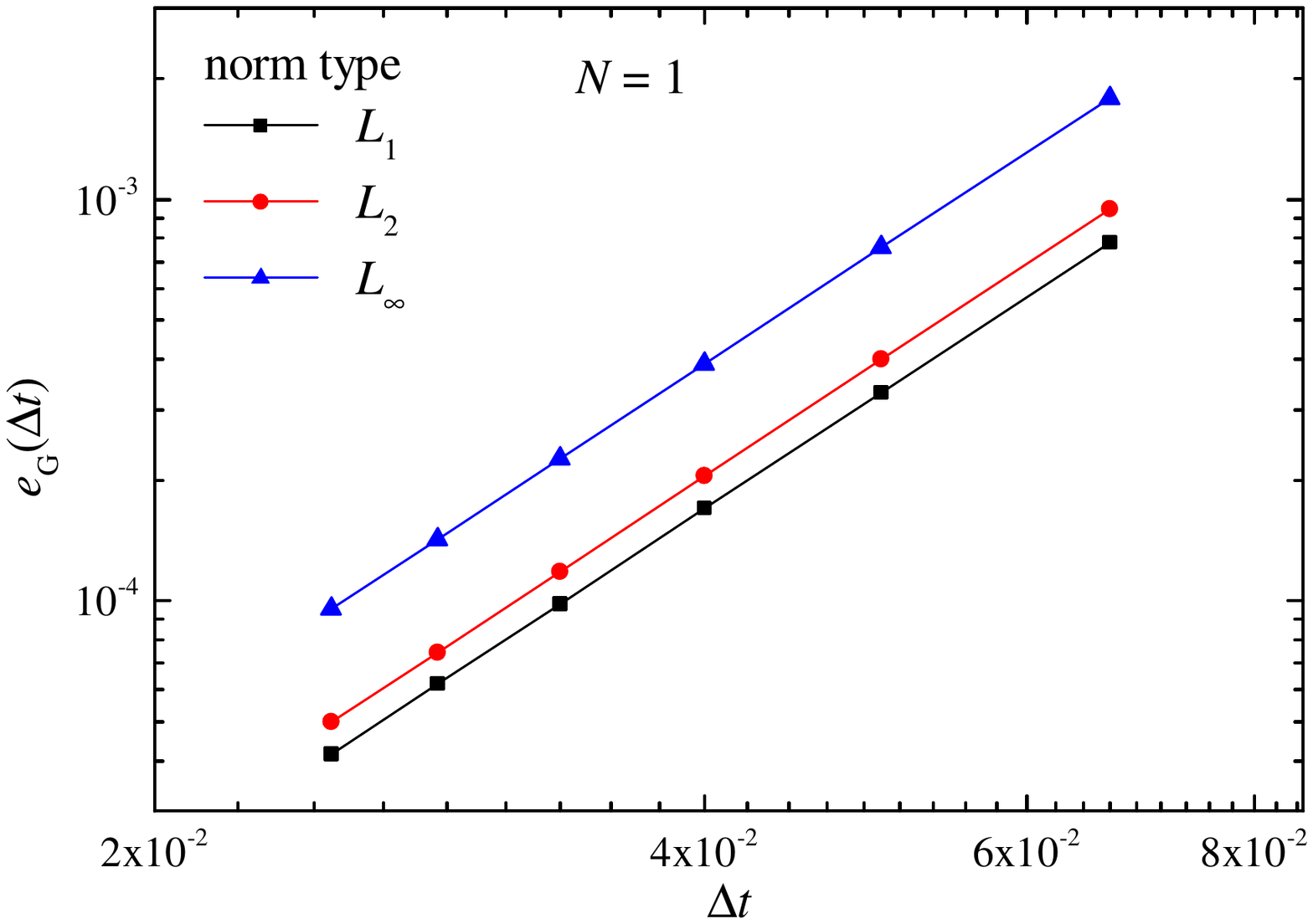}
\vspace{-8mm}\caption{\label{fig:third_order_ivp_ex_1_sols:c1}}
\end{subfigure}
\begin{subfigure}{0.23\textwidth}
\includegraphics[width=\textwidth]{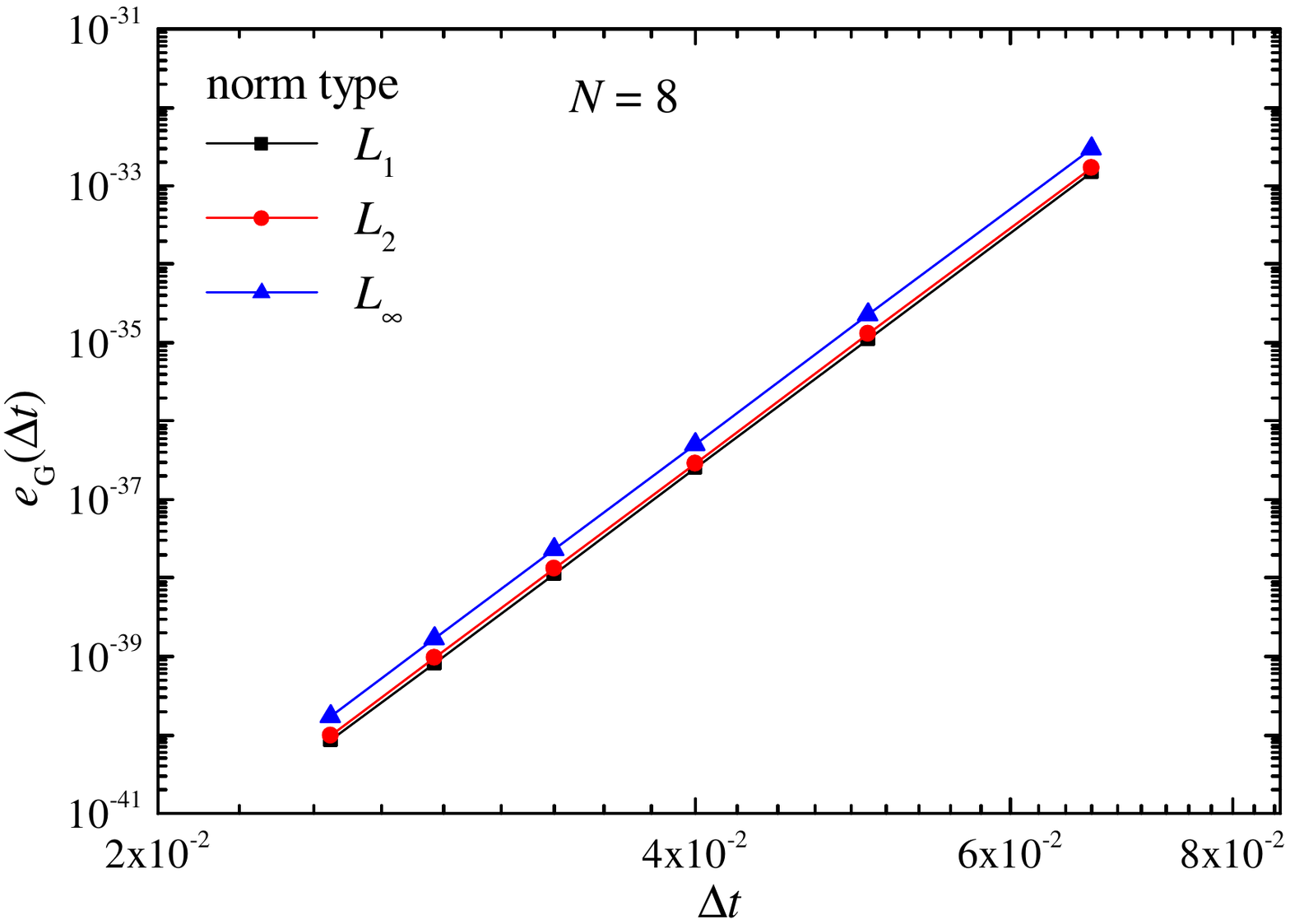}
\vspace{-8mm}\caption{\label{fig:third_order_ivp_ex_1_sols:c2}}
\end{subfigure}
\begin{subfigure}{0.23\textwidth}
\includegraphics[width=\textwidth]{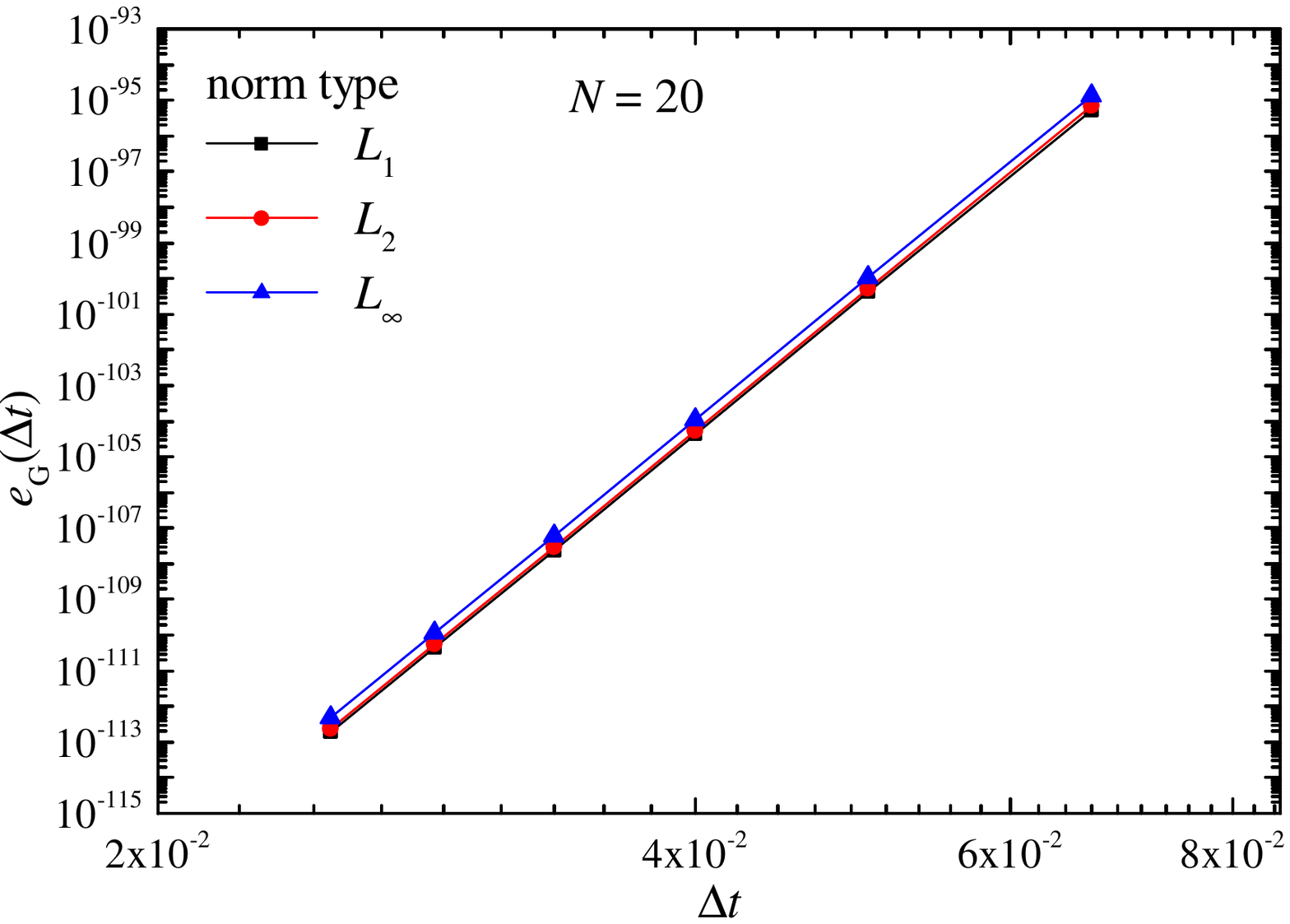}
\vspace{-8mm}\caption{\label{fig:third_order_ivp_ex_1_sols:c3}}
\end{subfigure}
\begin{subfigure}{0.23\textwidth}
\includegraphics[width=\textwidth]{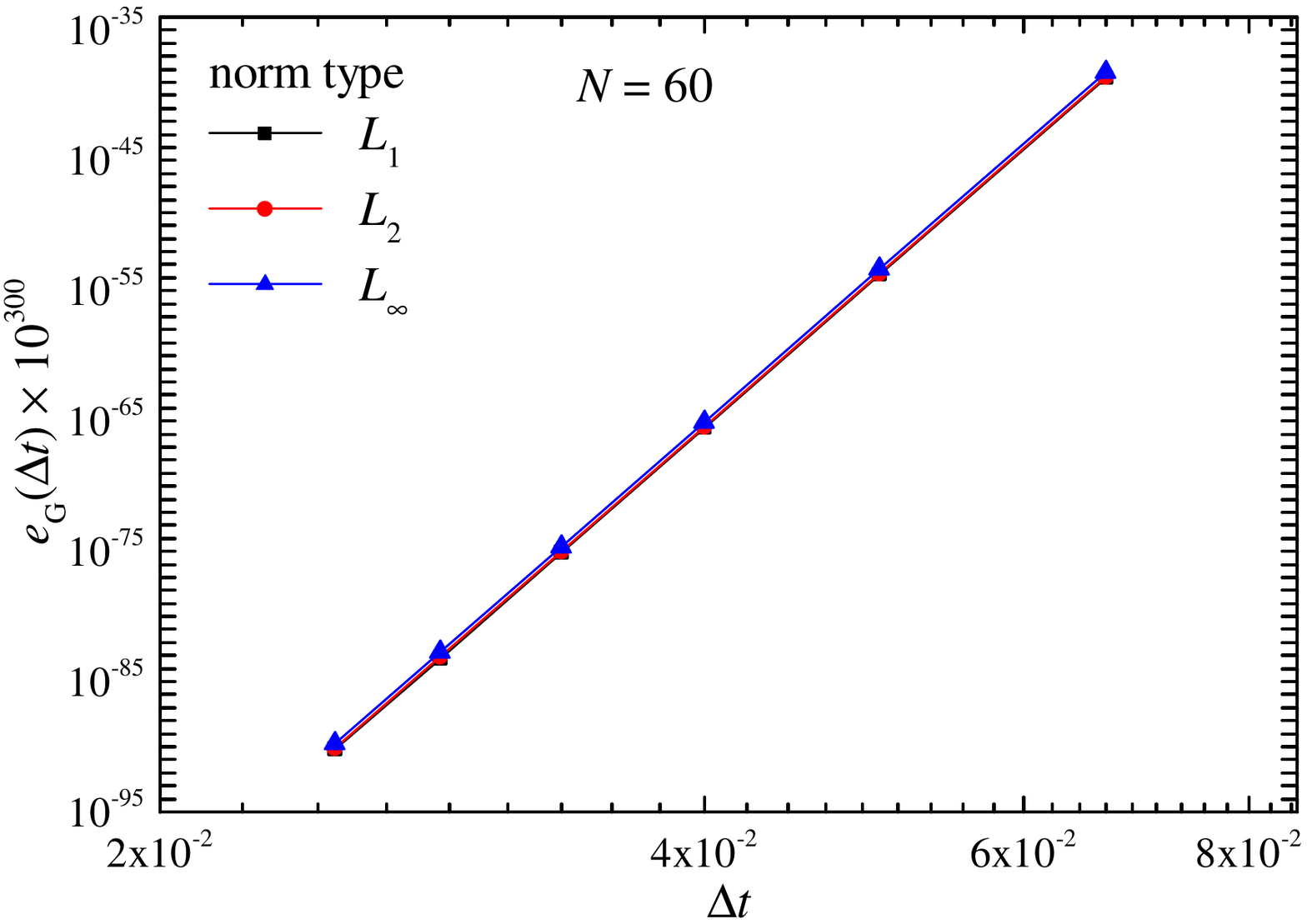}
\vspace{-8mm}\caption{\label{fig:third_order_ivp_ex_1_sols:c4}}
\end{subfigure}\\
\begin{subfigure}{0.23\textwidth}
\includegraphics[width=\textwidth]{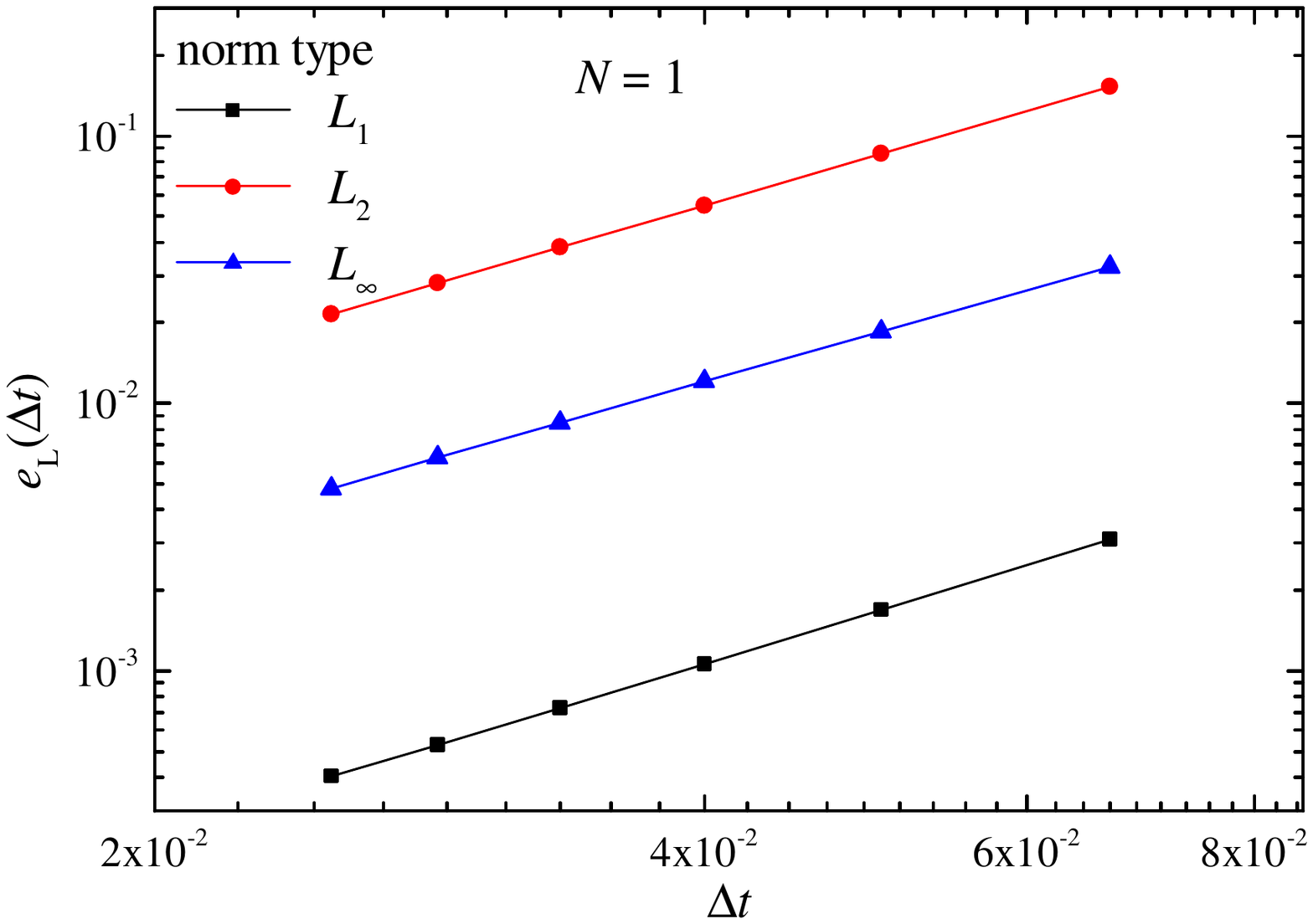}
\vspace{-8mm}\caption{\label{fig:third_order_ivp_ex_1_sols:d1}}
\end{subfigure}
\begin{subfigure}{0.23\textwidth}
\includegraphics[width=\textwidth]{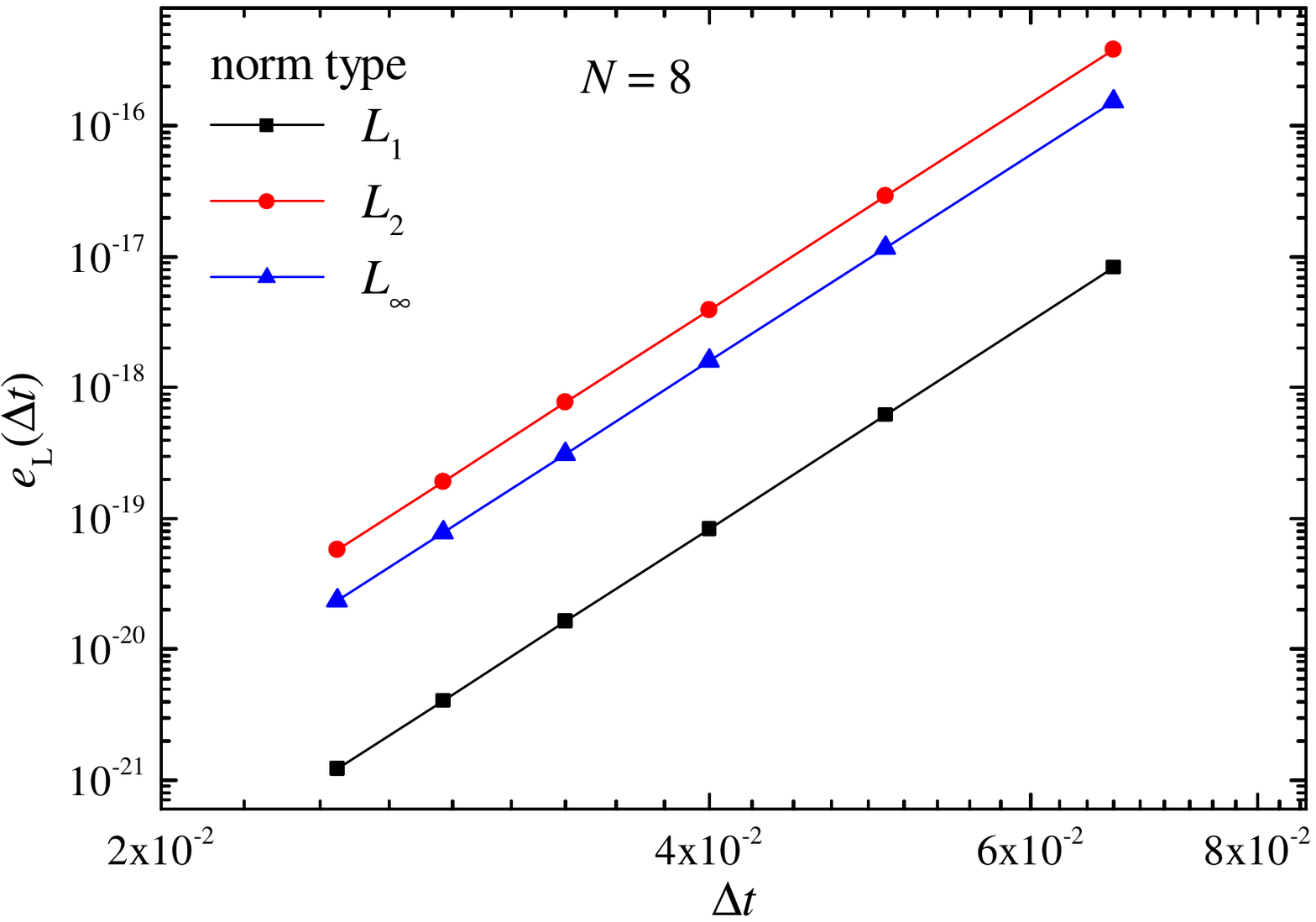}
\vspace{-8mm}\caption{\label{fig:third_order_ivp_ex_1_sols:d2}}
\end{subfigure}
\begin{subfigure}{0.23\textwidth}
\includegraphics[width=\textwidth]{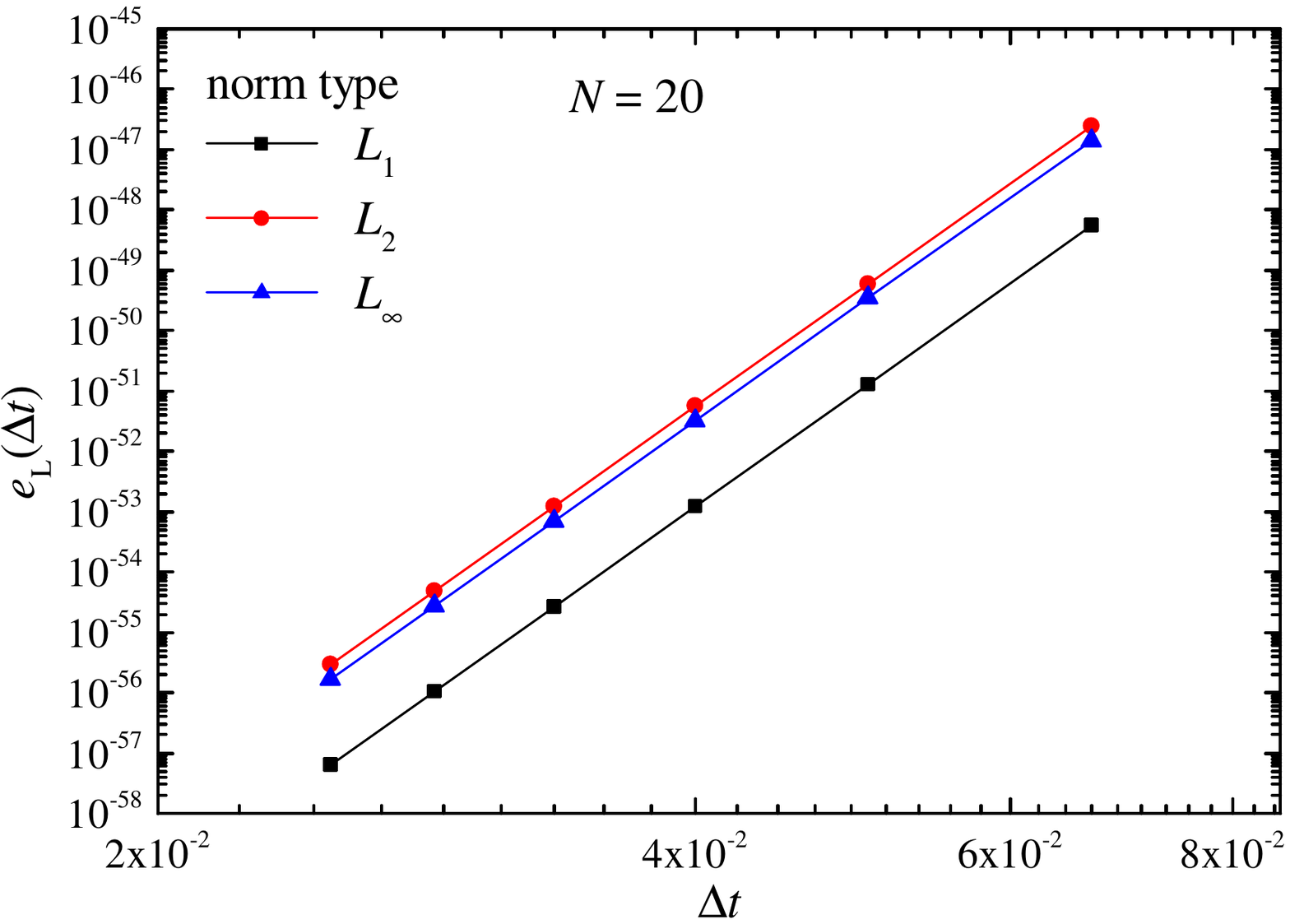}
\vspace{-8mm}\caption{\label{fig:third_order_ivp_ex_1_sols:d3}}
\end{subfigure}
\begin{subfigure}{0.23\textwidth}
\includegraphics[width=\textwidth]{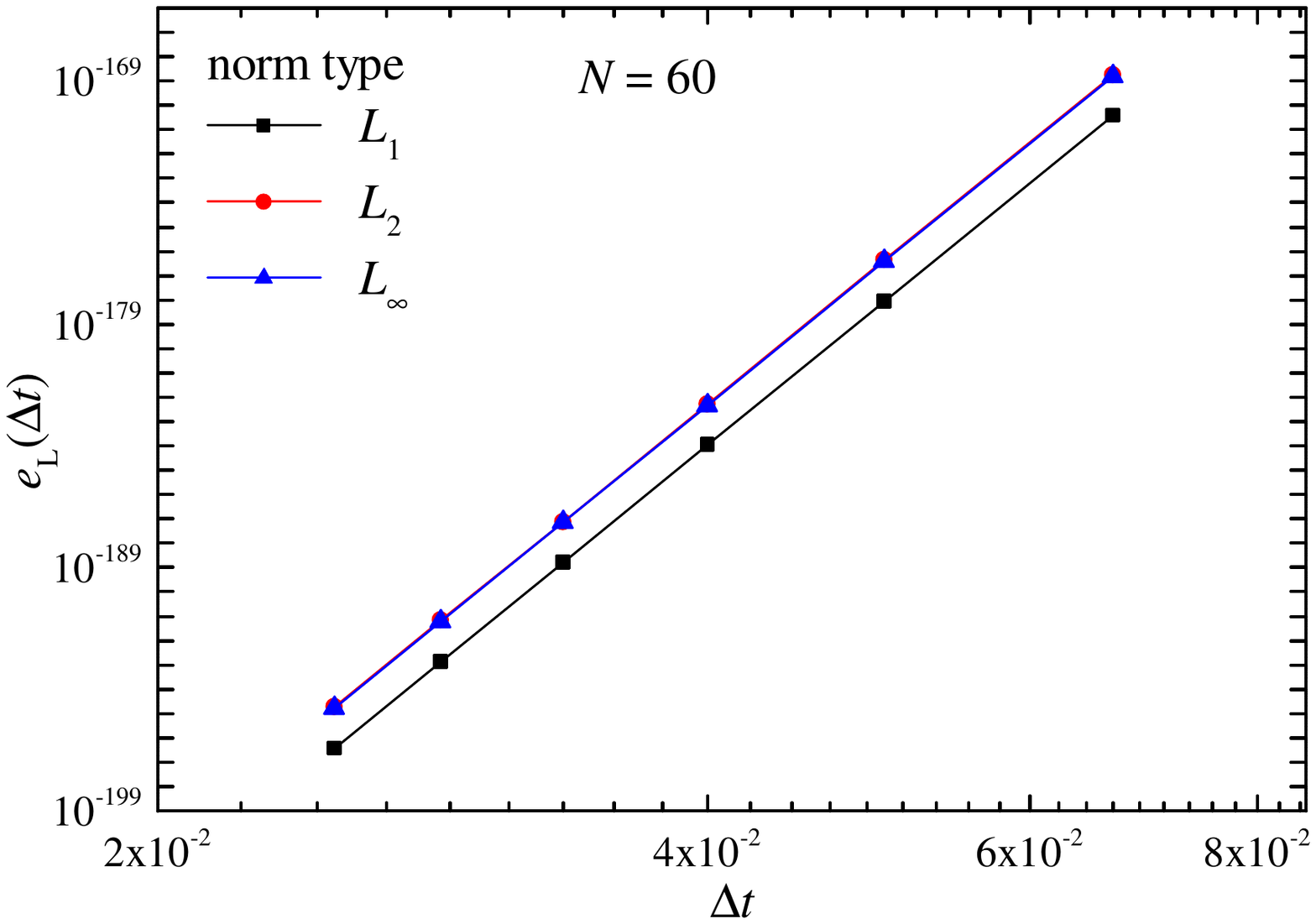}
\vspace{-8mm}\caption{\label{fig:third_order_ivp_ex_1_sols:d4}}
\end{subfigure}\\
\caption{%
Numerical solution of the problem (\ref{eq:third_order_ivp_ex_1}). Comparison of the solution at nodes $\mathbf{u}_{n}$, the local solution $\mathbf{u}_{L}(t)$ and the exact solution $\mathbf{u}^{\rm ex}(t)$ for components $u_{1}$ (\subref{fig:third_order_ivp_ex_1_sols:a1}, \subref{fig:third_order_ivp_ex_1_sols:b1}), $u_{2}$ (\subref{fig:third_order_ivp_ex_1_sols:a2}, \subref{fig:third_order_ivp_ex_1_sols:b2}) and $u_{3}$ (\subref{fig:third_order_ivp_ex_1_sols:a3}, \subref{fig:third_order_ivp_ex_1_sols:b3}), the error $\varepsilon(t)$ (\subref{fig:third_order_ivp_ex_1_sols:a4}, \subref{fig:third_order_ivp_ex_1_sols:b4}), obtained using polynomials with degrees $N = 1$ (\subref{fig:third_order_ivp_ex_1_sols:a1}, \subref{fig:third_order_ivp_ex_1_sols:a2}, \subref{fig:third_order_ivp_ex_1_sols:a3}, \subref{fig:third_order_ivp_ex_1_sols:a4}) and $N = 8$ (\subref{fig:third_order_ivp_ex_1_sols:b1}, \subref{fig:third_order_ivp_ex_1_sols:b2}, \subref{fig:third_order_ivp_ex_1_sols:b3}, \subref{fig:third_order_ivp_ex_1_sols:b4}). Dependence of the global error for the solution at nodes $e_{G}$ (\subref{fig:third_order_ivp_ex_1_sols:c1}, \subref{fig:third_order_ivp_ex_1_sols:c2}, \subref{fig:third_order_ivp_ex_1_sols:c3}, \subref{fig:third_order_ivp_ex_1_sols:c4}) and the local solution $e_{L}$ (\subref{fig:third_order_ivp_ex_1_sols:d1}, \subref{fig:third_order_ivp_ex_1_sols:d2}, \subref{fig:third_order_ivp_ex_1_sols:d3}, \subref{fig:third_order_ivp_ex_1_sols:d4}) on the discretization step $\mathrm{\Delta}t$, obtained in the norms $L_{1}$, $L_{2}$ and $L_{\infty}$, obtained using polynomials with degrees $N = 1$ (\subref{fig:third_order_ivp_ex_1_sols:c1}, \subref{fig:third_order_ivp_ex_1_sols:d1}), $N = 8$ (\subref{fig:third_order_ivp_ex_1_sols:c2}, \subref{fig:third_order_ivp_ex_1_sols:d2}), $N = 20$ (\subref{fig:third_order_ivp_ex_1_sols:c3}, \subref{fig:third_order_ivp_ex_1_sols:d3}) and $N = 60$ (\subref{fig:third_order_ivp_ex_1_sols:c4}, \subref{fig:third_order_ivp_ex_1_sols:d4}).
}
\label{fig:third_order_ivp_ex_1_sols}
\end{figure} 
\begin{table}[h!]
\centering
\caption{%
Orders of convergence $p_{L_{1}}$, $p_{L_{2}}$, $p_{L_{\infty}}$ of the ADER-DG numerical method with local DG predictor calculated in norms $L_{1}$, $L_{2}$, $L_{\infty}$ for the problem (\ref{eq:third_order_ivp_ex_1}); $N$ is the degree of the basis polynomials $\varphi_{p}$ used for the DG representations. The orders of convergence $p$ for the solution at nodes $\mathbf{u}_{n}$ are compared to the theoretical value $p_{\rm (G)} = 2N+1$, the orders of convergence $p$ for the local solution $\mathbf{u}_{L}$ represented between the nodes are compared to the theoretical value $p_{\rm (L)} = N+1$.
}
\label{tab:third_order_ivp_ex_1_sols}
\begin{tabular}{@{}lllllllll@{}}
\toprule
& \multicolumn{4}{c}{order $p$ for solution at nodes} & \multicolumn{4}{c}{order $p$ for local solution} \\
\midrule
$N$ & $p_{L_{1}}$ & $p_{L_{2}}$ & $p_{L_{\infty}}$ & theor. & $p_{L_{1}}$ & $p_{L_{2}}$ & $p_{L_{\infty}}$ & theor. \\
\midrule
$1$	&	$2.99$	&	$3.00$	&	$2.98$		&	$3$	&	$2.08$	&	$1.99$	&	$1.94$		&	$2$	\\
$2$	&	$4.96$	&	$4.97$	&	$4.96$		&	$5$	&	$2.99$	&	$2.97$	&	$2.94$		&	$3$	\\
$3$	&	$6.95$	&	$6.96$	&	$6.97$		&	$7$	&	$3.99$	&	$3.97$	&	$3.94$		&	$4$	\\
$4$	&	$9.03$	&	$9.06$	&	$9.00$		&	$9$	&	$4.99$	&	$4.97$	&	$4.94$		&	$5$	\\
$5$	&	$11.03$	&	$11.04$	&	$10.99$		&	$11$	&	$6.00$	&	$5.98$	&	$5.95$		&	$6$	\\
$6$	&	$13.01$	&	$13.02$	&	$12.99$		&	$13$	&	$7.00$	&	$6.98$	&	$6.95$		&	$7$	\\
$7$	&	$15.01$	&	$15.02$	&	$14.99$		&	$15$	&	$8.00$	&	$7.98$	&	$7.95$		&	$8$	\\
$8$	&	$16.99$	&	$17.00$	&	$16.99$		&	$17$	&	$9.00$	&	$8.98$	&	$8.95$		&	$9$	\\
$9$	&	$18.98$	&	$18.99$	&	$18.99$		&	$19$	&	$10.00$	&	$9.98$	&	$9.95$		&	$10$	\\
$10$	&	$20.96$	&	$20.98$	&	$20.98$		&	$21$	&	$11.00$	&	$10.98$	&	$10.95$		&	$11$	\\
\midrule
$15$	&	$30.99$	&	$31.00$	&	$30.99$		&	$31$	&	$16.00$	&	$15.98$	&	$15.95$		&	$16$	\\
$20$	&	$41.01$	&	$41.02$	&	$41.00$		&	$41$	&	$21.00$	&	$20.98$	&	$20.95$		&	$21$	\\
$25$	&	$51.03$	&	$51.04$	&	$51.00$		&	$51$	&	$26.00$	&	$25.98$	&	$25.95$		&	$26$	\\
$30$	&	$60.99$	&	$61.00$	&	$60.99$		&	$61$	&	$31.00$	&	$30.98$	&	$30.95$		&	$31$	\\
$35$	&	$71.00$	&	$71.01$	&	$71.00$		&	$71$	&	$36.00$	&	$35.98$	&	$35.96$		&	$36$	\\
$40$	&	$81.02$	&	$81.03$	&	$81.00$		&	$81$	&	$41.00$	&	$40.98$	&	$40.96$		&	$41$	\\
$45$	&	$91.04$	&	$91.06$	&	$91.00$		&	$91$	&	$46.00$	&	$45.98$	&	$45.96$		&	$46$	\\
$50$	&	$100.99$	&	$101.00$	&	$101.00$		&	$101$	&	$51.00$	&	$50.98$	&	$50.96$		&	$51$	\\
$55$	&	$111.01$	&	$111.02$	&	$111.00$		&	$111$	&	$56.00$	&	$55.98$	&	$55.96$		&	$56$	\\
$60$	&	$120.98$	&	$120.99$	&	$121.01$		&	$121$	&	$61.00$	&	$60.98$	&	$60.96$		&	$61$	\\
\bottomrule
\end{tabular}
\end{table} 

The first three examples (\ref{eq:harm_osc}), (\ref{eq:exp_diss}) and (\ref{eq:bratu_problem}) of selected problems were related to the solution of second-order differential equations, which reduced to the solution of the IVP for the ODE of a system of two first-order equations. In the following three sections, examples of problems for the third-order differential equations will be considered~\cite{dg_ivp_ode_5}, which reduced to the solution of the IVP for the ODE of a system of three first-order equations.

The first example for IVP third-order ODE was chosen in the following form:
\begin{equation}
\begin{split}
&\dddot{x} = 2\ddot{x} + 3\dot{x} - 10 x + (34 t - 16) \exp(-2t) - 10 t^{2} + 6t + 34;\\
&x(0) = 3;\quad \dot{x}(0) = 0;\quad \ddot{x}(0) = 0;\qquad 0 \leqslant t \leqslant 1;\\
\end{split}
\end{equation}
with exact analytical solution of this problem having the form $x(t) = t^{2}\exp(-2t) - t^{2} + 3$. The IVP for a first-order ODE system corresponding to this problem has the following form:
\begin{equation}\label{eq:third_order_ivp_ex_1}
\begin{split}
&\frac{du_{1}}{dt} = u_{2};\quad \frac{du_{2}}{dt} = u_{3};\\
&\dfrac{du_{3}}{dt} = 2 u_{3} + 3 u_{2} - 10 u_{1} + (34 t - 16) \exp(-2t) - 10 t^{2} + 6t + 34;\\
&u_{1}(0) = 1;\quad u_{2}(0) = 0;\quad u_{3}(0) = 0;\\
\end{split}
\end{equation}
where $\mathbf{u} = [u_{1},\, u_{2},\, u_{3}]^{T} = [x,\, \dot{x},\, \ddot{x}]^{T}$ is a desired three component vector function. The solution definition domain $\Omega$ was chosen as $0 \leqslant t \leqslant 1$. The full three component exact analytical solution of this problem has the form
\begin{equation}
\mathbf{u}^{\rm ex} = \left[
\begin{array}{c}
t^{2}\exp(-2t) - t^{2} + 3\\
2 t \left((1-t^{2})\exp(-2t) - 1\right)\\
2 \left((1-4t+2t^{2})\exp(-2t)-1\right)
\end{array}
\right].
\end{equation}

The study of the convergence of the numerical solution was carried out on the basis of the analysis of the solution for the set of $6$ different uniform discretizations of the domain of definition $\Omega$ of the desired function, with the number of grid nodes $L = 16$, $21$, $26$, $31$, $36$, $41$; note that the zero node, in which the initial condition was defined, also refers to grid nodes, so the discretization step was determined by the expression $\Delta t = 2\pi/(L-1)$.

The results of the study of the application of the numerical method for this problem are presented in Fig.~\ref{fig:third_order_ivp_ex_1_sols} and in Table~\ref{tab:third_order_ivp_ex_1_sols}. The results obtained, in terms of the accuracy of the numerical solution when compared with the exact analytical solution, demonstrate all the same properties that were revealed in the previous considered problems for second-order differential equations. The error $\varepsilon(t)$ for $N > 1$ shows that the error scales for the local solution and the solution at nodes differ by several orders of magnitude, reaching values of $10$--$15$ and $190$--$200$ orders of magnitude for a degree of $N = 8$ and $N = 60$, respectively.

The values of the global errors $e$ were calculated for the solution at nodes $\mathbf{u}_{n}$ and the local solution $\mathbf{u}_{L}$, in three different norms $L_{1}$, $L_{2}$, $L_{\infty}$ (\ref{eq:global_eps_def}), which are shown in Fig.~\ref{fig:third_order_ivp_ex_1_sols}. Based on this, the convergence orders $p_{L_{1}}$, $p_{L_{2}}$, $p_{L_{\infty}}$ were calculated, which were presented in Table~\ref{tab:third_order_ivp_ex_1_sols}. The convergence orders $p$ for the solution at nodes are compared to the theoretical value $p_{\rm (G)} = 2N+1$, the convergence orders $p$ for the local solution represented between the nodes are compared to the theoretical value $p_{\rm (L)} = N+1$, in accordance with (\ref{eq:expect_orders}). It is interesting to note that all the calculated values of the convergence orders $p$ correspond to the expected values $p_{\rm (G)}$ and $p_{\rm (L)}$. There are no significant deviations from the expected values.

As a result, it becomes clear that the solution at nodes $\mathbf{u}_{n}$ obtained by the ADER-DG numerical method with a local DG predictor demonstrates the classical superconvergence $2N+1$ expected from such numerical methods, at least in problem (\ref{eq:third_order_ivp_ex_1}). The local solution $\mathbf{u}_{L}$ obtained by the DG predictor demonstrates the classical convergence characteristic of DG methods. 

\subsubsection{Example 5}

\begin{figure}[h!]
\captionsetup[subfigure]{%
position=bottom,
textfont=normalfont,
singlelinecheck=off,
justification=raggedright,
font=tiny
}
\centering
\begin{subfigure}{0.23\textwidth}
\includegraphics[width=\textwidth]{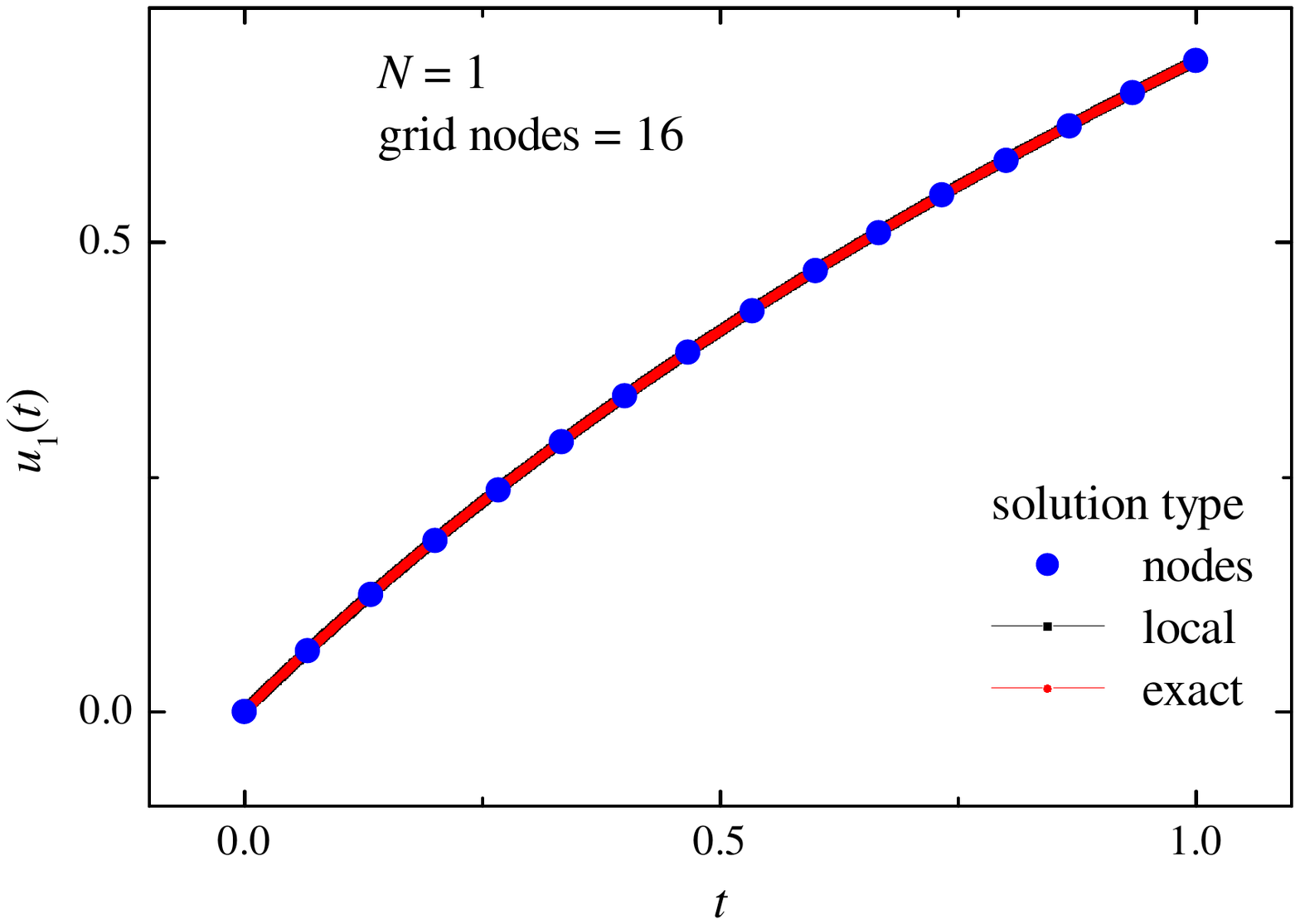}
\vspace{-8mm}\caption{\label{fig:third_order_ivp_ex_2_sols:a1}}
\end{subfigure}
\begin{subfigure}{0.23\textwidth}
\includegraphics[width=\textwidth]{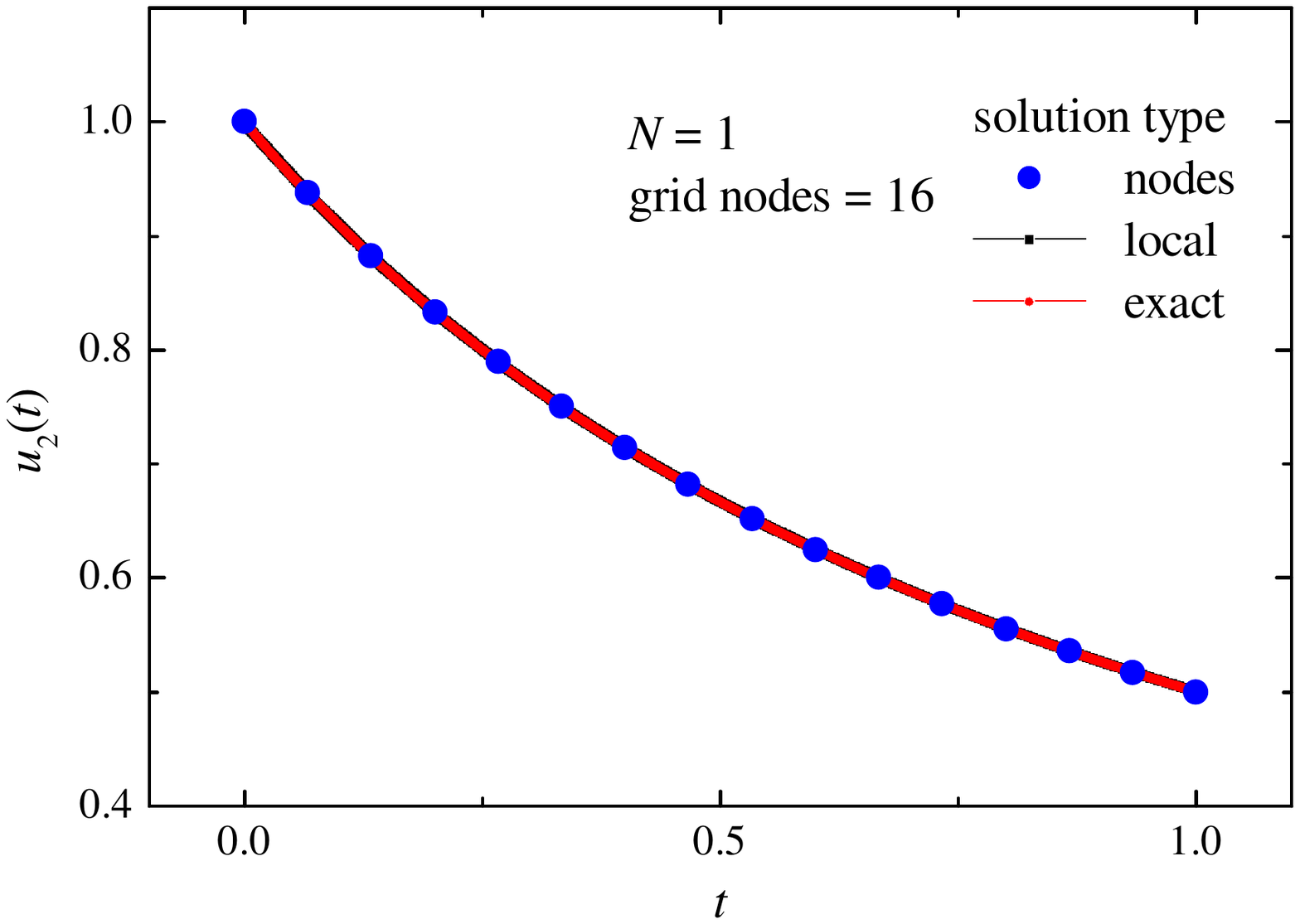}
\vspace{-8mm}\caption{\label{fig:third_order_ivp_ex_2_sols:a2}}
\end{subfigure}
\begin{subfigure}{0.23\textwidth}
\includegraphics[width=\textwidth]{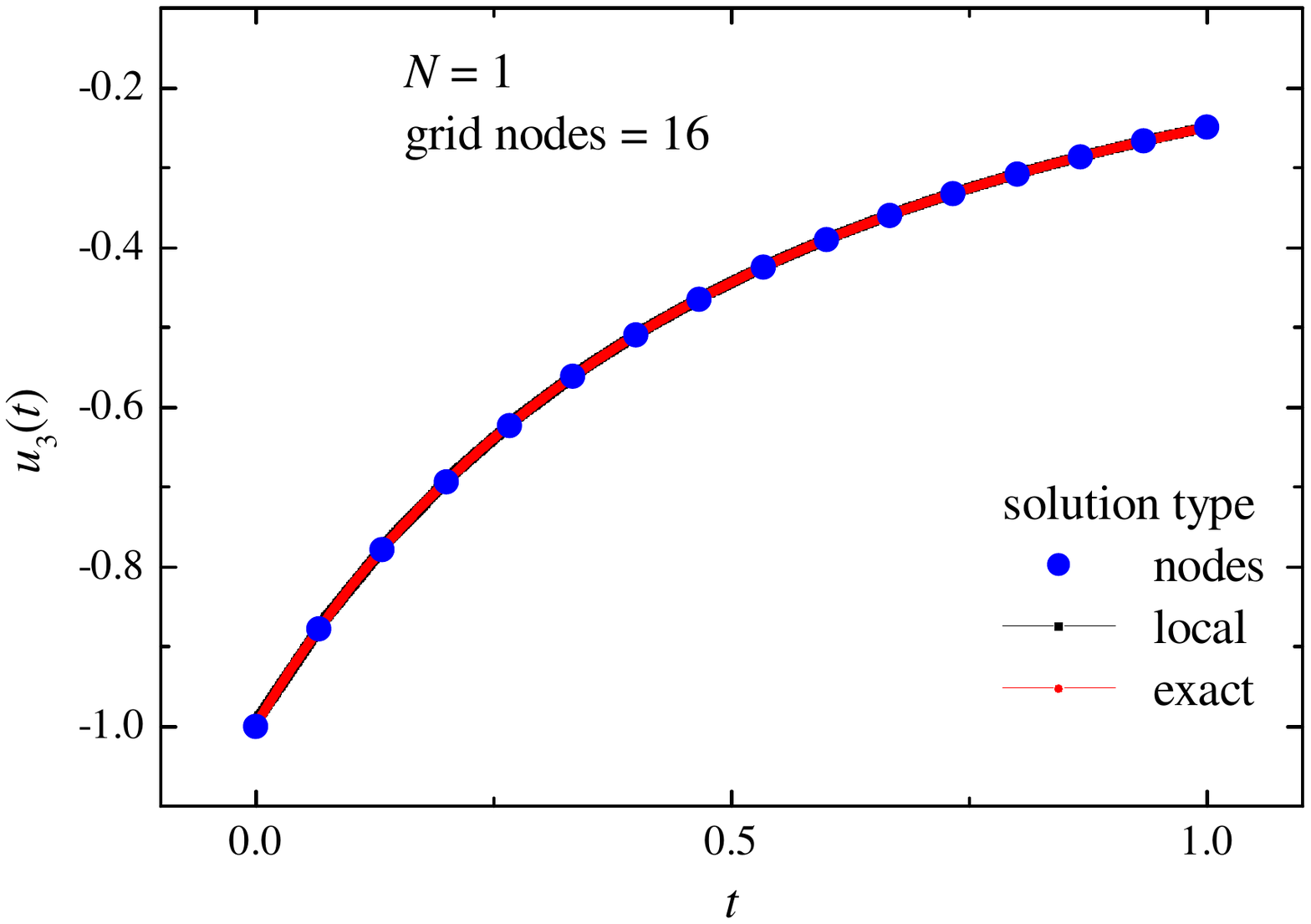}
\vspace{-8mm}\caption{\label{fig:third_order_ivp_ex_2_sols:a3}}
\end{subfigure}
\begin{subfigure}{0.23\textwidth}
\includegraphics[width=\textwidth]{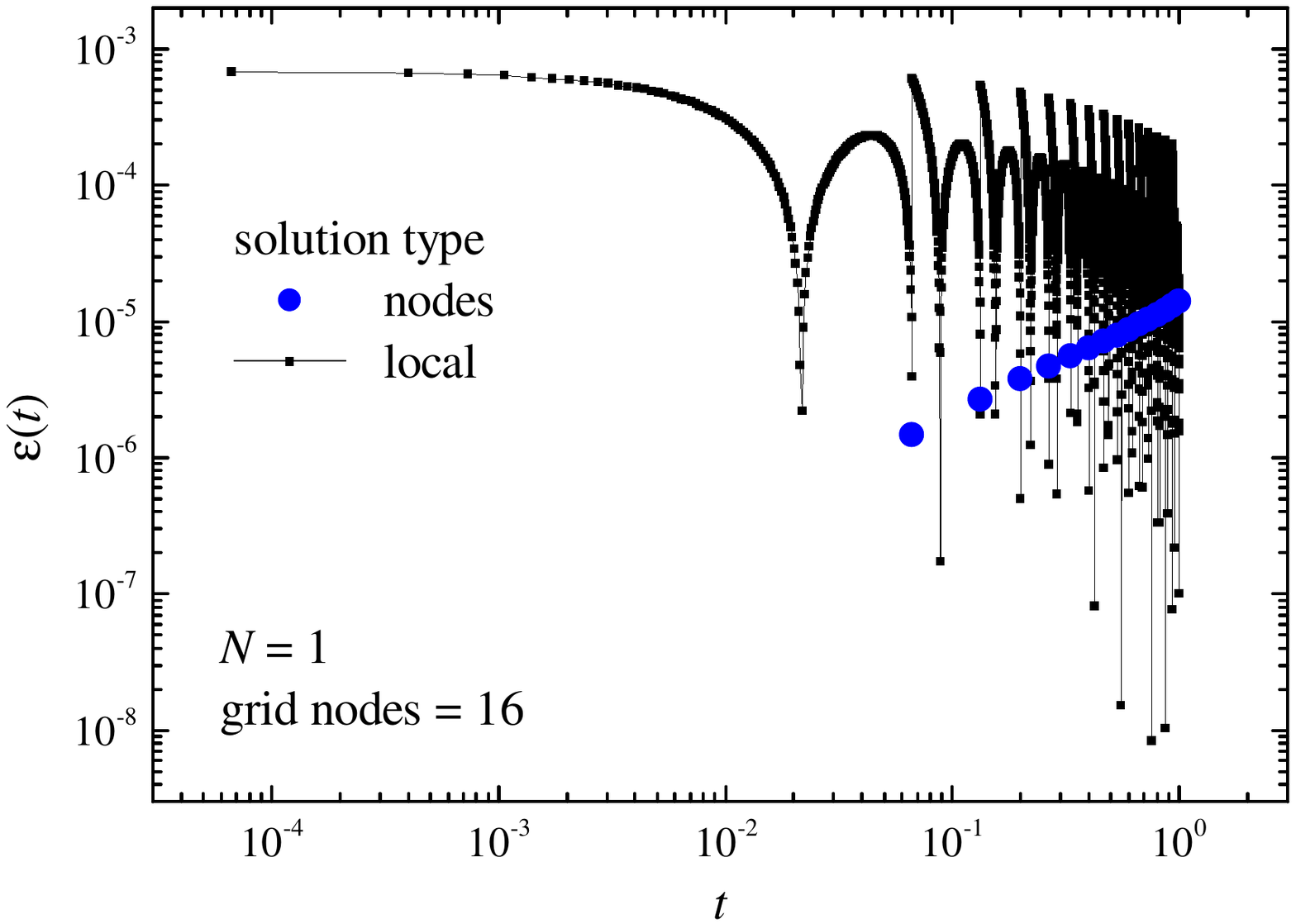}
\vspace{-8mm}\caption{\label{fig:third_order_ivp_ex_2_sols:a4}}
\end{subfigure}\\
\begin{subfigure}{0.23\textwidth}
\includegraphics[width=\textwidth]{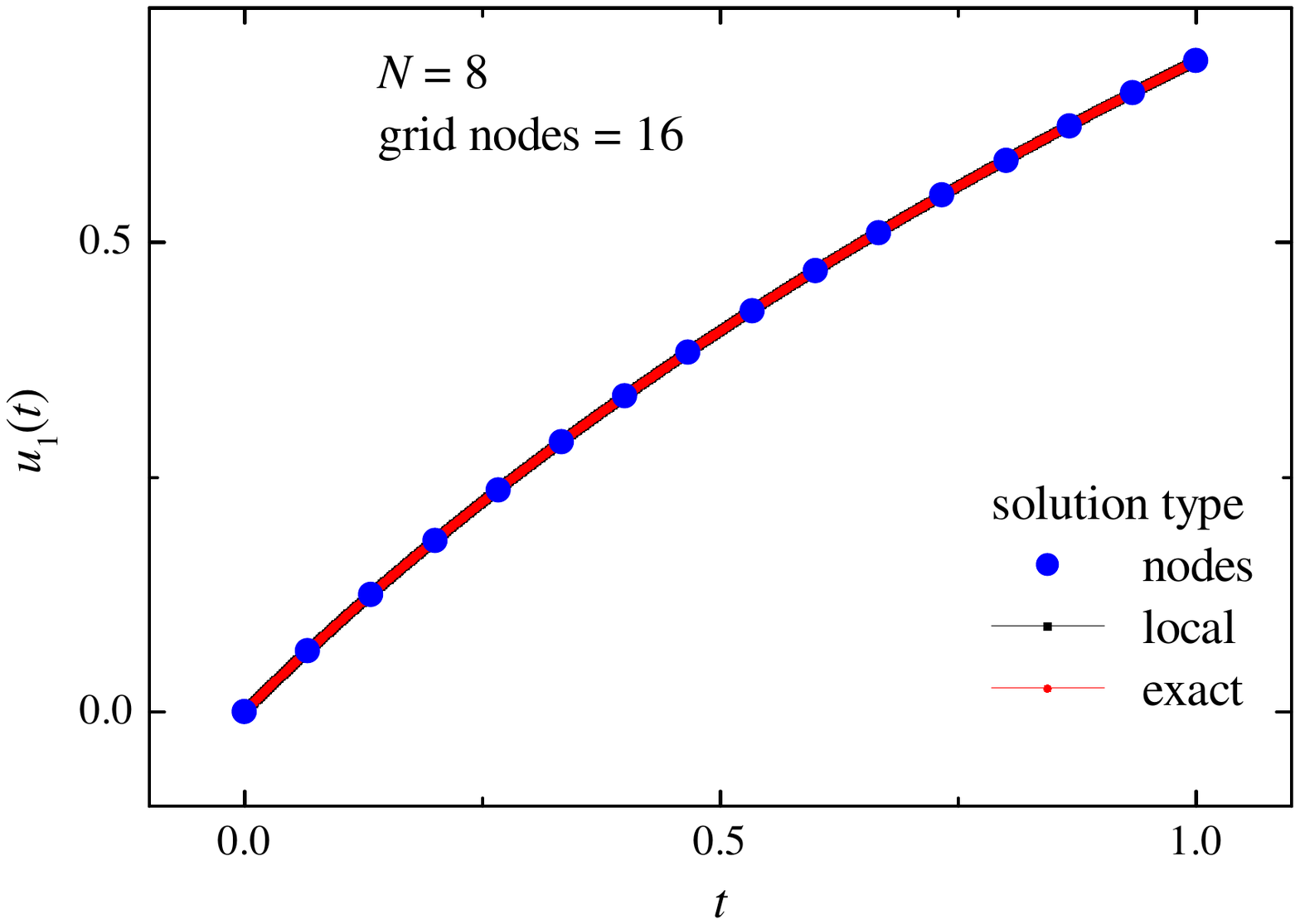}
\vspace{-8mm}\caption{\label{fig:third_order_ivp_ex_2_sols:b1}}
\end{subfigure}
\begin{subfigure}{0.23\textwidth}
\includegraphics[width=\textwidth]{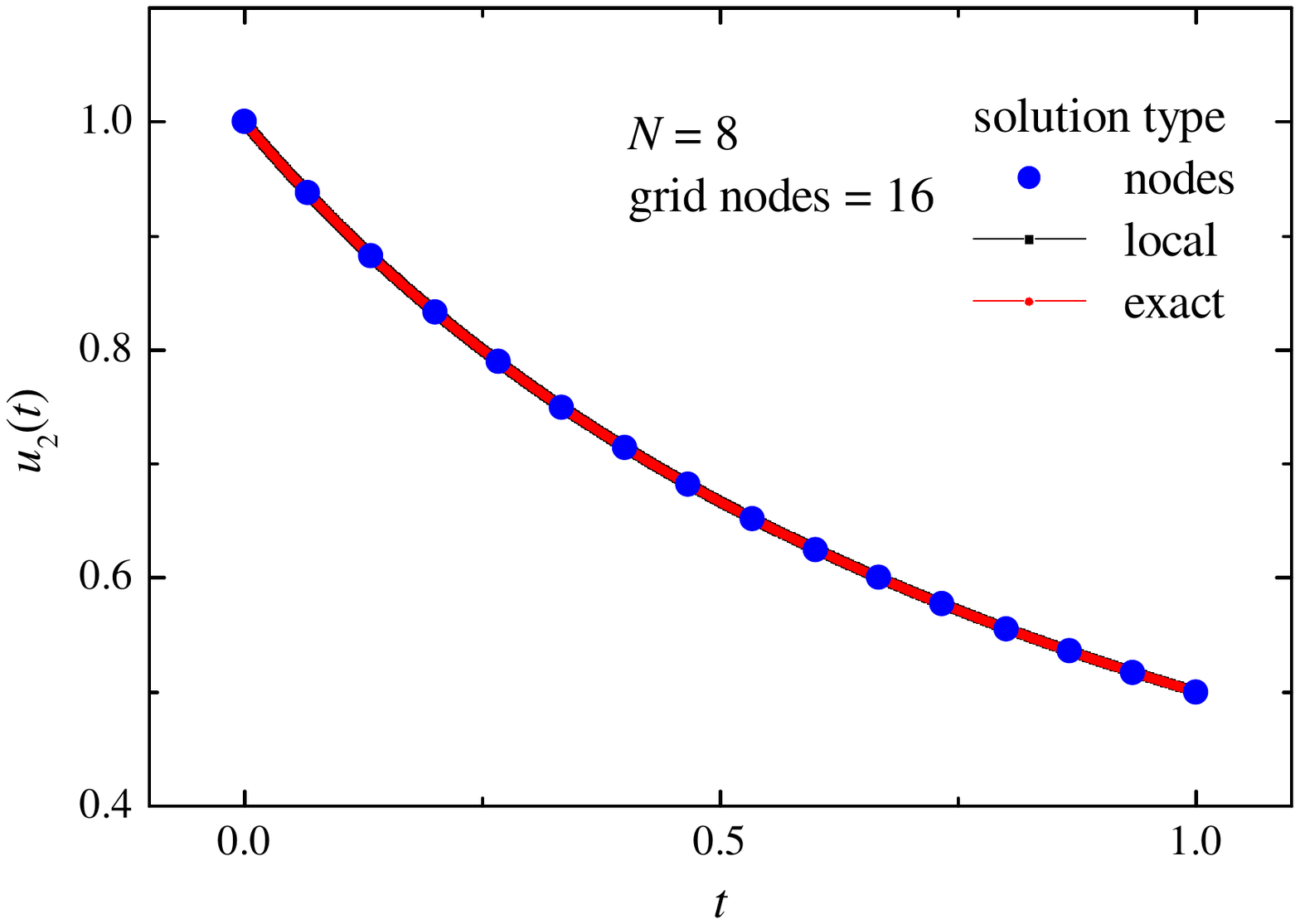}
\vspace{-8mm}\caption{\label{fig:third_order_ivp_ex_2_sols:b2}}
\end{subfigure}
\begin{subfigure}{0.23\textwidth}
\includegraphics[width=\textwidth]{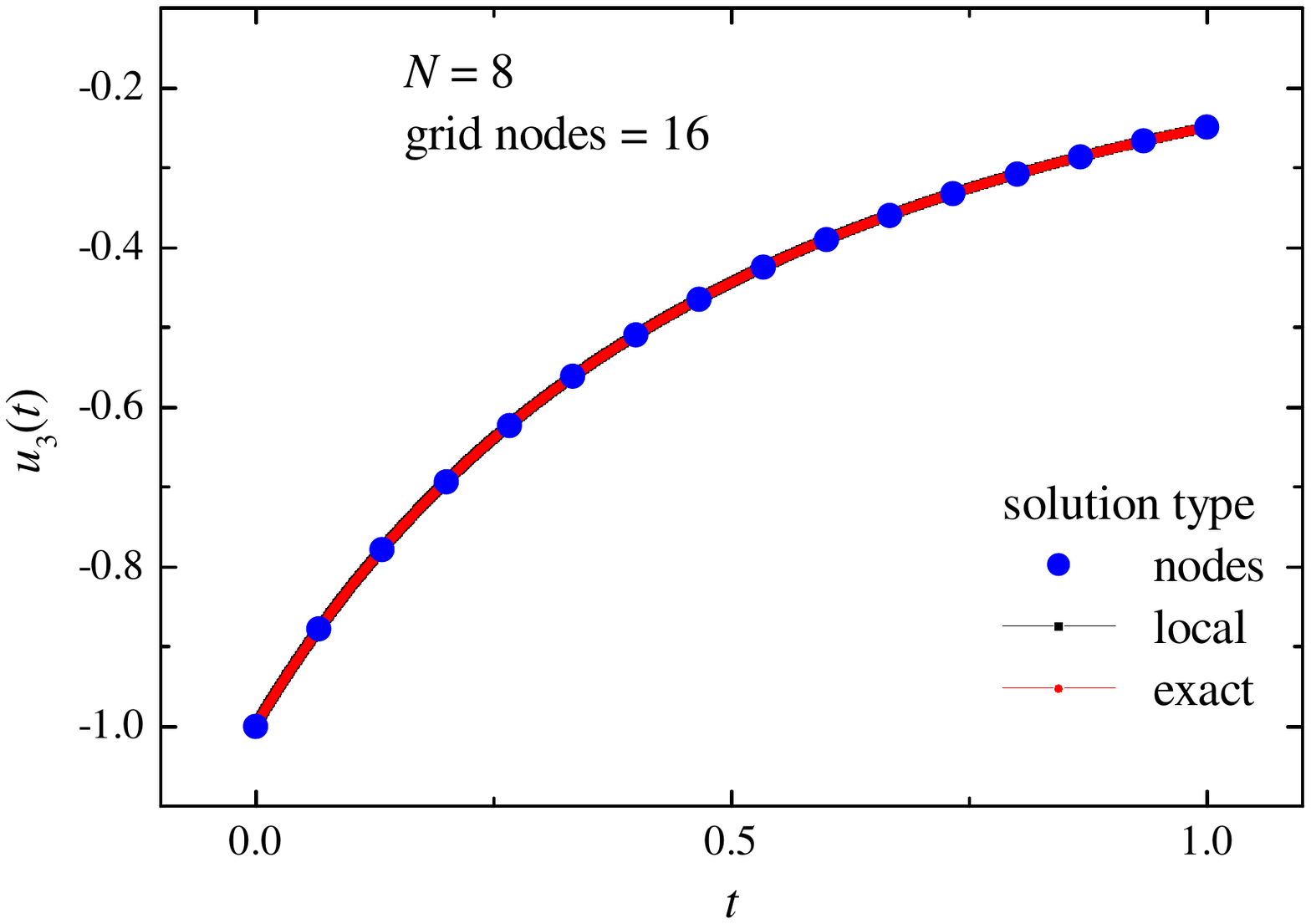}
\vspace{-8mm}\caption{\label{fig:third_order_ivp_ex_2_sols:b3}}
\end{subfigure}
\begin{subfigure}{0.23\textwidth}
\includegraphics[width=\textwidth]{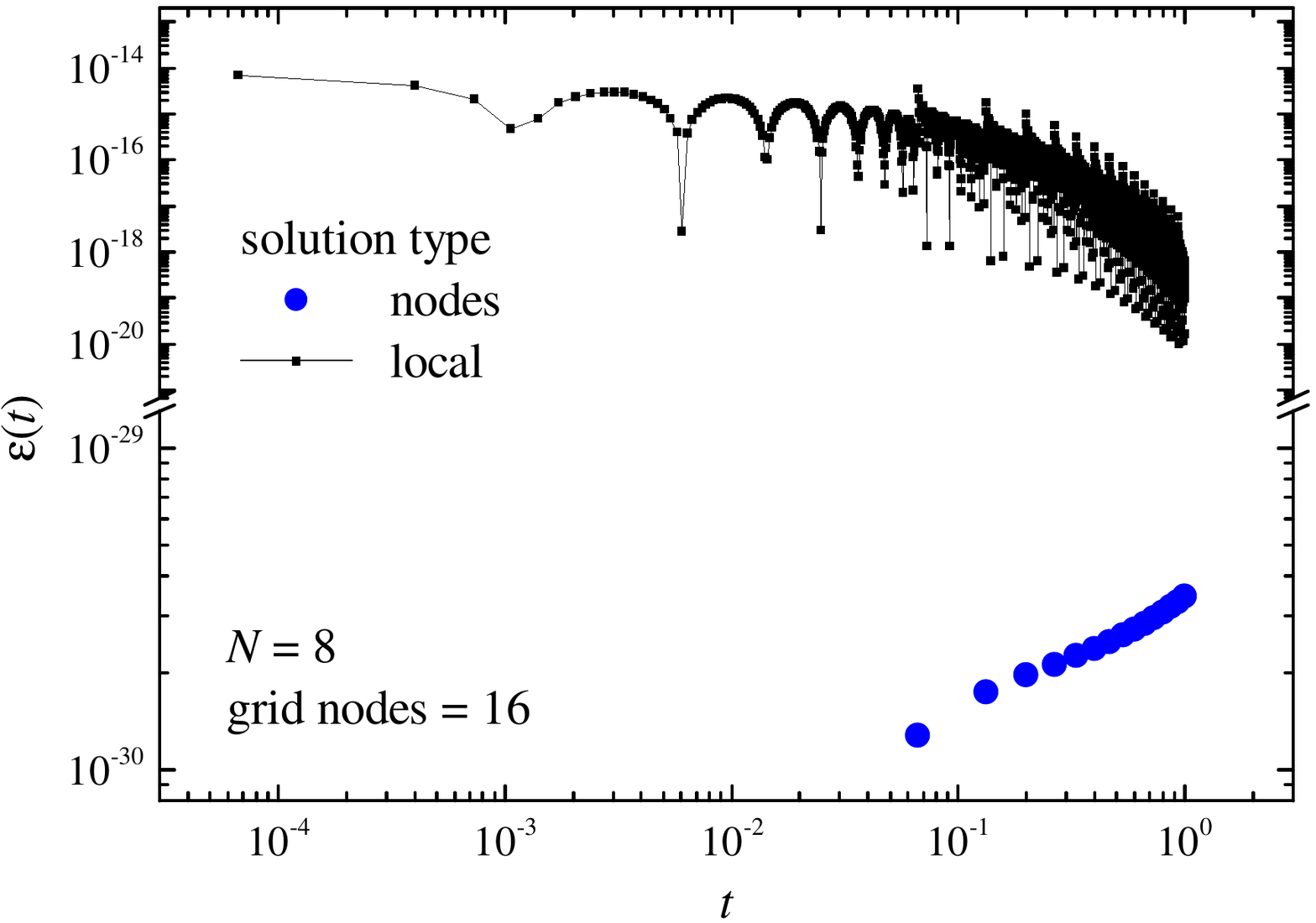}
\vspace{-8mm}\caption{\label{fig:third_order_ivp_ex_2_sols:b4}}
\end{subfigure}\\
\begin{subfigure}{0.23\textwidth}
\includegraphics[width=\textwidth]{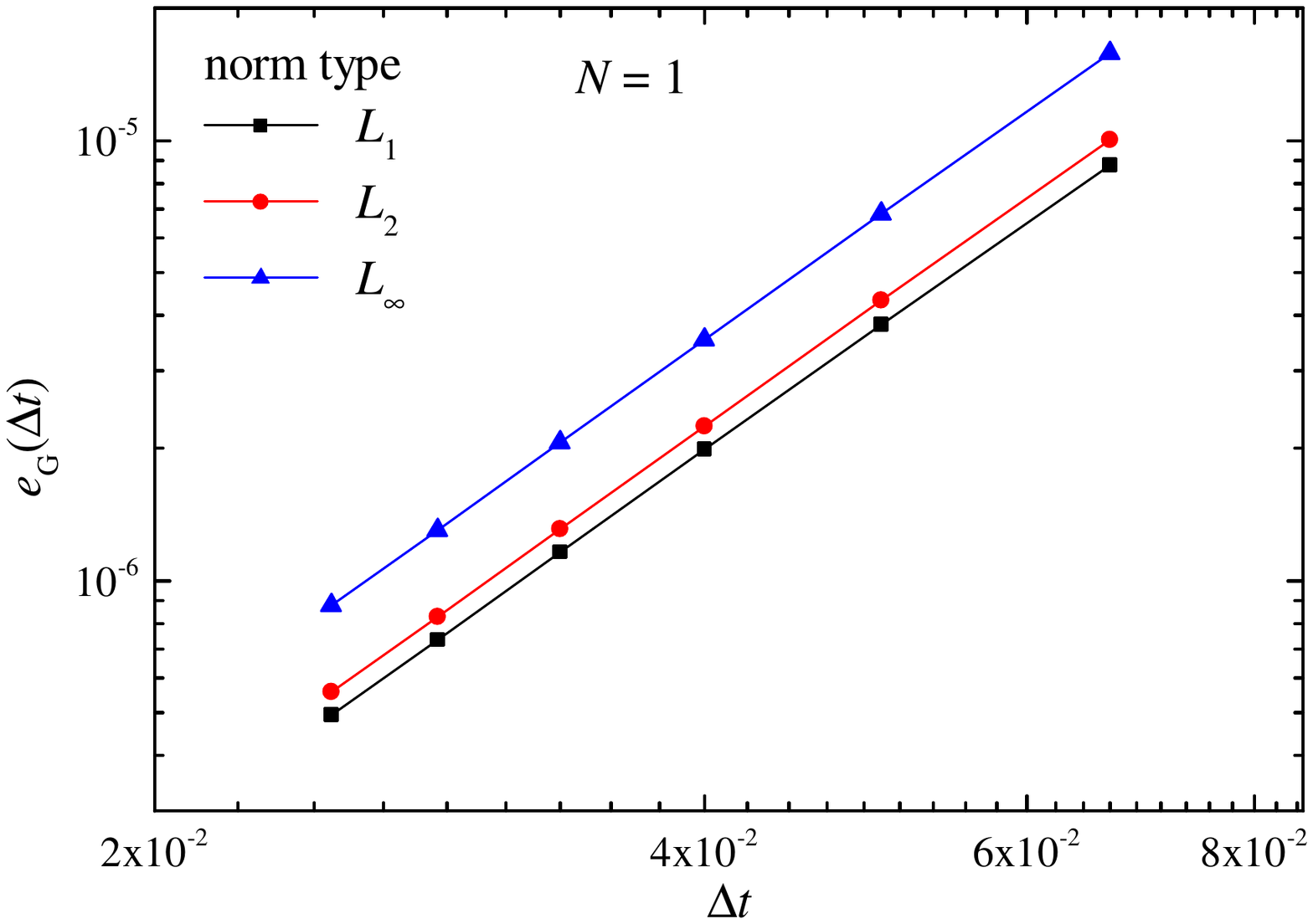}
\vspace{-8mm}\caption{\label{fig:third_order_ivp_ex_2_sols:c1}}
\end{subfigure}
\begin{subfigure}{0.23\textwidth}
\includegraphics[width=\textwidth]{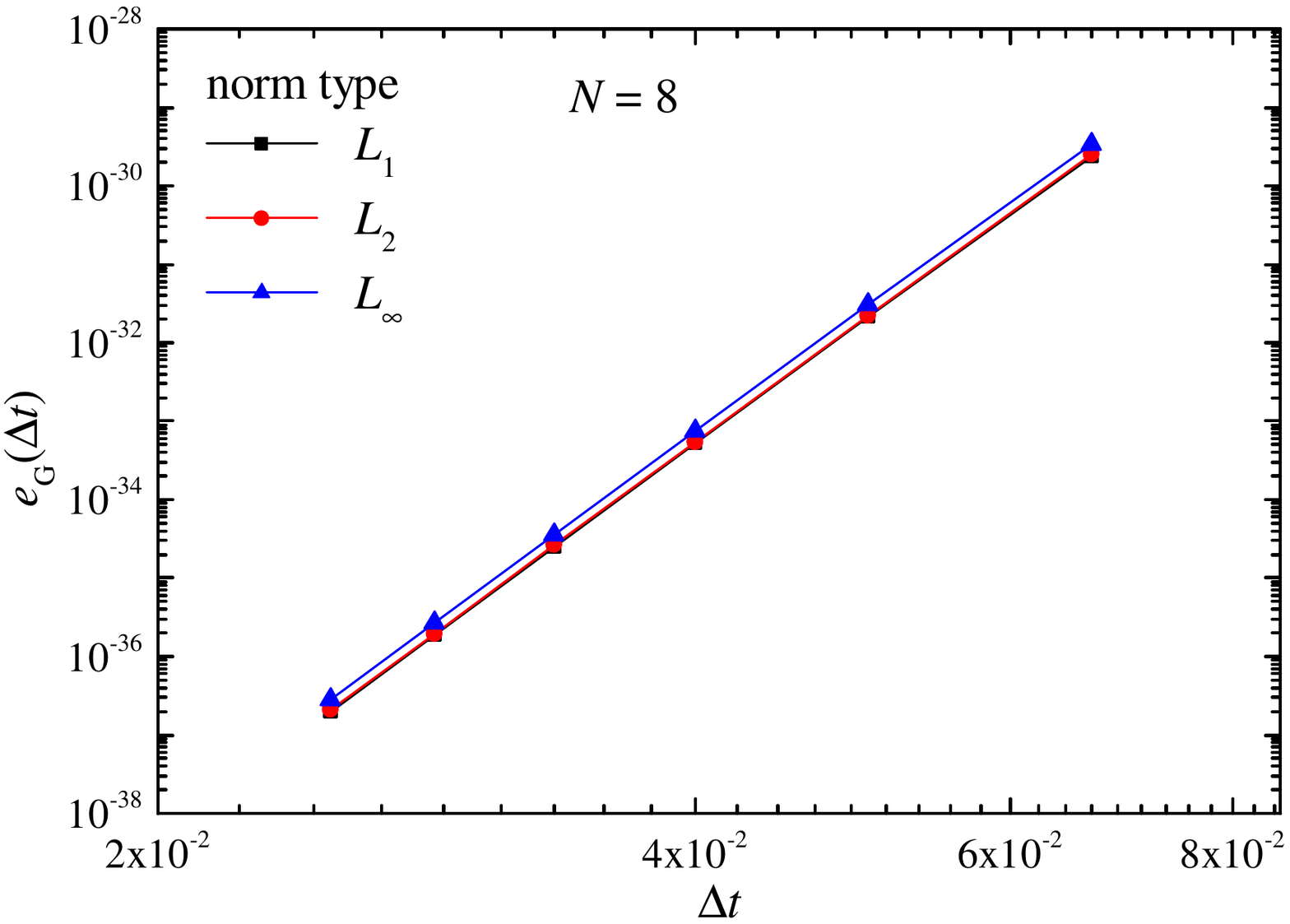}
\vspace{-8mm}\caption{\label{fig:third_order_ivp_ex_2_sols:c2}}
\end{subfigure}
\begin{subfigure}{0.23\textwidth}
\includegraphics[width=\textwidth]{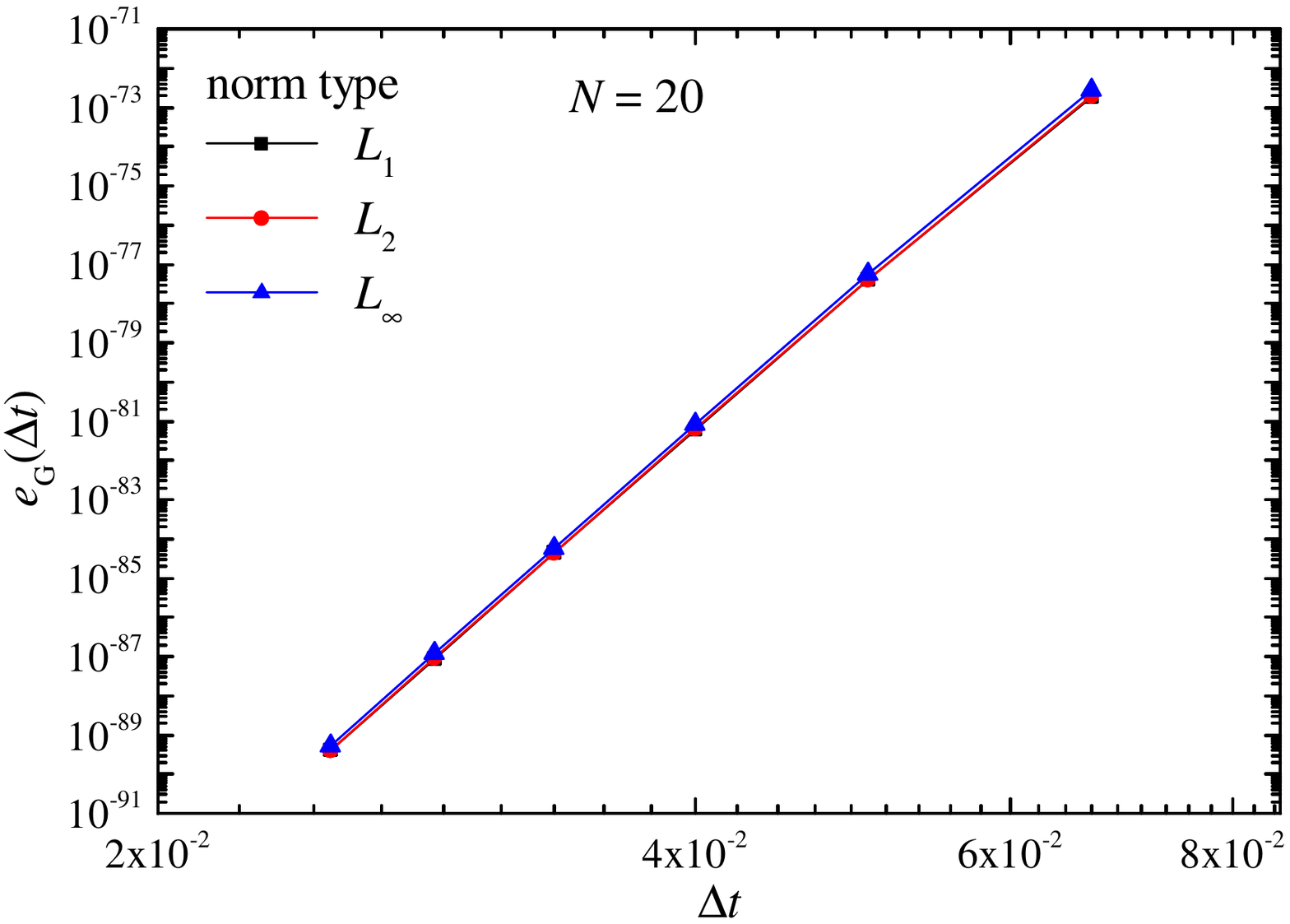}
\vspace{-8mm}\caption{\label{fig:third_order_ivp_ex_2_sols:c3}}
\end{subfigure}
\begin{subfigure}{0.23\textwidth}
\includegraphics[width=\textwidth]{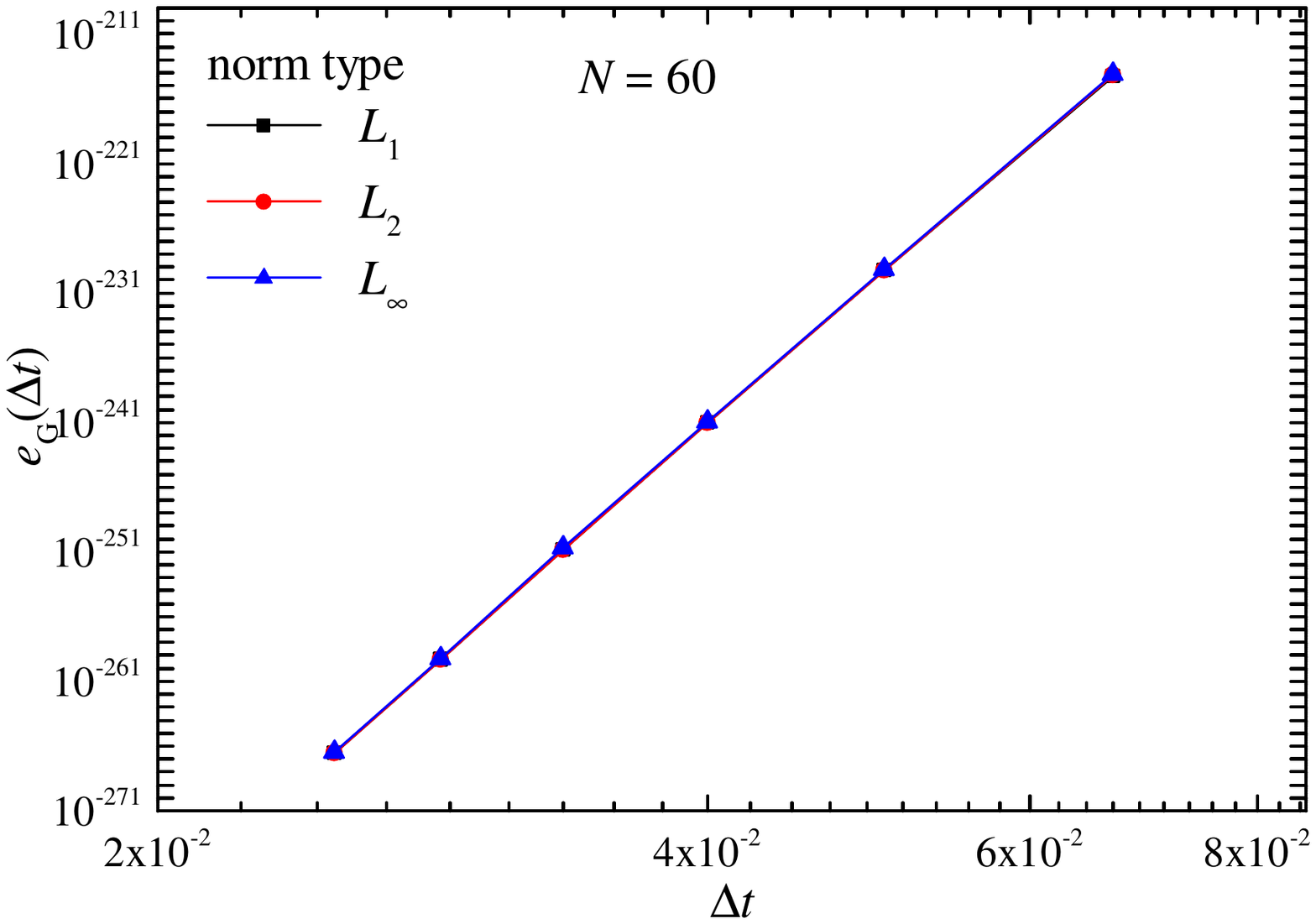}
\vspace{-8mm}\caption{\label{fig:third_order_ivp_ex_2_sols:c4}}
\end{subfigure}\\
\begin{subfigure}{0.23\textwidth}
\includegraphics[width=\textwidth]{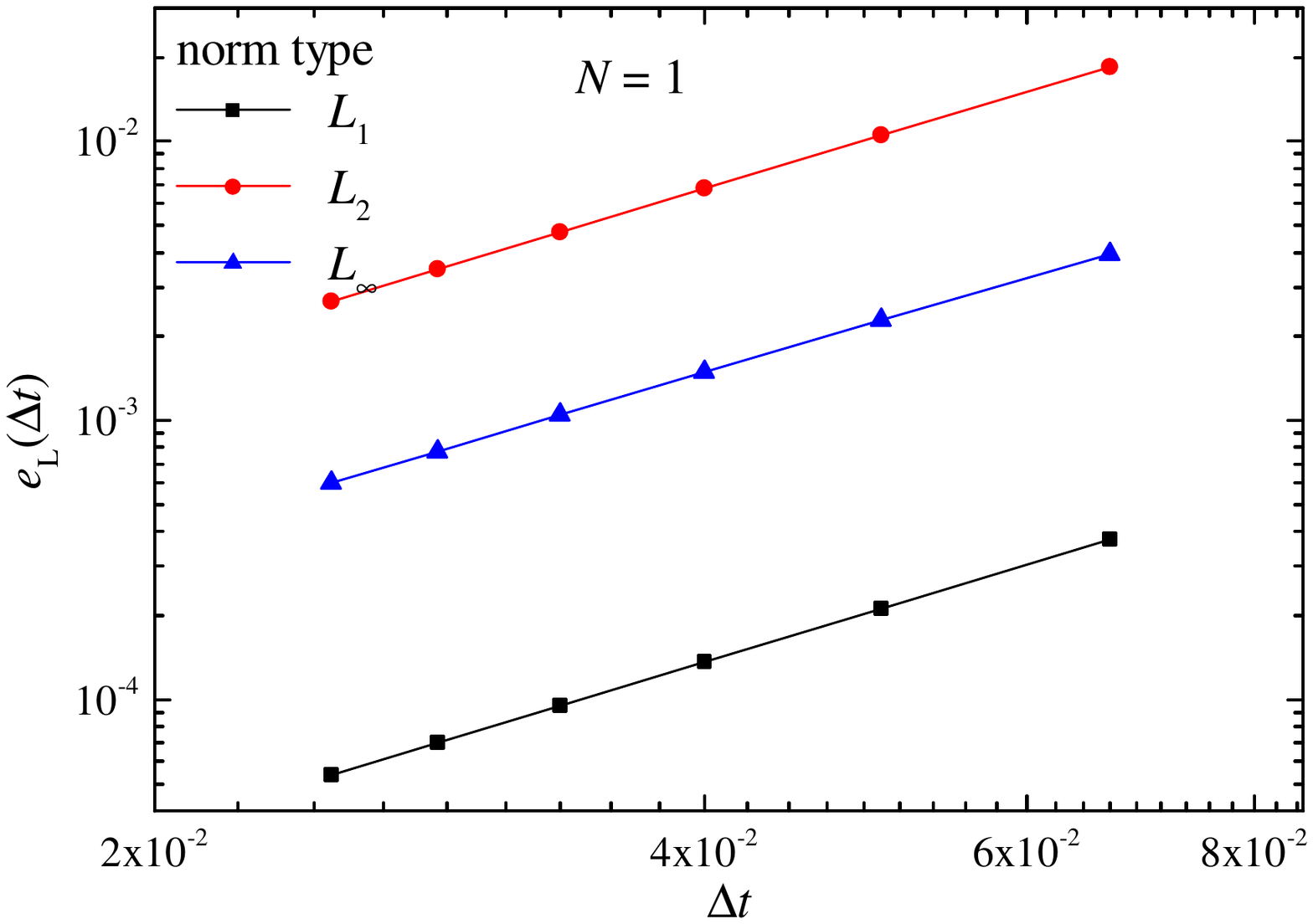}
\vspace{-8mm}\caption{\label{fig:third_order_ivp_ex_2_sols:d1}}
\end{subfigure}
\begin{subfigure}{0.23\textwidth}
\includegraphics[width=\textwidth]{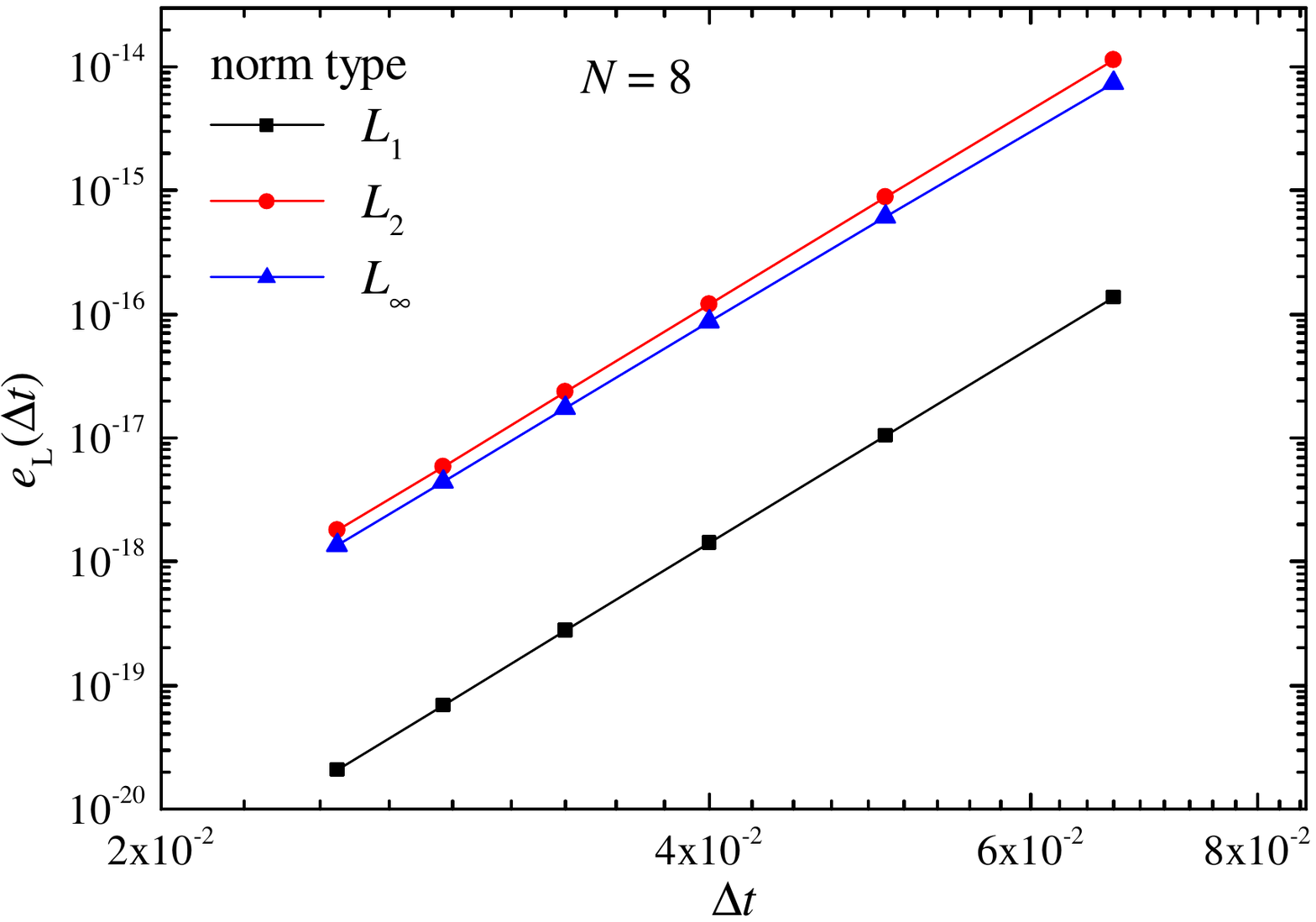}
\vspace{-8mm}\caption{\label{fig:third_order_ivp_ex_2_sols:d2}}
\end{subfigure}
\begin{subfigure}{0.23\textwidth}
\includegraphics[width=\textwidth]{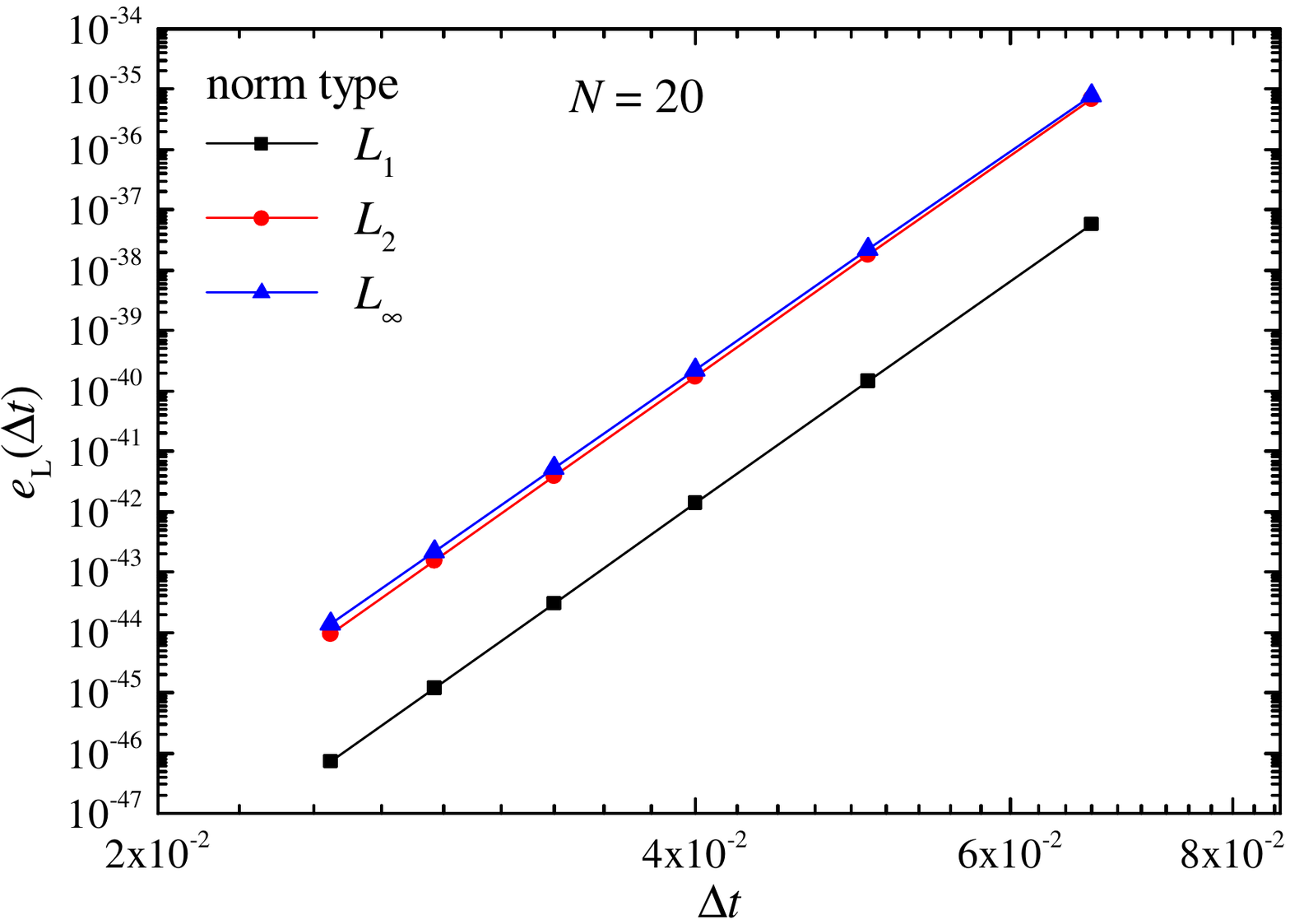}
\vspace{-8mm}\caption{\label{fig:third_order_ivp_ex_2_sols:d3}}
\end{subfigure}
\begin{subfigure}{0.23\textwidth}
\includegraphics[width=\textwidth]{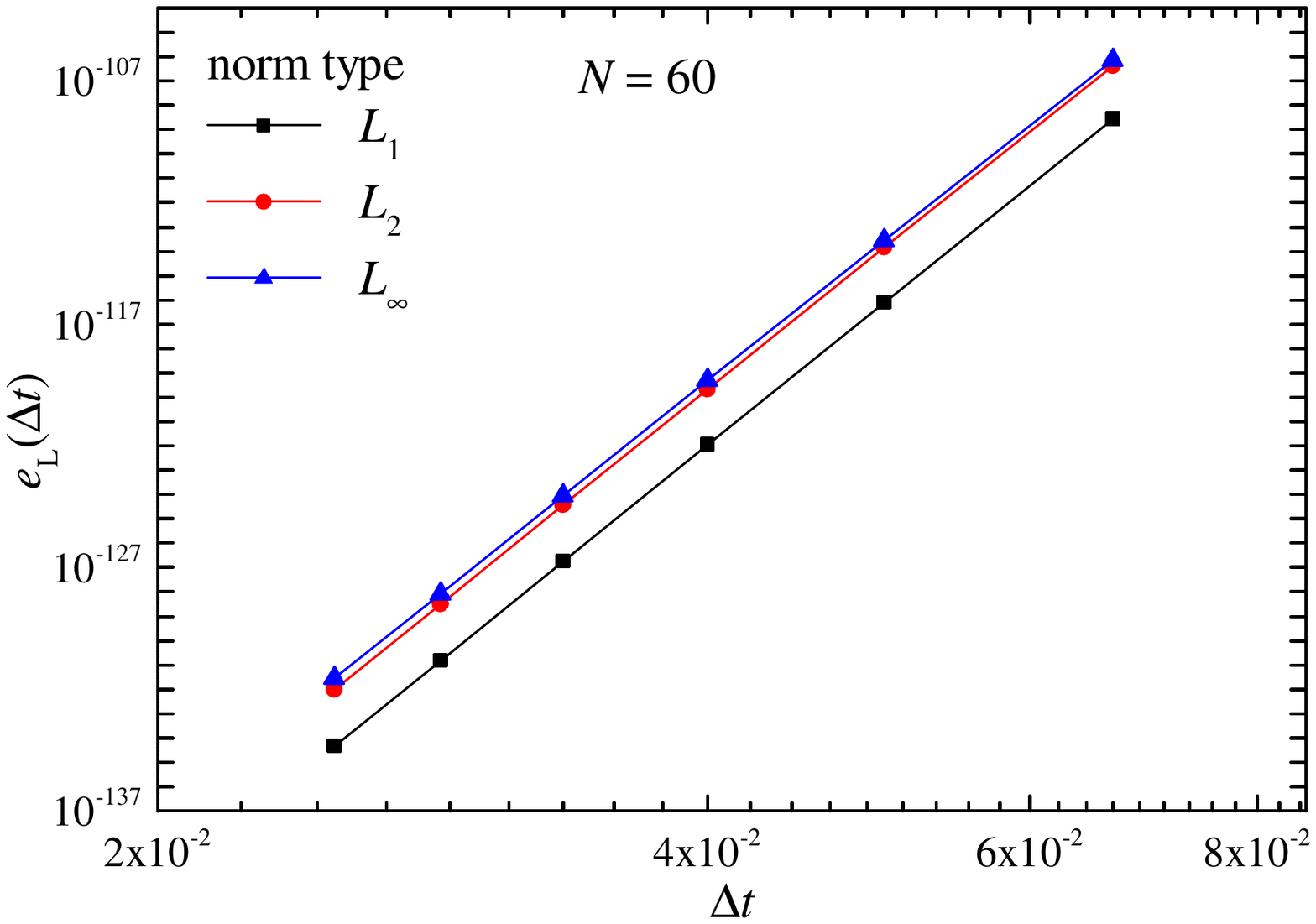}
\vspace{-8mm}\caption{\label{fig:third_order_ivp_ex_2_sols:d4}}
\end{subfigure}\\
\caption{%
Numerical solution of the problem (\ref{eq:third_order_ivp_ex_2}). Comparison of the solution at nodes $\mathbf{u}_{n}$, the local solution $\mathbf{u}_{L}(t)$ and the exact solution $\mathbf{u}^{\rm ex}(t)$ for components $u_{1}$ (\subref{fig:third_order_ivp_ex_2_sols:a1}, \subref{fig:third_order_ivp_ex_2_sols:b1}), $u_{2}$ (\subref{fig:third_order_ivp_ex_2_sols:a2}, \subref{fig:third_order_ivp_ex_2_sols:b2}) and $u_{3}$ (\subref{fig:third_order_ivp_ex_2_sols:a3}, \subref{fig:third_order_ivp_ex_2_sols:b3}), the error $\varepsilon(t)$ (\subref{fig:third_order_ivp_ex_2_sols:a4}, \subref{fig:third_order_ivp_ex_2_sols:b4}), obtained using polynomials with degrees $N = 1$ (\subref{fig:third_order_ivp_ex_2_sols:a1}, \subref{fig:third_order_ivp_ex_2_sols:a2}, \subref{fig:third_order_ivp_ex_2_sols:a3}, \subref{fig:third_order_ivp_ex_2_sols:a4}) and $N = 8$ (\subref{fig:third_order_ivp_ex_2_sols:b1}, \subref{fig:third_order_ivp_ex_2_sols:b2}, \subref{fig:third_order_ivp_ex_2_sols:b3}, \subref{fig:third_order_ivp_ex_2_sols:b4}). Dependence of the global error for the solution at nodes $e_{G}$ (\subref{fig:third_order_ivp_ex_2_sols:c1}, \subref{fig:third_order_ivp_ex_2_sols:c2}, \subref{fig:third_order_ivp_ex_2_sols:c3}, \subref{fig:third_order_ivp_ex_2_sols:c4}) and the local solution $e_{L}$ (\subref{fig:third_order_ivp_ex_2_sols:d1}, \subref{fig:third_order_ivp_ex_2_sols:d2}, \subref{fig:third_order_ivp_ex_2_sols:d3}, \subref{fig:third_order_ivp_ex_2_sols:d4}) on the discretization step $\mathrm{\Delta}t$, obtained in the norms $L_{1}$, $L_{2}$ and $L_{\infty}$, obtained using polynomials with degrees $N = 1$ (\subref{fig:third_order_ivp_ex_2_sols:c1}, \subref{fig:third_order_ivp_ex_2_sols:d1}), $N = 8$ (\subref{fig:third_order_ivp_ex_2_sols:c2}, \subref{fig:third_order_ivp_ex_2_sols:d2}), $N = 20$ (\subref{fig:third_order_ivp_ex_2_sols:c3}, \subref{fig:third_order_ivp_ex_2_sols:d3}) and $N = 60$ (\subref{fig:third_order_ivp_ex_2_sols:c4}, \subref{fig:third_order_ivp_ex_2_sols:d4}).
}
\label{fig:third_order_ivp_ex_2_sols}
\end{figure} 
\begin{table}[h!]
\centering
\caption{%
Orders of convergence $p_{L_{1}}$, $p_{L_{2}}$, $p_{L_{\infty}}$ of the ADER-DG numerical method with local DG predictor calculated in norms $L_{1}$, $L_{2}$, $L_{\infty}$ for the problem (\ref{eq:third_order_ivp_ex_2}); $N$ is the degree of the basis polynomials $\varphi_{p}$ used for the DG representations. The orders of convergence $p$ for the solution at nodes $\mathbf{u}_{n}$ are compared to the theoretical value $p_{\rm (G)} = 2N+1$, the orders of convergence $p$ for the local solution $\mathbf{u}_{L}$ represented between the nodes are compared to the theoretical value $p_{\rm (L)} = N+1$.
}
\label{tab:third_order_ivp_ex_2_sols}
\begin{tabular}{@{}lllllllll@{}}
\toprule
& \multicolumn{4}{c}{order $p$ for solution at nodes} & \multicolumn{4}{c}{order $p$ for local solution} \\
\midrule
$N$ & $p_{L_{1}}$ & $p_{L_{2}}$ & $p_{L_{\infty}}$ & theor. & $p_{L_{1}}$ & $p_{L_{2}}$ & $p_{L_{\infty}}$ & theor. \\
\midrule
$1$	&	$2.93$	&	$2.95$	&	$2.95$		&	$3$	&	$1.99$	&	$1.97$	&	$1.92$		&	$2$	\\
$2$	&	$4.91$	&	$4.92$	&	$4.92$		&	$5$	&	$2.99$	&	$2.97$	&	$2.90$		&	$3$	\\
$3$	&	$6.87$	&	$6.88$	&	$6.88$		&	$7$	&	$3.99$	&	$3.96$	&	$3.88$		&	$4$	\\
$4$	&	$8.83$	&	$8.84$	&	$8.84$		&	$9$	&	$4.99$	&	$4.96$	&	$4.86$		&	$5$	\\
$5$	&	$10.78$	&	$10.79$	&	$10.80$		&	$11$	&	$5.98$	&	$5.96$	&	$5.84$		&	$6$	\\
$6$	&	$12.73$	&	$12.75$	&	$12.75$		&	$13$	&	$6.98$	&	$6.95$	&	$6.82$		&	$7$	\\
$7$	&	$14.68$	&	$14.70$	&	$14.70$		&	$15$	&	$7.98$	&	$7.95$	&	$7.80$		&	$8$	\\
$8$	&	$16.63$	&	$16.64$	&	$16.65$		&	$17$	&	$8.98$	&	$8.94$	&	$8.78$		&	$9$	\\
$9$	&	$18.57$	&	$18.58$	&	$18.59$		&	$19$	&	$9.97$	&	$9.93$	&	$9.76$		&	$10$	\\
$10$	&	$20.50$	&	$20.52$	&	$20.52$		&	$21$	&	$10.97$	&	$10.93$	&	$10.74$		&	$11$	\\
\midrule
$15$	&	$30.07$	&	$30.08$	&	$30.09$		&	$31$	&	$15.95$	&	$15.89$	&	$15.63$		&	$16$	\\
$20$	&	$39.25$	&	$39.27$	&	$39.27$		&	$41$	&	$20.92$	&	$20.83$	&	$20.53$		&	$21$	\\
$25$	&	$48.24$	&	$48.26$	&	$48.26$		&	$51$	&	$25.88$	&	$25.77$	&	$25.43$		&	$26$	\\
$30$	&	$59.00$	&	$59.02$	&	$59.02$		&	$61$	&	$30.84$	&	$30.70$	&	$30.33$		&	$31$	\\
$35$	&	$69.88$	&	$69.90$	&	$69.90$		&	$71$	&	$35.79$	&	$35.62$	&	$35.23$		&	$36$	\\
$40$	&	$80.60$	&	$80.62$	&	$80.63$		&	$81$	&	$40.74$	&	$40.54$	&	$40.12$		&	$41$	\\
$45$	&	$91.56$	&	$91.57$	&	$91.58$		&	$91$	&	$45.68$	&	$45.45$	&	$45.02$		&	$46$	\\
$50$	&	$102.18$	&	$102.20$	&	$102.21$		&	$101$	&	$50.62$	&	$50.36$	&	$49.92$		&	$51$	\\
$55$	&	$113.97$	&	$113.99$	&	$114.00$		&	$111$	&	$55.55$	&	$55.26$	&	$54.82$		&	$56$	\\
$60$	&	$122.93$	&	$122.95$	&	$122.96$		&	$121$	&	$60.48$	&	$60.17$	&	$59.72$		&	$61$	\\
\bottomrule
\end{tabular}
\end{table} 

The previous example (\ref{eq:third_order_ivp_ex_1}) was related to solving the IVP for a third order linear differential equation. This and the next examples will be related to the solution of third-order non-linear equations~\cite{dg_ivp_ode_5}. This example for IVP third-order ODE was chosen in the following form:
\begin{equation}
\dddot{x} = \frac{4}{(1+t)^{3}} -2\exp(-3x);\qquad
x(0) = 0;\ \dot{x}(0) = 1;\ \ddot{x}(0) = -1;\quad 0 \leqslant t \leqslant 1;
\end{equation}
with exact analytical solution of this problem having the form $x(t) = \ln(1+t)$. The IVP for a first-order non-linear ODE system corresponding to this problem has the following form:
\begin{equation}\label{eq:third_order_ivp_ex_2}
\begin{split}
&\frac{du_{1}}{dt} = u_{2};\quad \frac{du_{2}}{dt} = u_{3};\quad
\frac{du_{3}}{dt} = \frac{4}{(1+t)^{3}} -2\exp(-3u_{1});\\
&u_{1}(0) = 0;\quad u_{2}(0) = 1;\quad u_{3}(0) = -1;
\end{split}
\end{equation}
where $\mathbf{u} = [u_{1},\, u_{2},\, u_{3}]^{T} = [x,\, \dot{x},\, \ddot{x}]^{T}$ is a desired three component vector function. The solution definition domain $\Omega$ was chosen as $0 \leqslant t \leqslant 1$. The full three component exact analytical solution of this problem has the form
\begin{equation}
\mathbf{u}^{\rm ex} = \left[
\ln(1+t),\quad
\dfrac{1}{1+t},\quad
-\dfrac{1}{(1+t)^{2}}
\right]^{T}.
\end{equation}

The study of the convergence of the numerical solution was carried out on the basis of the analysis of the solution for the set of $6$ different uniform discretizations of the domain of definition $\Omega$ of the desired function, with the number of grid nodes $L = 16$, $21$, $26$, $31$, $36$, $41$; note that the zero node, in which the initial condition was defined, also refers to grid nodes, so the discretization step was determined by the expression $\Delta t = 2\pi/(L-1)$.

The results of the study of the application of the numerical method for this problem are presented in Fig.~\ref{fig:third_order_ivp_ex_2_sols} and in Table~\ref{tab:third_order_ivp_ex_2_sols}. The error $\varepsilon(t)$ for $N > 1$ shows that the error scales for the local solution and the solution at nodes differ by several orders of magnitude, reaching values of $10$--$12$ and $100$--$120$ orders of magnitude for a degree of $N = 8$ and $N = 60$, respectively. The local error $\varepsilon(t)$ of the solution at nodes shows an approximately linear growth with $t$. The calculated convergence orders $p_{L_{1}}$, $p_{L_{2}}$, $p_{L_{\infty}}$ for various functional norms are presented in Table~\ref{tab:third_order_ivp_ex_2_sols}. The convergence orders $p$ for the solution at nodes are compared to the theoretical value $p_{\rm (G)} = 2N+1$, the convergence orders $p$ for the local solution represented between the nodes are compared to the theoretical value $p_{\rm (L)} = N+1$, in accordance with (\ref{eq:expect_orders}). The main part of the obtained convergence orders $p$ sufficiently correspond to the expected theoretical values $p_{\rm (G)}$ and $p_{\rm (L)}$. In the range of polynomial degrees $N = 50, \ldots, 60$, the empirical values of the convergence orders $p$ for the solution at nodes exceed the expected values $p_{\rm (G)}$ by more than $1$. In the case of a local solution $\mathbf{u}_{L}$, the convergence orders $p$ correspond to the expected values $p_{\rm (L)}$.

As a result, it becomes clear that the solution at nodes $\mathbf{u}_{n}$ obtained by the ADER-DG numerical method with a local DG predictor demonstrates the classical superconvergence $2N+1$ expected from such numerical methods, at least in problem (\ref{eq:third_order_ivp_ex_2}). The local solution $\mathbf{u}_{L}$ obtained by the DG predictor demonstrates the classical convergence characteristic of DG methods. 

\subsubsection{Example 6}

\begin{figure}[h!]
\captionsetup[subfigure]{%
position=bottom,
textfont=normalfont,
singlelinecheck=off,
justification=raggedright,
font=tiny
}
\centering
\begin{subfigure}{0.23\textwidth}
\includegraphics[width=\textwidth]{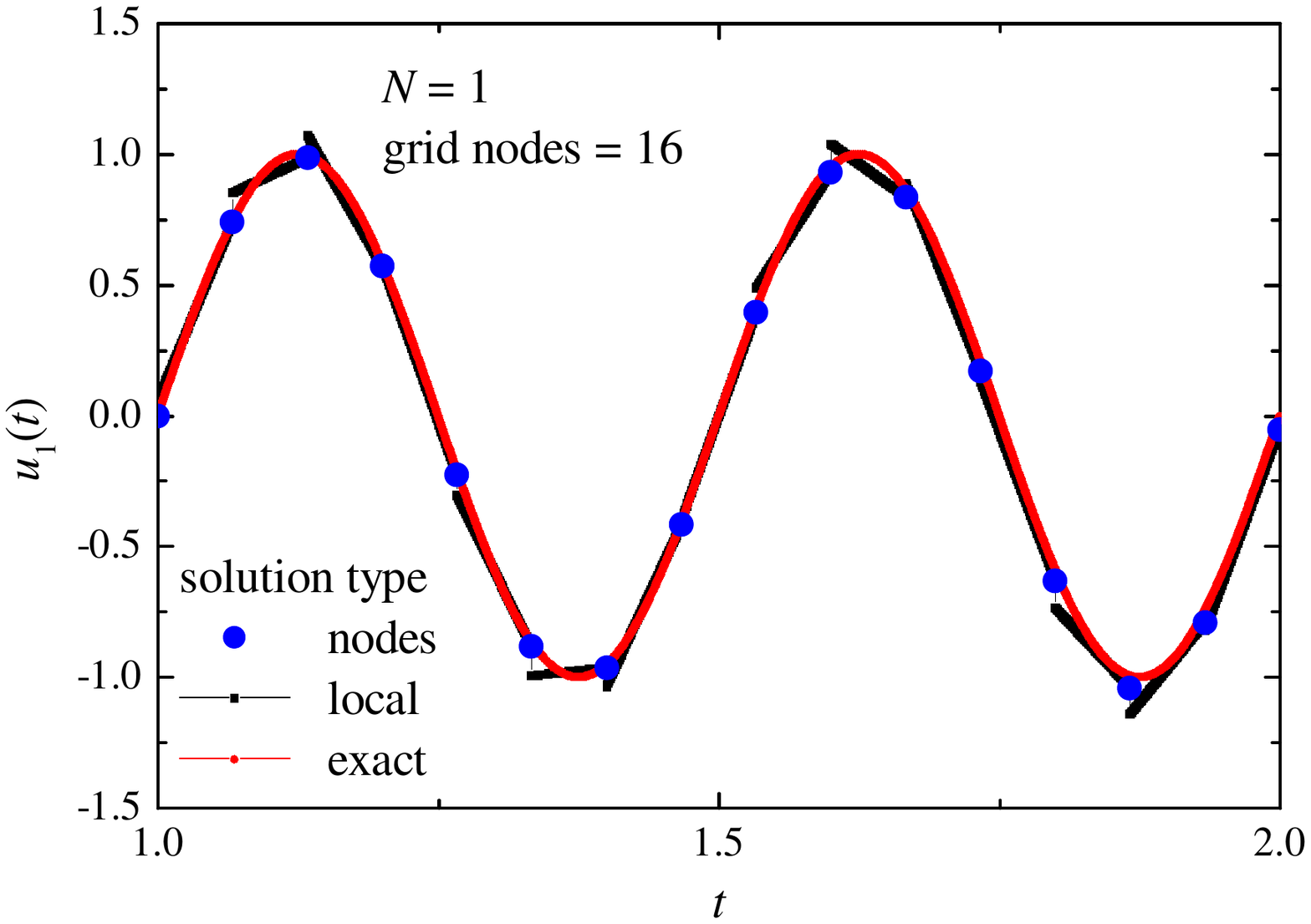}
\vspace{-8mm}\caption{\label{fig:third_order_ivp_ex_3_sols:a1}}
\end{subfigure}
\begin{subfigure}{0.23\textwidth}
\includegraphics[width=\textwidth]{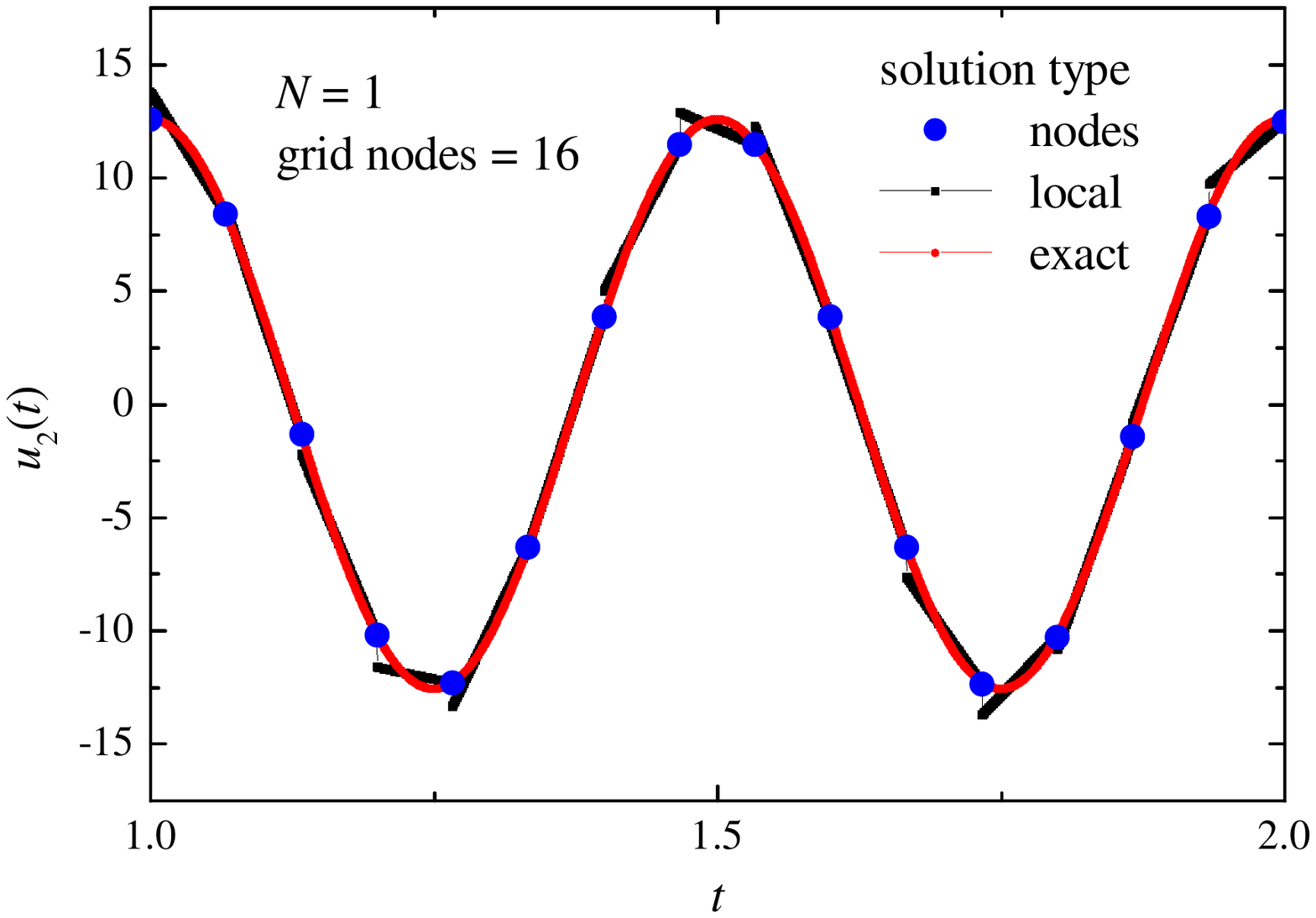}
\vspace{-8mm}\caption{\label{fig:third_order_ivp_ex_3_sols:a2}}
\end{subfigure}
\begin{subfigure}{0.23\textwidth}
\includegraphics[width=\textwidth]{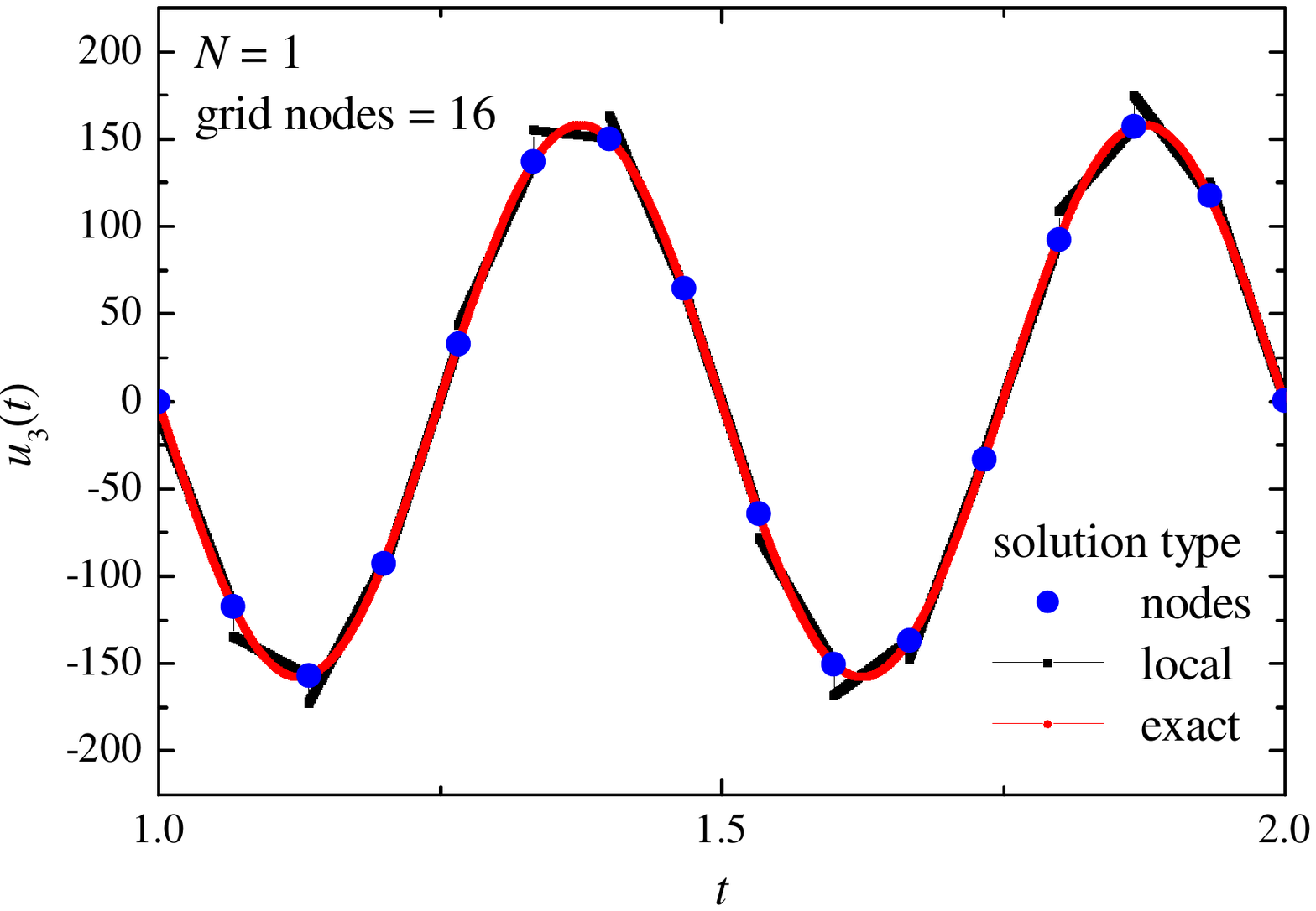}
\vspace{-8mm}\caption{\label{fig:third_order_ivp_ex_3_sols:a3}}
\end{subfigure}
\begin{subfigure}{0.23\textwidth}
\includegraphics[width=\textwidth]{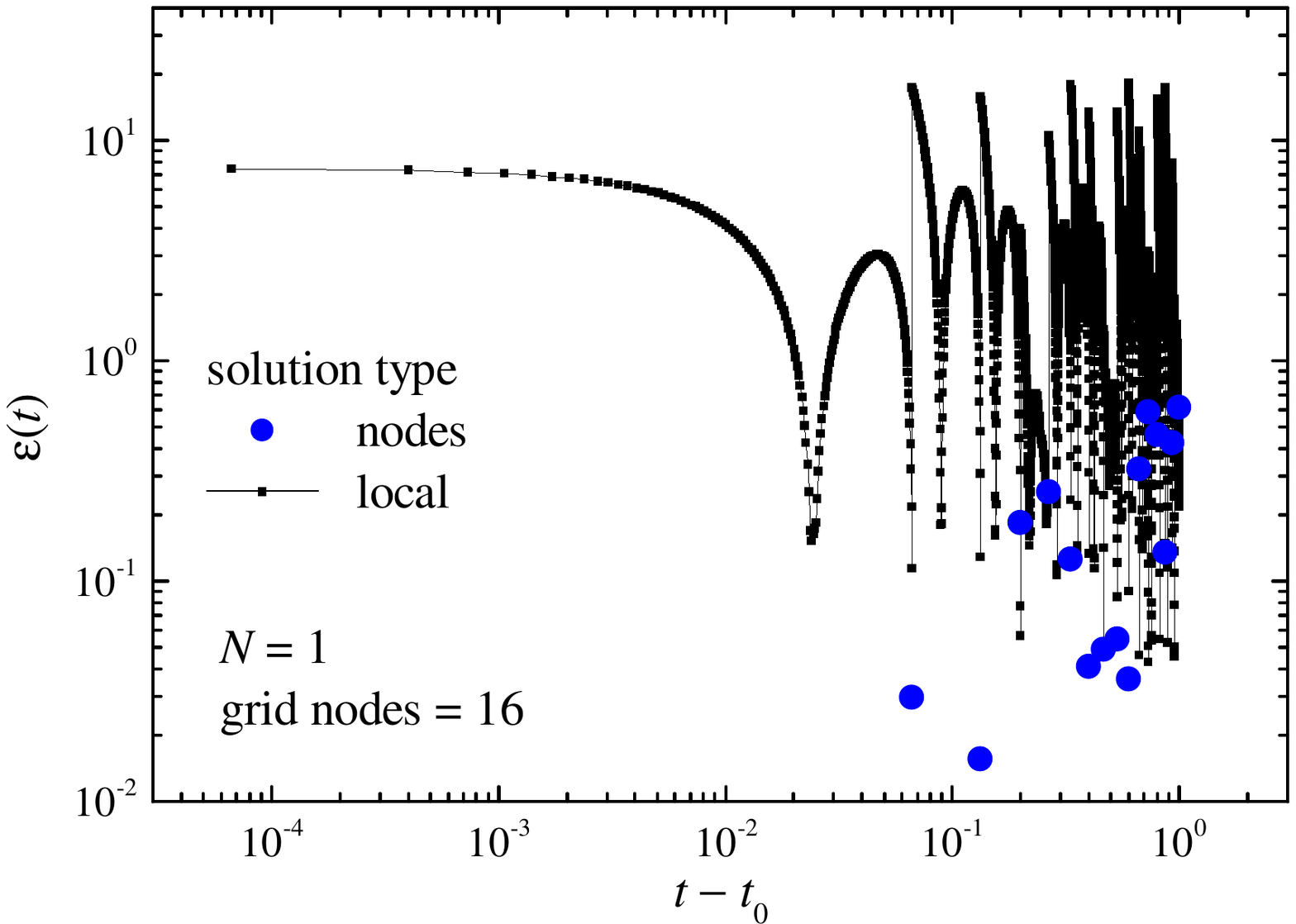}
\vspace{-8mm}\caption{\label{fig:third_order_ivp_ex_3_sols:a4}}
\end{subfigure}\\
\begin{subfigure}{0.23\textwidth}
\includegraphics[width=\textwidth]{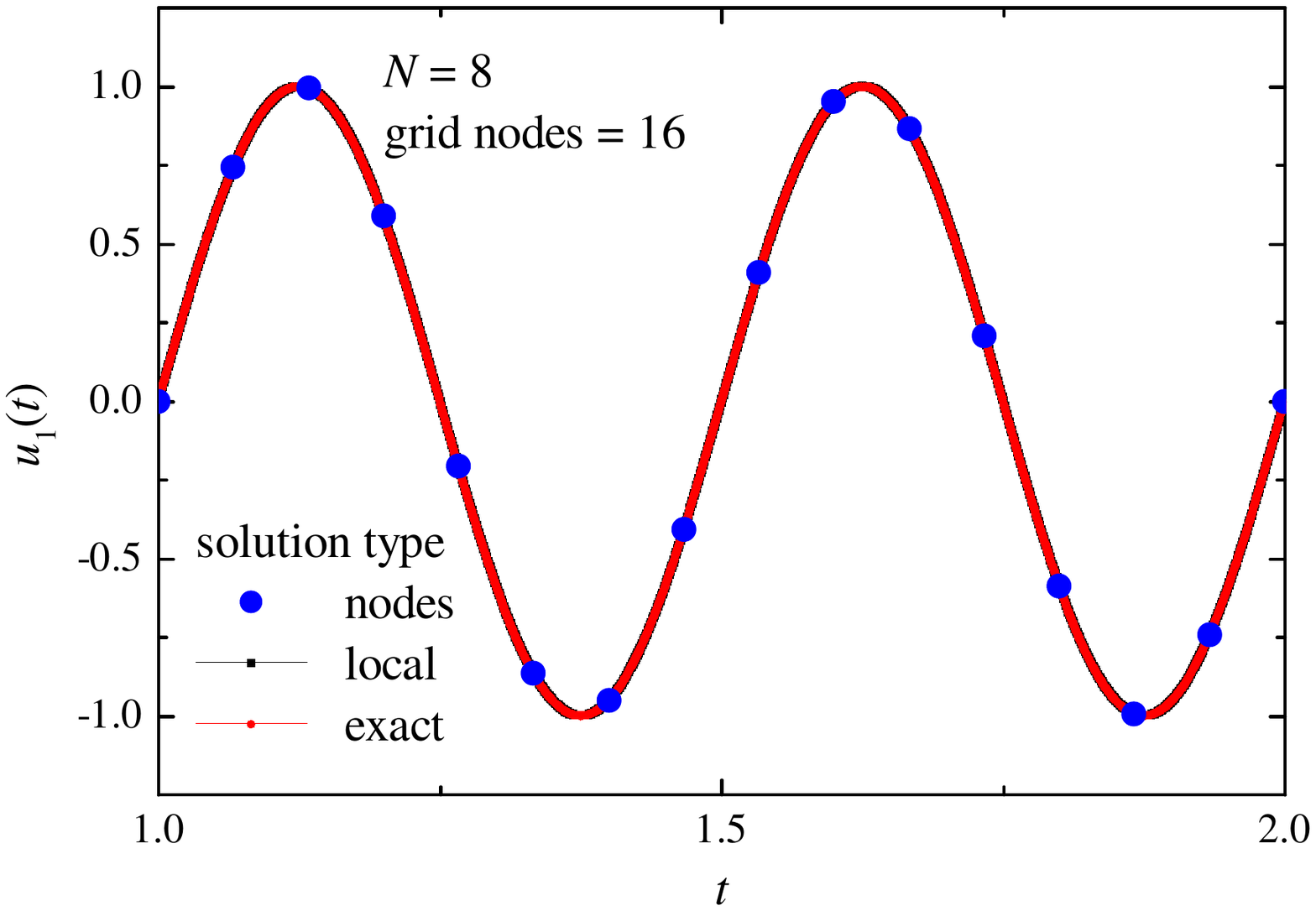}
\vspace{-8mm}\caption{\label{fig:third_order_ivp_ex_3_sols:b1}}
\end{subfigure}
\begin{subfigure}{0.23\textwidth}
\includegraphics[width=\textwidth]{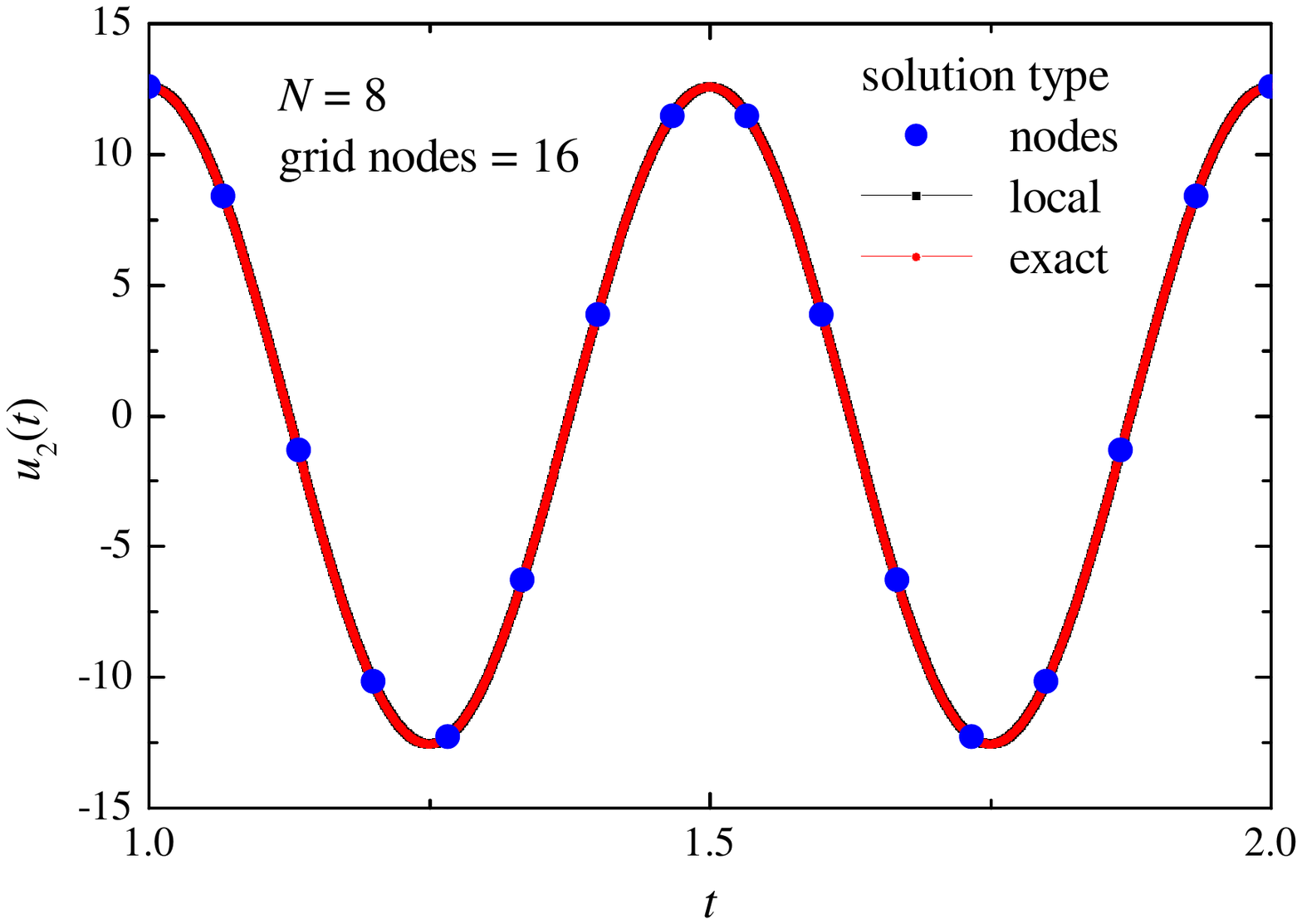}
\vspace{-8mm}\caption{\label{fig:third_order_ivp_ex_3_sols:b2}}
\end{subfigure}
\begin{subfigure}{0.23\textwidth}
\includegraphics[width=\textwidth]{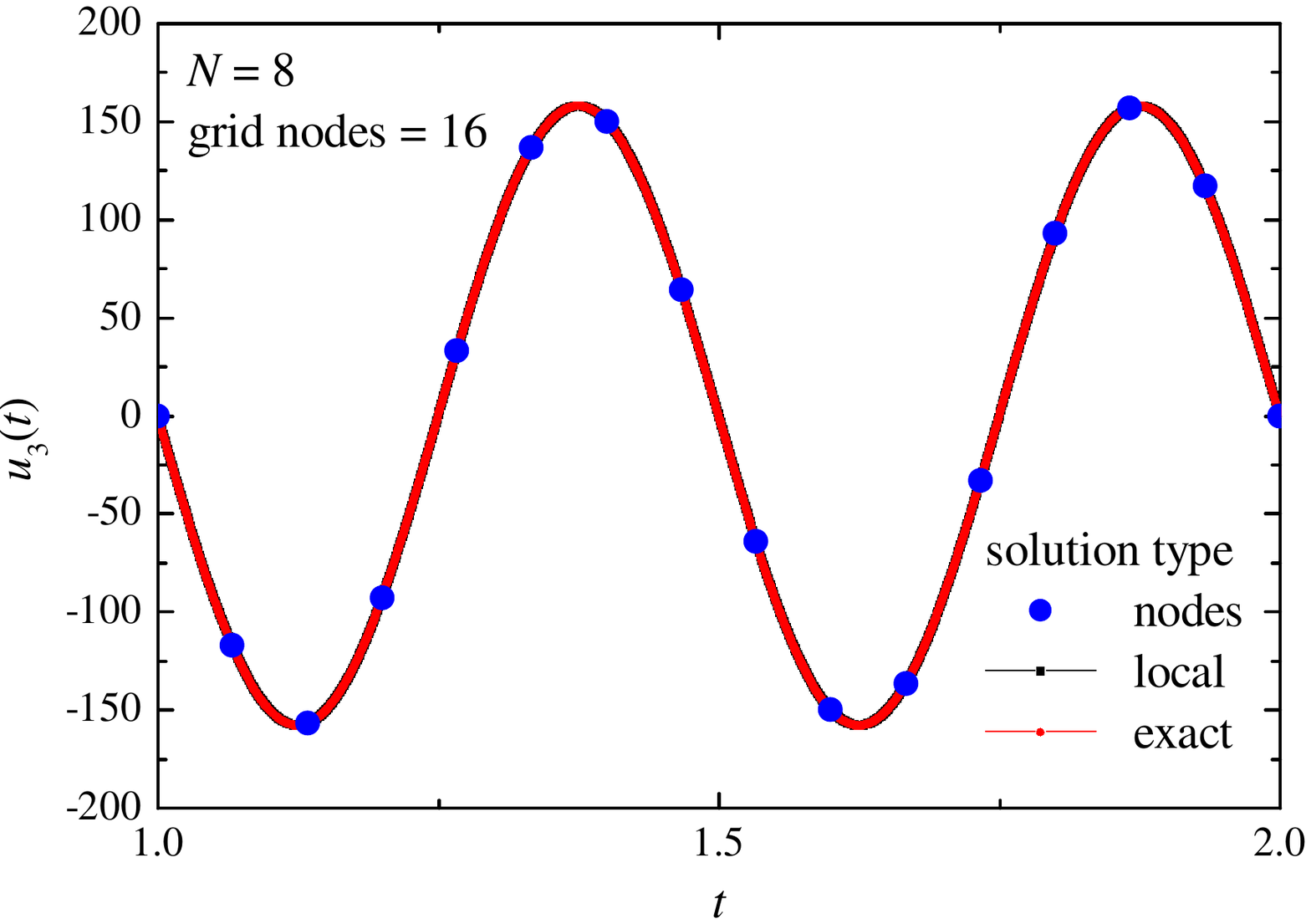}
\vspace{-8mm}\caption{\label{fig:third_order_ivp_ex_3_sols:b3}}
\end{subfigure}
\begin{subfigure}{0.23\textwidth}
\includegraphics[width=\textwidth]{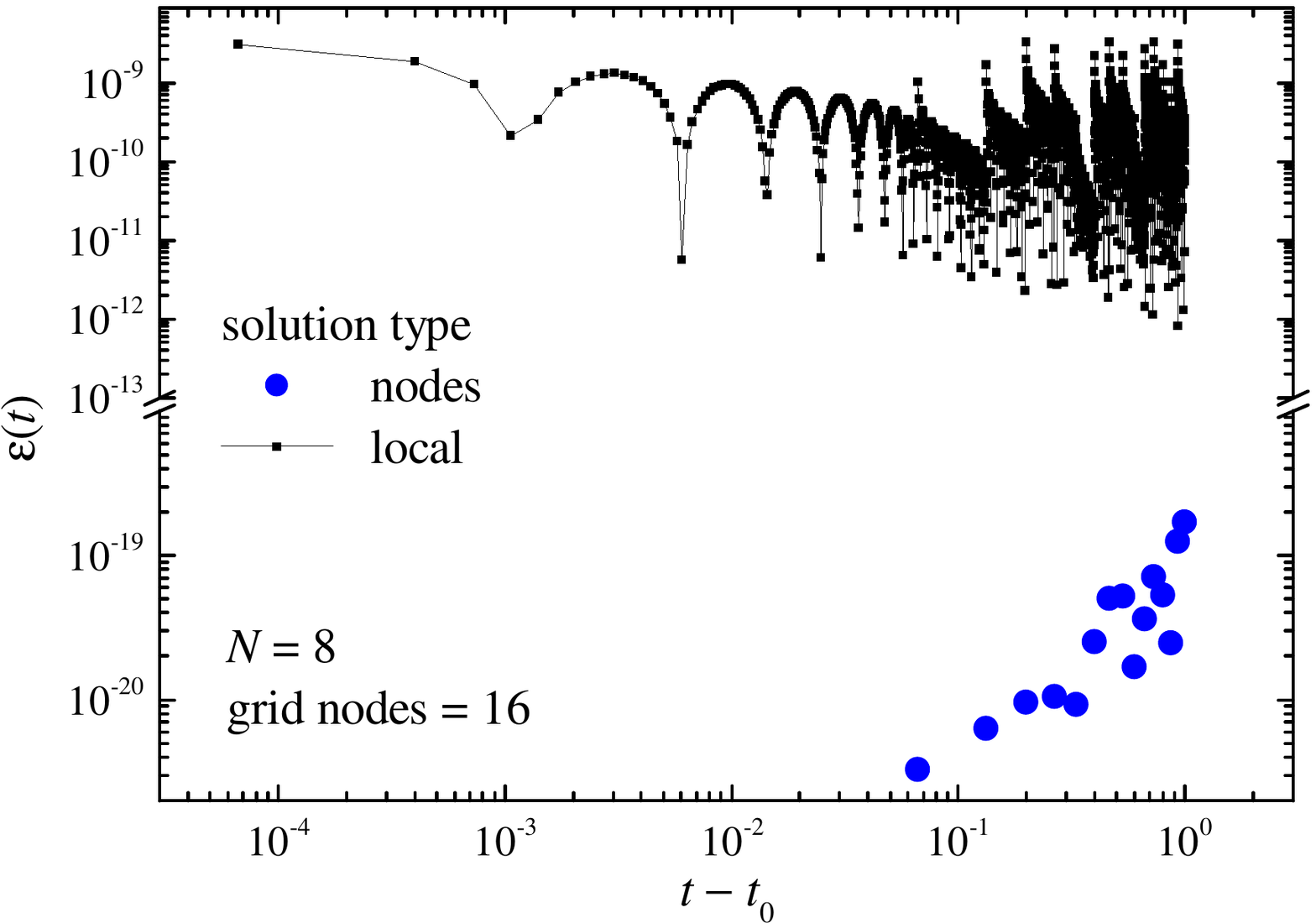}
\vspace{-8mm}\caption{\label{fig:third_order_ivp_ex_3_sols:b4}}
\end{subfigure}\\
\begin{subfigure}{0.23\textwidth}
\includegraphics[width=\textwidth]{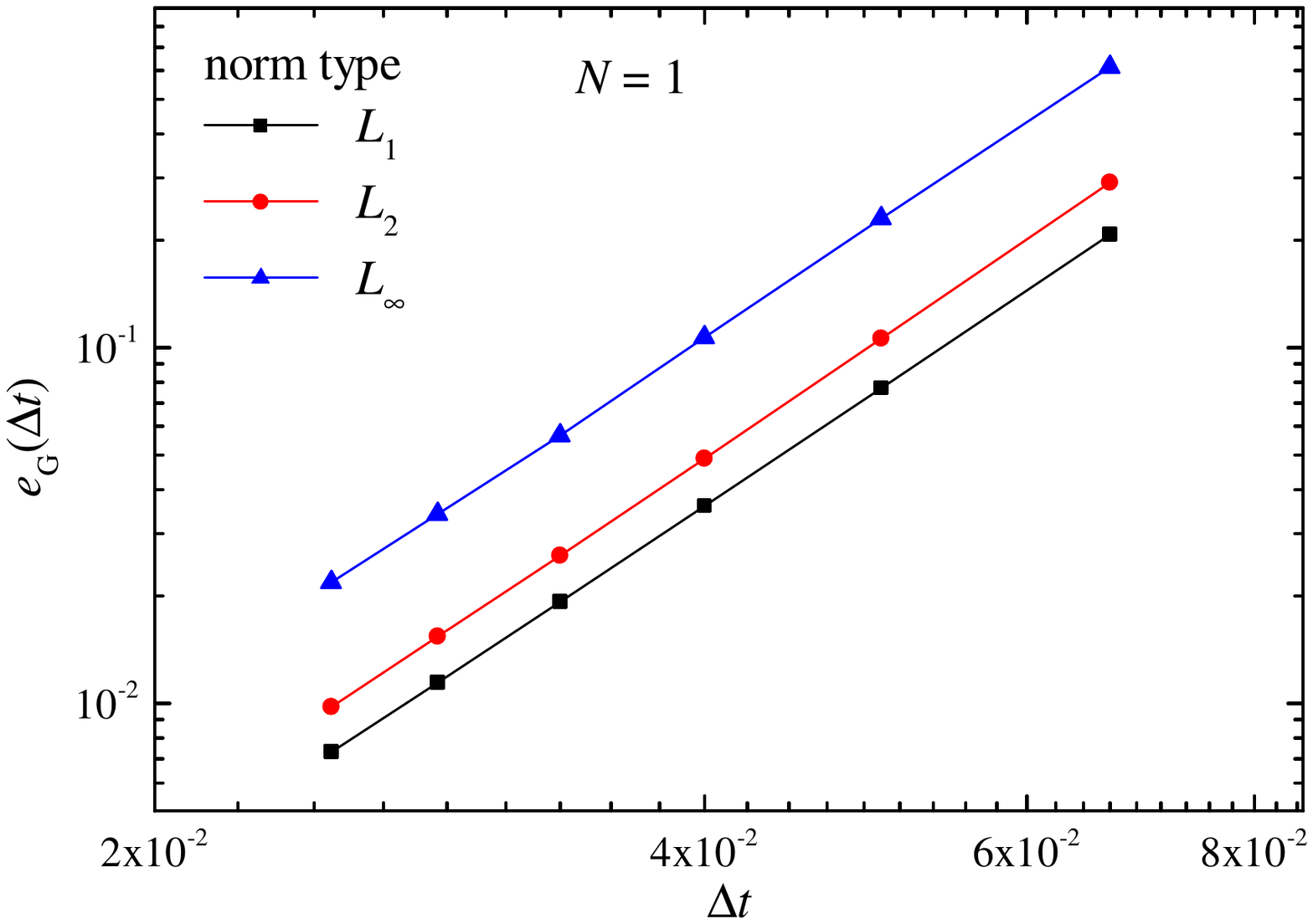}
\vspace{-8mm}\caption{\label{fig:third_order_ivp_ex_3_sols:c1}}
\end{subfigure}
\begin{subfigure}{0.23\textwidth}
\includegraphics[width=\textwidth]{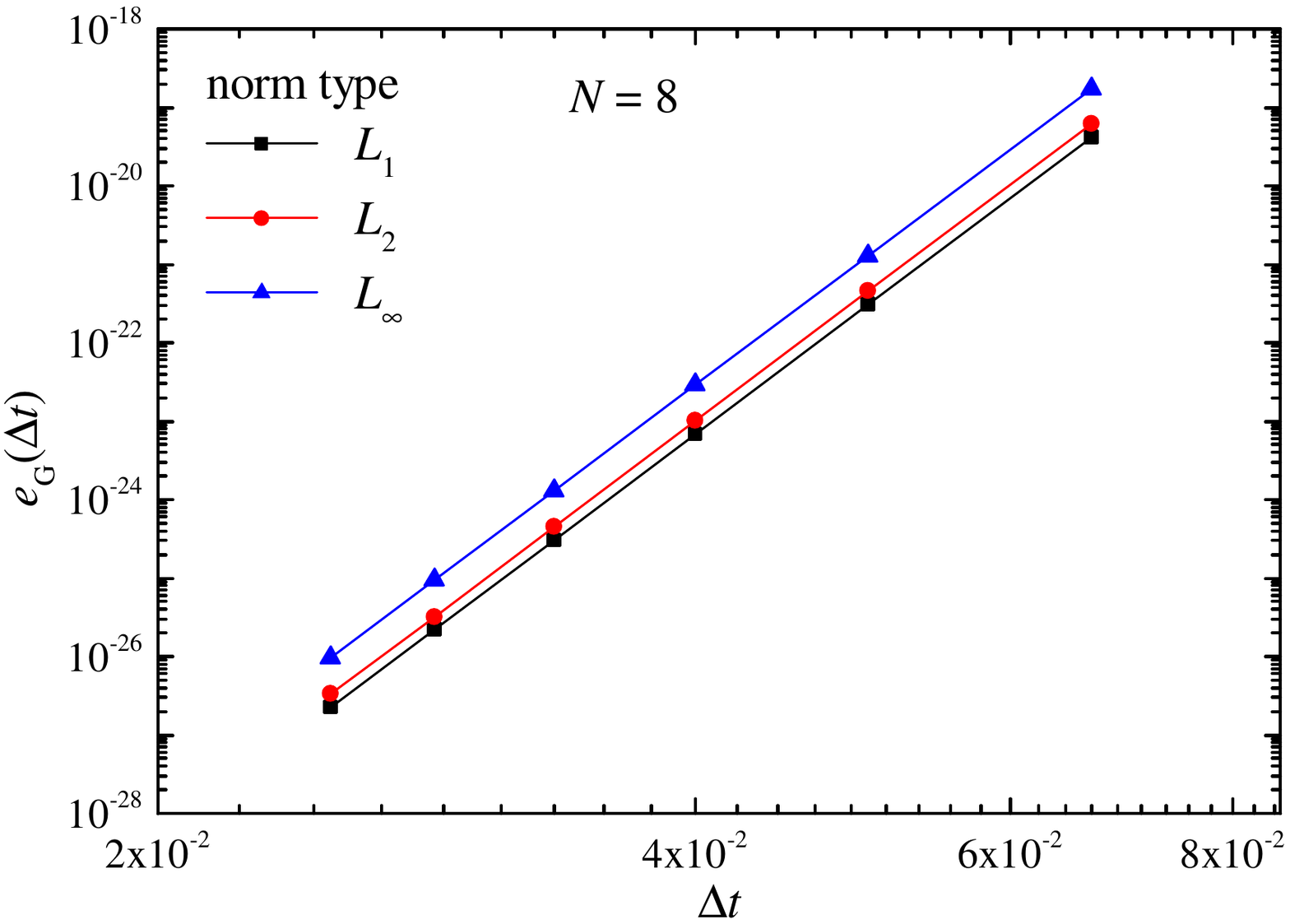}
\vspace{-8mm}\caption{\label{fig:third_order_ivp_ex_3_sols:c2}}
\end{subfigure}
\begin{subfigure}{0.23\textwidth}
\includegraphics[width=\textwidth]{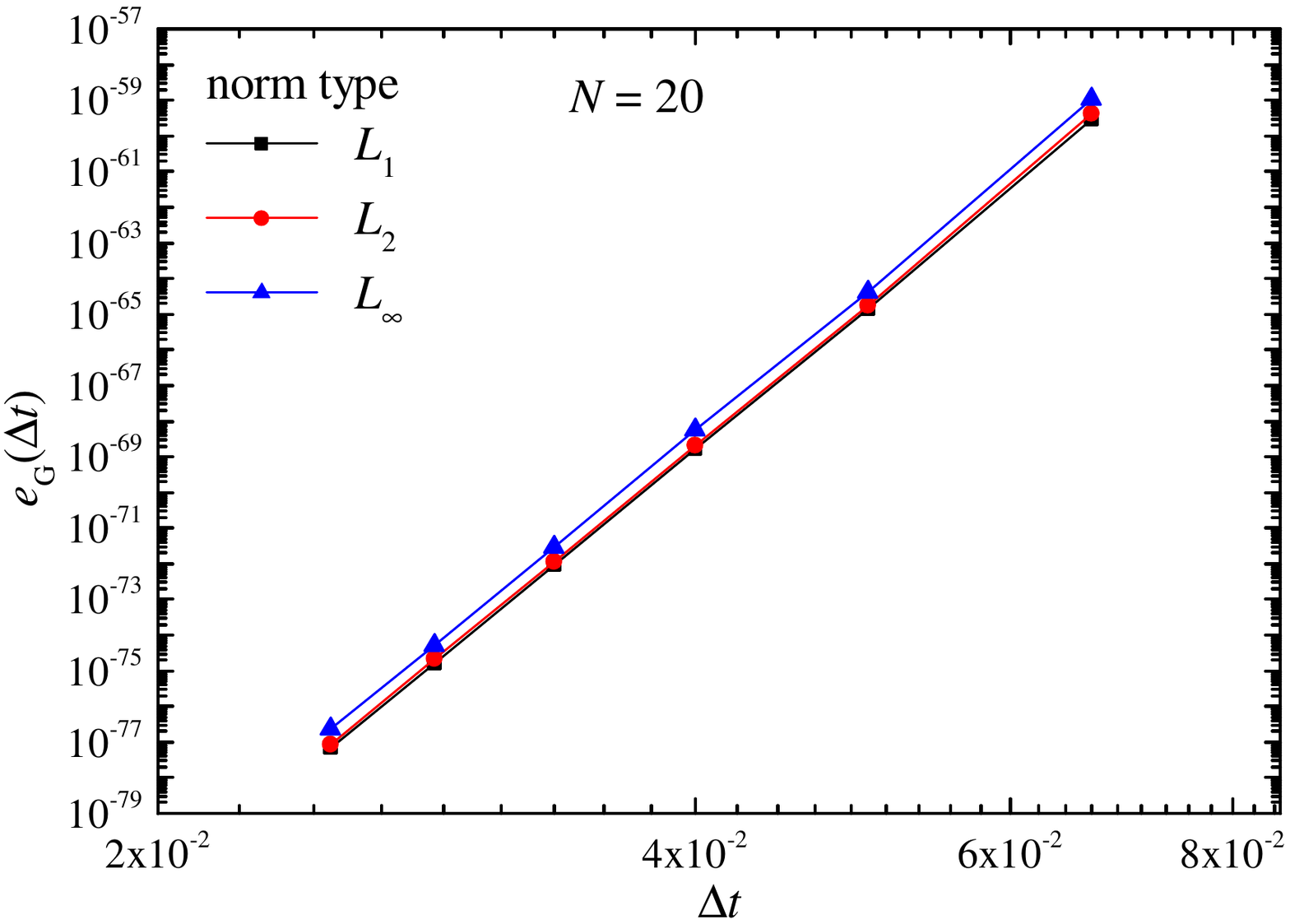}
\vspace{-8mm}\caption{\label{fig:third_order_ivp_ex_3_sols:c3}}
\end{subfigure}
\begin{subfigure}{0.23\textwidth}
\includegraphics[width=\textwidth]{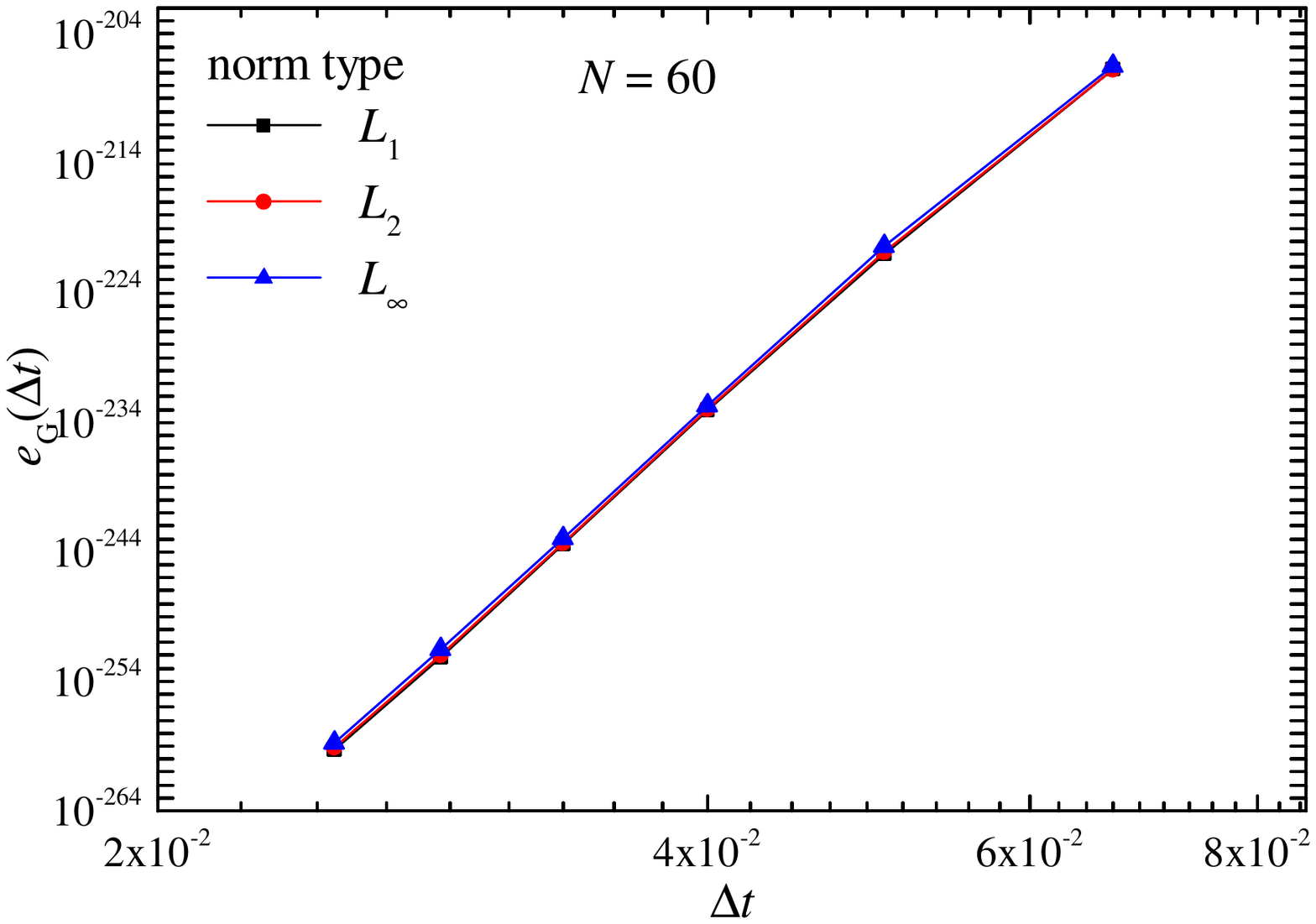}
\vspace{-8mm}\caption{\label{fig:third_order_ivp_ex_3_sols:c4}}
\end{subfigure}\\
\begin{subfigure}{0.23\textwidth}
\includegraphics[width=\textwidth]{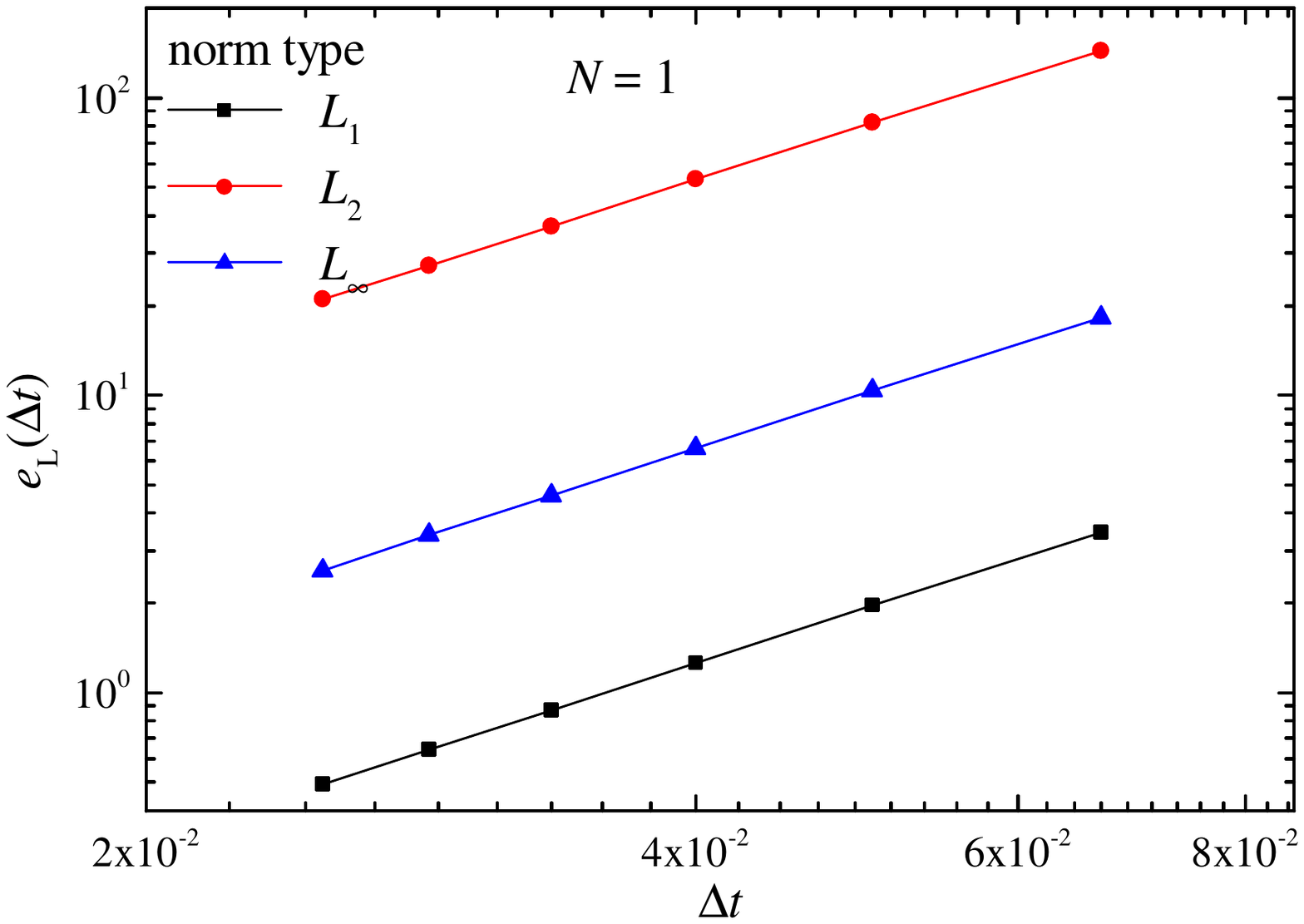}
\vspace{-8mm}\caption{\label{fig:third_order_ivp_ex_3_sols:d1}}
\end{subfigure}
\begin{subfigure}{0.23\textwidth}
\includegraphics[width=\textwidth]{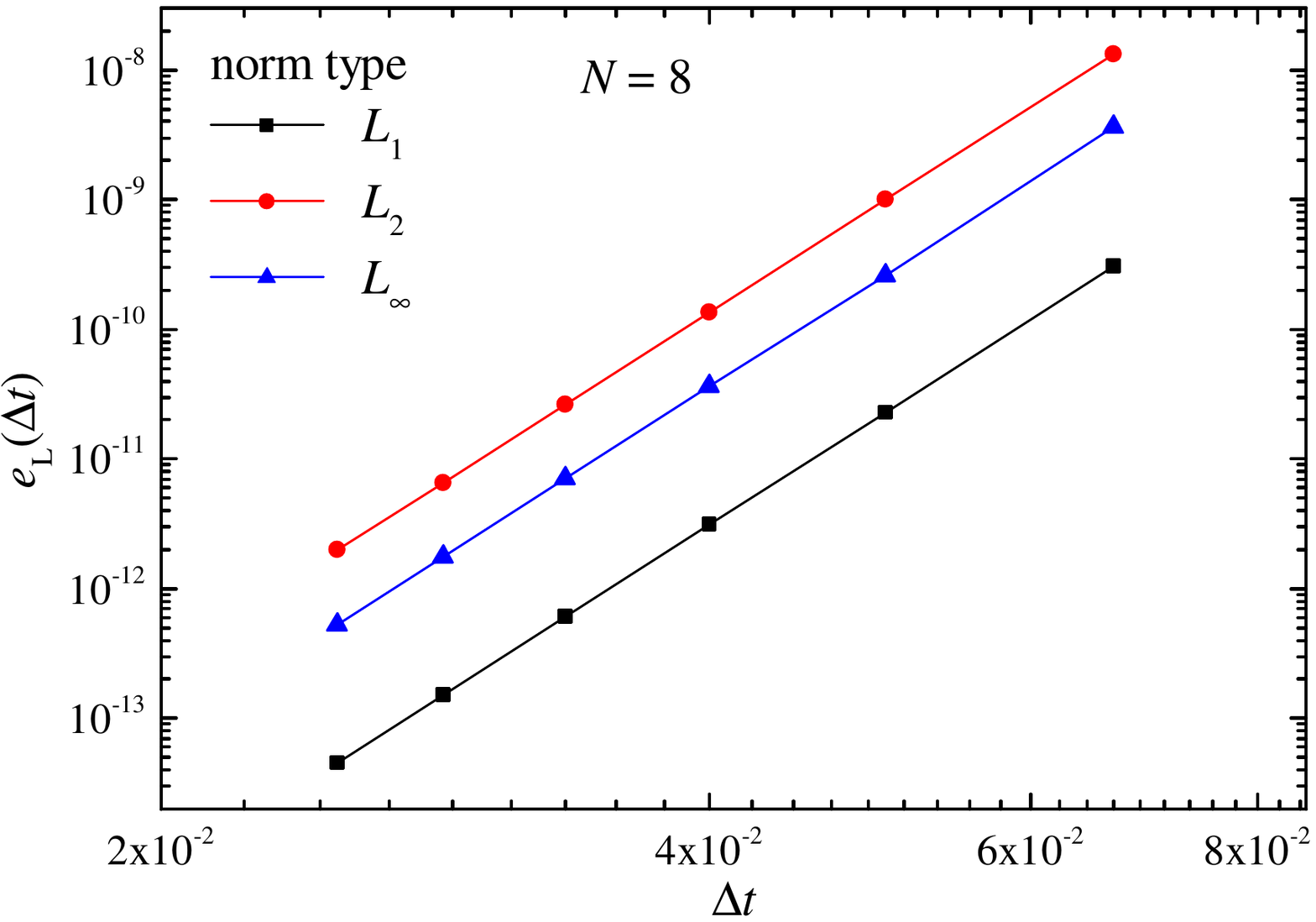}
\vspace{-8mm}\caption{\label{fig:third_order_ivp_ex_3_sols:d2}}
\end{subfigure}
\begin{subfigure}{0.23\textwidth}
\includegraphics[width=\textwidth]{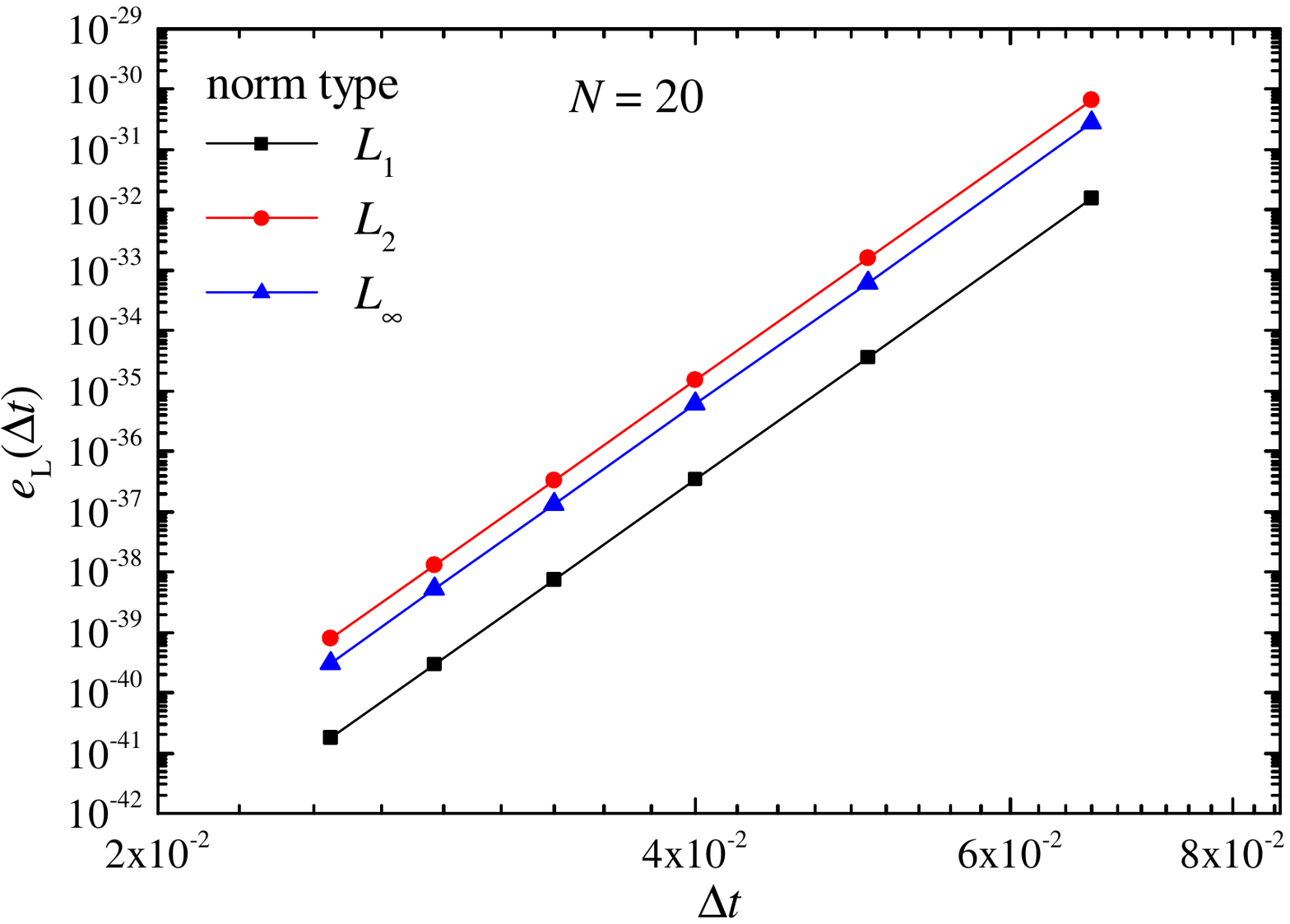}
\vspace{-8mm}\caption{\label{fig:third_order_ivp_ex_3_sols:d3}}
\end{subfigure}
\begin{subfigure}{0.23\textwidth}
\includegraphics[width=\textwidth]{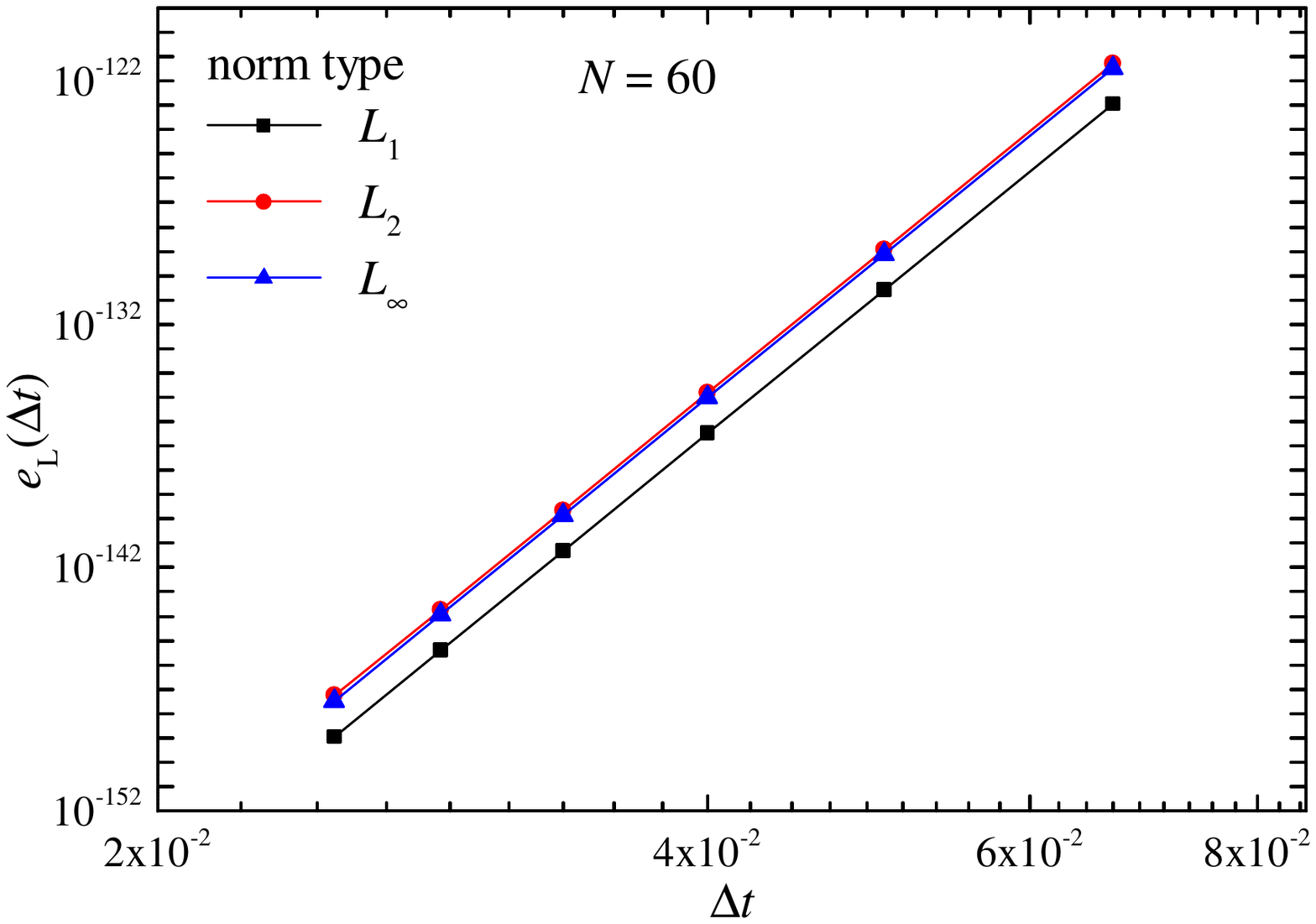}
\vspace{-8mm}\caption{\label{fig:third_order_ivp_ex_3_sols:d4}}
\end{subfigure}\\
\caption{%
Numerical solution of the problem (\ref{eq:third_order_ivp_ex_3}). Comparison of the solution at nodes $\mathbf{u}_{n}$, the local solution $\mathbf{u}_{L}(t)$ and the exact solution $\mathbf{u}^{\rm ex}(t)$ for components $u_{1}$ (\subref{fig:third_order_ivp_ex_3_sols:a1}, \subref{fig:third_order_ivp_ex_3_sols:b1}), $u_{2}$ (\subref{fig:third_order_ivp_ex_3_sols:a2}, \subref{fig:third_order_ivp_ex_3_sols:b2}) and $u_{3}$ (\subref{fig:third_order_ivp_ex_3_sols:a3}, \subref{fig:third_order_ivp_ex_3_sols:b3}), the error $\varepsilon(t)$ (\subref{fig:third_order_ivp_ex_3_sols:a4}, \subref{fig:third_order_ivp_ex_3_sols:b4}), obtained using polynomials with degrees $N = 1$ (\subref{fig:third_order_ivp_ex_3_sols:a1}, \subref{fig:third_order_ivp_ex_3_sols:a2}, \subref{fig:third_order_ivp_ex_3_sols:a3}, \subref{fig:third_order_ivp_ex_3_sols:a4}) and $N = 8$ (\subref{fig:third_order_ivp_ex_3_sols:b1}, \subref{fig:third_order_ivp_ex_3_sols:b2}, \subref{fig:third_order_ivp_ex_3_sols:b3}, \subref{fig:third_order_ivp_ex_3_sols:b4}). Dependence of the global error for the solution at nodes $e_{G}$ (\subref{fig:third_order_ivp_ex_3_sols:c1}, \subref{fig:third_order_ivp_ex_3_sols:c2}, \subref{fig:third_order_ivp_ex_3_sols:c3}, \subref{fig:third_order_ivp_ex_3_sols:c4}) and the local solution $e_{L}$ (\subref{fig:third_order_ivp_ex_3_sols:d1}, \subref{fig:third_order_ivp_ex_3_sols:d2}, \subref{fig:third_order_ivp_ex_3_sols:d3}, \subref{fig:third_order_ivp_ex_3_sols:d4}) on the discretization step $\mathrm{\Delta}t$, obtained in the norms $L_{1}$, $L_{2}$ and $L_{\infty}$, obtained using polynomials with degrees $N = 1$ (\subref{fig:third_order_ivp_ex_3_sols:c1}, \subref{fig:third_order_ivp_ex_3_sols:d1}), $N = 8$ (\subref{fig:third_order_ivp_ex_3_sols:c2}, \subref{fig:third_order_ivp_ex_3_sols:d2}), $N = 20$ (\subref{fig:third_order_ivp_ex_3_sols:c3}, \subref{fig:third_order_ivp_ex_3_sols:d3}) and $N = 60$ (\subref{fig:third_order_ivp_ex_3_sols:c4}, \subref{fig:third_order_ivp_ex_3_sols:d4}).
}
\label{fig:third_order_ivp_ex_3_sols}
\end{figure} 
\begin{table}[h!]
\centering
\caption{%
Orders of convergence $p_{L_{1}}$, $p_{L_{2}}$, $p_{L_{\infty}}$ of the ADER-DG numerical method with local DG predictor calculated in norms $L_{1}$, $L_{2}$, $L_{\infty}$ for the problem (\ref{eq:third_order_ivp_ex_3}); $N$ is the degree of the basis polynomials $\varphi_{p}$ used for the DG representations. The orders of convergence $p$ for the solution at nodes $\mathbf{u}_{n}$ are compared to the theoretical value $p_{\rm (G)} = 2N+1$, the orders of convergence $p$ for the local solution $\mathbf{u}_{L}$ represented between the nodes are compared to the theoretical value $p_{\rm (L)} = N+1$.
}
\label{tab:third_order_ivp_ex_3_sols}
\begin{tabular}{@{}lllllllll@{}}
\toprule
& \multicolumn{4}{c}{order $p$ for solution at nodes} & \multicolumn{4}{c}{order $p$ for local solution} \\
\midrule
$N$ & $p_{L_{1}}$ & $p_{L_{2}}$ & $p_{L_{\infty}}$ & theor. & $p_{L_{1}}$ & $p_{L_{2}}$ & $p_{L_{\infty}}$ & theor. \\
\midrule
$1$	&	$3.41$	&	$3.46$	&	$3.41$		&	$3$	&	$1.99$	&	$1.97$	&	$2.00$		&	$2$	\\
$2$	&	$5.43$	&	$5.17$	&	$4.48$		&	$5$	&	$2.99$	&	$2.97$	&	$2.98$		&	$3$	\\
$3$	&	$6.90$	&	$6.92$	&	$6.83$		&	$7$	&	$4.00$	&	$3.97$	&	$4.00$		&	$4$	\\
$4$	&	$8.99$	&	$9.00$	&	$8.93$		&	$9$	&	$4.99$	&	$4.98$	&	$4.98$		&	$5$	\\
$5$	&	$11.01$	&	$11.03$	&	$10.96$		&	$11$	&	$6.00$	&	$5.98$	&	$6.00$		&	$6$	\\
$6$	&	$13.03$	&	$13.04$	&	$12.98$		&	$13$	&	$6.99$	&	$6.98$	&	$6.98$		&	$7$	\\
$7$	&	$15.04$	&	$15.06$	&	$14.99$		&	$15$	&	$8.00$	&	$7.98$	&	$8.00$		&	$8$	\\
$8$	&	$17.05$	&	$17.06$	&	$17.00$		&	$17$	&	$8.99$	&	$8.98$	&	$8.98$		&	$9$	\\
$9$	&	$19.06$	&	$19.08$	&	$19.01$		&	$19$	&	$10.00$	&	$9.98$	&	$10.00$		&	$10$	\\
$10$	&	$21.08$	&	$21.09$	&	$21.02$		&	$21$	&	$10.99$	&	$10.98$	&	$10.98$		&	$11$	\\
\midrule
$15$	&	$31.23$	&	$31.32$	&	$31.31$		&	$31$	&	$16.00$	&	$15.98$	&	$16.00$		&	$16$	\\
$20$	&	$41.29$	&	$41.38$	&	$41.34$		&	$41$	&	$20.99$	&	$20.98$	&	$20.99$		&	$21$	\\
$25$	&	$52.50$	&	$52.47$	&	$52.21$		&	$51$	&	$26.00$	&	$25.98$	&	$26.00$		&	$26$	\\
$30$	&	$63.39$	&	$63.31$	&	$62.87$		&	$61$	&	$30.99$	&	$30.98$	&	$30.99$		&	$31$	\\
$35$	&	$73.57$	&	$73.45$	&	$73.07$		&	$71$	&	$36.00$	&	$35.98$	&	$36.00$		&	$36$	\\
$40$	&	$84.73$	&	$84.71$	&	$84.47$		&	$81$	&	$40.99$	&	$40.98$	&	$40.99$		&	$41$	\\
$45$	&	$94.24$	&	$94.12$	&	$93.61$		&	$91$	&	$46.01$	&	$45.98$	&	$46.00$		&	$46$	\\
$50$	&	$105.70$	&	$105.63$	&	$105.52$		&	$101$	&	$50.99$	&	$50.98$	&	$50.99$		&	$51$	\\
$55$	&	$113.40$	&	$113.21$	&	$112.58$		&	$111$	&	$56.01$	&	$55.98$	&	$56.00$		&	$56$	\\
$60$	&	$124.25$	&	$124.04$	&	$123.74$		&	$121$	&	$60.99$	&	$60.98$	&	$60.99$		&	$61$	\\
\bottomrule
\end{tabular}
\end{table} 

The third example for IVP third-order ODE was chosen in the following form~\cite{dg_ivp_ode_5}:
\begin{equation}
\begin{split}
&\dddot{x} = x\ddot{x} - \dfrac{2}{t}\dot{x} + 16\pi^{2}x^{2} + \left(\dfrac{8\pi}{t}-64\pi^{3}\right)\cos(4\pi t);\\
&x(0) = 0;\ \dot{x}(0) = 4\pi;\ \ddot{x}(0) = 0;\quad 1 \leqslant t \leqslant 2;
\end{split}
\end{equation}
with exact analytical solution of this problem having the form $x(t) = \sin(4\pi t)$. The IVP for a first-order non-linear ODE system corresponding to this problem has the following form:
\begin{equation}\label{eq:third_order_ivp_ex_3}
\begin{split}
&\frac{du_{1}}{dt} = u_{2};\quad \frac{du_{2}}{dt} = u_{3};\\
&\frac{du_{3}}{dt} = u_{1}u_{3} - \dfrac{2}{t}u_{2} + 16\pi^{2}u_{1}^{2} + \left(\dfrac{8\pi}{t}-64\pi^{3}\right)\cos(4\pi t);\\
&u_{1}(0) = 0;\quad u_{2}(0) = 4\pi;\quad u_{3}(0) = 0;
\end{split}
\end{equation}
where $\mathbf{u} = [u_{1},\, u_{2},\, u_{3}]^{T} = [x,\, \dot{x},\, \ddot{x}]^{T}$ is a desired three component vector function. The solution definition domain $\Omega$ was chosen as $1 \leqslant t \leqslant 2$. The full three component exact analytical solution of this problem has the form
\begin{equation}
\mathbf{u}^{\rm ex} = \left[
\sin(4\pi t),\quad
4\pi\cos(4\pi t),\quad
-(4\pi)^{2}\sin(4\pi t)
\right]^{T}.
\end{equation}

The study of the convergence of the numerical solution was carried out on the basis of the analysis of the solution for the set of $6$ different uniform discretizations of the domain of definition $\Omega$ of the desired function, with the number of grid nodes $L = 16$, $21$, $26$, $31$, $36$, $41$; note that the zero node, in which the initial condition was defined, also refers to grid nodes, so the discretization step was determined by the expression $\Delta t = 2\pi/(L-1)$.

The results of the study of the application of the numerical method for this problem are presented in Fig.~\ref{fig:third_order_ivp_ex_3_sols} and in Table~\ref{tab:third_order_ivp_ex_3_sols}. The results obtained, in terms of the accuracy of the numerical solution when compared with the exact analytical solution, demonstrate all the same properties that were revealed in the previous considered problems. The error $\varepsilon(t)$ for $N > 1$ shows that the error scales for the local solution and the solution at nodes differ by several orders of magnitude, reaching values of $8$--$10$ and $80$--$100$ orders of magnitude for a degree of $N = 8$ and $N = 60$, respectively.

The calculated convergence orders $p_{L_{1}}$, $p_{L_{2}}$, $p_{L_{\infty}}$ for various functional norms are presented in Table~\ref{tab:third_order_ivp_ex_3_sols}. The convergence orders $p$ for the solution at nodes are compared to the theoretical value $p_{\rm (G)} = 2N+1$, the convergence orders $p$ for the local solution represented between the nodes are compared to the theoretical value $p_{\rm (L)} = N+1$, in accordance with (\ref{eq:expect_orders}). The main part of the obtained convergence orders $p$ sufficiently correspond to the expected theoretical values $p_{\rm (G)}$ and $p_{\rm (L)}$. In the region of large values of the degrees of polynomials $N \geqslant 25$, convergence orders $p$ are observed that exceed the expected values $p_{\rm (G)}$, which was an unexpected result. In the range of polynomial degrees $N = 50, \ldots, 60$, the empirical values of the convergence orders $p$ for the solution at nodes exceed the expected values $p_{\rm (G)}$ by more than $2$--$4$. In the case of a local solution $\mathbf{u}_{L}$, the convergence orders $p$ correspond to the expected values $p_{\rm (L)}$.

As a result, it becomes clear that the solution at nodes $\mathbf{u}_{n}$ obtained by the ADER-DG numerical method with a local DG predictor demonstrates the classical superconvergence $2N+1$ expected from such numerical methods, at least in problem (\ref{eq:third_order_ivp_ex_3}). The local solution $\mathbf{u}_{L}$ obtained by the DG predictor demonstrates the classical convergence characteristic of DG methods.

\subsection{Stiff test}

\begin{figure}[h!]
\captionsetup[subfigure]{%
position=bottom,
textfont=normalfont,
singlelinecheck=off,
justification=raggedright,
font=tiny
}
\centering
\begin{subfigure}{0.29\textwidth}
\includegraphics[width=\textwidth]{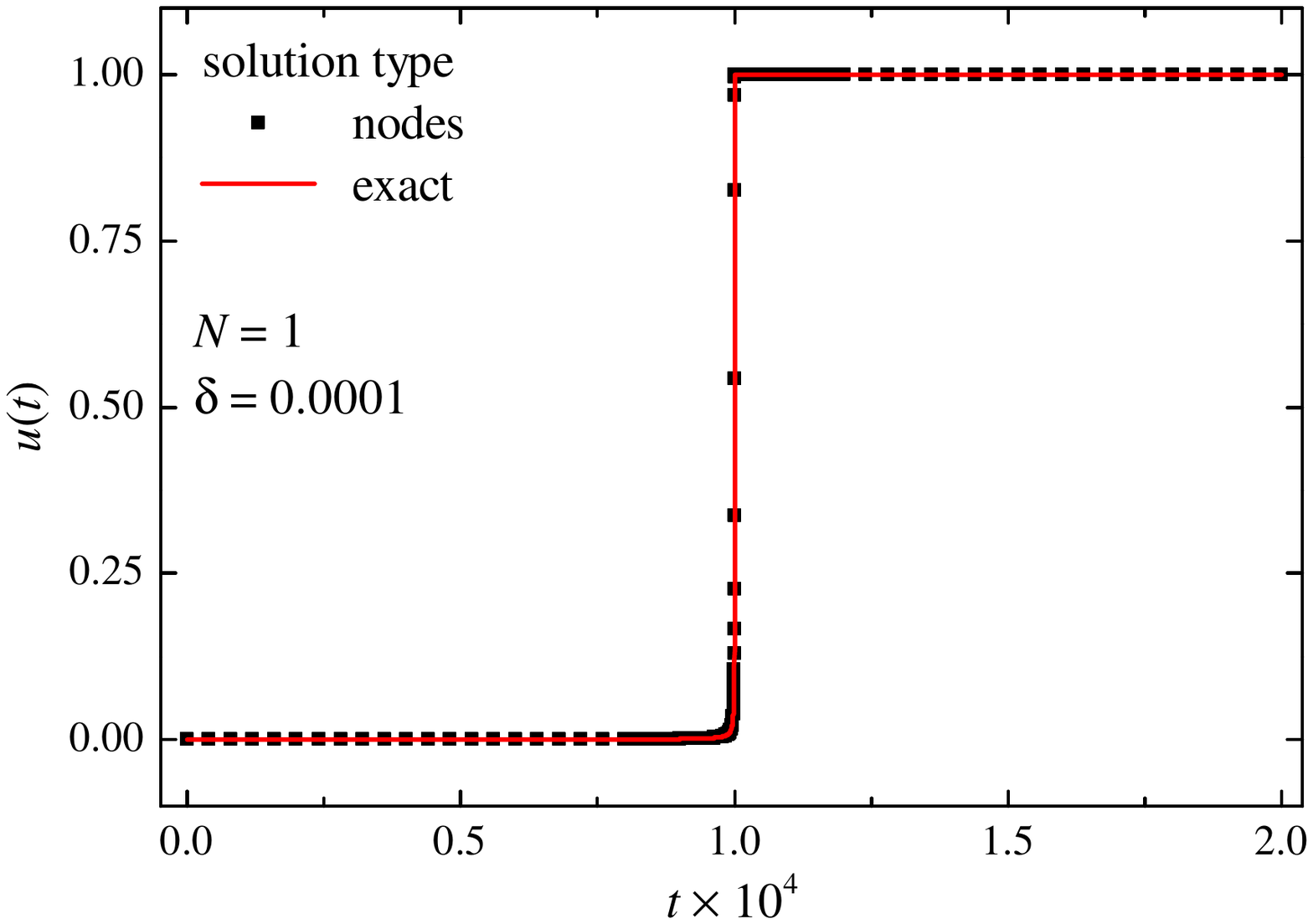}
\vspace{-9mm}\caption{\label{fig:shampine_flame_d_4:a1}}
\end{subfigure}\hspace{5mm}
\begin{subfigure}{0.29\textwidth}
\includegraphics[width=\textwidth]{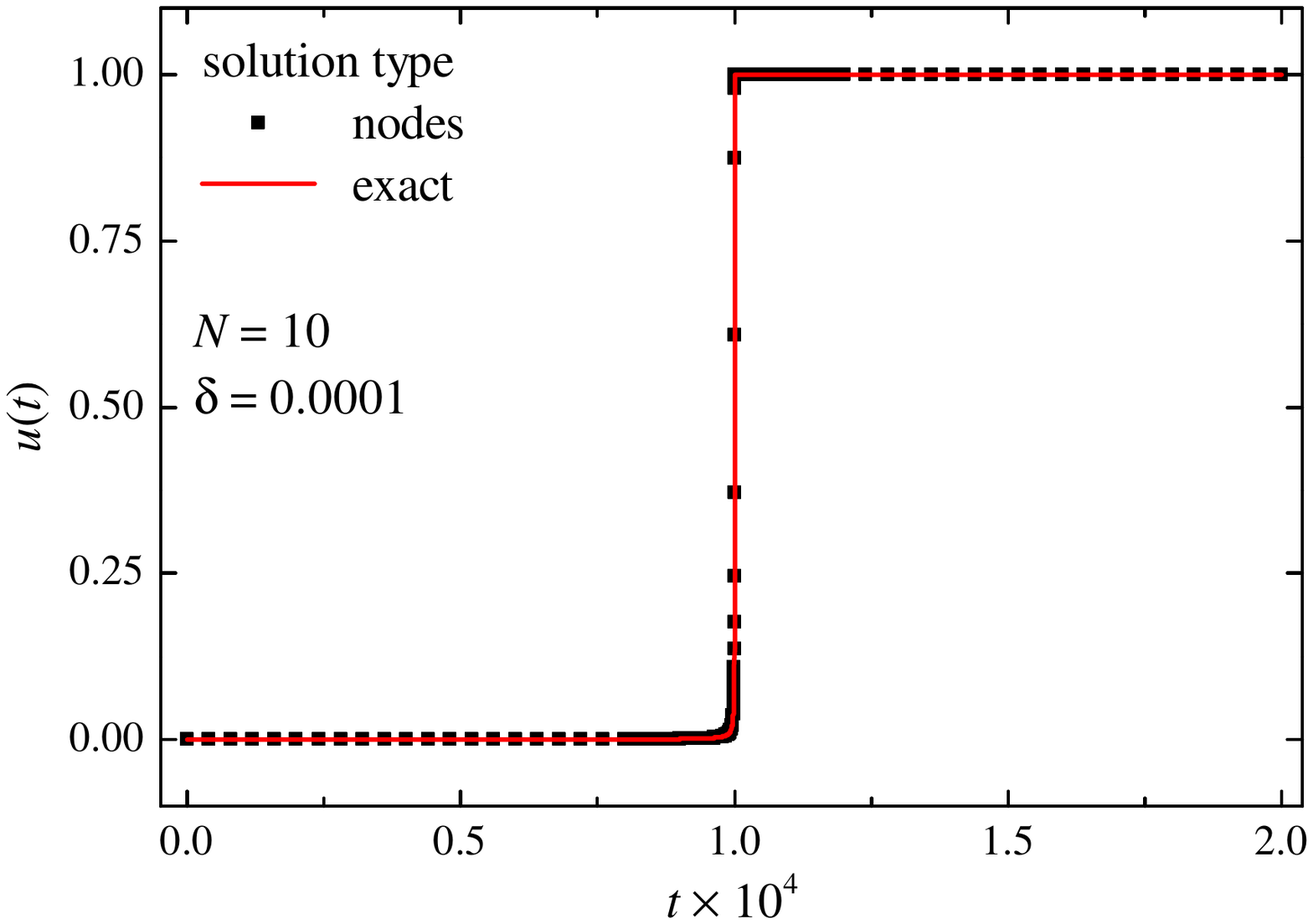}
\vspace{-9mm}\caption{\label{fig:shampine_flame_d_4:a2}}
\end{subfigure}\hspace{5mm}
\begin{subfigure}{0.29\textwidth}
\includegraphics[width=\textwidth]{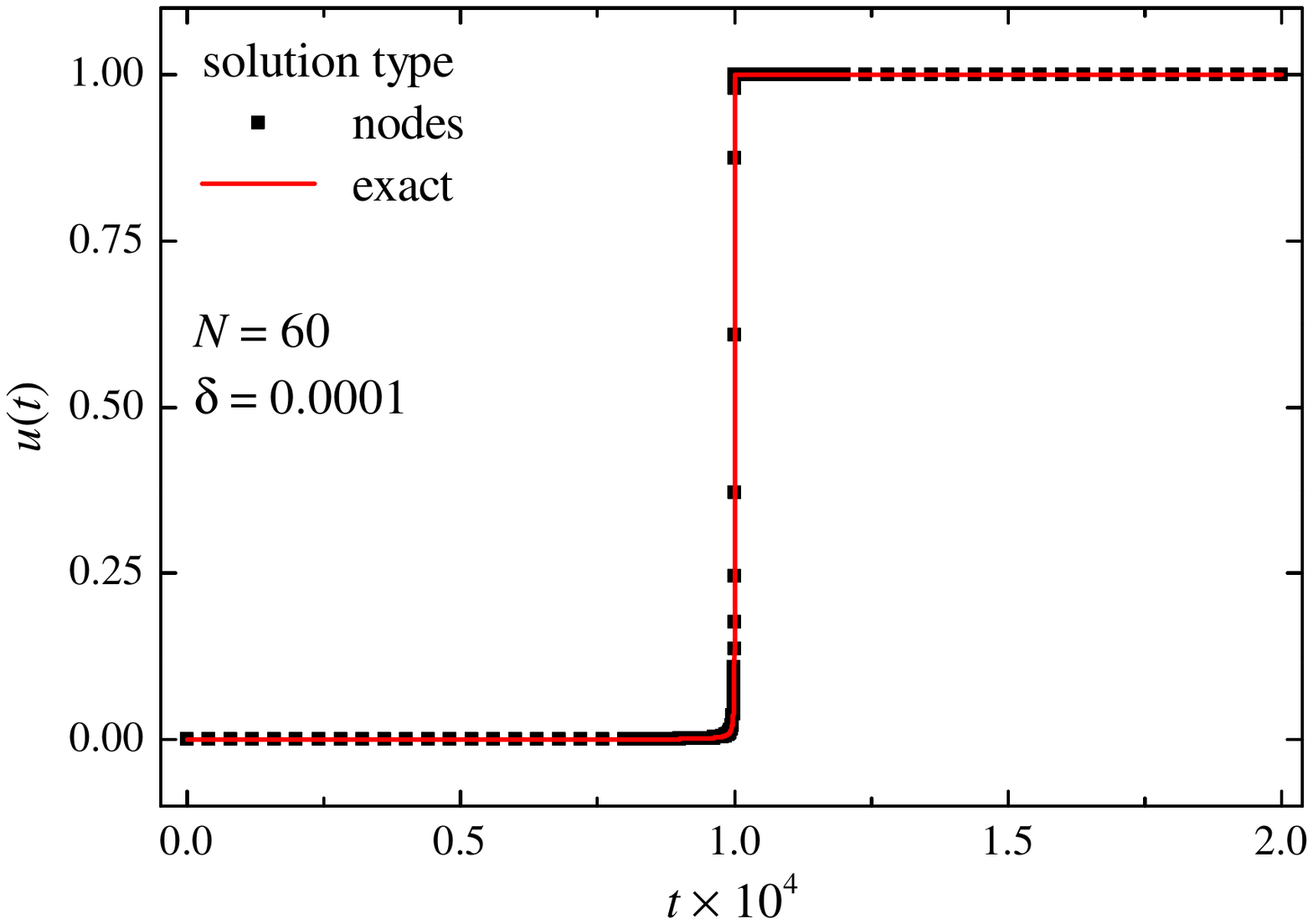}
\vspace{-9mm}\caption{\label{fig:shampine_flame_d_4:a3}}
\end{subfigure}\\
\begin{subfigure}{0.29\textwidth}
\includegraphics[width=\textwidth]{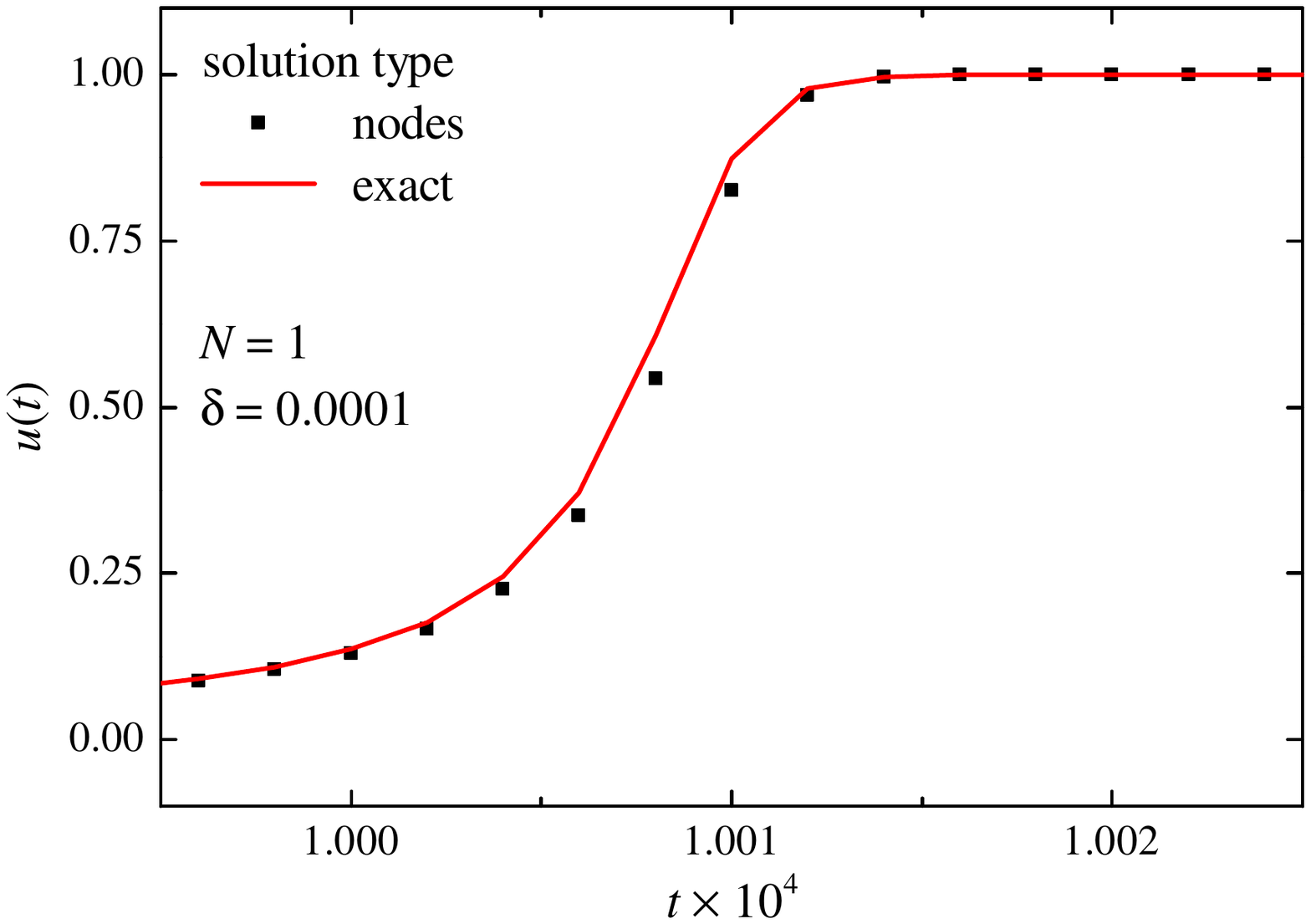}
\vspace{-9mm}\caption{\label{fig:shampine_flame_d_4:b1}}
\end{subfigure}\hspace{5mm}
\begin{subfigure}{0.29\textwidth}
\includegraphics[width=\textwidth]{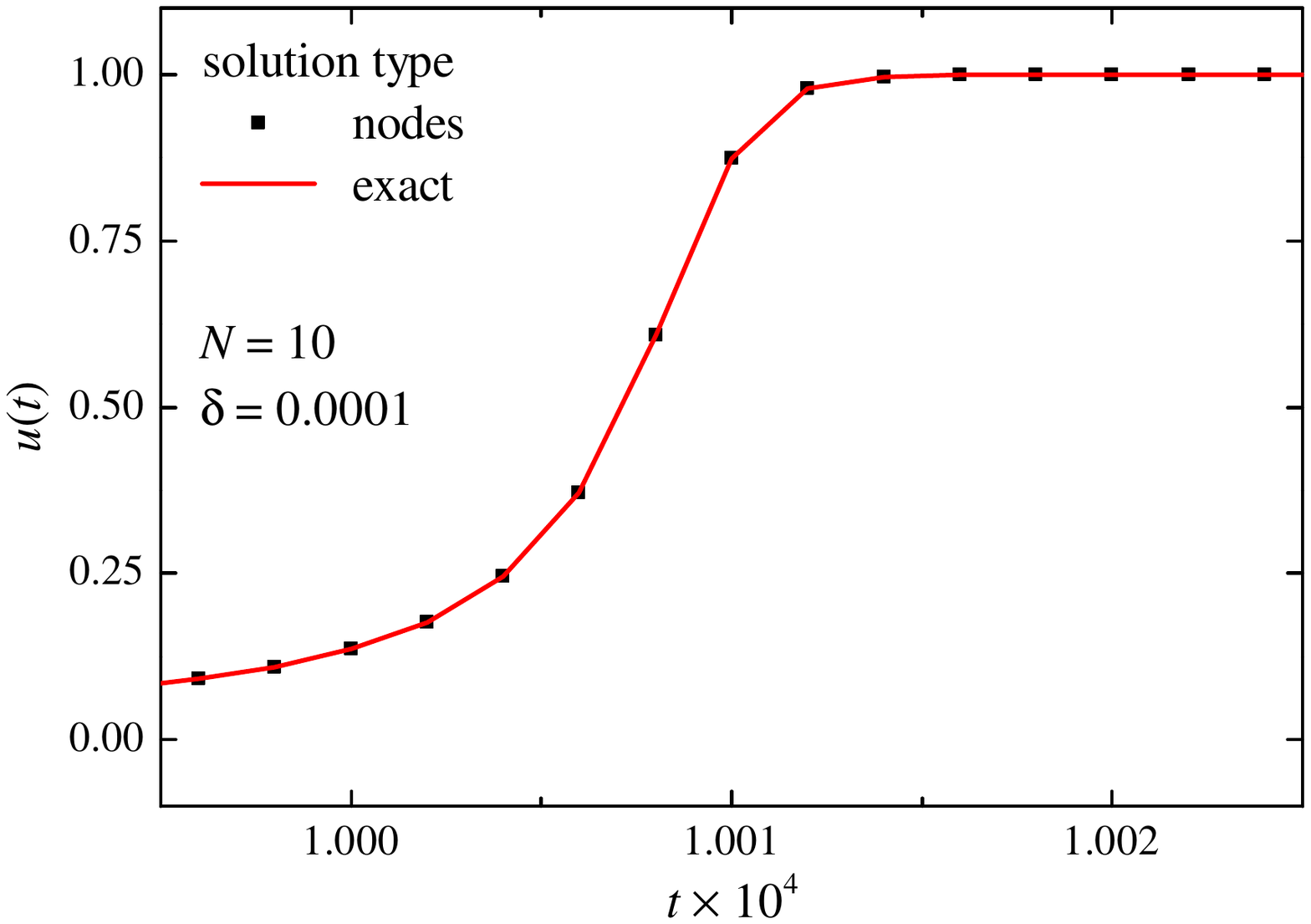}
\vspace{-9mm}\caption{\label{fig:shampine_flame_d_4:b2}}
\end{subfigure}\hspace{5mm}
\begin{subfigure}{0.29\textwidth}
\includegraphics[width=\textwidth]{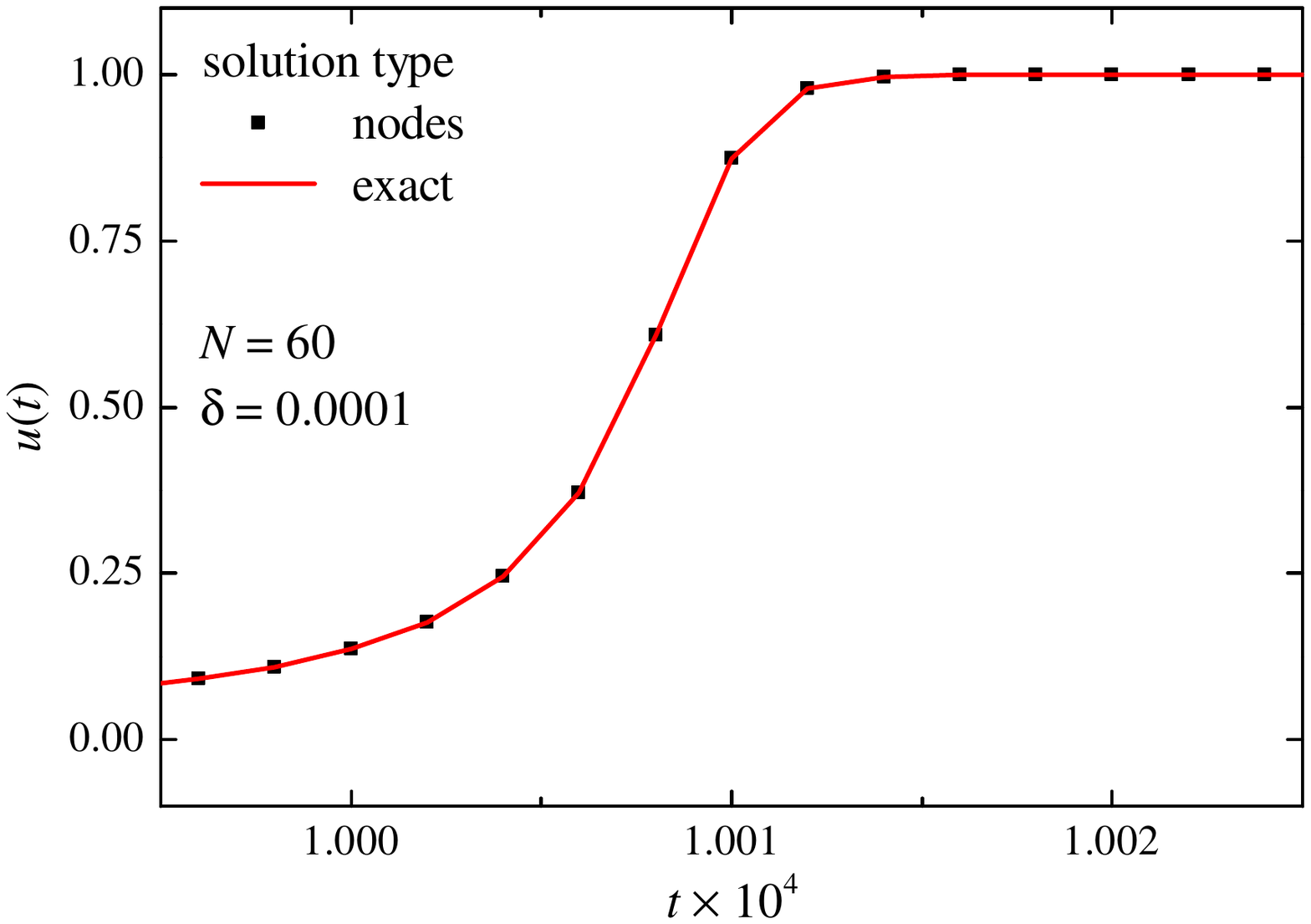}
\vspace{-9mm}\caption{\label{fig:shampine_flame_d_4:b3}}
\end{subfigure}\\
\begin{subfigure}{0.29\textwidth}
\includegraphics[width=\textwidth]{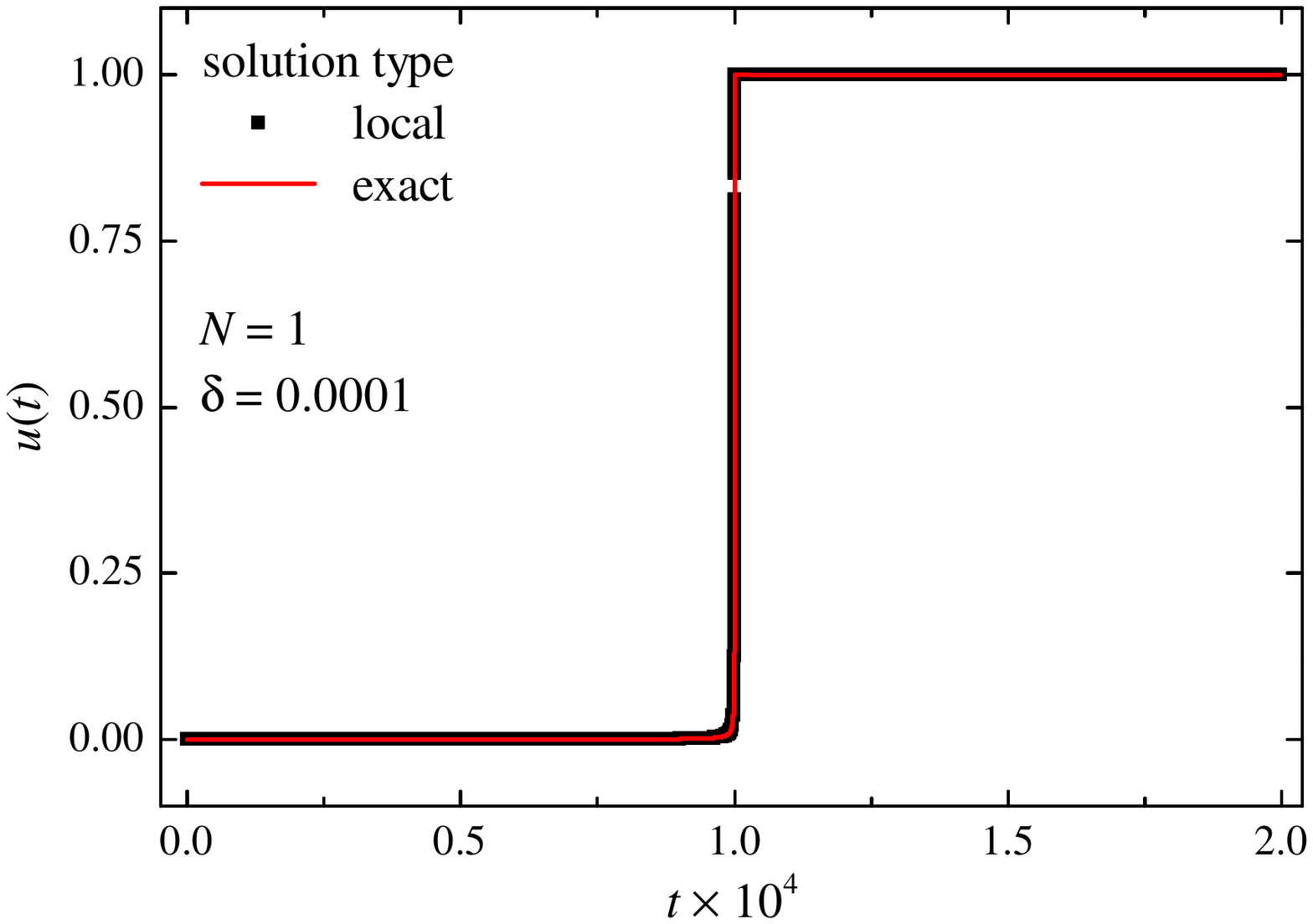}
\vspace{-9mm}\caption{\label{fig:shampine_flame_d_4:c1}}
\end{subfigure}\hspace{5mm}
\begin{subfigure}{0.29\textwidth}
\includegraphics[width=\textwidth]{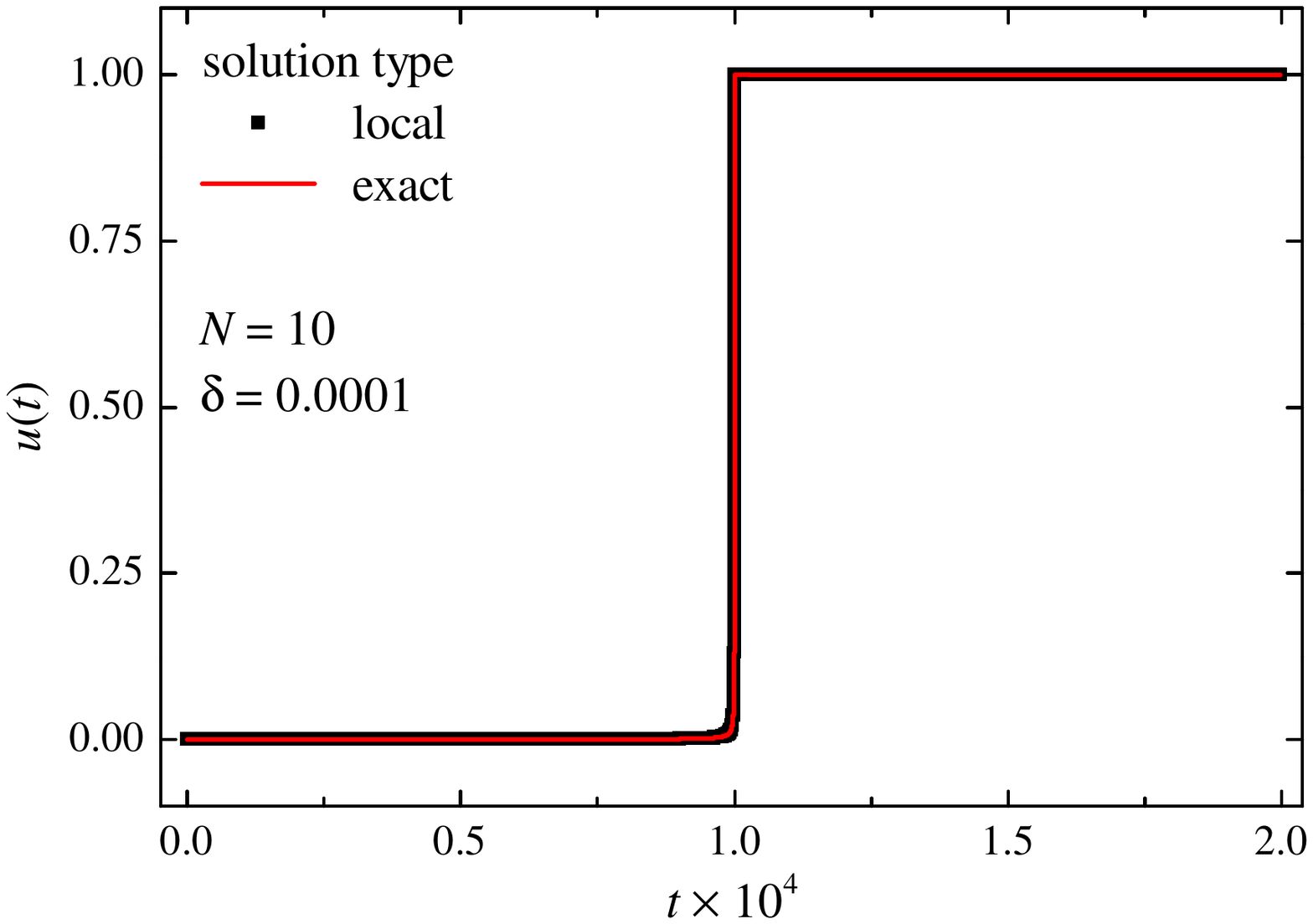}
\vspace{-9mm}\caption{\label{fig:shampine_flame_d_4:c2}}
\end{subfigure}\hspace{5mm}
\begin{subfigure}{0.29\textwidth}
\includegraphics[width=\textwidth]{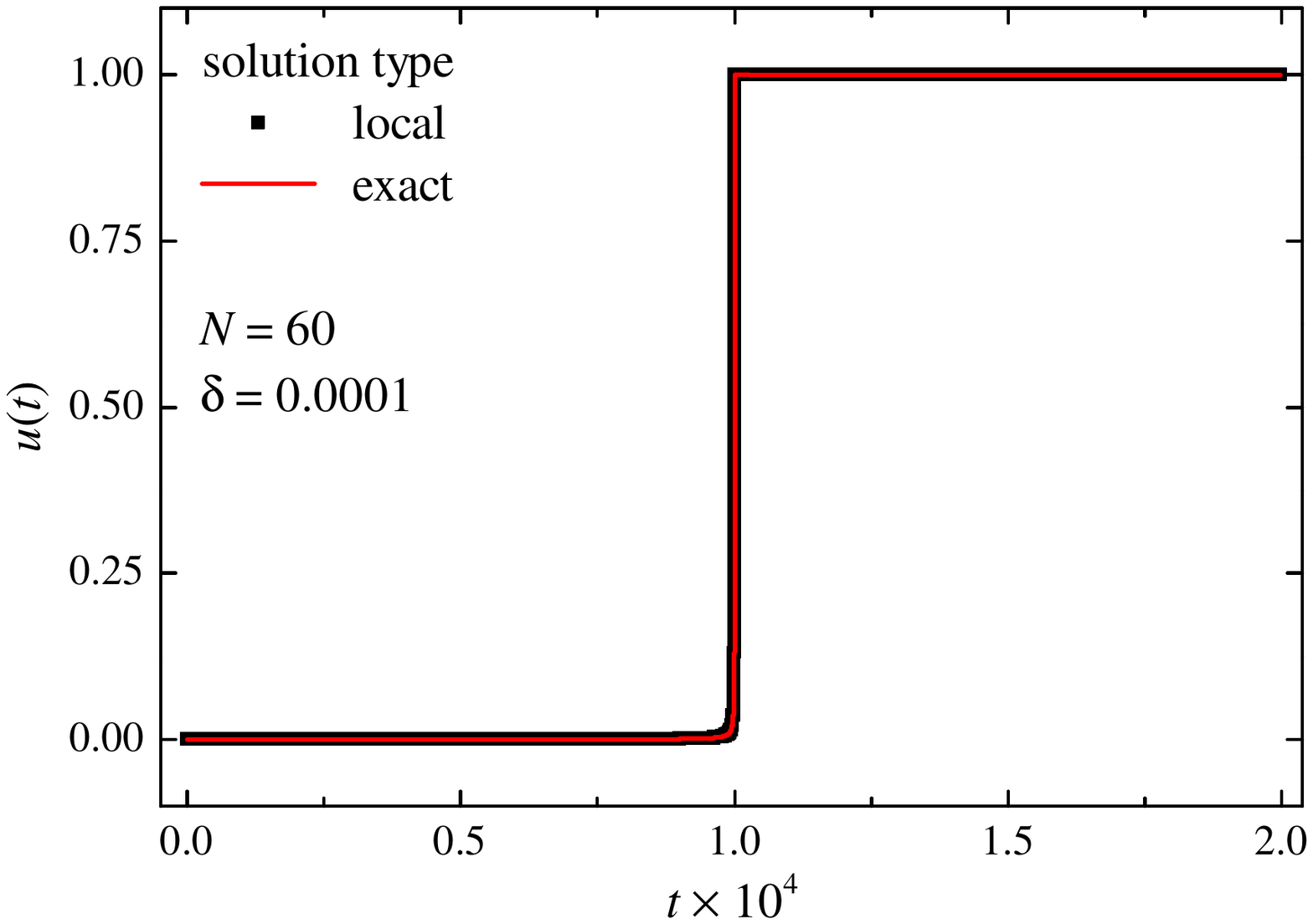}
\vspace{-9mm}\caption{\label{fig:shampine_flame_d_4:c3}}
\end{subfigure}\\
\begin{subfigure}{0.29\textwidth}
\includegraphics[width=\textwidth]{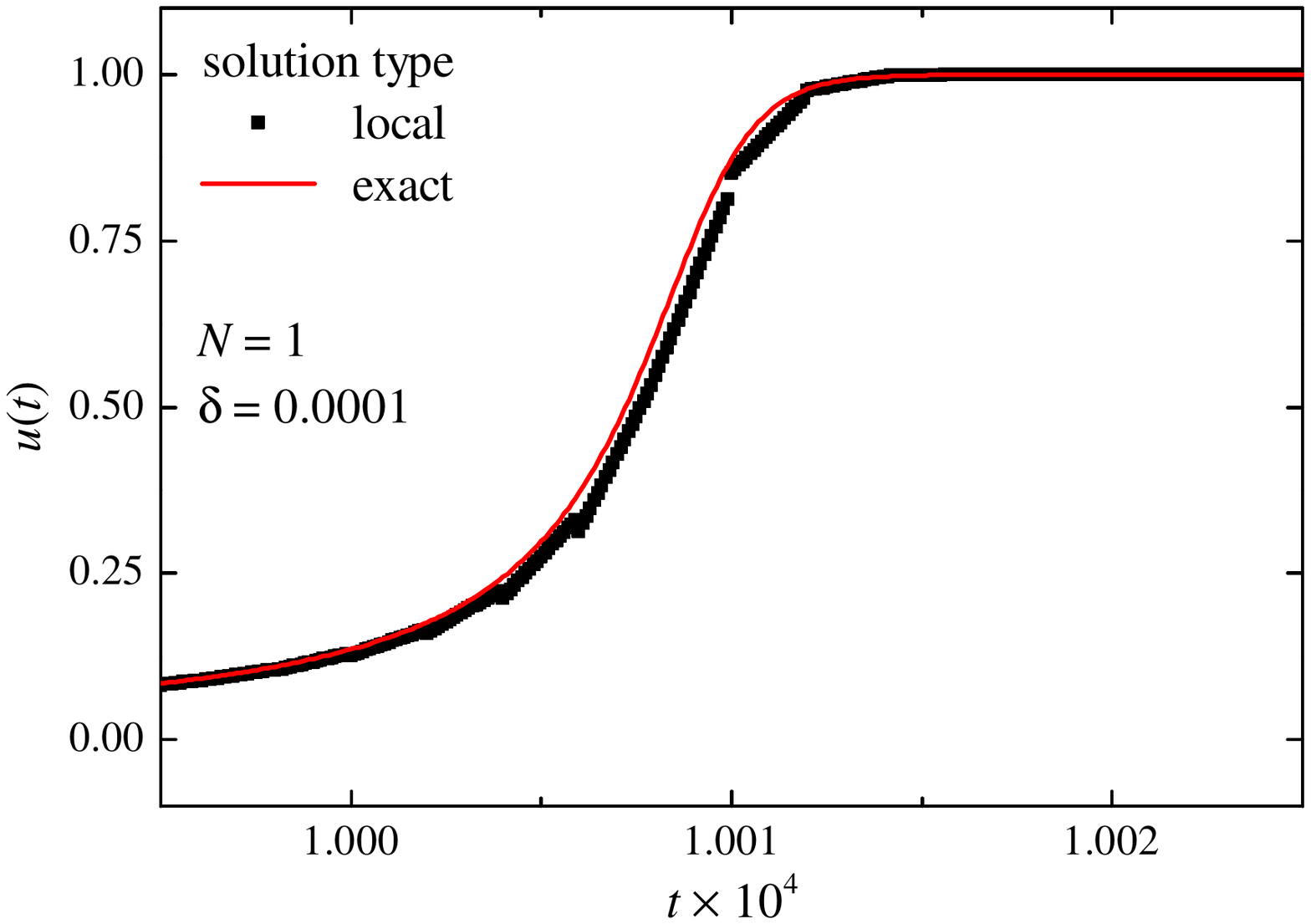}
\vspace{-9mm}\caption{\label{fig:shampine_flame_d_4:d1}}
\end{subfigure}\hspace{5mm}
\begin{subfigure}{0.29\textwidth}
\includegraphics[width=\textwidth]{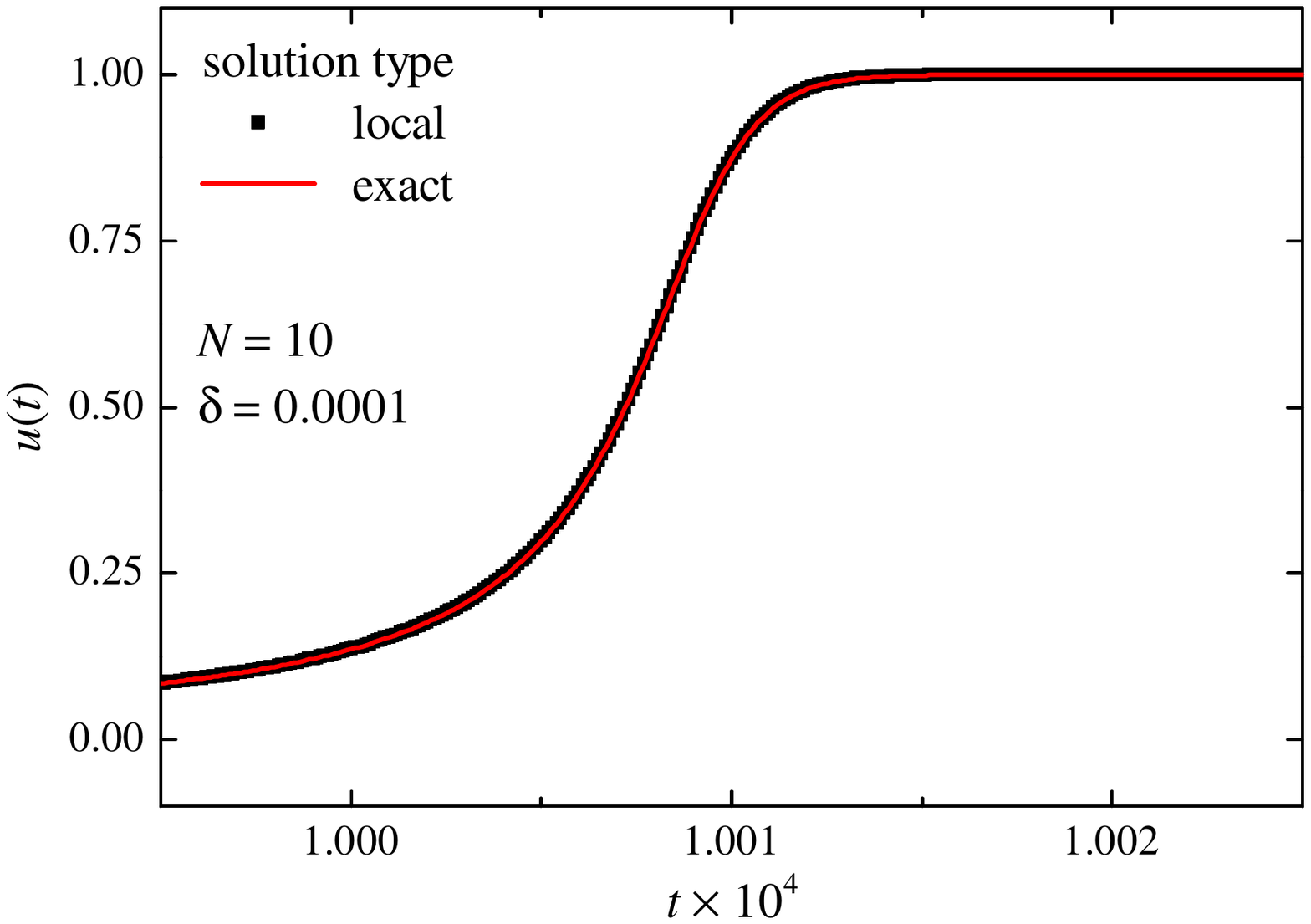}
\vspace{-9mm}\caption{\label{fig:shampine_flame_d_4:d2}}
\end{subfigure}\hspace{5mm}
\begin{subfigure}{0.29\textwidth}
\includegraphics[width=\textwidth]{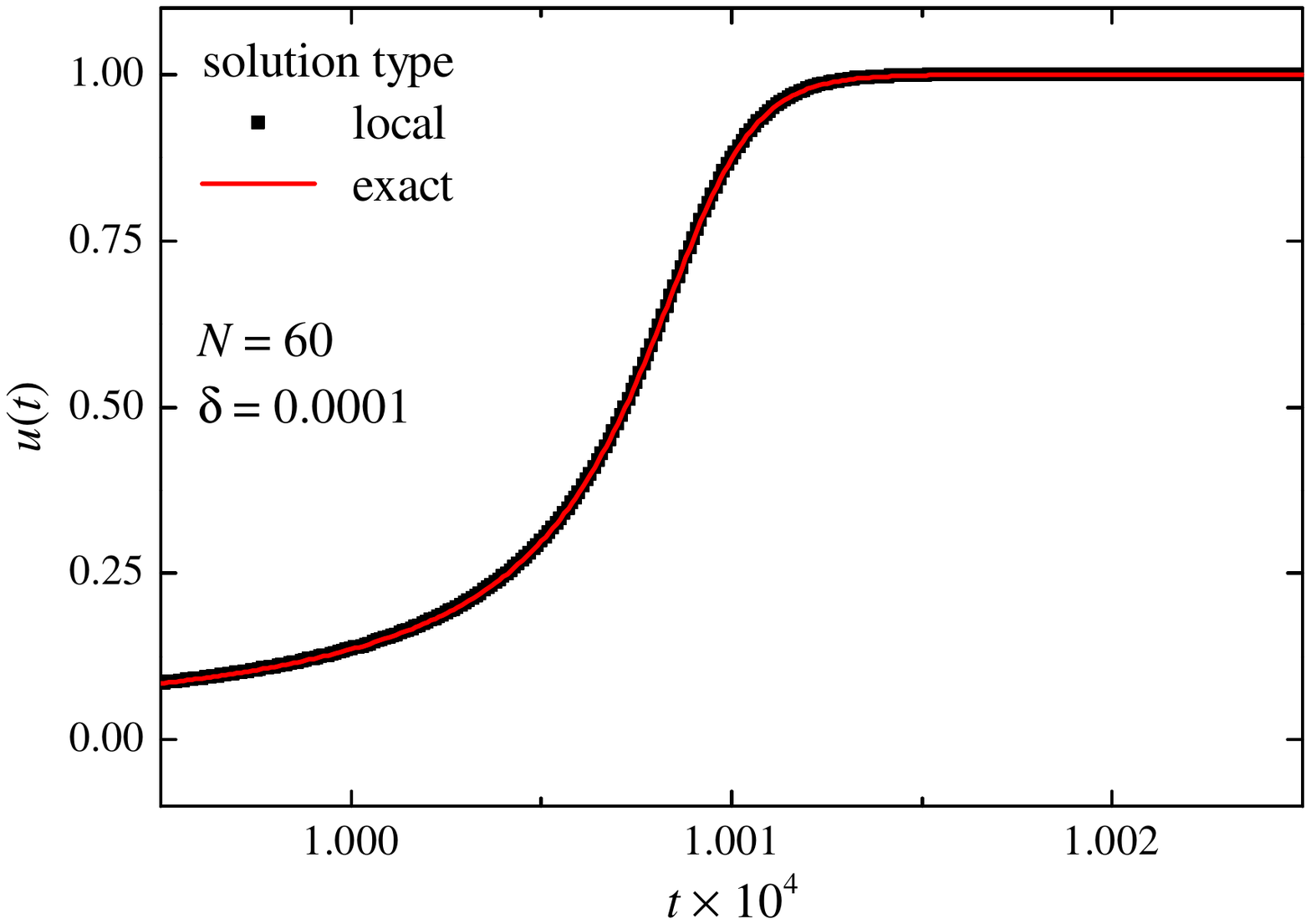}
\vspace{-9mm}\caption{\label{fig:shampine_flame_d_4:d3}}
\end{subfigure}\\
\begin{subfigure}{0.29\textwidth}
\includegraphics[width=\textwidth]{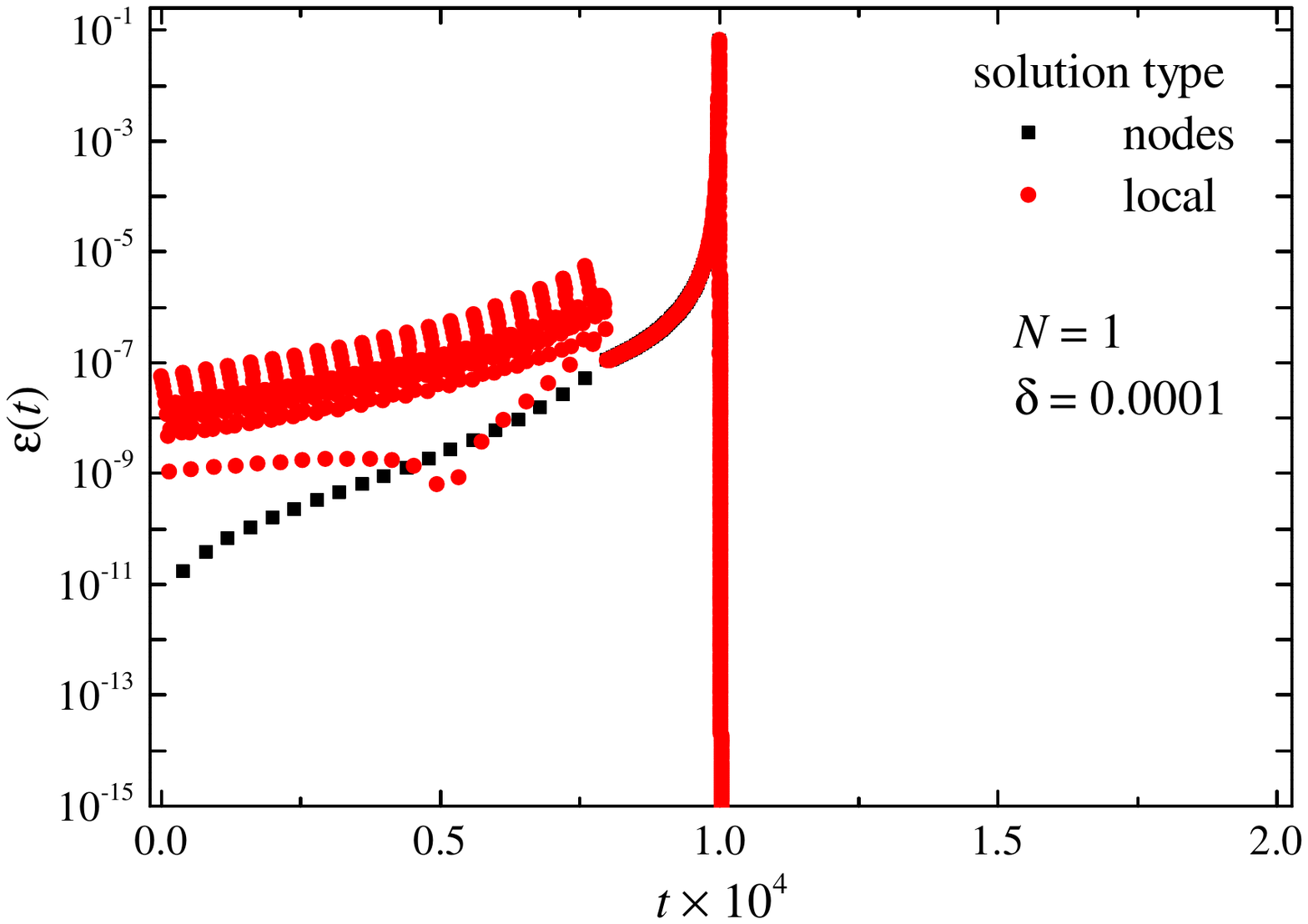}
\vspace{-9mm}\caption{\label{fig:shampine_flame_d_4:e1}}
\end{subfigure}\hspace{5mm}
\begin{subfigure}{0.29\textwidth}
\includegraphics[width=\textwidth]{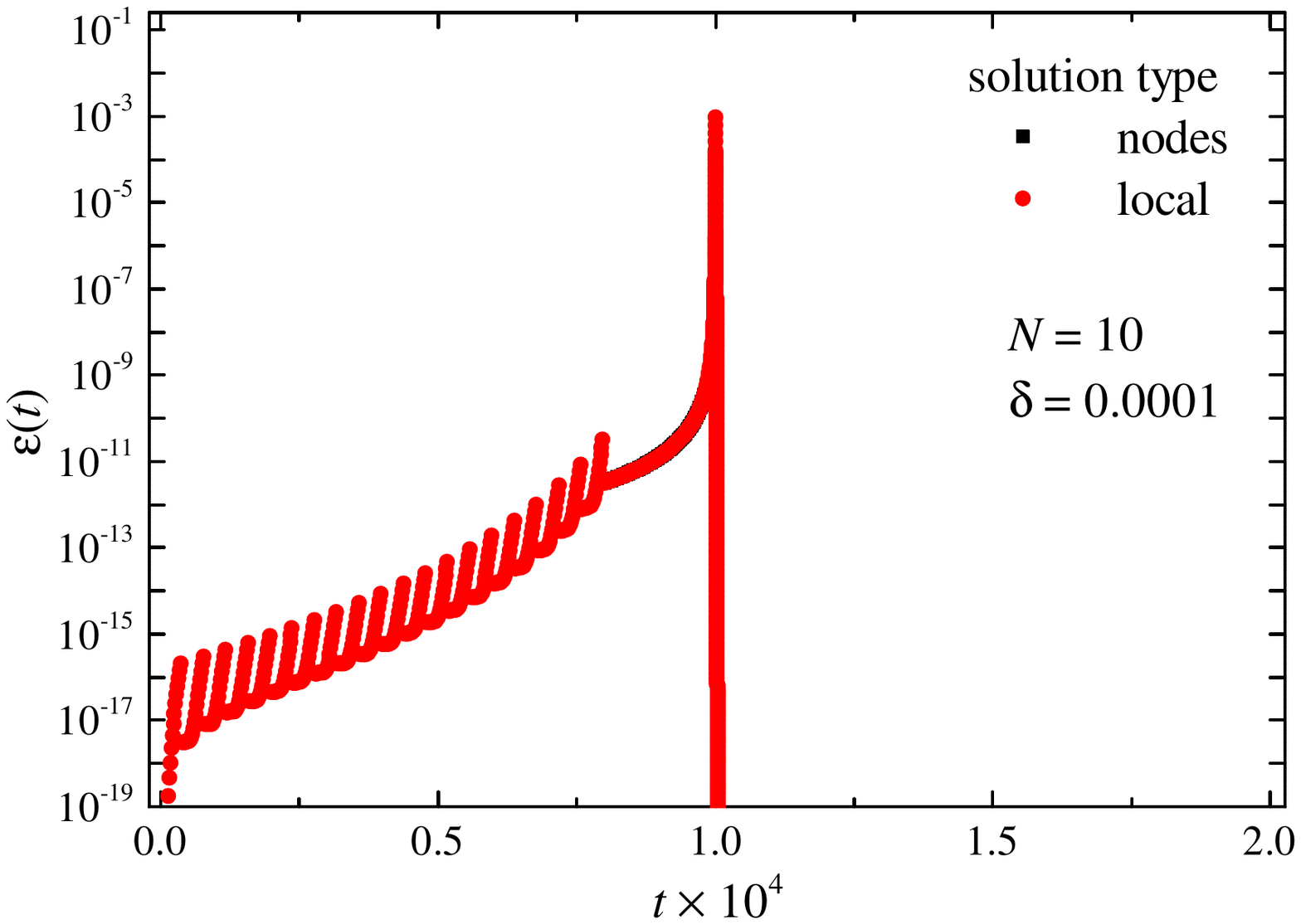}
\vspace{-9mm}\caption{\label{fig:shampine_flame_d_4:e2}}
\end{subfigure}\hspace{5mm}
\begin{subfigure}{0.29\textwidth}
\includegraphics[width=\textwidth]{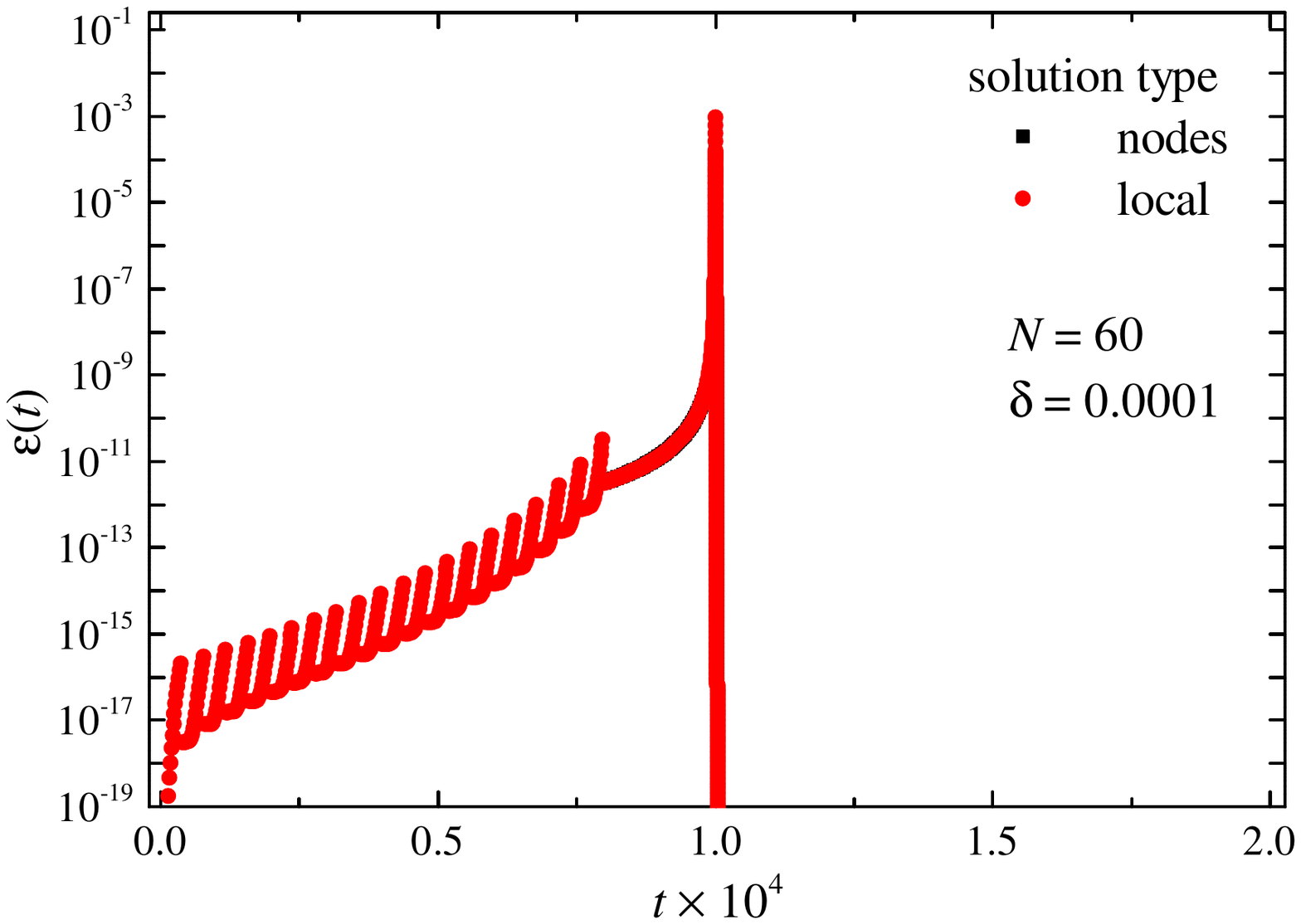}
\vspace{-9mm}\caption{\label{fig:shampine_flame_d_4:e3}}
\end{subfigure}\\
\caption{%
Numerical solution of the problem (\ref{eq:shampine_flame}) with $\delta = 10^{-4}$ in the domain $0 \leqslant t \leqslant 2/\delta$ with $2040$ grid cells: $20$ each in subdomains $t \in [0, 0.4/\delta]$ and $t \in [0.6/\delta, 2/\delta]$ and $2000$ in subdomain $t \in [0.4/\delta, 0.6/\delta]$, using the ADER-DG numerical method with a local DG predictor with degrees of polynomials $N = 1$ (left column: \subref{fig:shampine_flame_d_4:a1}, \subref{fig:shampine_flame_d_4:b1}, \subref{fig:shampine_flame_d_4:c1}, \subref{fig:shampine_flame_d_4:d1}, \subref{fig:shampine_flame_d_4:e1}), $N = 10$ (center column: \subref{fig:shampine_flame_d_4:a2}, \subref{fig:shampine_flame_d_4:b2}, \subref{fig:shampine_flame_d_4:c2}, \subref{fig:shampine_flame_d_4:d2}, \subref{fig:shampine_flame_d_4:e2}) and $N = 60$ (right column: \subref{fig:shampine_flame_d_4:a3}, \subref{fig:shampine_flame_d_4:b3}, \subref{fig:shampine_flame_d_4:c3}, \subref{fig:shampine_flame_d_4:d3}, \subref{fig:shampine_flame_d_4:e3}). Presented data: numerical solution $u(t)$ at nodes in the whole domain $t \in [0, 2/\delta]$ (\subref{fig:shampine_flame_d_4:a1}, \subref{fig:shampine_flame_d_4:a2}, \subref{fig:shampine_flame_d_4:a3}) and in the subdomain $t \in [0.9995, 1.0025]$ (\subref{fig:shampine_flame_d_4:b1}, \subref{fig:shampine_flame_d_4:b2}, \subref{fig:shampine_flame_d_4:b3}); local solution $u(t)$ in the whole domain $t \in [0, 2/\delta]$ (\subref{fig:shampine_flame_d_4:c1}, \subref{fig:shampine_flame_d_4:c2}, \subref{fig:shampine_flame_d_4:c3}) and in the subdomain $t \in [0.9995, 1.0025]$ (\subref{fig:shampine_flame_d_4:d1}, \subref{fig:shampine_flame_d_4:d2}, \subref{fig:shampine_flame_d_4:d3}); dependence of the point-wise error $\varepsilon(t)$ of the numerical solution for the solution at the nodes and the local solution (\subref{fig:shampine_flame_d_4:e1}, \subref{fig:shampine_flame_d_4:e2}, \subref{fig:shampine_flame_d_4:e3}); red line is the exact solution (\ref{eq:shampine_flame:exact_solution}).
}
\label{fig:shampine_flame_d_4}
\end{figure}

\begin{figure}[h!]
\captionsetup[subfigure]{%
position=bottom,
textfont=normalfont,
singlelinecheck=off,
justification=raggedright,
font=tiny
}
\centering
\begin{subfigure}{0.29\textwidth}
\includegraphics[width=\textwidth]{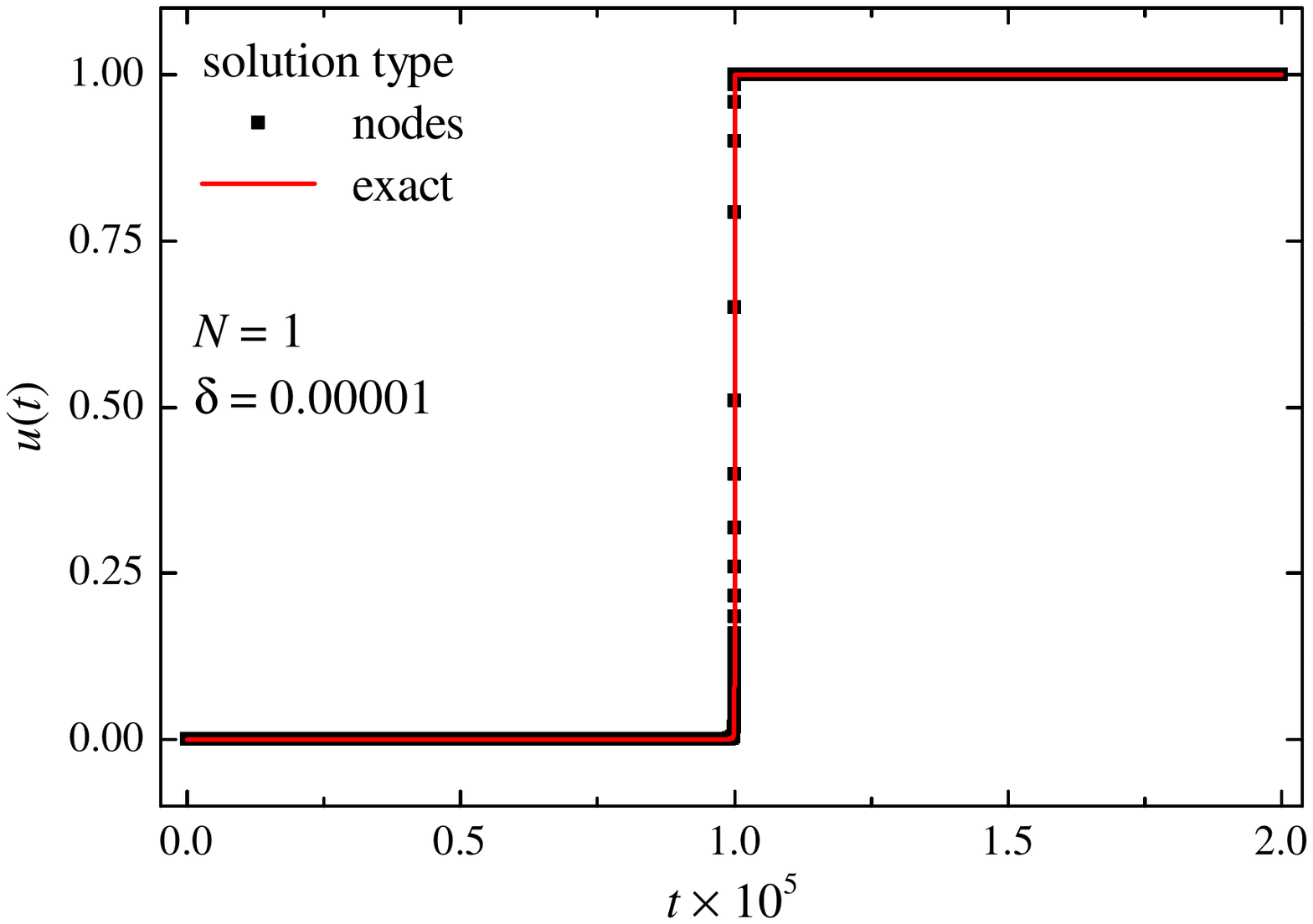}
\vspace{-9mm}\caption{\label{fig:shampine_flame_d_5:a1}}
\end{subfigure}\hspace{5mm}
\begin{subfigure}{0.29\textwidth}
\includegraphics[width=\textwidth]{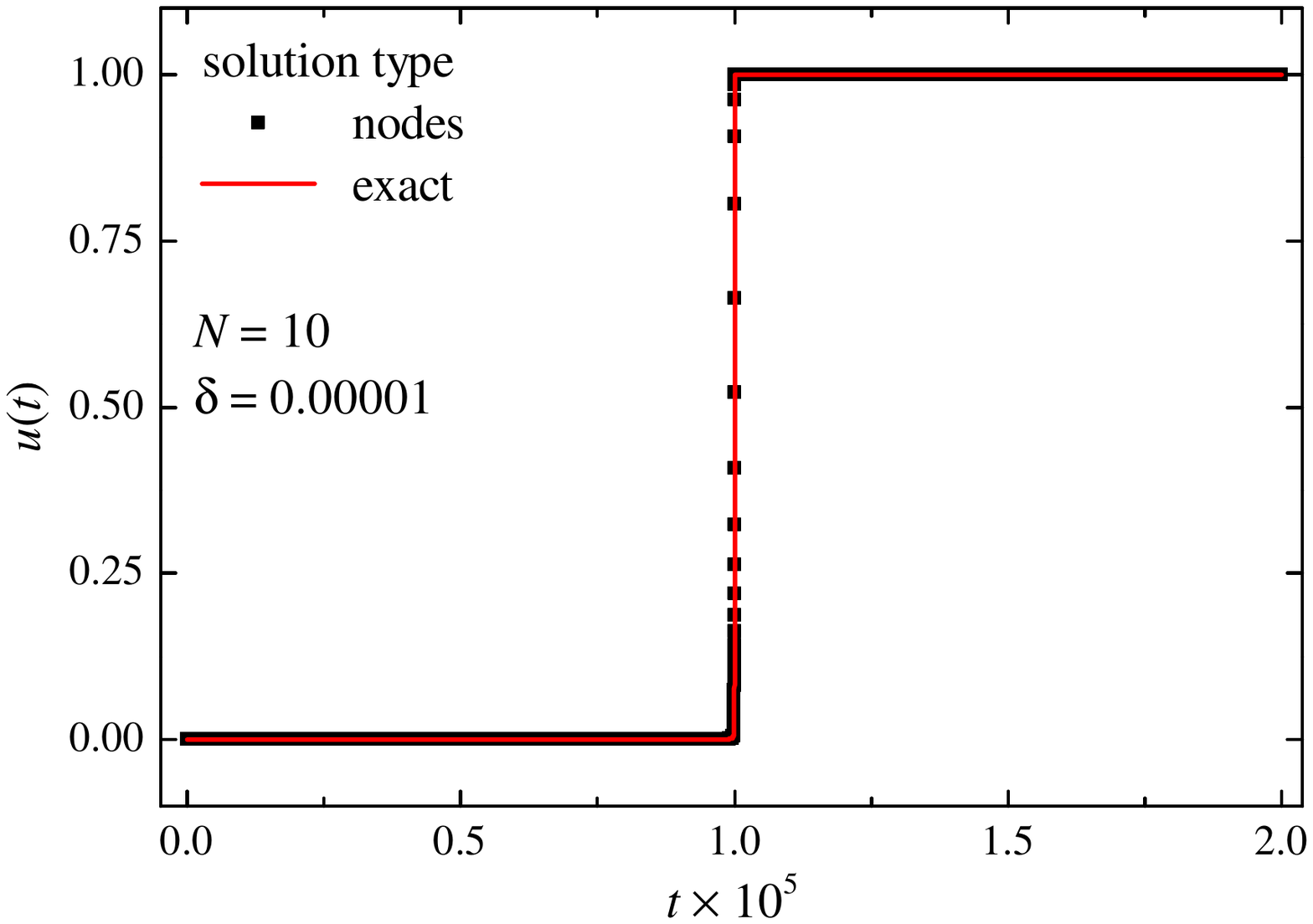}
\vspace{-9mm}\caption{\label{fig:shampine_flame_d_5:a2}}
\end{subfigure}\hspace{5mm}
\begin{subfigure}{0.29\textwidth}
\includegraphics[width=\textwidth]{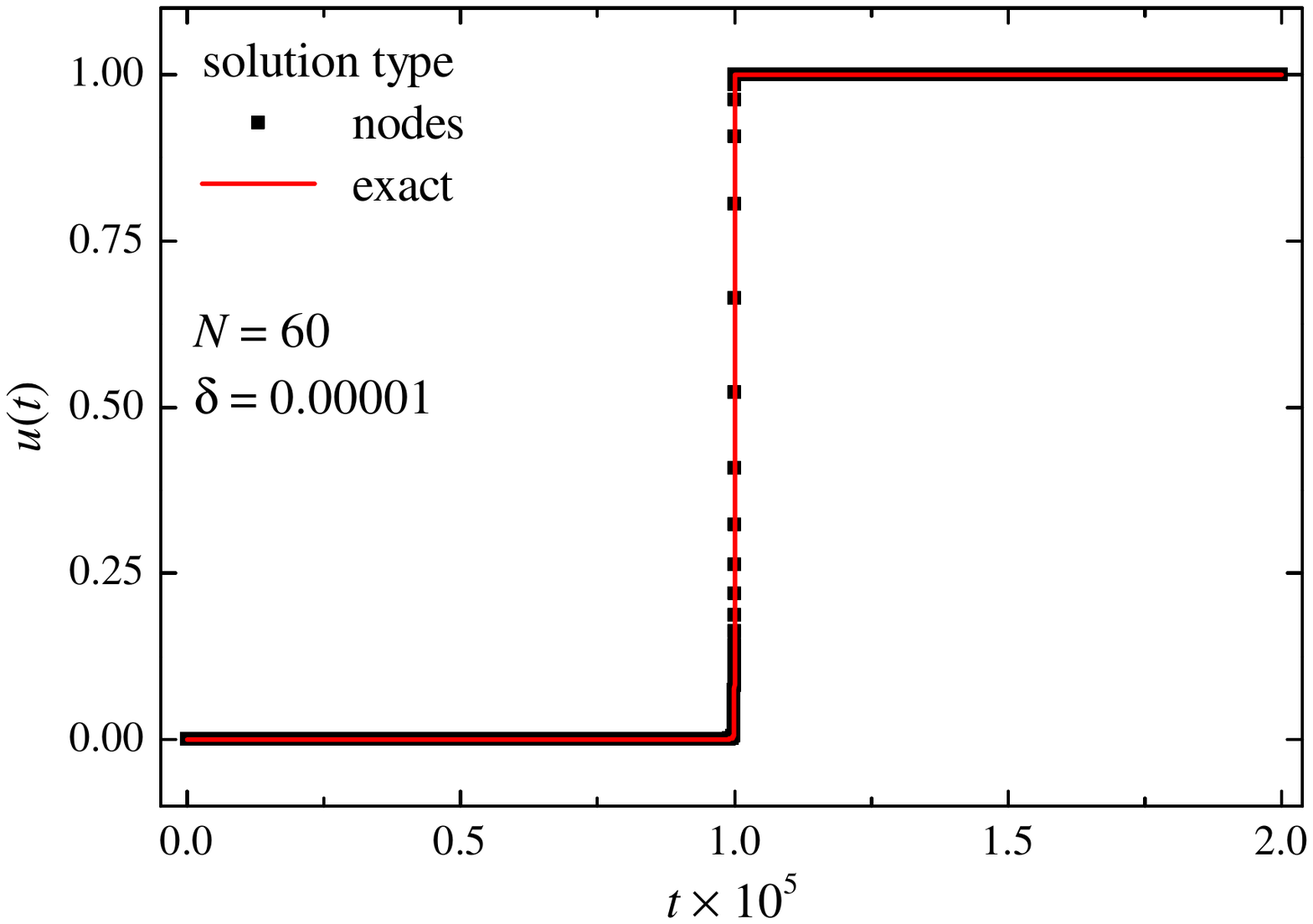}
\vspace{-9mm}\caption{\label{fig:shampine_flame_d_5:a3}}
\end{subfigure}\\
\begin{subfigure}{0.29\textwidth}
\includegraphics[width=\textwidth]{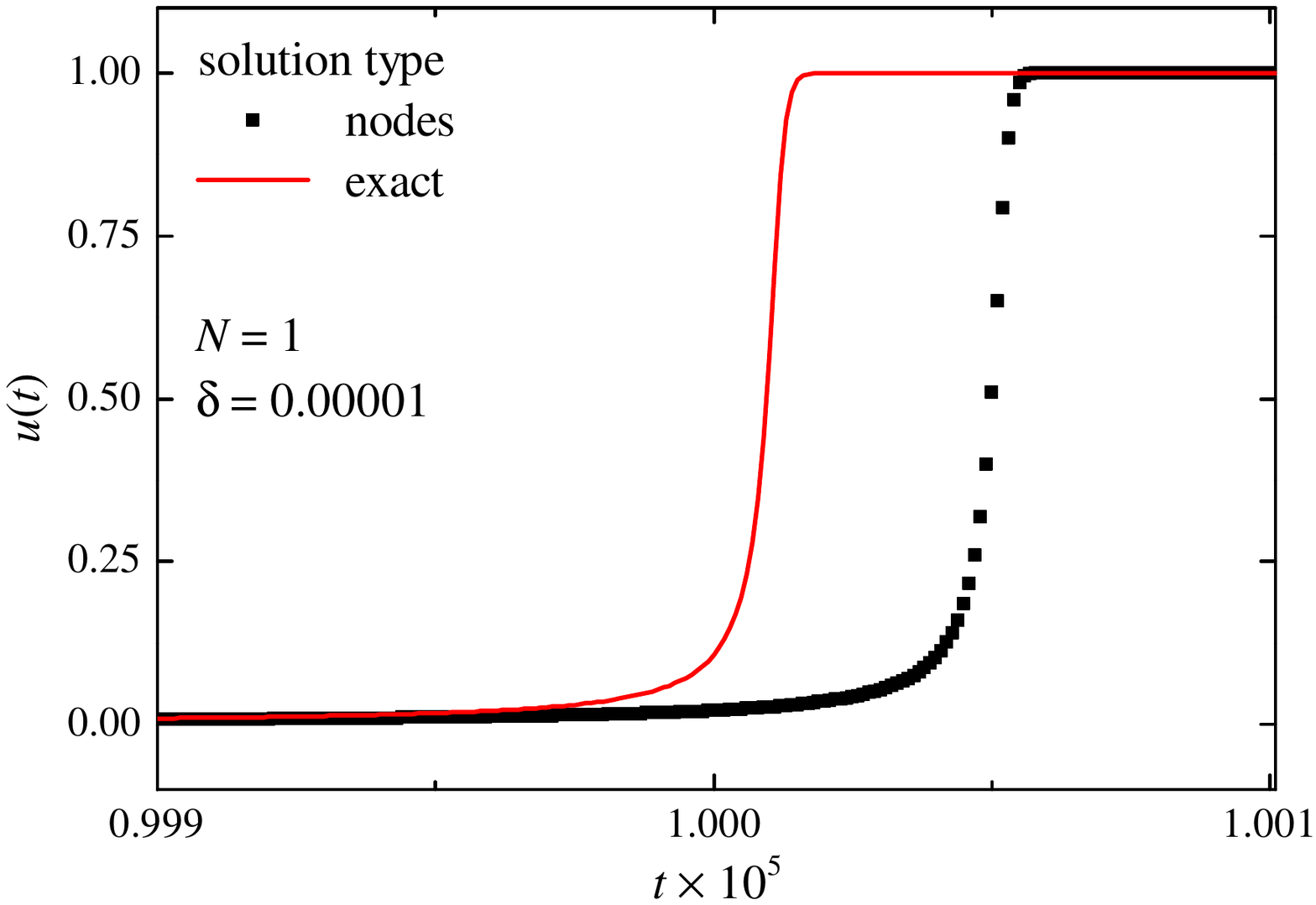}
\vspace{-9mm}\caption{\label{fig:shampine_flame_d_5:b1}}
\end{subfigure}\hspace{5mm}
\begin{subfigure}{0.29\textwidth}
\includegraphics[width=\textwidth]{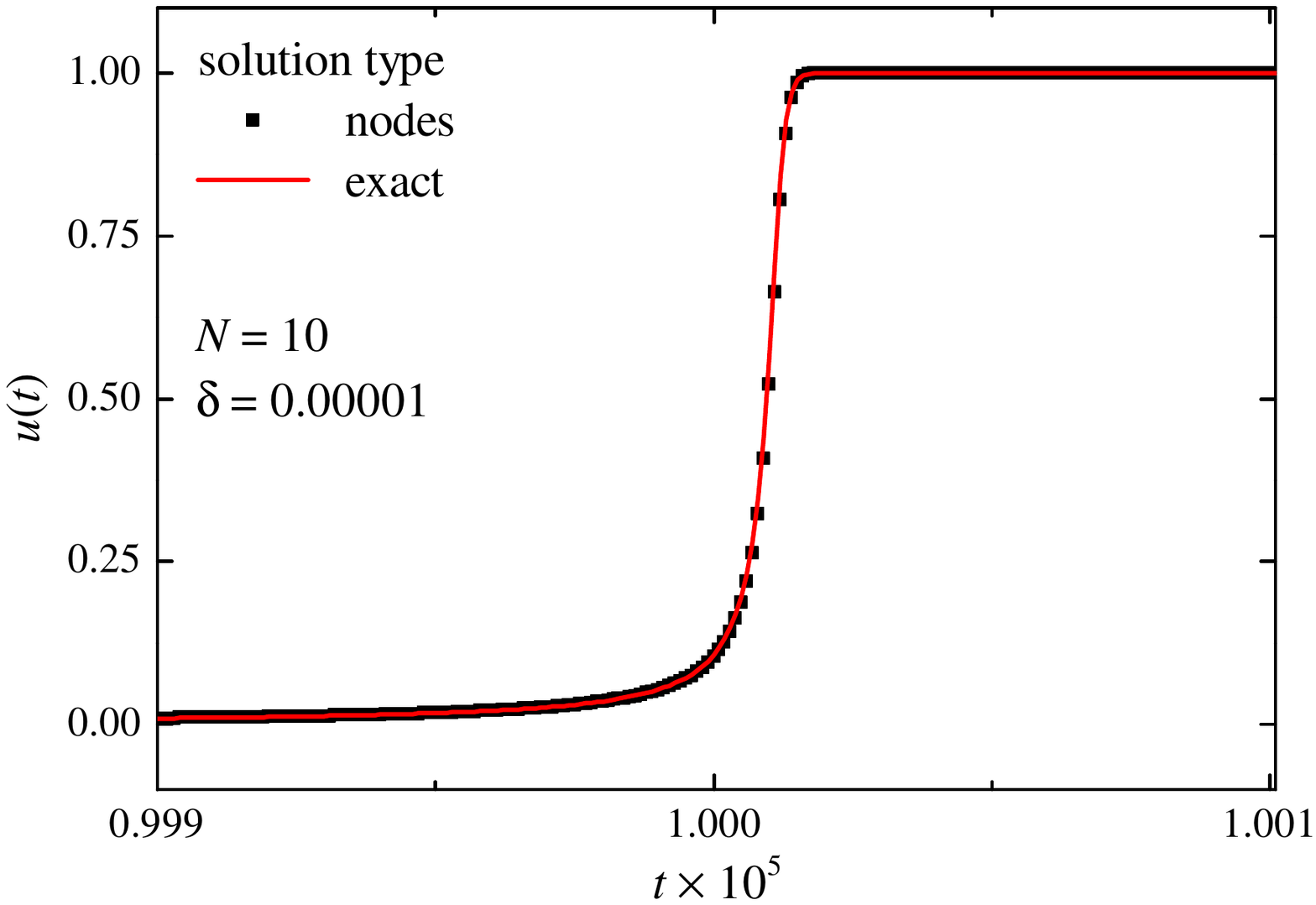}
\vspace{-9mm}\caption{\label{fig:shampine_flame_d_5:b2}}
\end{subfigure}\hspace{5mm}
\begin{subfigure}{0.29\textwidth}
\includegraphics[width=\textwidth]{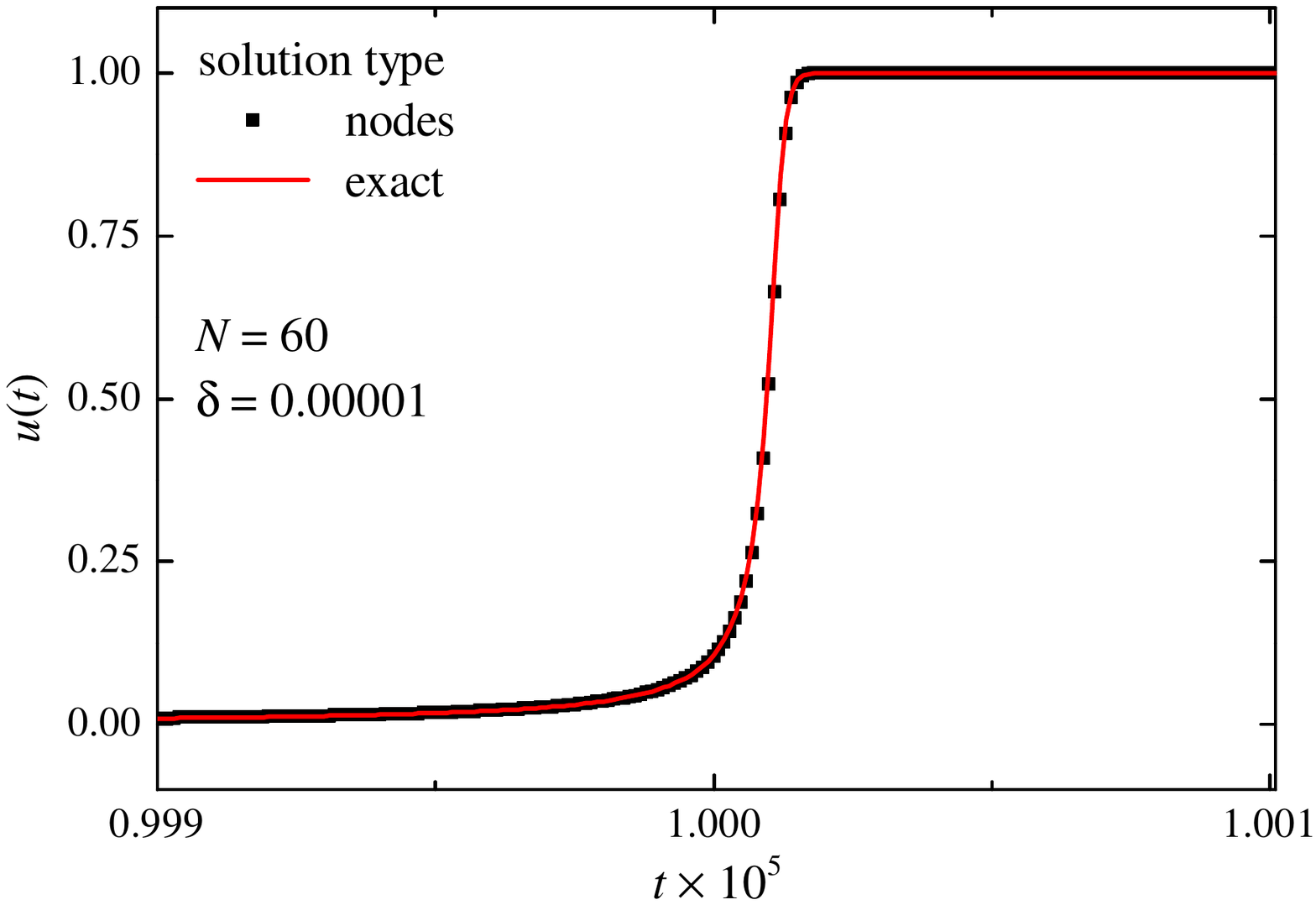}
\vspace{-9mm}\caption{\label{fig:shampine_flame_d_5:b3}}
\end{subfigure}\\
\begin{subfigure}{0.29\textwidth}
\includegraphics[width=\textwidth]{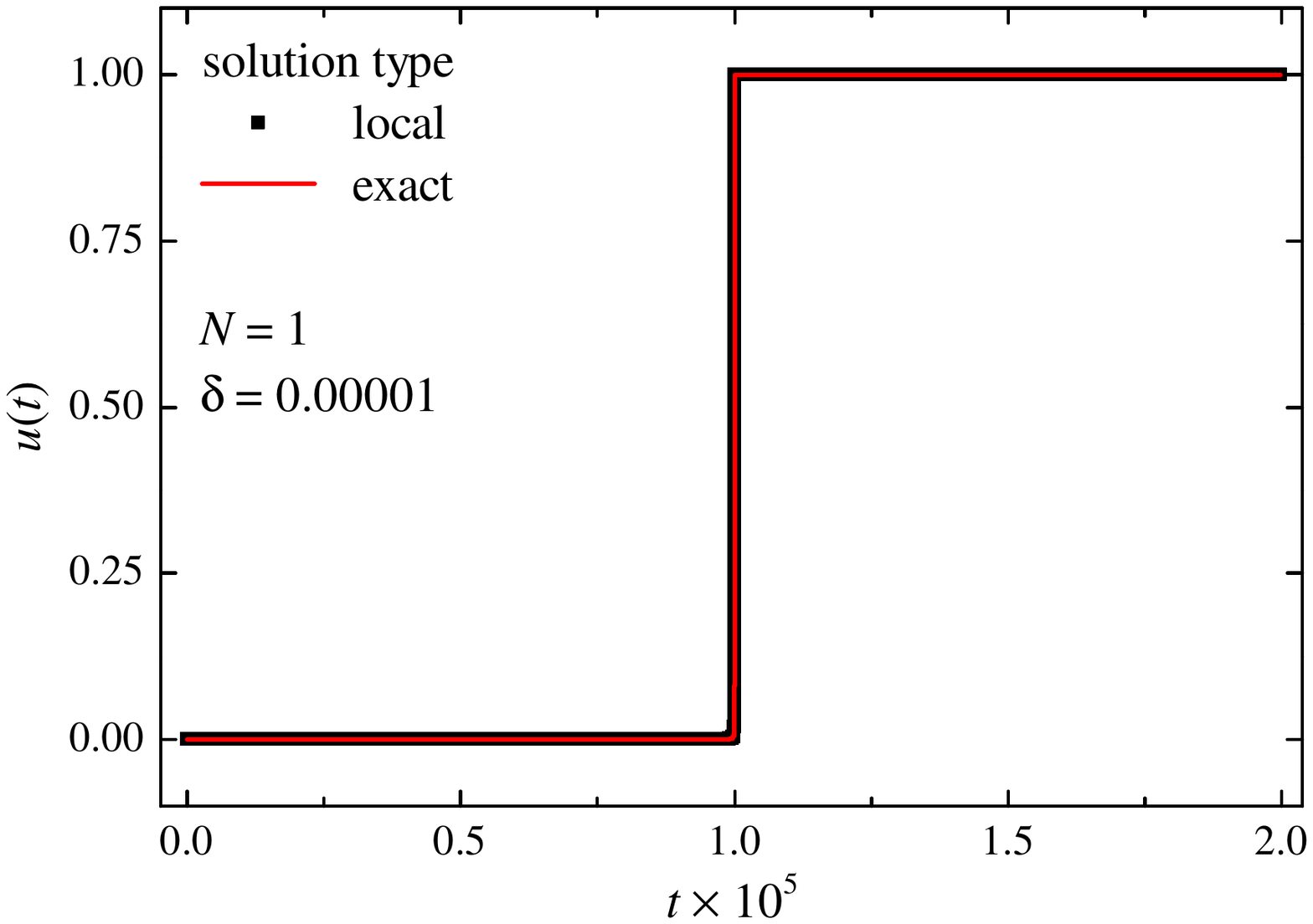}
\vspace{-9mm}\caption{\label{fig:shampine_flame_d_5:c1}}
\end{subfigure}\hspace{5mm}
\begin{subfigure}{0.29\textwidth}
\includegraphics[width=\textwidth]{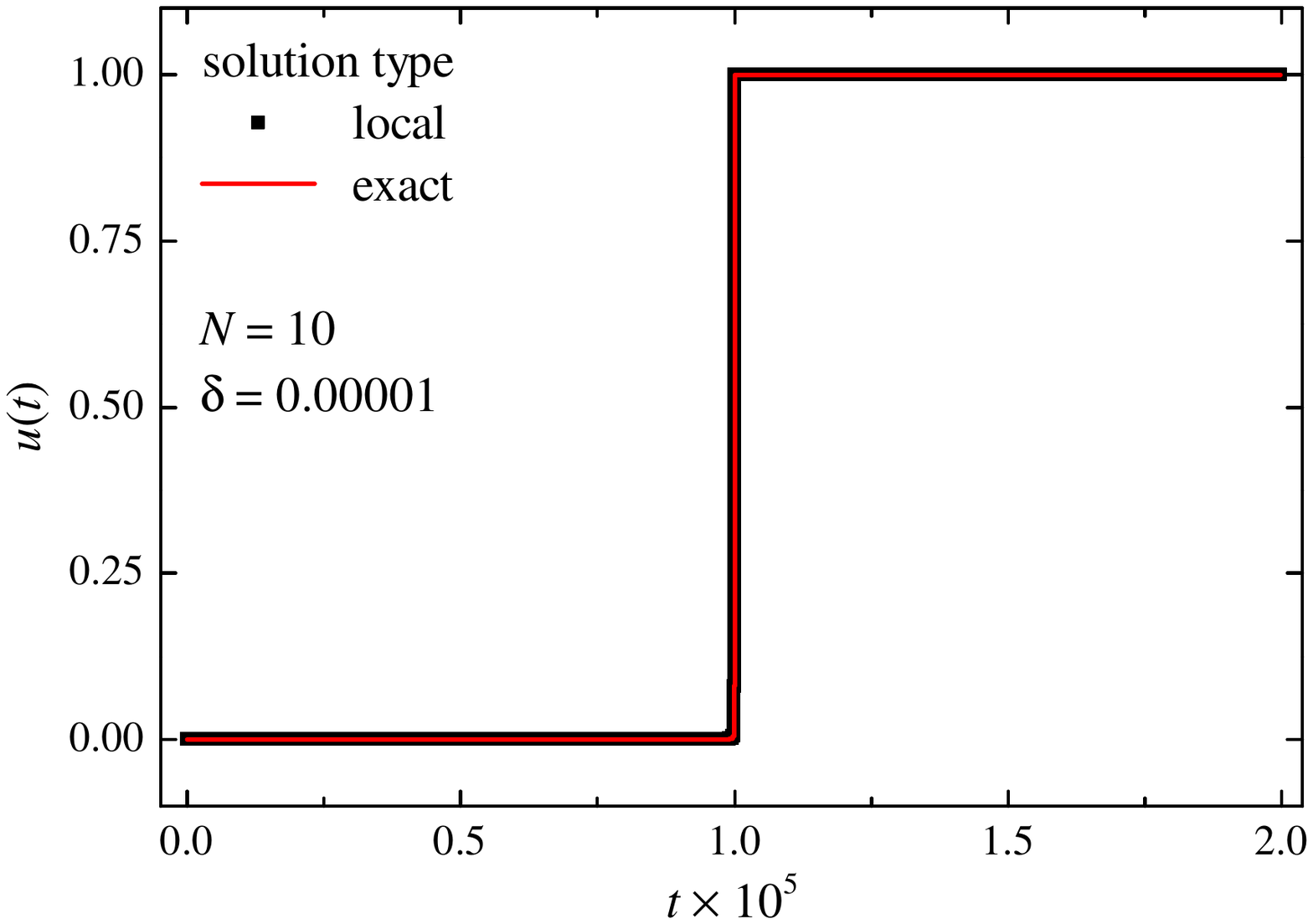}
\vspace{-9mm}\caption{\label{fig:shampine_flame_d_5:c2}}
\end{subfigure}\hspace{5mm}
\begin{subfigure}{0.29\textwidth}
\includegraphics[width=\textwidth]{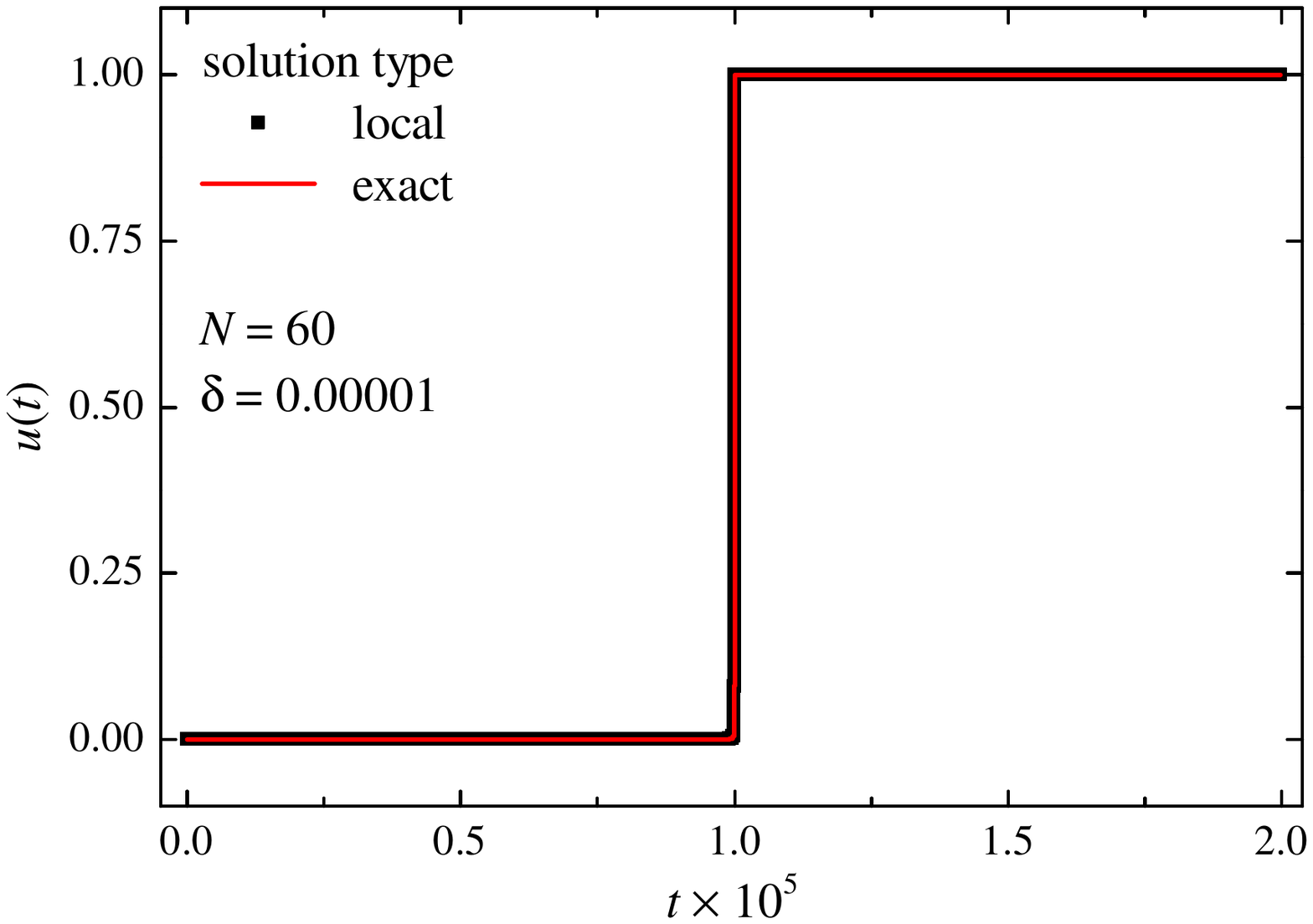}
\vspace{-9mm}\caption{\label{fig:shampine_flame_d_5:c3}}
\end{subfigure}\\
\begin{subfigure}{0.29\textwidth}
\includegraphics[width=\textwidth]{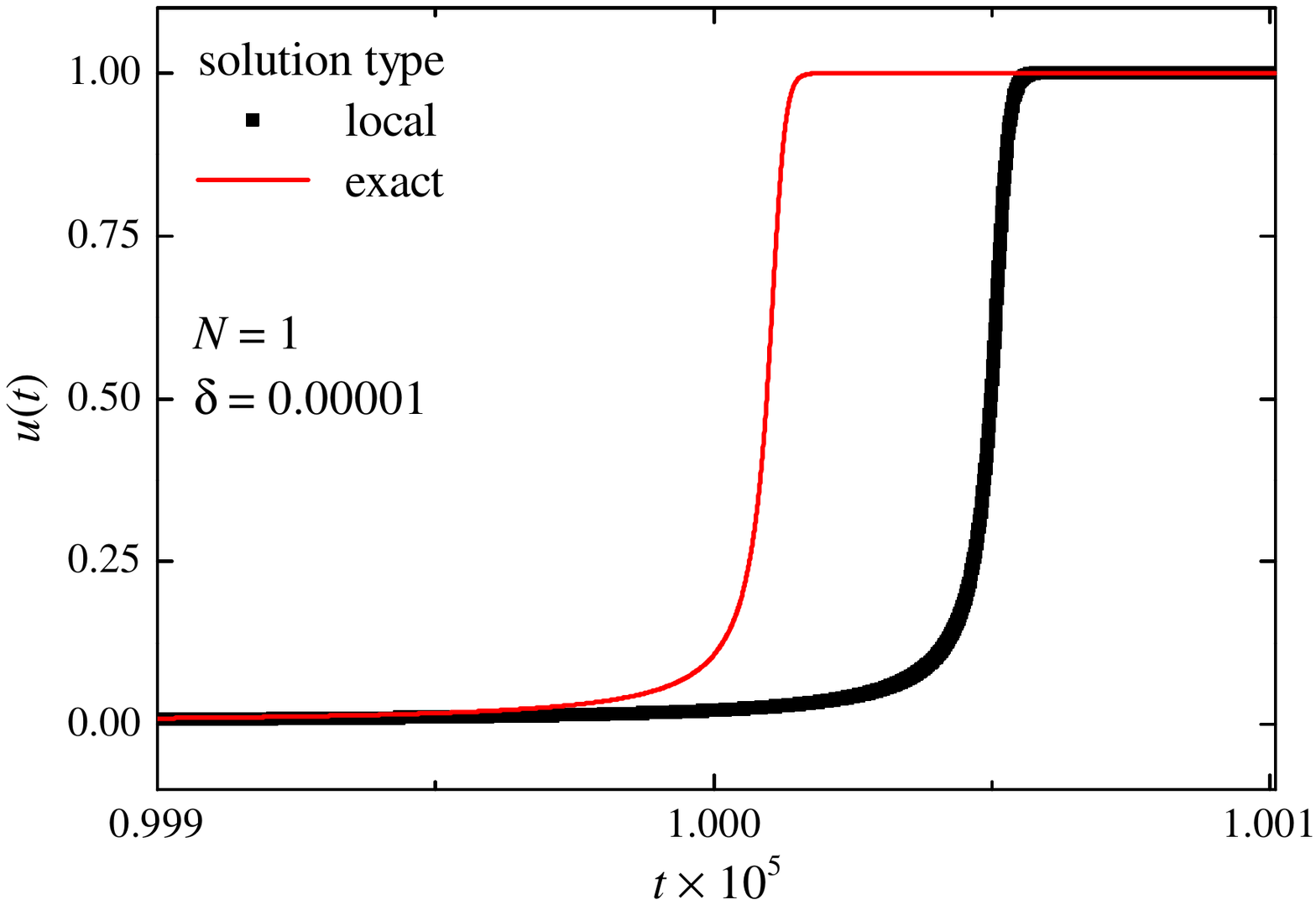}
\vspace{-9mm}\caption{\label{fig:shampine_flame_d_5:d1}}
\end{subfigure}\hspace{5mm}
\begin{subfigure}{0.29\textwidth}
\includegraphics[width=\textwidth]{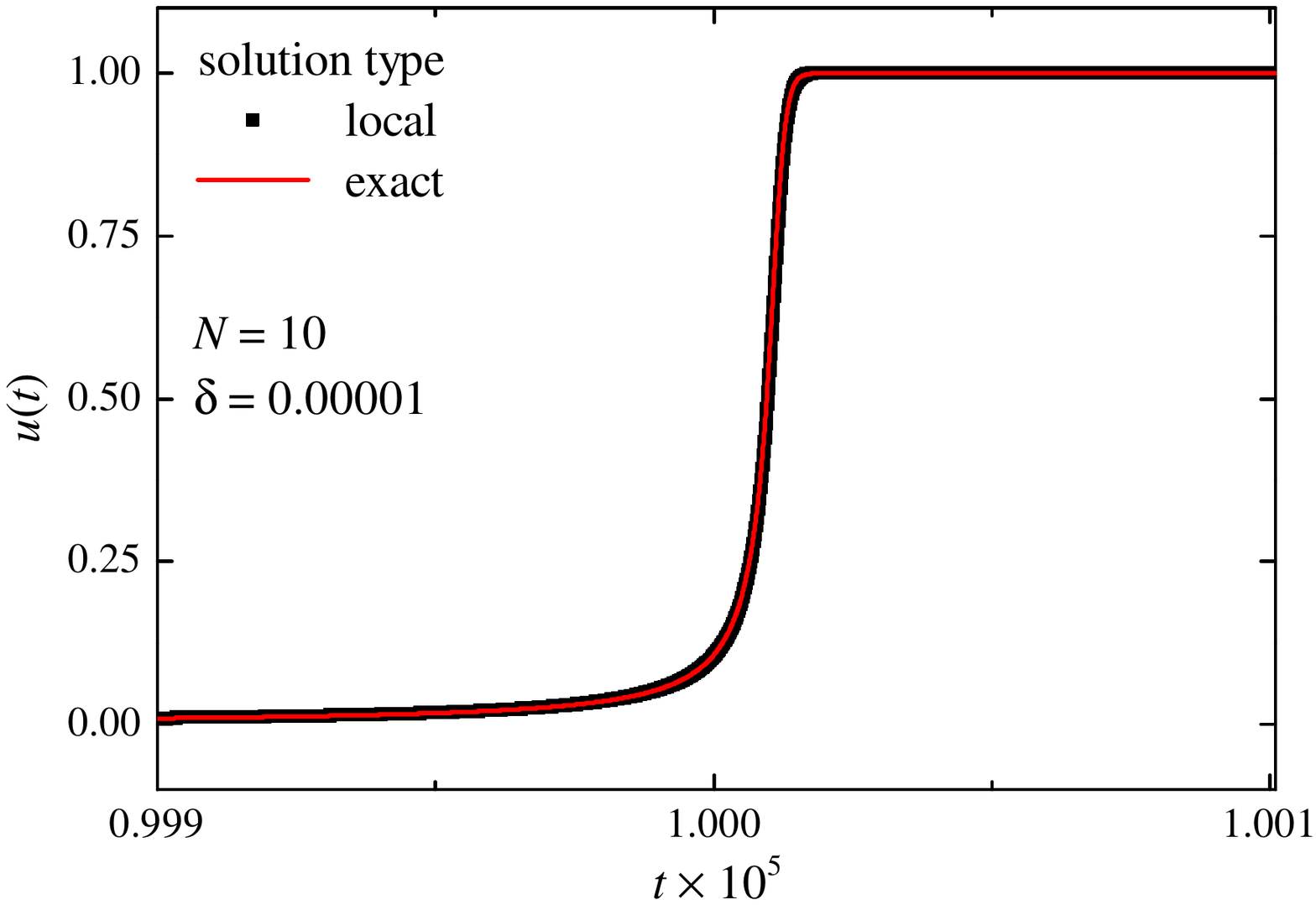}
\vspace{-9mm}\caption{\label{fig:shampine_flame_d_5:d2}}
\end{subfigure}\hspace{5mm}
\begin{subfigure}{0.29\textwidth}
\includegraphics[width=\textwidth]{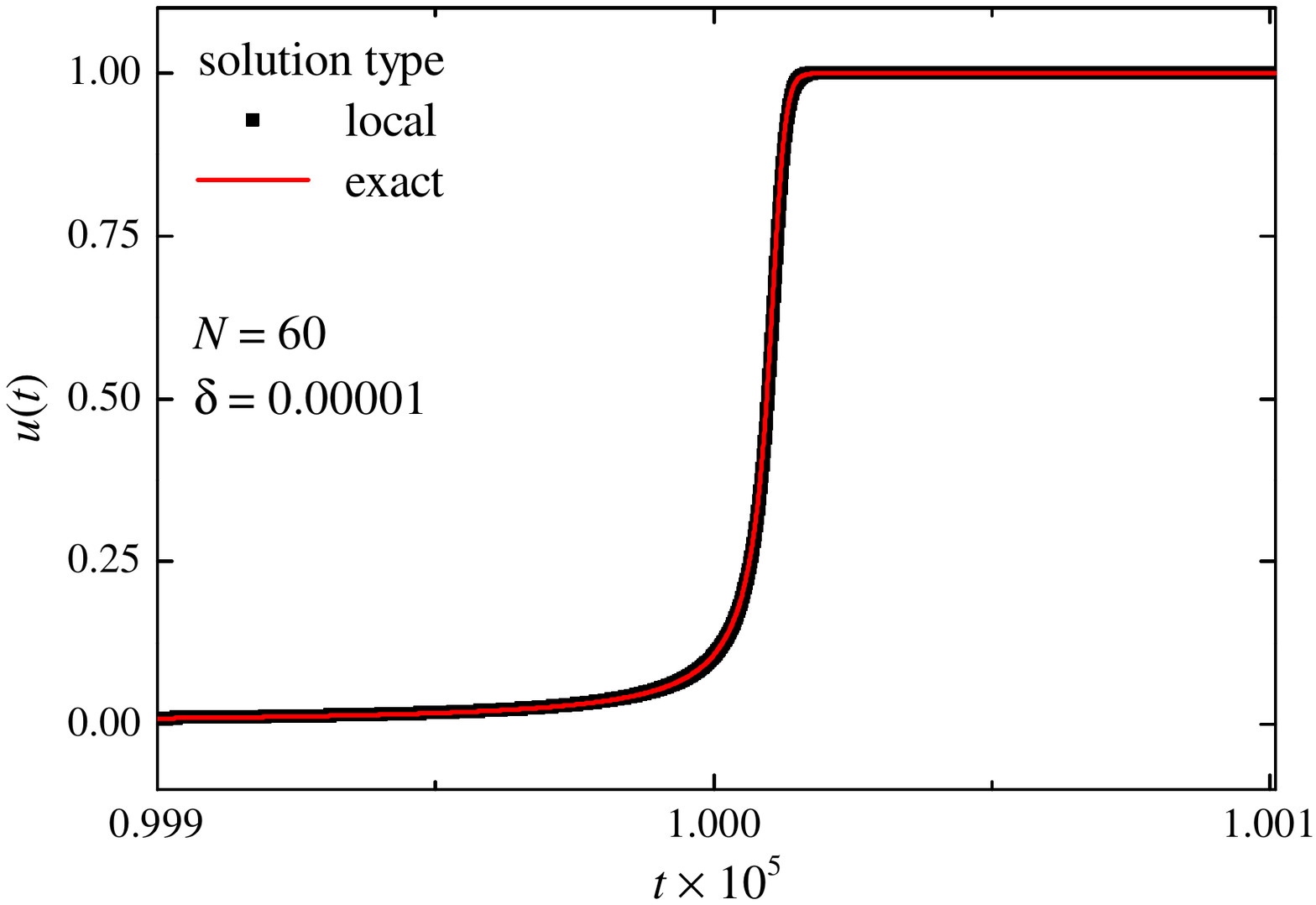}
\vspace{-9mm}\caption{\label{fig:shampine_flame_d_5:d3}}
\end{subfigure}\\
\begin{subfigure}{0.29\textwidth}
\includegraphics[width=\textwidth]{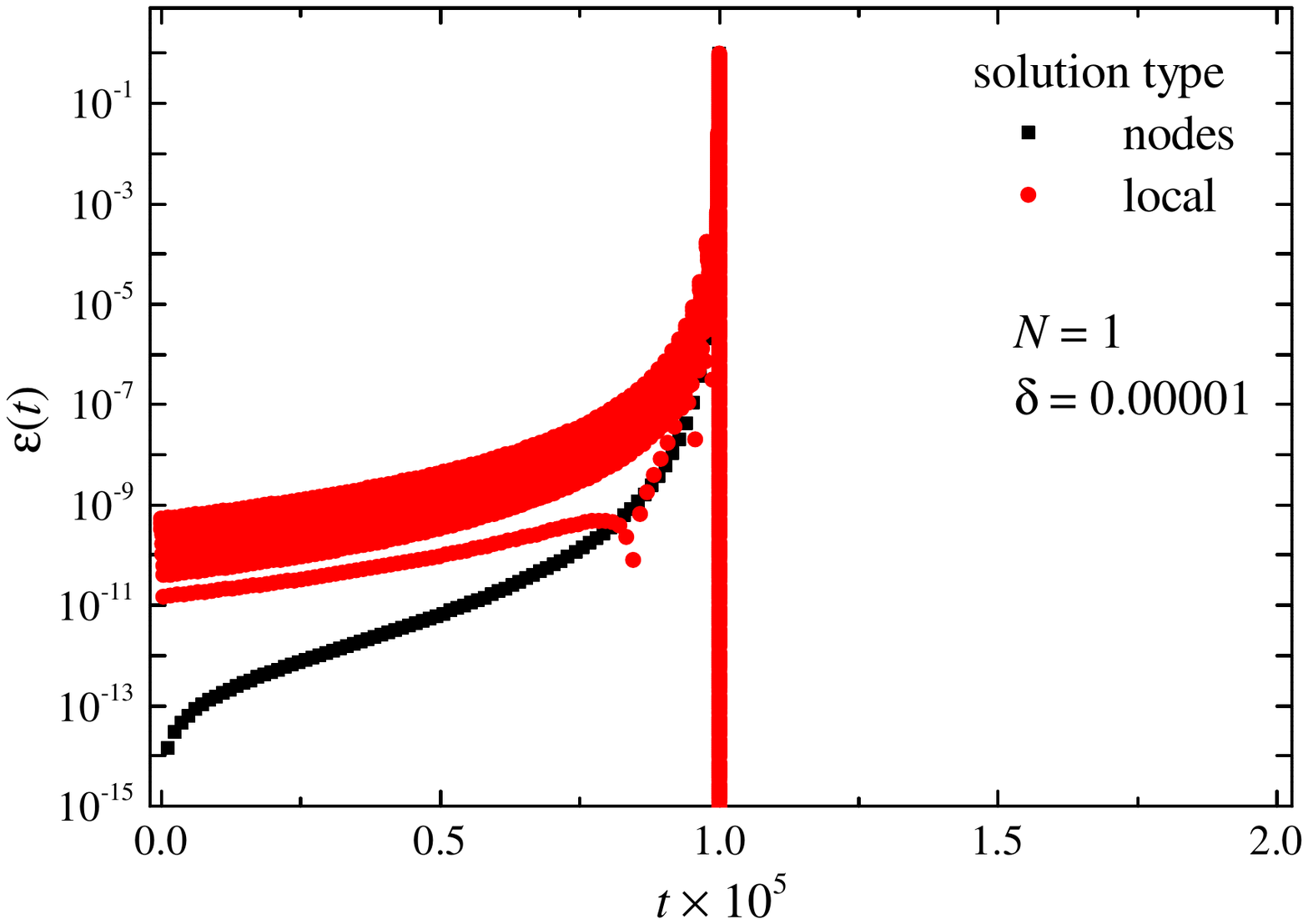}
\vspace{-9mm}\caption{\label{fig:shampine_flame_d_5:e1}}
\end{subfigure}\hspace{5mm}
\begin{subfigure}{0.29\textwidth}
\includegraphics[width=\textwidth]{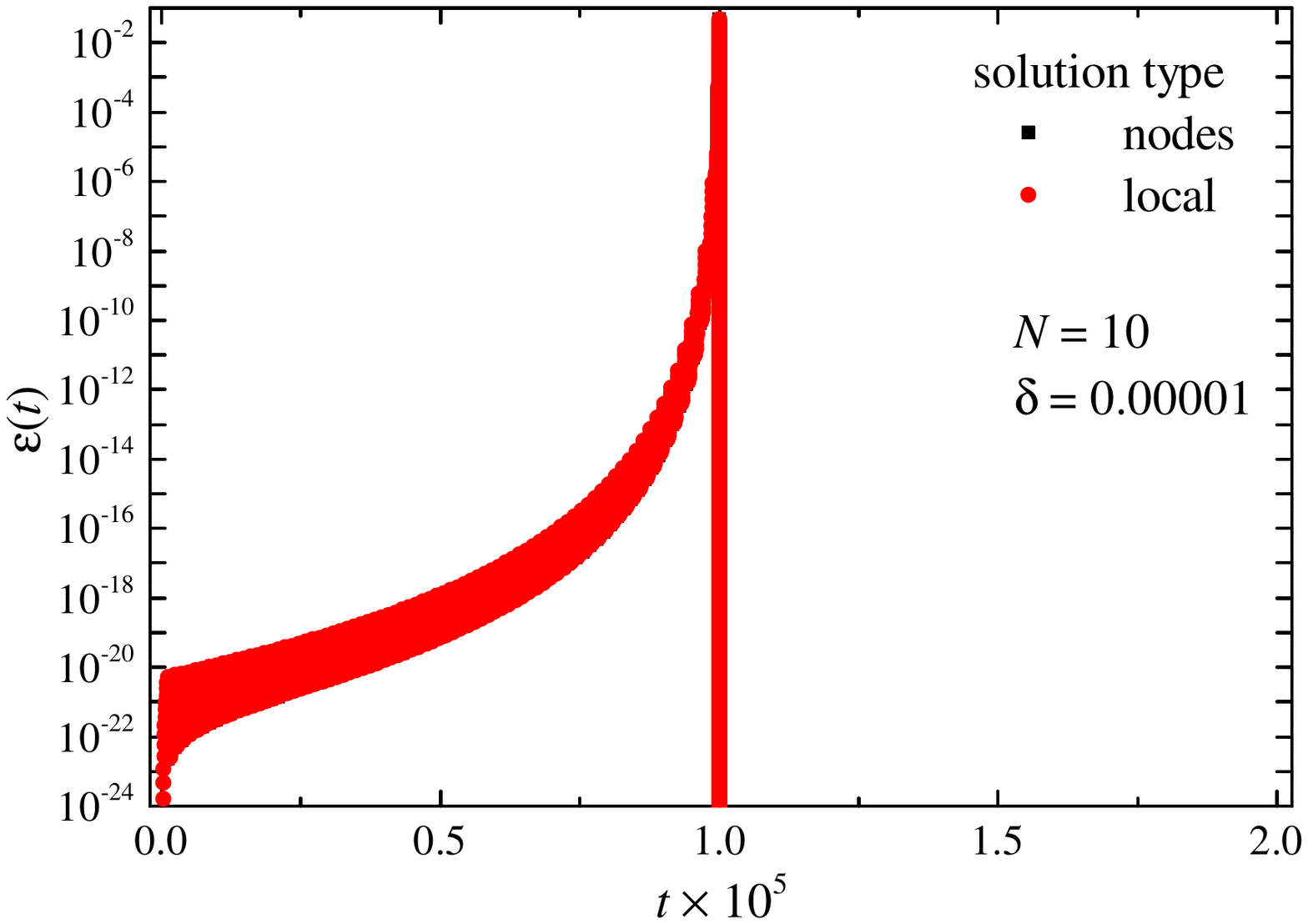}
\vspace{-9mm}\caption{\label{fig:shampine_flame_d_5:e2}}
\end{subfigure}\hspace{5mm}
\begin{subfigure}{0.29\textwidth}
\includegraphics[width=\textwidth]{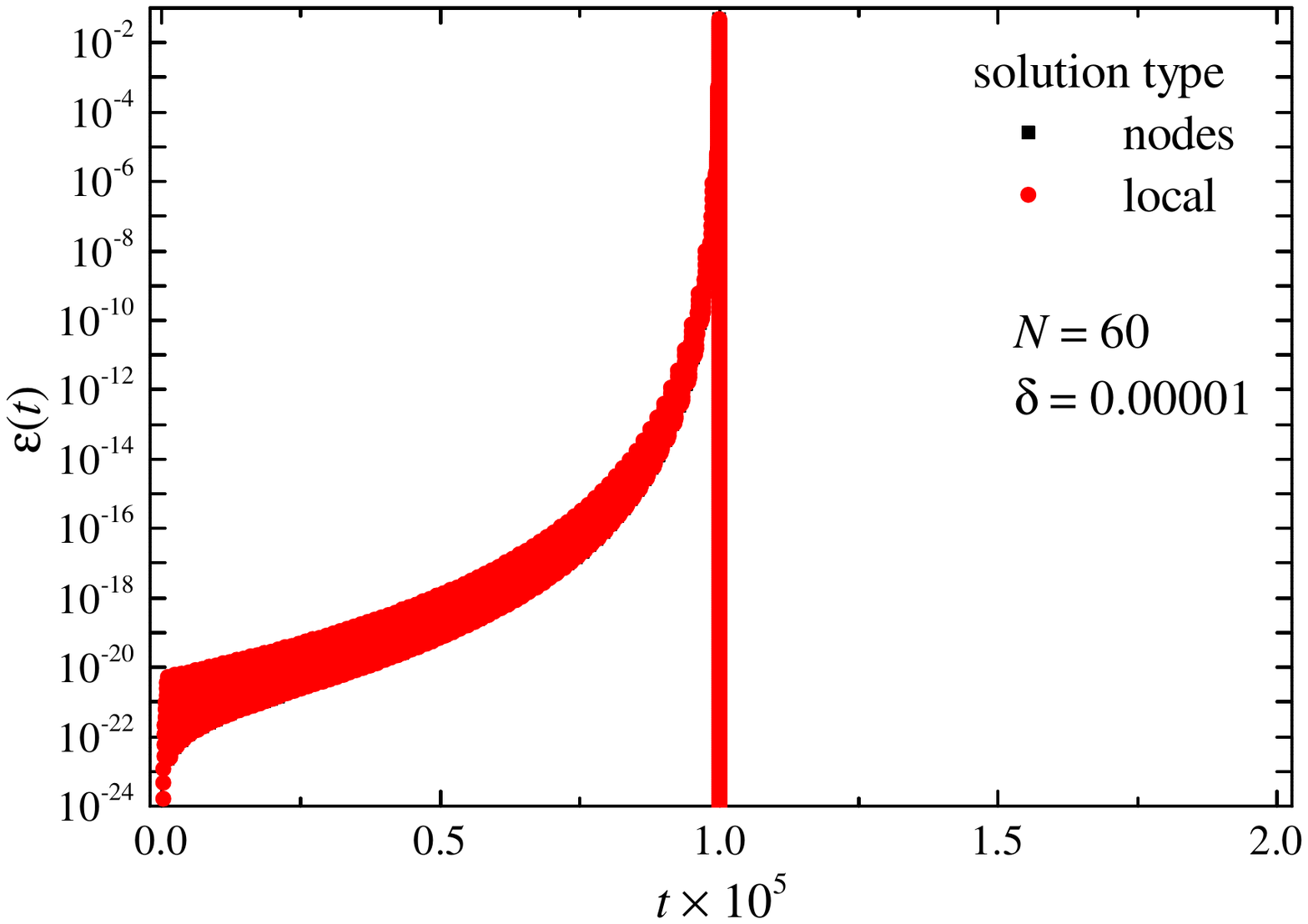}
\vspace{-9mm}\caption{\label{fig:shampine_flame_d_5:e3}}
\end{subfigure}\\
\caption{%
Numerical solution of the problem (\ref{eq:shampine_flame}) with $\delta = 10^{-5}$ in the domain $0 \leqslant t \leqslant 2/\delta$ with $2160$ grid cells: $80$ each in subdomains $t \in [0, 0.495/\delta]$ and $t \in [0.505/\delta, 2/\delta]$ and $2000$ in subdomain $t \in [0.495/\delta, 0.505/\delta]$, using the ADER-DG numerical method with a local DG predictor with degrees of polynomials $N = 1$ (left column: \subref{fig:shampine_flame_d_5:a1}, \subref{fig:shampine_flame_d_5:b1}, \subref{fig:shampine_flame_d_5:c1}, \subref{fig:shampine_flame_d_5:d1}, \subref{fig:shampine_flame_d_5:e1}), $N = 10$ (center column: \subref{fig:shampine_flame_d_5:a2}, \subref{fig:shampine_flame_d_5:b2}, \subref{fig:shampine_flame_d_5:c2}, \subref{fig:shampine_flame_d_5:d2}, \subref{fig:shampine_flame_d_5:e2}) and $N = 60$ (right column: \subref{fig:shampine_flame_d_5:a3}, \subref{fig:shampine_flame_d_5:b3}, \subref{fig:shampine_flame_d_5:c3}, \subref{fig:shampine_flame_d_5:d3}, \subref{fig:shampine_flame_d_5:e3}). Presented data: numerical solution $u(t)$ at nodes in the whole domain $t \in [0, 2/\delta]$ (\subref{fig:shampine_flame_d_5:a1}, \subref{fig:shampine_flame_d_5:a2}, \subref{fig:shampine_flame_d_5:a3}) and in the subdomain $t \in [0.999, 1.001]$ (\subref{fig:shampine_flame_d_5:b1}, \subref{fig:shampine_flame_d_5:b2}, \subref{fig:shampine_flame_d_5:b3}); local solution $u(t)$ in the whole domain $t \in [0, 2/\delta]$ (\subref{fig:shampine_flame_d_5:c1}, \subref{fig:shampine_flame_d_5:c2}, \subref{fig:shampine_flame_d_5:c3}) and in the subdomain $t \in [0.999, 1.001]$ (\subref{fig:shampine_flame_d_5:d1}, \subref{fig:shampine_flame_d_5:d2}, \subref{fig:shampine_flame_d_5:d3}); dependence of the point-wise error $\varepsilon(t)$ of the numerical solution for the solution at the nodes and the local solution (\subref{fig:shampine_flame_d_5:e1}, \subref{fig:shampine_flame_d_5:e2}, \subref{fig:shampine_flame_d_5:e3}); red line is the exact solution (\ref{eq:shampine_flame:exact_solution}).
}
\label{fig:shampine_flame_d_5}
\end{figure}

Numerical methods of ADER-DG are known not only for their high accuracy and the ability to obtain an arbitrarily high order, but also for the ability to solve systems of equations with extremely high stiffness~\cite{ader_stiff_1, ader_stiff_2, ader_stiff_3}. In the works~\cite{ader_stiff_1, ader_stiff_2}, numerical methods based on the ADER paradigm were used to solve stiff problems associated in particular with combustion and detonation. In this work~\cite{ader_stiff_3}, the finite element numerical method ADER-DG with a posteriori correction of the solution by the finite-volume ADER-WENO method was also successfully used to simulate the development and propagation of a one-dimensional detonation wave. Therefore, it is of interest to study the capabilities of the numerical ADER-DG method for solving the initial value problem for a stiff ODE system. In this work, the initial value problem for the equation was chosen as the stiff ODE:
\begin{equation}\label{eq:shampine_flame}
\frac{du}{dt} = u^{2} - u^{3};\qquad u_{1}(0) = \delta;\qquad 0 \leqslant t \leqslant \frac{2}{\delta};
\end{equation} 
where $\delta$ is a small constant. The equation is a very simple model of the fireball balances in a flammable medium~\cite{stiff_flame_1}: the quadratic term $u^{2}$ is related to the area of the combustion domain boundary to which oxygen is supplied, and the cubic term $u^{3}$ determines a value proportional to the volume of the fireball. This problem is well known among stiff ODEs~\cite{stiff_flame_1, stiff_flame_2, stiff_flame_3, stiff_flame_4, stiff_flame_5}. The exact analytical solution of the problem (\ref{eq:shampine_flame}) can be presented in the following form:
\begin{equation}\label{eq:shampine_flame:exact_solution}
u(t) = \frac{1}{W(a\cdot\exp(a - t)) + 1};\quad a = \frac{1}{\delta} - 1;
\end{equation} 
where $W=W(x)$ is the Lambert $W$ function that determines the solution to the transcendental algebraic equation $W\exp(W)=x$. The investigations of the solution to the problem (\ref{eq:shampine_flame}) presented in the works~\cite{stiff_flame_1, stiff_flame_3, stiff_flame_4} and the work~\cite{stiff_flame_2} shows that in the case of small values of parameter $\delta \ll 1$, the equation demonstrate high stiffness in the domain near point $1/\delta$, while case $\delta = 10^{-4}$ can be called stiff, and it is usually used to demonstrate~\cite{stiff_flame_1, stiff_flame_2, stiff_flame_3, stiff_flame_4} the capabilities of the numerical method for solving stiff ODE systems, and case $\delta = 10^{-5}$ is characterized by extreme stiffness. It is known that the solution $u(t)$ has asymptotic behavior $u(t \gg 1/\delta) \simeq 1$, and in the domain near point $1/\delta$ there is a very sharp change in the function from a small value $u \simeq \delta$ to an asymptotic one $u(t \rightarrow \infty) \rightarrow 1$.

A numerical solution to the problem (\ref{eq:shampine_flame}) was obtained using the numerical method ADER-DG using degrees of polynomials $N = 1,\, 2,\, \ldots,\, 60$, which is similar to the numerical solutions of non-stiff problems presented above, and for two parameter $\delta$ values: $\delta = 10^{-4}$ and $\delta = 10^{-5}$. Due to the fact that the areas of occurrence of significant stiffness of the solution is localized in the vicinity of point $1/\delta$, a variable grid step $\Delta t$ was chosen. In case $\delta = 10^{-4}$, the decision definition domain $\Omega = [0, 2/\delta]$ was divided into three subdomains: $\Omega_{1} = [0, 0.4/\delta]$, $\Omega_{2} = [0.4/\delta, 0.6/\delta]$ and $\Omega_{3} = [0.6/\delta, 2/\delta]$, in the $\Omega_{1}$ and $\Omega_{3}$ subdomains cases of $10$, $12$, $15$ and $20$ grid nodes were selected, in the $\Omega_{2}$ subdomain cases of $1000$, $1200$, $1500$ and $2000$ grid nodes were selected. In case $\delta = 10^{-5}$, the decision definition domain $\Omega = [0, 2/\delta]$ was divided into three subdomains: $\Omega_{1} = [0, 0.495/\delta]$, $\Omega_{2} = [0.495/\delta, 0.505/\delta]$ and $\Omega_{3} = [0.505/\delta, 2/\delta]$, in the $\Omega_{1}$ and $\Omega_{3}$ subdomains cases of $40$, $48$, $60$ and $80$ grid nodes were selected, in the $\Omega_{2}$ subdomain cases of $1000$, $1200$, $1500$ and $2000$ grid nodes were selected. The results of solving the problem in case $\delta = 10^{-4}$ are presented in Fig.~\ref{fig:shampine_flame_d_4}, the results of solving the problem in case 2 are presented in Fig.~\ref{fig:shampine_flame_d_5}.

\begin{table}[h!]
\centering
\caption{%
Orders of convergence $p_{L_{1}}$, $p_{L_{2}}$, $p_{L_{\infty}}$ of the ADER-DG numerical method with local DG predictor calculated in norms $L_{1}$, $L_{2}$, $L_{\infty}$ for the stiff problem (\ref{eq:shampine_flame}) with $\delta = 10^{-4}$ and $10^{-5}$; $N$ is the degree of the basis polynomials $\varphi_{p}$ used for the DG representations. The orders of convergence $p$ for the solution at nodes $\mathbf{u}_{n}$ are compared to the theoretical value $p_{\rm (G)} = 2N+1$, the orders of convergence $p$ for the local solution $\mathbf{u}_{L}$ represented between the nodes are compared to the theoretical value $p_{\rm (L)} = N+1$.
}
\label{tab:shampine_flame_d_4}
\begin{tabular}{@{}lllllllll@{}}
\toprule
& \multicolumn{8}{c}{$\delta = 10^{-4}$} \\
\midrule
& \multicolumn{4}{c}{order $p$ for solution at nodes} & \multicolumn{4}{c}{order $p$ for local solution} \\
\midrule
$N$ & $p_{L_{1}}$ & $p_{L_{2}}$ & $p_{L_{\infty}}$ & theor. & $p_{L_{1}}$ & $p_{L_{2}}$ & $p_{L_{\infty}}$ & theor. \\
\midrule
$1$	&	$3.05$	&	$2.91$	&	$2.50$		&	$3$	&	$2.71$	&	$2.21$	&	$1.74$		&	$2$	\\
$2$	&	$7.99$	&	$8.10$	&	$8.03$		&	$5$	&	$3.02$	&	$2.90$	&	$2.59$		&	$3$	\\
$3$	&	$8.15$	&	$8.02$	&	$7.65$		&	$7$	&	$4.21$	&	$3.81$	&	$3.45$		&	$4$	\\
$4$	&	$8.15$	&	$8.02$	&	$7.65$		&	$9$	&	$5.94$	&	$4.83$	&	$4.31$		&	$5$	\\
$5$	&	$8.15$	&	$8.02$	&	$7.65$		&	$11$	&	$6.48$	&	$6.40$	&	$5.62$		&	$6$	\\
$6$	&	$8.15$	&	$8.02$	&	$7.65$		&	$13$	&	$6.50$	&	$6.45$	&	$5.93$		&	$7$	\\
$7$	&	$8.15$	&	$8.02$	&	$7.65$		&	$15$	&	$6.50$	&	$6.45$	&	$5.93$		&	$8$	\\
$8$	&	$8.15$	&	$8.02$	&	$7.65$		&	$17$	&	$6.50$	&	$6.45$	&	$5.93$		&	$9$	\\
$9$	&	$8.15$	&	$8.02$	&	$7.65$		&	$19$	&	$6.50$	&	$6.45$	&	$5.93$		&	$10$	\\
\midrule
& \multicolumn{8}{c}{$\delta = 10^{-5}$} \\
\midrule
& \multicolumn{4}{c}{order $p$ for solution at nodes} & \multicolumn{4}{c}{order $p$ for local solution} \\
\midrule
$N$ & $p_{L_{1}}$ & $p_{L_{2}}$ & $p_{L_{\infty}}$ & theor. & $p_{L_{1}}$ & $p_{L_{2}}$ & $p_{L_{\infty}}$ & theor. \\
\midrule
$1$	&	$2.00$	&	$1.10$	&	$0.04$		&	$3$	&	$1.85$	&	$1.81$	&	$0.75$		&	$2$	\\
$2$	&	$3.35$	&	$3.36$	&	$3.32$		&	$5$	&	$2.58$	&	$2.05$	&	$1.55$		&	$3$	\\
$3$	&	$4.14$	&	$4.02$	&	$3.76$		&	$7$	&	$3.99$	&	$3.40$	&	$2.01$		&	$4$	\\
$4$	&	$4.20$	&	$4.07$	&	$3.77$		&	$9$	&	$4.19$	&	$4.15$	&	$2.48$		&	$5$	\\
$5$	&	$4.21$	&	$4.08$	&	$3.78$		&	$11$	&	$4.22$	&	$4.20$	&	$3.42$		&	$6$	\\
$6$	&	$4.21$	&	$4.08$	&	$3.78$		&	$13$	&	$4.22$	&	$4.20$	&	$3.40$		&	$7$	\\
$7$	&	$4.21$	&	$4.08$	&	$3.78$		&	$15$	&	$4.22$	&	$4.20$	&	$3.40$		&	$8$	\\
$8$	&	$4.21$	&	$4.08$	&	$3.78$		&	$17$	&	$4.22$	&	$4.20$	&	$3.40$		&	$9$	\\
$9$	&	$4.21$	&	$4.08$	&	$3.78$		&	$19$	&	$4.22$	&	$4.20$	&	$3.40$		&	$10$	\\
\bottomrule
\end{tabular}
\end{table}

It should be noted that in the case of non-stiff problems presented above, both the Picard iteration method and the Newton method could be used to solve a system of nonlinear algebraic equations (\ref{dg_pred}) of the local DG predictor, however, in the case of a stiff problem, the use of Picard iterations leads to a divergence of iterations for large grid steps $\Delta t$ and the use of simple initial iteration values~\cite{Jackson_2017}, therefore it is preferable to use Newton’s method for solving systems of nonlinear algebraic equations (\ref{dg_pred}) of the local DG predictor. It should also be noted that the work~\cite{ader_dg_axioms} proposes effective ways to select initial iterations for a system of equations (\ref{dg_pred}).

The results for case $\delta = 10^{-4}$ presented in Fig.~\ref{fig:shampine_flame_d_4} demonstrate the high efficiency of the numerical ADER-DG method for solving the stiff problem (\ref{eq:shampine_flame}). A sharp change in the solution $u(t)$ in the vicinity of a point $1/\delta$ is well resolved both in the solution at nodes and in the local solution. Artifacts of the numerical solution, as presented and described in the works~\cite{stiff_flame_1, stiff_flame_2}, are not observed in this case. The numerical solution adequately corresponds to the exact analytical solution in all considered cases of polynomial degrees $N$. To the right of the stiffness localization region $t \simeq 1/\delta$, the solution quickly grows to the asymptotic value $u(t \rightarrow \infty) \rightarrow 1$, therefore the local error $\varepsilon(t)$ reaches such small values of $\sim 10^{-128}$-$10^{-248}$ compared to the local $\varepsilon(t)$ error on the left, and is not observed in Fig.~\ref{fig:shampine_flame_d_4:e1},~\ref{fig:shampine_flame_d_4:e2} and~\ref{fig:shampine_flame_d_4:e3}. The presented results also demonstrate that as the degree of polynomials $N$ increases, the local error $\varepsilon(t)$ does not decrease significantly, and in the presented cases of degrees $10$ and $60$, differences are not visually observed.

In the case of extreme stiffness $\delta = 10^{-5}$, the results for which are presented in Fig.~\ref{fig:shampine_flame_d_5}, one can see the discrepancy between the numerical solution and the exact analytical solution in the stiffness localization region $t \simeq 1/\delta$ in the case of polynomials of degree $N = 1$. With an increase in the degree of polynomials $N$ or in the case of decrease of the grid step $\Delta t$, this effect, when the region of a sharp gradient of the solution $u(t)$ is to the right of the exact position, is completely disappears. This artifacts was observed in works~\cite{stiff_flame_1, stiff_flame_2} for other numerical methods, and can also be observed for the Radau IIA numerical method, which is included in the \texttt{scipy} module of the programming language \texttt{python}. The presented results for local error $\varepsilon(t)$ demonstrate properties similar to case $\delta = 10^{-4}$. The presented artifacts of the numerical solution are expected~\cite{Butcher_book_2016, Hairer_book_2, stiff_flame_1, stiff_flame_2, stiff_flame_3, stiff_flame_4, stiff_flame_5} for cases of extreme stiffness.

It should be noted that the calculated empirical convergence orders $p$ in this case of the stiff problem (\ref{eq:shampine_flame}) do not have such large values as in the case of non-stiff problems  presented above. Empirical orders of convergence $p$ are presented in Table~\ref{tab:shampine_flame_d_4}. The presented results $p$ clearly demonstrate that in case $\delta = 10^{-4}$ there is a correspondence between the empirical $p$ and theoretical $p_{\rm (L)}$ and $p_{\rm (G)}$ values of order convergence up to polynomial degree $N \leqslant 3$ for solution in nodes and $N \leqslant 4$ for local solution; and in case $\delta = 10^{-5}$ there is no correspondence between the empirical and theoretical values. This is due to the general nature of the solution to the problem -- the solution changes too sharply in the vicinity of the point $1/\delta$ where stiffness occurs.

As a result of the analysis of the obtained numerical solution to the stiff problem, it may be concluded that the numerical method ADER-DG is very well suited for solving stiff and extremely stiff initial value problems for ODE systems.

\subsection{Computational costs}
\label{sec:comp_costs}

The ADER-DG numerical method makes it possible to obtain an arbitrarily high order, which is well achieved for linear and nonlinear non-stiff ODE systems. Of course, it is of some interest to determine the computational costs that accompany obtaining a numerical solution by this method, and their comparison with the computational costs of standard ODE solvers, such as the classical Runge-Kutta methods, which are widely represented in existing numerical software libraries.

Preliminary estimates of computational costs were carried out from the position of estimating the number of calculations (so-called evolutions) of the function $\mathbf{F}(\mathbf{u}, t)$ on the right side of the system of equations (\ref{eq:ivp_chosen_form}). Calculating the local solution $\mathbf{q}(\tau)$ in the ADER-DG method with polynomial degree $N$ as a result of solving a generally nonlinear system of algebraic equations (\ref{dg_pred}) requires at least $(N+1)^{2}$ function calculations for each Picard iteration. The same number of function $\mathbf{F}$ calculations is necessary in the case of each iteration of Newton's method, but in this case, $(N+1)^{2}$ calculations of the Jacobian matrix $\mathrm{J} = \partial\mathbf{F}/\partial\mathbf{u}$ of the function $\mathbf{F}$ are also added, which are combined into a joint $[(N+1) \cdot K] \times [(N+1) \cdot K]$ matrix, with which it is necessary to directly solve a SLAE, calculate the LU decomposition, or calculate the inverse matrix. Calculating a solution at grid nodes also requires $N+1$ function calculations (\ref{eq:sol_in_node}), however, these function values can be calculated and stored during the iteration process when solving a system of algebraic predictor equations (\ref{dg_pred}). In the case of explicit $s$-stage Runge-Kutta methods, the function $\mathbf{F}$ of the right-hand side of the ODE of the system is calculated no more than $s$ times. In the case of $s \leqslant 4$, there are explicit Runge-Kutta methods, the order $p$ of which is equal to the number of stages $s$, but in cases of very high orders the number of stages increases faster than the order of the numerical method~\cite{Butcher_book_2016, Hairer_book_1} -- $8$th order methods are $11$-stage. It is clear that it is not possible to predict the minimum number of stages for explicit Runge-Kutta methods of $121$st order, equal to the $p_{\rm G}$ order of the numerical method ADER-DG, the results for which are presented in this article. However, it is clear that in addition to fairly small values of the degree of polynomials $N$, the computational costs of the ADER-DG method in terms of function $\mathbf{F}$ calculations will be noticeably higher than similar costs for explicit Runge-Kutta methods, which is however difficult to reasonably conclude for very high orders, for which the minimum number of stages is unknown~\cite{Butcher_book_2016}. Also, the use of explicit Runge-Kutta methods leads to the difficulties of solving stiff ODE systems, which are not difficulties for the ADER-DG method; and although there are explicit Runge-Kutta methods with a tunable stability region that can be applied to solving stiff problems, these methods do not have sufficient universality and robustness. Implicit Runge-Kutta methods~\cite{Butcher_book_2016, Hairer_book_2}, which are well used for solving stiff ODE systems, and for which such difficulties do not arise with an increase in the number of stages with increasing order of the method, are conceptually closer in terms of this analysis to ADER-DG methods. Calculating a solution using the implicit $s$-stage Runge-Kutta method in the general case requires obtaining a solution to a nonlinear system of algebraic equations, which requires at least $s$ function $\mathbf{F}$ calculations for each Picard iteration. The same number of function calculations is necessary in the case of each iteration of Newton's method, but in this case $s$ calculations of the Jacobian matrix $\mathrm{J}$ of the function $\mathbf{F}$ are also added, which are arranged into a joint $s K \times s K$ matrix $I \otimes I - \Delta t \cdot A \otimes \mathrm{J}$, with which it is necessary to directly solve a system of linear algebraic equations or calculate the LU decomposition. It is clear that in specific implicit Runge-Kutta methods, these estimates for Newton's method may be incorrect -- DIRK, SDIRK, QSDIRK, EDIRK, ESDIRK, QESDIRK, FSAL-DIRK and others methods require less computations of the function $\mathbf{F}$, Jacobian matrix or LU decomposition compared to the general case~\cite{Butcher_book_2016, Hairer_book_1, Hairer_book_2}. Therefore, the computational costs of the ADER-DG method in terms of function $\mathbf{F}$ and Jacobian matrix $\mathrm{J}$ calculations will be higher than similar costs for implicit Runge-Kutta methods. The work~\cite{ader_dg_axioms} proposes an effective method for calculating the initial iteration for the Picard iterations, based on the use of continuous extension Runge-Kutta (CERK) schemes~\cite{cerk_1992, cerk_2011}, which allows to significantly reduce the number of iterations in this case, however, this will not affect the final asymptotic estimates. Therefore, the computational costs of the ADER-DG method in terms of function $\mathbf{F}$ and Jacobian matrix $\mathrm{J}$ calculations will be higher than similar costs for the classical Runge-Kutta methods. A similar situation will occur when compared with $s$-stage Rosenbrock methods~\cite{Hairer_book_2}, which in this regard can be analyzed similarly.

\begin{figure}[h!]
\captionsetup[subfigure]{%
position=bottom,
textfont=normalfont,
singlelinecheck=off,
justification=raggedright,
font=small
}
\centering
\begin{subfigure}{0.49\textwidth}
\includegraphics[width=\textwidth]{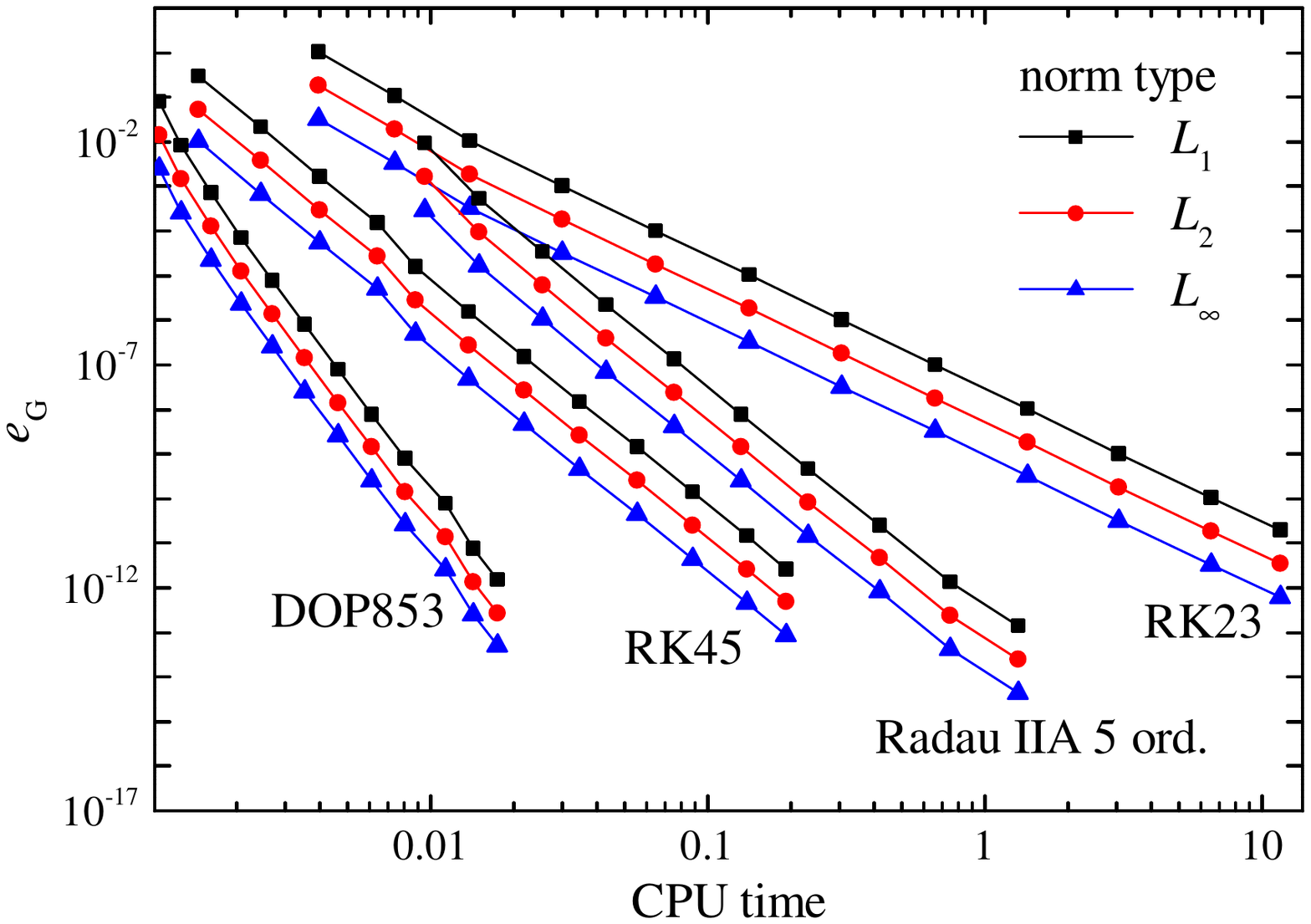}
\vspace{-9mm}\caption{\label{fig:comp_costs:a}}
\end{subfigure}
\begin{subfigure}{0.49\textwidth}
\includegraphics[width=\textwidth]{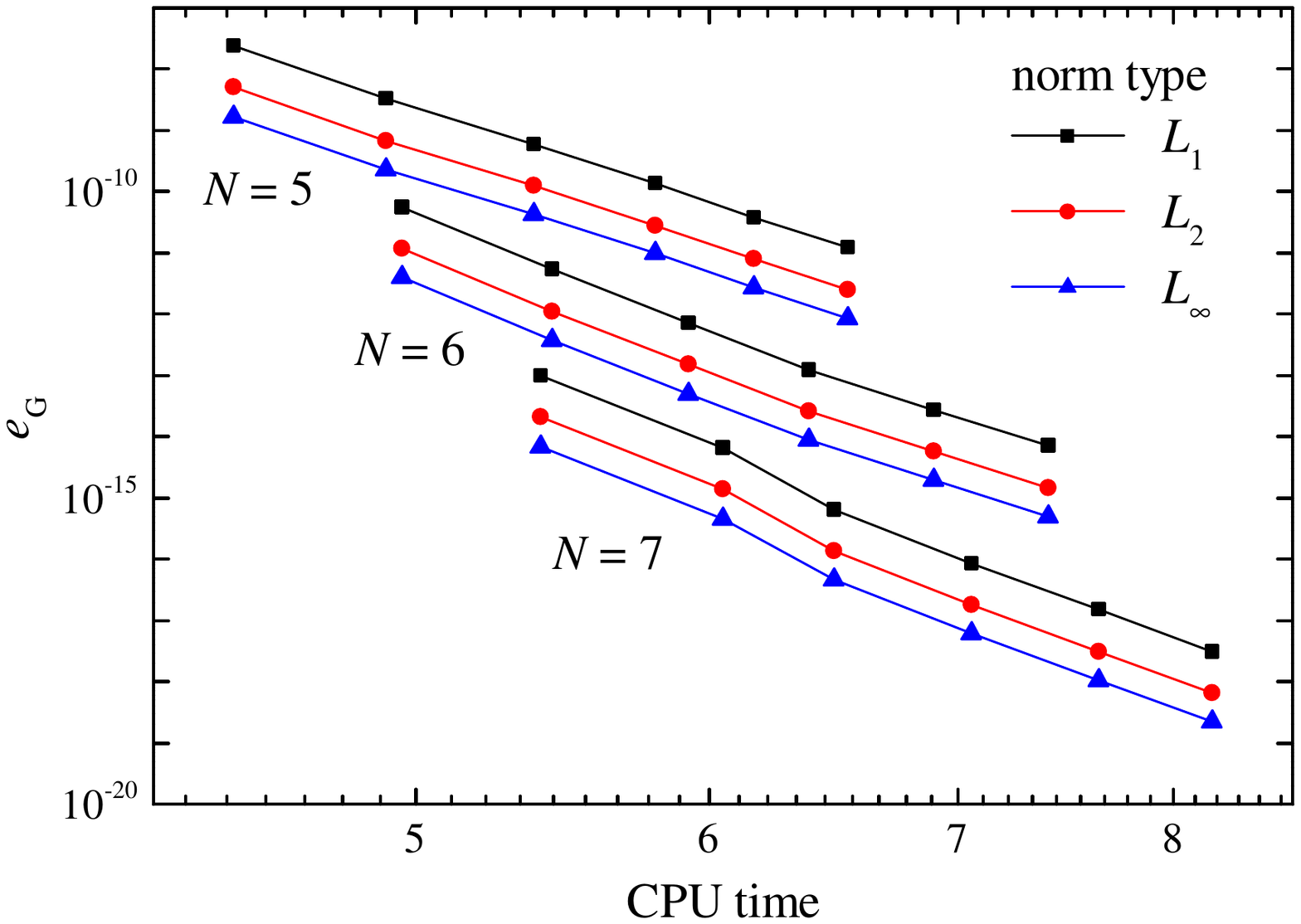}
\vspace{-9mm}\caption{\label{fig:comp_costs:b}}
\end{subfigure}\\
\begin{subfigure}{0.49\textwidth}
\includegraphics[width=\textwidth]{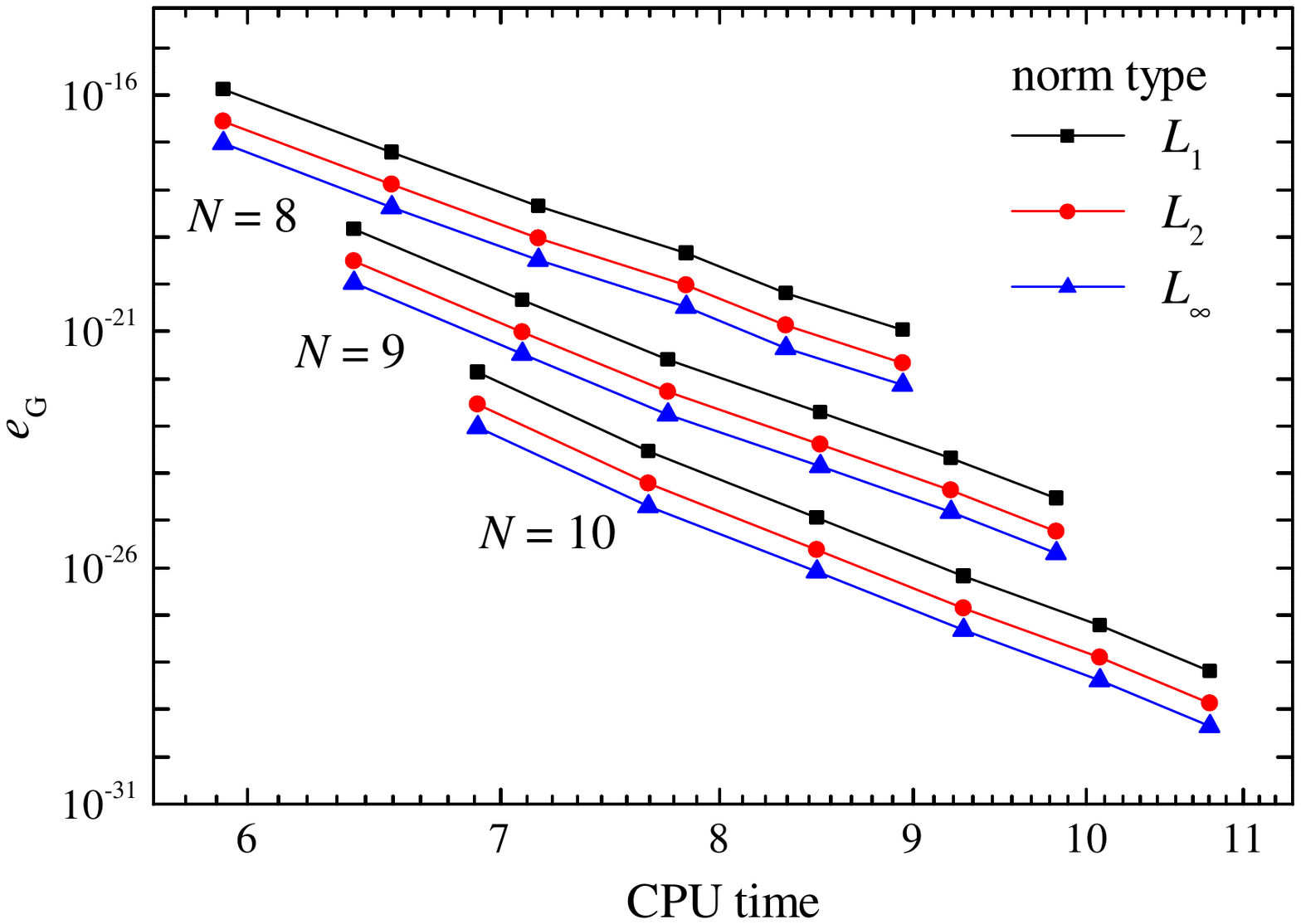}
\vspace{-9mm}\caption{\label{fig:comp_costs:c}}
\end{subfigure}
\begin{subfigure}{0.49\textwidth}
\includegraphics[width=\textwidth]{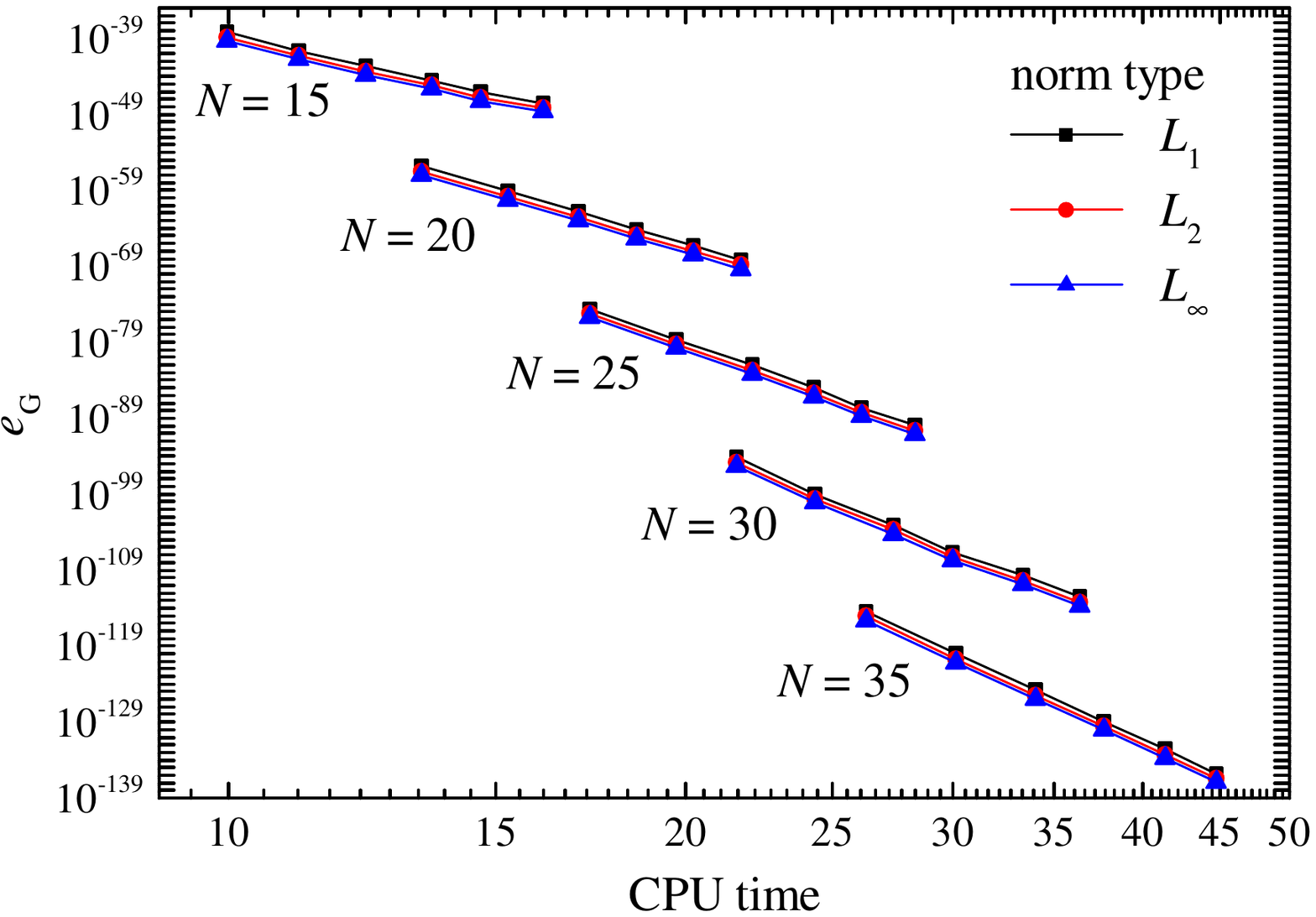}
\vspace{-9mm}\caption{\label{fig:comp_costs:d}}
\end{subfigure}\\
\begin{subfigure}{0.49\textwidth}
\includegraphics[width=\textwidth]{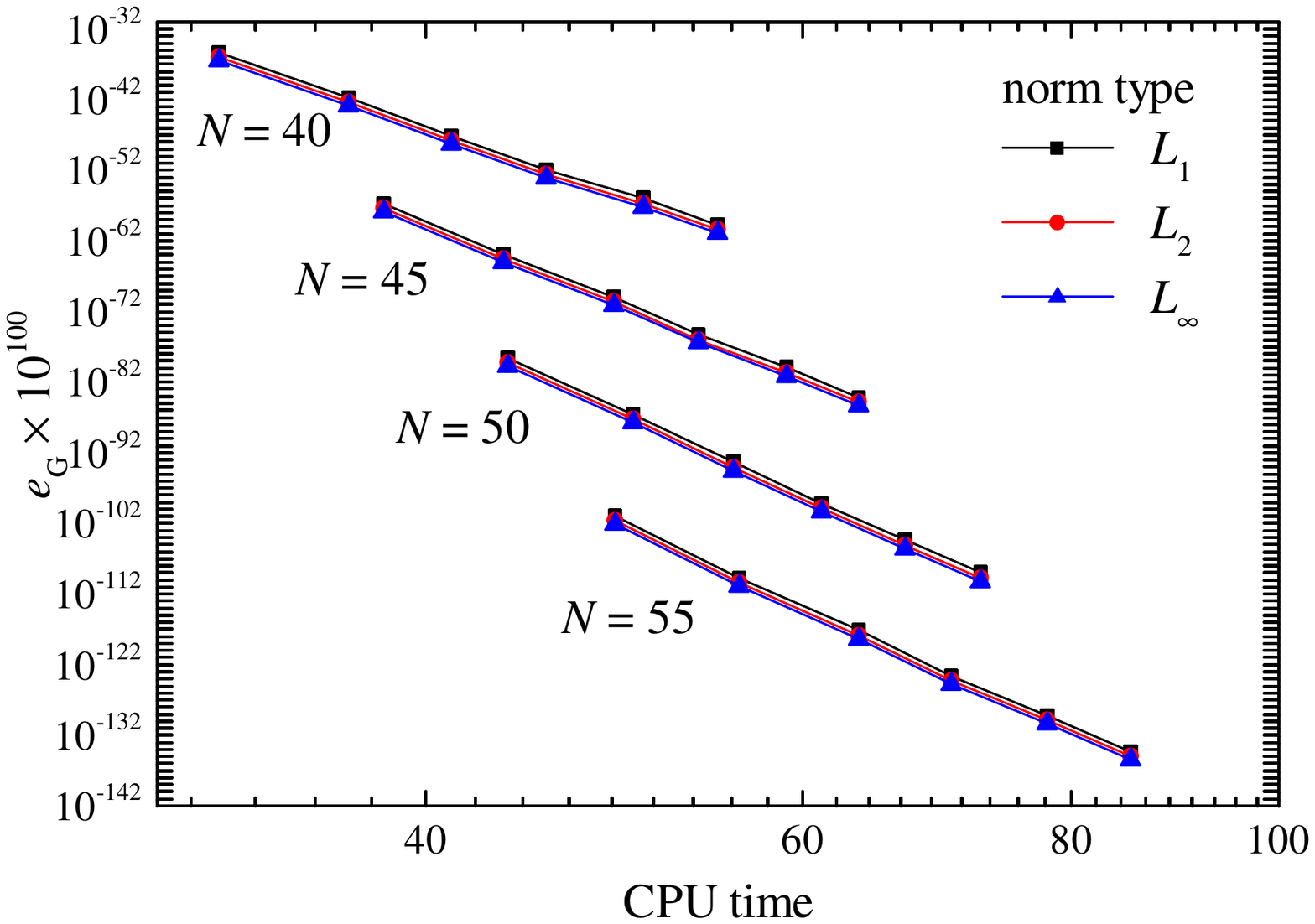}
\vspace{-9mm}\caption{\label{fig:comp_costs:e}}
\end{subfigure}
\begin{subfigure}{0.49\textwidth}
\includegraphics[width=\textwidth]{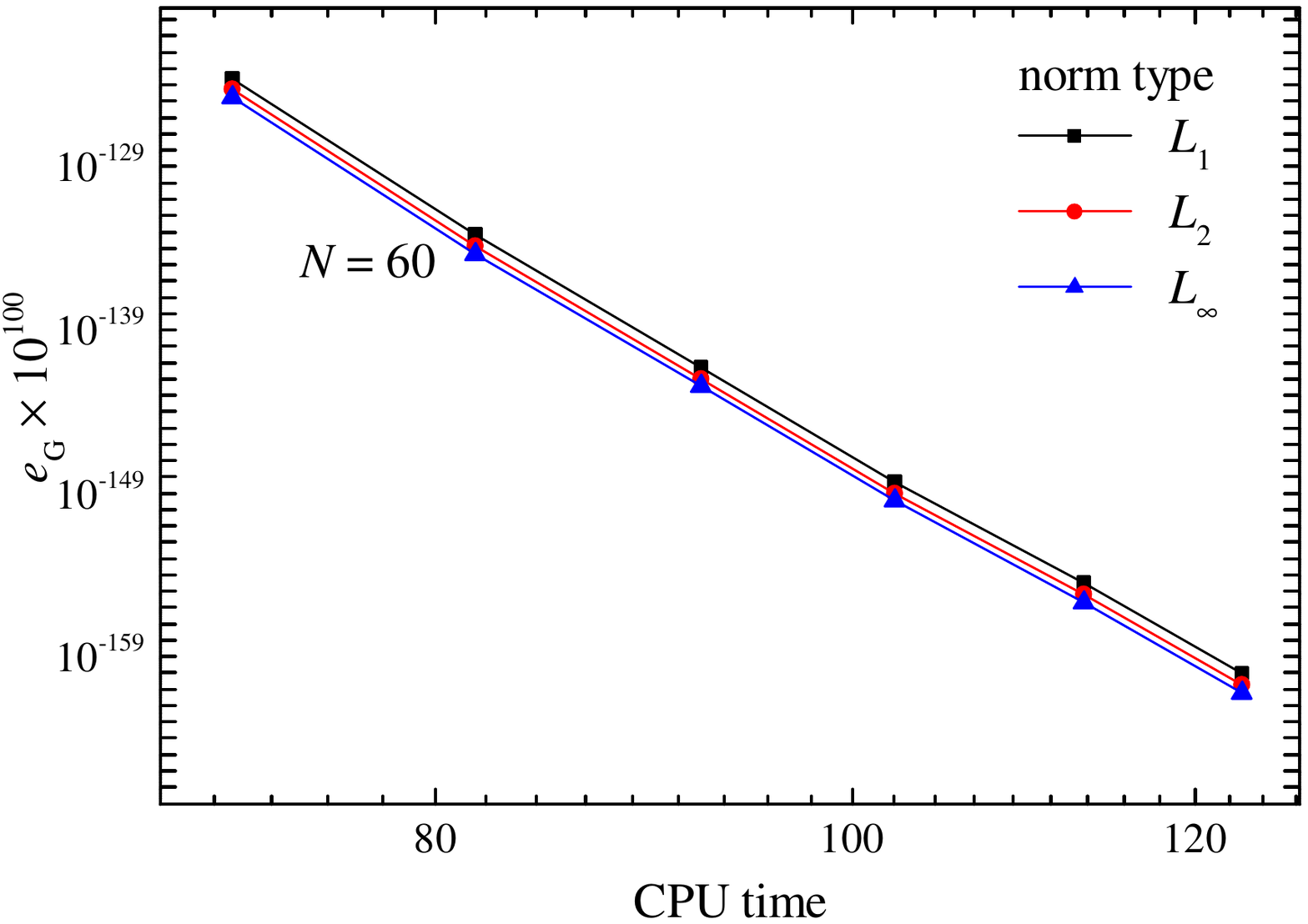}
\vspace{-9mm}\caption{\label{fig:comp_costs:f}}
\end{subfigure}\\
\caption{%
Computational costs, measured in CPU time, of the numerical ADDER-DG method with a local DG predictor for solving the IVP for ODE (\ref{eq:harm_osc_prec}) in domain $t \in [0, 10\pi]$ depending on the error values $e_{\rm G}$ of the numerical solution and comparison with the computational costs of standard Runge-Kutta methods from the implementation in \texttt{scipy}. Presented data: computational costs for the Runge-Kutta numerical methods RK23, RK45, DOP853 and Radau IIA family of order 5 (\subref{fig:comp_costs:a}), computational costs of the ADER-DG method for degrees of polynomials $N = 5$, $6$, $7$ (\subref{fig:comp_costs:b}), $8$, $9$, $10$ (\ref{fig:comp_costs:c}), $15$, $20$, $25$, $30$, $35$ (\subref{fig:comp_costs:d}), $40$, $45$, $50$, $55$ (\subref{fig:comp_costs:e}) and $60$ (\subref{fig:comp_costs:f}). CPU time is measured in seconds. In Figs.~(\subref{fig:comp_costs:e}) and~(\subref{fig:comp_costs:f}) the rescaled error $e_{\rm G} \times 10^{100}$ is presented.
}
\label{fig:comp_costs}
\end{figure}

The computation times of solving the initial value problem for the ODE system (\ref{eq:harm_osc_prec}) in domain $\Omega = [0, 10\pi]$ were obtained. The dependence of the computation time, which was measured as CPU time, for the ADER-DG method on the error $e_{\rm G}$ of the numerical solution is presented in Fig.~\ref{fig:comp_costs} for several values of the degrees of polynomials $N = 1,\, 2,\, \ldots,\, 60$. To compare the results obtained, three classical explicit Runge-Kutta methods -- embedded methods RK23~\cite{rkm_3} of order 3(2) and RK45~\cite{rkm_1, rkm_2} of order 5(4), DOP853~\cite{Hairer_book_1} of order 8, and one implicit Runge-Kutta method -- Radau IIA family method of order 5~\cite{Hairer_book_2}, which are presented in the \texttt{scipy} module of the \texttt{python} programming language, were selected. The dependence of CPU time on the error $e_{\rm G}$ of the numerical solution for selected reference classical numerical methods is also presented in Fig.~\ref{fig:comp_costs}. The presented results for the ADER-DG method were obtained using a software implementation on the \texttt{python} programming language using the \texttt{mpmath} module for using high-precision floating point numbers. The implementation of the classical Runge-Kutta methods in the \texttt{scipy} module of the \texttt{python} programming language was carried out using a binding of high-performance implementations in the \texttt{Fortan} and \texttt{C} languages. It is clear that the high-performance implementation in the core of the \texttt{scipy} module is more efficient compared to a direct implementation in \texttt{python} programming language, however, in terms of comparison in logarithmic scale, this only leads to a shift of the curves in Fig.~\ref{fig:comp_costs} to the right.

A direct comparison of the computational costs obtained shows that the dependence of the CPU time on the error $e_{\rm G}$ has a pronounced power-law dependence, which is characterized by the convergence orders of the methods. The numerical values of the CPU time for the ADER-DG method are greater than the similar values for the classical Runge-Kutta methods with comparable orders $p$. However, increasing the order $p$ of convergence for the ADER-DG method leads to the fact that no adequately small step for classical Runge-Kutta methods can achieve such a low error $e_{\rm G}$ -- in cases of polynomial degrees $N \geqslant 30$-$35$, the order of convergence for the solution at nodes is $p_{\rm (G)} = 2N+1 = 61$-$71$, so the choice of even a coarse grid and large step $\Delta t$ allows one to obtain a solution with such a low error $e_{\rm G}$, which is unattainable for classical Runge-Kutta methods (of course, also using high-precision floating-point numbers), which is due precisely to the order of convergence.

A direct comparison of computational costs in terms of computation time is not a representative result in this particular case. This is due to the fact that the ADER-DG numerical methods make it possible to obtain a numerical solution even on very coarse grids with a large step $\Delta t$, while the local solution $\mathbf{q}(\tau)$ obtained by the local DG predictor allows one to obtain an accurate ``subgrid'' resolution, while classical numerical methods, such as the Runge-Kutta methods, make it possible to obtain a numerical solution at nodes $t_{n}$, and to restore the numerical solution between nodes it is necessary to use some type of reconstruction of the solution -- for example, interpolation, however, in the case of a very coarse grid and large step $\Delta t$, and complex properties of the solution, this does not allow a sufficiently accurate reconstruction of the solution between grid nodes $t_{n}$.

\section*{Conclusion}
\label{sec:conclusion}

In conclusion, it should be noted that in this work an adaptation of the arbitrary high order ADER-DG numerical method with local DG predictor for solving the IVP for a first-order non-linear ODE system is proposed. The implementation of the method is based on the use of the nodal basis of Lagrange interpolation polynomials with nodal points at the nodes of the Gauss-Legendre quadrature formula. This made it possible to use the properties of the nodal basis associated with the lack of calculation of the integrals of the function $\mathbf{F}$ on the right side of the ODE system. An important feature of this variant of the DG method is the significant simplicity of its algorithmic and software implementations, which in no way reduces the accuracy and efficiency of the method. The proposed numerical method is fully one-step ODE solver, with uniform steps in terms of formula implementation. It was shown that the proposed version of the ADER-DG numerical method is $A$-stable and $L$-stable. 

The calculated empirical orders of convergence $p$ of the numerical method for the polynomials degrees $1 \leqslant N \leqslant 60$ are presented. It is shown that the numerical method demonstrates superconvergence with convergence order $p_{\rm (G)} = 2N+1$ for the solution $\mathbf{u}_{n}$ at grid nodes, while the local solution $\mathbf{u}_{L}(t)$ obtained using the local DG predictor has convergence order $p_{\rm (L)} = N+1$, which is usually expected from classical implementations of DG methods for solving ODE, which has been demonstrated in numerical examples. Among the empirical orders of convergence, there are cases when the empirical values significantly exceed the expected values.

It was demonstrated that an important applied feature of this implementation of the numerical method is the possibility of using the local solution as a solution with a subgrid resolution, which makes it possible to obtain a detailed solution even on very coarse coordinate grids. The scale of the error of the local solution $\varepsilon(t)$, when calculating using standard representations of single or double precision floating point numbers, using large values of the degree $N$, practically does not differ from the error of the solution at the grid nodes. It should be noted that the local solution is not some kind of reconstruction of the solution from the set of values at the nodes, it is precisely the representation of the solution in the form of a finite element.

The capabilities of the ADER-DG method for solving stiff ODE systems characterized by extreme stiffness are demonstrated. It was shown that the ADER-DG numerical method is capable of obtaining a numerical solution for the stiff ODE. In this case, a significant decrease in the orders of convergence is observed, which is associated with the sharp change in the solution in the domain of high stiffness. It has been demonstrated that the ADER-DG numerical method is capable of obtaining a numerical solution for a problem with extremely high stiffness.

Estimates of the computational costs of the ADER-DG numerical method are obtained. It was shown that the computational costs of the ADER-DG method in terms of function on the right side of the system of equations and Jacobian matrix calculations will be higher than similar costs for the classical Runge-Kutta methods. A direct empirical computational costs obtained from measurements of CPU time is not a representative result in this particular case -- this is due to the fact that the ADER-DG numerical methods make it possible to obtain a numerical solution even on very coarse grids with a large step, while the local solution allows one to obtain an accurate ``subgrid'' resolution without some type of reconstruction of the solution.

\bmhead{Acknowledgments}
The reported study was supported by the Russian Science Foundation grant No. 21-71-00118 \url{https://rscf.ru/en/project/21-71-00118/}.

The author would like to thank the anonymous reviewers for their encouraging comments and remarks that helped to improve the quality and readability of this paper. The author would like to thank Popova A.P. for help in correcting the English text.

\section*{Declarations}

\subsection*{Data Availability}
The datasets generated during and/or analysed during the current study are available from the corresponding author on reasonable request.

\subsection*{Declaration of competing interest}
The author declares that he has no conflict of interest.





\end{document}